\documentclass[a4paper, abstract=on, sbibliography=totoc,parskip=half, 11pt] {scrreprt}

\usepackage{tikz}
\usepackage[german]{babel}
\usepackage[latin1]{inputenc}
\usepackage[T1]{fontenc}
\usepackage{tocstyle}
\usepackage{amssymb}
\usepackage{amsmath}
\usepackage{amscd}
\usepackage[amsmath, thmmarks]{ntheorem}
\usepackage[all, cmtip]{xy}
\usepackage{enumerate}
\usepackage{cleveref}
\usepackage{array}
\usepackage{chngcntr}
\usepackage{graphicx}
\usepackage{url}
\usepackage{multicol}

\addto{\captionsgerman}{} 

\typearea{12}
\addtokomafont{disposition}{\rmfamily}
\usetocstyle{allwithdot}

\newtheorem{zaehler}{}[chapter]
\newtheorem{satz}[zaehler]{Satz}
\newtheorem{lemma}[zaehler]{Lemma}
\newtheorem{kor}[zaehler]{Korollar}
\newtheorem{prop}[zaehler]{Proposition}
\theorembodyfont{\normalfont}
\newtheorem{bem}[zaehler]{Bemerkung}
\newtheorem{bspl}[zaehler]{Beispiel}
\newtheorem{defn}[zaehler]{Definition}
\theoremsymbol{\ensuremath{\square}}
\theoremseparator{.}
\newtheorem*{bew}{Beweis}

\crefformat{lemma}{Lemma #1}
\crefformat{prop}{Proposition #1}
\crefformat{kor}{Korollar #1}
\crefformat{bem}{Bemerkung #1}
\crefformat{bspl}{Beispiel #1}
\crefformat{defn}{Definition #1}

\DeclareMathOperator{\Rad}{Rad}
\DeclareMathOperator{\Dec}{Dec}

\DeclareMathOperator{\Bild}{Bild}
\DeclareMathOperator{\Ker}{Kern}
\DeclareMathOperator{\Coker}{Kokern}
\DeclareMathOperator{\Ext}{Ext}
\DeclareMathOperator{\End}{End}
\DeclareMathOperator{\Hom}{Hom}
\DeclareMathOperator{\Ht}{ht}
\DeclareMathOperator{\dimv}{\underline{dim}}
\newcommand{\fmod}[1]{\ensuremath{#1\text{\normalfont-mod}}}
\newcommand{\fdar}[1]{\ensuremath{#1\text{\normalfont-dar}}}

\title{\begin{large}Diplomarbeit\end{large}\\\textbf{Über erreichbare und baumartige unzerlegbare Darstellungen von Köchern}
}
\author{Wolfgang Peternell\\Bergische Universität Wuppertal\\ \\Betreuer:\\Prof. Dr. Klaus Bongartz}
\date{Februar 2012}

\begin{document}
\pagenumbering{roman}
\maketitle
\begin{abstract}
Let $A$ be a finite-dimensional algebra over an algebraically closed field. The problem of constructing indecomposable $A$-modules inductively from simple ones by means of exact sequences - called accessibility - is the starting point of the present diploma-thesis. It has lead us to the consideration of exceptional and indecomposable tree-representations of finite quivers. Following Ringel, we prove his result that exceptional representations are tree-representations. We give a detailed description of the various aspects of the Schofield-Induction which plays an important role in the proof. Moreover we introduce a functor (strong hypotheses being given) which enables us to construct indecomposable modules of an algebra from indecomposable representations of a certain bipartite quiver. We also give a proof of Ringel's result that each exceptional representation of dimension $d>1$ of a generalized Kronecker quiver has an indecomposable factor- or subrepresentation of dimension $d-1$. The thesis is concluded by some calculations showing the accessibility of representations of the 3-Kronecker-quiver in small dimensions.
\end{abstract}
\setcounter{page}{1}
\tableofcontents
\newpage

\setcounter{page}{1}
\pagenumbering{arabic}
\chapter*{Einleitung}
\addcontentsline{toc}{chapter}{Einleitung}
\setcounter{zaehler}{0}
\counterwithout{zaehler}{chapter}

Sei $A$ eine endlichdimensionale Algebra über einem algebraisch abgeschlossenen Körper $k$. Alle Moduln seien endlichdimensional. Zwei Resultate von Bongartz über die Beziehung zwischen unzerlegbaren $A$-Moduln der Länge $n$ und solchen der Länge $n-1$ bilden den Ausgangspunkt dieser Arbeit. Der erste Satz gibt -- für spezielle Algebren -- eine sehr genaue Auskunft über diese Beziehung:
\begin{satz}[\protect{\cite{Bong4} und \cite[Theorem 6.1]{Bong2}}]\label{einl:satz1} $A$ sei zahm-verkleidet oder darstellungsendlich. Sei $M$ ein unzerlegbarer, nicht-einfacher $A$-Modul der Länge $n$. Dann besitzt $M$ einen unzerlegbaren Unter- oder Faktormodul der Länge $n-1$.
\end{satz}
In der allgemeinen Situation gilt folgende Aussage:
\begin{satz}[\protect{\cite{Bong1}}]\label{einl:satz2} Sei $M$ ein unzerlegbarer, nicht-einfacher $A$-Modul der Länge $n$. Dann existiert ein unzerlegbarer $A$-Modul der Länge $n-1$.
\end{satz}
Die Frage, ob sich \cref{einl:satz1} sogar auf beliebige Algebren verallgemeinern läßt, ist bis heute offen. Ringel führt dazu in \cite{Ri1} den Begriff der Erreichbarkeit eines $A$-Moduls ein: Demnach sind alle einfachen Moduln erreichbar; ein nicht-einfacher Modul ist genau dann erreichbar, wenn er unzerlegbar ist und einen erreichbaren Unter- oder Faktormodul der Länge $l(M)-1$ besitzt. Es stellt sich also die Frage, ob alle unzerlegbaren $A$-Moduln erreichbar sind. Als eine erste Annäherung an dieses Problem gibt Ringel die folgende Verallgemeinerung von \cref{einl:satz2} an:
\begin{satz}[\protect{\cite{Ri1}}]\label{einl:satz3} Sei $M$ ein unzerlegbarer $A$-Modul der Länge $n$. Dann existiert auch ein erreichbarer $A$-Modul der Länge $n$.
\end{satz}
Ringel beweist diese Aussage für nicht-distributive Algebren, indem er die Argumente von Bongartz in \cite[Abschnitt 1]{Bong1} für den elementaren nicht-distributiven Fall von \cref{einl:satz2} verallgemeinert. Für den distributiven Fall verweist Ringel auf den Hauptteil von \cite{Bong1}. Wir wollen ein letztes allgemeines Resultat nicht verschweigen: 
\begin{satz}[\protect{\cite[Corollary 2.2]{Ri6}}]\label{einl:satz4} Sei $M$ ein unzerlegbarer, nicht-einfacher $A$-Modul der Länge $n$. Dann existiert eine exakte Sequenz
\[0\longrightarrow U_1\longrightarrow M\longrightarrow U_2\longrightarrow 0\]
von $A$-Moduln, wobei die $U_i$ unzerlegbar sind.
\end{satz}\enlargethispage{\baselineskip}
Der Beweis von Ringel ist verblüffend einfach und verwendet den Begriff des Gabriel-Roiter-Maßes. In dieser Terminologie ist $U_1$ ein echter Untermodul von $M$ von maximalem Gabriel-Roiter-Maß. Wir konnten diesen vielversprechenden Ansatz leider nicht weiterverfolgen.\newpage

Es liegt nahe, die Frage der Erreichbarkeit für den Fall einer Köcheralgebra $A=kK$ eines endlichen, zykellosen, wilden Köchers $K$ zu untersuchen. Wir nennen eine Wurzel $\underline{d}\in\mathbb{Z}^{K_0}$ erreichbar, falls eine erreichbare $K$-Darstellung $M$ mit $\dimv M=\underline{d}$ existiert und stellen die Frage, ob jede Wurzel von $K$ erreichbar ist. Vom geometrischen Standpunkt aus betrachtet sollten dabei zunächst die unzerlegbaren $K$-Darstellungen ohne Selbsterweiterungen von besonderem Interesse sein, da sie in den entsprechenden Darstellungsvarietäten generisch sind. Solche Darstellungen heißen auch exzeptionell; die Dimensionsvektoren exzeptioneller Darstellungen sind gerade die reellen Schur-Wurzeln von K.

Um das Problem weiter einzuschränken, betrachten wir in dieser Arbeit den scheinbar übersichtlichsten Fall wilder Köcher, nämlich den der $n$-Kronecker-Köcher $Q_n$ mit zwei Punkten $x_1$ und $x_2$ und $n$ Pfeilen von $x_1$ nach $x_2$ für $n\ge3$. In diesem Fall ist jede reelle Wurzel schon eine Schur-Wurzel und die exzeptionellen Darstellungen sind gerade die präprojektiven und präinjektiven Moduln über der zugehörigen Köcheralgebra $\Lambda_n$. Wir haben unsere Untersuchungen mit einigen konkreten Rechnungen und geometrischen Betrachtungen für unzerlegbare Darstellungen des $3$-Kronecker-Köchers der Dimension $\leq6$ begonnen, die in Kapitel \ref{rech} zusammengefaßt sind. 

Leider haben sich diese Rechnungen als wenig erhellend erwiesen, weswegen wir uns in einem zweiten Schritt den exzeptionellen $Q_n$-Darstellungen zugewandt haben. Diese Darstellungen kann man stets zu Darstellungen der universellen Überlagerung $\widetilde{Q}_n$ liften. Außerdem lassen sich sowohl die exzeptionellen $Q_n$-Darstellungen als auch ihre Liftungen mit Hilfe von Spiegelungsfunktoren aus einfachen Darstellungen gewinnen. Um die Erreichbarkeit einer exzeptionellen Darstellung zu verifizieren, genügt es, die Erreichbarkeit ihrer Liftung zu überprüfen. Der Vorteil der Liftung besteht darin, daß sie einen größeren Träger hat und damit in gewissem Sinne dünner ist als die zugehörige $Q_n$-Darstellung. Durch einfache kombinatorische Betrachtungen kann man so einige interessante exakte Sequenzen für exzeptionelle $Q_n$-Darstellungen konstruieren; von einem Beweis der Erreichbarkeit ist man dann allerdings noch weit entfernt.

Diese Vorgehensweise wird in den Abschnitten \ref{baum:exk} und \ref{baum:kron} genau beschrieben. Dort findet der Leser auch die Konstruktion der eben erwähnten exakten Sequenzen für exzeptionelle $n$-Kronecker-Darstellungen (vgl. \cref{baum:kron:satz1}). Dieses Ergebnis stammt von Ringel und ist die einzige allgemeine Aussage der vorliegenden Arbeit, die in Richtung der Erreichbarkeit exzeptioneller $n$-Kronecker-Darstellungen zielt. Genauer besagt sie, daß jeder präprojektive, nicht-einfache Modul der Dimension $l$ über der $n$-Kronecker-Algebra $\Lambda_n$ einen $l-1$-dimensionalen, unzerlegbaren Faktormodul besitzt. Wir konnten mit einiger Mühe verifizieren, daß die ersten vier Moduln der präprojektiven Komponente der $3$-Kronecker-Algebra $\Lambda_3$ auch erreichbar sind (vgl. Abschnitt \ref{baum:erreich}), haben durch unsere Rechnungen aber kein allgemeines Verfahren gefunden, um die Erreichbarkeit aller präprojektiven $\Lambda_3$-Moduln zu zeigen. Nach unserer Ansicht sind diese Rechnungen für ein Verständnis des Problems aber aufschlußreicher als unsere Analysen in Kapitel \ref{rech}, und wir erhalten immerhin ein Beispiel für eine $76$-dimensionale erreichbare Darstellung.

Das eben erwähnte Ergebnis von Ringel ist ein Nebenprodukt seiner Untersuchungen von Baumdarstellungen von Köchern; dies sind solche Darstellungen, die einen baumartigen Koeffizientenköcher besitzen (vgl. Abschnitt \ref{baum:beg} für die Definitionen). Unser Interesse für die reellen Schur-Wurzeln eines Köchers -- bzw. für die zugehörigen exzeptionellen Darstellungen -- haben uns motiviert, den folgenden Satz von Ringel in diese Arbeit aufzunehmen und ausführlich zu beweisen -- wir folgen in Abschnitt \ref{baum:hauptres} dem eleganten Beweis von Ringel in \cite{Ri4}.
\begin{satz}[\protect{\cite{Ri3}}]\label{einl:satz5} Sei $K$ ein endlicher Köcher und sei $M$ eine exzeptionelle $K$-Darstellung. Dann ist $M$ eine Baumdarstellung.
\end{satz}
Dabei stehen weniger das Resultat selbst als vielmehr die im Beweis verwendeten Methoden im Zentrum unseres Interesses. Das Verfahren der Schofield-Induktion zur induktiven Beschreibung von exzeptionellen Köcherdarstellungen ist ein Kernstück des Beweises. Der allgemeine Satz besagt, daß für jede exzeptionelle $K$-Darstellung $M$ eine volle Unterkategorie $\mathcal{C}$ von $\fdar{K}$ mit folgenden Eigenschaften existiert: $\mathcal{C}$ ist abgeschlossen unter Kernen, Kokernen und Erweiterungen, $M$ ist nicht einfach in $\mathcal{C}$ und $\mathcal{C}$ ist äquivalent zur Darstellungskategorie eines $n$-Kronecker-Köchers (vgl. \cref{scho:scho:satz1}). Man kann die einfachen Objekte von $\mathcal{C}$ überraschend leicht direkt angeben (vgl. \cref{scho:scho:bem4}), allerdings lassen sich die geforderten Eigenschaften nicht aus der Konstruktion ablesen. Da die verschiedenen erforderlichen Resultate in der Literatur ein wenig verstreut sind (vgl. \cite{CB2}, \cite{Ri2}, \cite{Scho}), geben wir in Kapitel \ref{scho} einen ausführlichen Überblick zu diesem Thema.

Im Beweis von \cref{einl:satz5} müssen außerdem Koeffizientenköcher von Mitteltermen von exakten Sequenzen aus den Koeffizientenköchern der Randterme zusammengeklebt werden. In Kapitel \ref{unz} konstruieren wir dazu einen geeigneten Funktor für endliche Köcher $K$, mit dessen Hilfe unter speziellen Voraussetzungen aus unzerlegbaren Baumdarstellungen von $n$-Kronecker-Köchern unzerlegbare Baumdarstellungen von $K$ konstruiert werden können (vgl. dazu auch \cref{baum:hauptres:lemma3} und \cref{baum:hauptres:kor1}). Die Idee zu diesem Kapitel stammt von Bongartz und scheint in dieser Form noch nicht in der Literatur beschrieben worden zu sein.

Es sollte betont werden, daß im Beweis des Satzes von Ringel an keiner Stelle Koeffizientenköcher von Darstellungen konkret konstruiert werden. In \cite{Ri3} gibt Ringel allerdings Matrixpräsentationen der präprojektiven Darstellungen der $n$-Kronecker-Köcher an, für die man induktiv Dimensionsvektoren von präprojektiven Darstellungen berechnen muß. Für unsere Untersuchungen von exzeptionellen $Q_3$-Darstellungen konnten wir diese Matrixpräsentationen jedoch nicht verwenden. Auch bei Kenntnis einer Matrixpräsentation oder eines baumartigen Koeffizientenköchers einer Baumdarstellung hat man im allgemeinen keine Möglichkeit, die Erreichbarkeit der Darstellung zu zeigen (etwa durch sukzessive Eliminierung aller Blätter des Koeffizientenköchers).

Neben der bereits angesprochenen Konstruktion von exakten Sequenzen für exzeptionelle $n$-Kronecker-Darstellungen wollen wir noch einen weiteren Schnittpunkt des Problems der Erreichbarkeit mit der Theorie der Baummoduln erwähnen. Eine Vermutung von Ringel über Baummoduln besagt, daß zu jeder Wurzel $\underline{d}$ eines Köchers $K$ eine Baumdarstellung $M$ mit $\dimv M=\underline{d}$ existiert. Für den Fall der $n$-Kronecker Köcher wurde diese Vermutung zuerst von Weist in \cite{Weist2} bewiesen. Ringel gibt für diesen Fall in \cite{Ri4} einen alternativen Beweis an. Wieder basiert unser Interesse nicht in erster Linie auf dem Resultat, sondern auf der Vorgehensweise im Beweis, weswegen wir kurz auf die Details eingehen: Sei $\widetilde{Q}_n$ die universelle Überlagerung von $Q_n$ und $F_\lambda:\fdar{Q_n}\longrightarrow \fdar{\widetilde{Q}_n}$ der Push-Down-Funktor. Betrachte die folgenden Klassen von $Q_n$-Darstellungen:
\begin{align*}
\mathcal{D}&:=\{M\mid M=F_\lambda(\widetilde{M})\text{ für eine dünne, unzerlegbare $\widetilde{Q}_n$-Darstellung $\widetilde{M}$}\},\\
\mathcal{E}&:=\{M\mid M=F_\lambda(\widetilde{M})\text{ für eine exzeptionelle $\widetilde{Q}_n$-Darstellung $\widetilde{M}$}\}.
\end{align*}
Man sieht leicht, daß $\mathcal{D}\subseteq\mathcal{E}$ ist und daß die Klasse $\mathcal{E}$ abgeschlossen ist unter der Anwendung der Spiegelungsfunktoren von $Q_n$ (vgl. dazu etwa die elementare Bemerkung \ref{baum:exk:bem2}). Außerdem folgt mit \cref{einl:satz5}, daß alle Moduln in $\mathcal{E}$ Baummoduln sind. Ringel berechnet nun in \cite{Ri4} die Menge $\underline{\mathcal{D}}$ der Dimensionsvektoren aller $Q_n$-Darstellungen $M$ mit $M\in\mathcal{D}$. $\underline{\mathcal{D}}$ enthält einen Fundamentalbereich der Operation der Weyl-Gruppe auf den imaginären Wurzeln von $Q_n$ (vgl. etwa \cite[S. 75]{Kac}). Also existiert zu jeder imaginären Wurzel $\underline{d}$ eine Darstellung $M\in\mathcal{E}$ mit $\dimv M=\underline{d}$ und $M$ ist eine Baumdarstellung. Für reelle Wurzeln ist die Aussage klar, da im Fall der Kronecker-Algebren $Q_n$ die reellen Wurzeln gerade die Dimensionsvektoren exzeptioneller Darstellungen sind.

Es liegt nahe, für einen Beweis der Erreichbarkeit aller Wurzeln des $n$-Kronecker-Köchers analog vorzugehen. Genauer reduziert sich unsere Ausgangsfrage auf das folgende Problem: Zeige, daß die Klasse der erreichbaren Moduln abgeschlossen ist unter der Anwendung von Spiegelungsfunktoren. Aus dieser Aussage folgt die Erreichbarkeit: Klar ist, daß dann alle reellen Wurzeln von $Q_n$ erreichbar sind. Außerdem ist klar, daß alle Wurzeln in $\underline{\mathcal{D}}$ erreichbar sind; ist die Klasse der erreichbaren Moduln also abgeschlossen unter Anwendung der Spiegelungsfunktoren, so sind schon alle imaginären Wurzeln erreichbar.

Unsere Rechnungen in Abschnitt \ref{baum:erreich} sind von der Frage, ob erreichbare Darstellungen unter Spiegelungen erhalten bleiben, inspiriert. \cref{baum:erreich:lemma3} behandelt einen Spezialfall, der unter starken Voraussetzungen an die dort gegebene Situation angepaßt ist.

Die Ringel-Vermutung wurde von Crawley-Boevey verschärft. Sei dazu ein endlicher Köcher $K$ gegeben und sei $\underline{d}\in\mathbb{Z}^{K_0}$ eine Wurzel von $K$. Sei weiter $X:=\operatorname{Rep}_K(\underline{d})$ die Varietät der Darstellungen zum Dimensionsvektor $\underline{d}$. Für alle $i\geq1$ ist die Menge $U^{(i)}\subseteq X$ aller unzerlegbaren $K$-Darstellungen $V$ mit $\dim\End_K(V)=i$ lokal abgeschlossen in $X$. Die irreduziblen Komponenten der $U^{(i)}$ nennen wir die Schichten von Unzerlegbaren in $X$. Die Menge $U=\bigcup_{i\geq1}U^{(i)}$ ist dann konstruierbar. Man sieht leicht, daß die Menge $\mathcal{B}\subseteq X$ der unzerlegbaren Baumdarstellungen in $X$ Vereinigung von endlich vielen Bahnen in $X$ ist. $\mathcal{B}$ ist also ebenfalls konstruierbar. Die Ringel-Vermutung besagt, daß $\mathcal{B}\cap U\neq\emptyset$ ist. Die naheliegende Präzisierung lautet dann: Für alle $i\geq1$ und für jede irreduzible Komponente $S$ von $U^{(i)}$ ist $\mathcal{B}\cap S\neq\emptyset$. Die Berechnungen der Varietäten $U^{(i)}$ für den $3$-Kronecker-Köcher und den Dimensionsvektor $(3,3)$ in Kapitel \ref{rech} waren von dieser Frage motiviert; in diesem Fall sind alle $U^{(i)}$ irreduzibel und wir konnten für jede Schicht von Unzerlegbaren in Abschnitt \ref{rech:baum} eine Baumdarstellung angeben.\enlargethispage{10mm}

Ich danke an dieser Stelle meinem Betreuer Herrn Prof. Dr. Klaus Bongartz für seine zahlreichen Hinweise, Ermutigungen und für seine Geduld.
\chapter{Bezeichnungen}\label{bez}
\setcounter{zaehler}{0}
\numberwithin{zaehler}{chapter}
Für eine Algebra $A$ über einem Körper $k$ bezeichnen wir mit $\fmod{A}$ die Kategorie aller endlichdimensionalen $A$-Linksmoduln. Mit $A^{op}$ bezeichnen wir die Gegenalgebra von $A$ und mit $J(A)$ das Jacobson-Radikal von $A$. Alle betrachteten Moduln sind Linksmoduln. Ist $M$ ein $A$-Modul, so bezeichnen wir für alle $a\in A$ mit $M(a)\in \End_k(M)$ die Linksmultiplikation mit $a$.

Ist $\mathcal{A}$ eine additive Kategorie, so sei $\mathcal{R}_\mathcal{A}$ das Radikal von $\mathcal{A}$. Für eine Algebra $A$ schreiben wir kurz $\mathcal{R}_A:=\mathcal{R}_{\fmod{A}}$. Für ein Objekt $X$ in $\mathcal{A}$ ist $\operatorname{add}(X)$ die volle Unterkategorie von $\mathcal{A}$ bestehend aus allen direkten Summanden aller Potenzen von $X$.

Sei $K$ ein Köcher. Wir bezeichnen mit $K_0$ die Menge aller Punkte und mit $K_1$ die Menge aller Pfeile in $K$. Ist $\alpha\in K_1$, so ist $n\alpha\in K_0$ der Start von $\alpha$ und $s\alpha\in K_0$ die Spitze von $\alpha$. Wir schreiben dafür auch $\alpha:n\alpha\longrightarrow s\alpha$. $K$ heißt endlich, falls $K_0$ und $K_1$ endlich sind. Ist $K$ endlich, so sei $kK$ die Köcheralgebra von $K$ über $k$. Für $x\in K_0$ sei $e_x$ das Idempotent zum Punkt $x$. Schließlich sei $K^{op}$ der Gegenköcher von $K$ und $\left|K\right|$ der (ungerichtete) Graph zu $K$. Ein Baum ist ein zusammenhängender zykelloser (ungerichteter) Graph; ein Köcher heißt baumartig, falls sein zugrundeliegender Graph ein Baum ist.

Für einen Körper $k$ bezeichnen wir mit $K$-$_k$dar die Kategorie aller $K$-Darstellungen $M$ über $k$, so daß $\sum_{x\in K_0}\dim M(x)<\infty$ gilt und nennen eine solche Darstellung endlichdimensional. Wir zeichnen für jeden Punkt $x\in K_0$ einige Darstellungen aus: $P_x$ sei der Projektive, $I_x$ der Injektive und $E_x$ der Einfache zum Punkt $x$. $\dimv M$ sei der Dimensionsvektor von $M$. Mit $Tr(M)$ bezeichnen wir den Träger von $M$, also die Menge aller $x\in K_0$, so daß $M(x)\neq 0$ ist. Meist fassen wir $Tr(M)$ als vollen Unterköcher von $K$ auf. 

Für eine natürliche Zahl $n\geq0$ bezeichnen wir mit $Q_n$ den $n$-Kronecker-Köcher; also den Köcher mit zwei Punkten $x_1,x_2$ und $n$ Pfeilen von $x_1$ nach $x_2$. Die zugehörige $n$-Kronecker-Algebra sei $\Lambda_n:=kQ_n$.

Eine Operation einer Gruppe $G$ auf einer Menge $X$ ist bei uns immer eine Linksoperation. Mit $\operatorname{Iso}_G(x)$ bezeichnen wir die Isotropiegruppe eines Elementes $x\in X$. Mit $\operatorname{GL}_n(k)$ bezeichnen wir die allgemeine lineare Gruppe vom Grad $n$ über $k$.
\chapter{Grundlagen}
\setcounter{zaehler}{0}
\numberwithin{zaehler}{chapter}

\section{Varietäten von Köcherdarstellungen}\label{grund:var}
\setcounter{zaehler}{0}
\numberwithin{zaehler}{section}
Sei $k$ ein algebraisch abgeschlossener Körper. Alle topologischen Aussagen betreffen die Zariski-Topologie. Das Material dieses Abschnitts stammt hauptsächlich aus \cite{Bong5}. Unsere Resultate für Varietäten von Köcherdarstellungen sind dort für Varietäten von Moduln endlich erzeugter Algebren angegeben. Beide Konzepte sind über assoziierte Faserbündel eng miteinander verwandt (vgl. \cite[Example 5.18]{Bong5}) und alle Aussagen lassen sich auf unsere Situation übertragen.

Sei $G$ eine algebraische Gruppe und $X$ eine $G$-Varietät. Der Vollständigkeit halber fassen wir einige grundlegende Ergebnisse zusammen:
\begin{prop}\label{grund:var:prop1}\   
\begin{enumerate}
\item Für alle $x\in X$ ist $Gx$ lokal abgeschlossen in $X$, der Abschluß $\overline{Gx}$ von $Gx$ in $X$ ist $G$-stabil und die $G$-Varietät $\overline{Gx}-Gx$ besteht aus lauter Orbiten $Gy$ mit $\dim Gy< \dim Gx$. Es gilt die Formel
\[\dim Gx=\dim G-\dim\operatorname{Iso}_G(x).\]
\item Für alle natürlichen Zahlen $i\geq0$ ist die Menge aller $x\in X$ mit $\dim Gx\geq i$ offen in $X$.
\end{enumerate}
\end{prop}
Wegen Teil (2) der Proposition ist für alle $i\geq0$ die Menge $X_i$ aller Punkte $x$ von $X$ mit $\dim Gx=i$ lokal abgeschlossen in $X$. Die irreduziblen Komponenten der $X_i$ heißen die Schichten von $X$. Ist $X$ irreduzibel, so ist die Vereinigung aller $G$-Orbiten maximaler Dimension die eindeutig bestimmte offene Schicht in $X$. 

Wir schreiben für Elemente $x,y\in X$
\[x\preceq_{deg} y :\Leftrightarrow Gy\subseteq \overline{Gx}\]
und sagen dann, daß $x$ in $y$ entartet. Diese Relation induziert eine partielle Ordnung auf der Menge der $G$-Bahnen von $X$, die wir wir ebenfalls mit $\preceq_{deg}$ bezeichnen.

Sei $K$ ein endlicher Köcher und sei $\underline{d}=(d(x))_{x\in K_0}\in\mathbb{Z}^{K_0}$ ein Dimensionsvektor. Alle betrachteten $K$-Darstellungen seien endlichdimensional.
\begin{defn}\label{grund:var:def1}\ 
\begin{enumerate}\enlargethispage{\baselineskip}
\item Die Varietät $\operatorname{Rep}_K(\underline{d})$ der $K$-Darstellungen mit Dimensionsvektor $\underline{d}$ über $k$ ist definiert als der affine Raum
\[\operatorname{Rep}_K(\underline{d}):=\prod_{\alpha\in K_1}k^{d(s\alpha)\times d(n\alpha)}.\]\newpage
\item Die lineare algebraische Gruppe $G_K(\underline{d})$ ist gegeben durch
\[G_K(\underline{d}):=\prod_{x\in K_0}\operatorname{GL}_{d(x)}.\]
\item $G_K(\underline{d})$ operiert algebraisch auf $\operatorname{Rep}_K(\underline{d})$ via der Vorschrift
\[g*m:=\bigl(g_{s\alpha}\cdot m_\alpha\cdot g_{n\alpha}^{-1}\bigr)_{\alpha\in K_1}\text{ für alle } g\in G_K(\underline{d}),\ m\in\operatorname{Rep}_K(\underline{d}).\]
\end{enumerate}
\end{defn}
Wir bezeichnen für Elemente $m\in\operatorname{Rep}_K(\underline{d})$ die zugehörige $K$-Darstellung mit $\operatorname{rep}m$. Die Zuordnung $m\mapsto \operatorname{rep} m$ induziert eine Bijektion zwischen der Menge der $G_K(\underline{d})$-Orbiten von $\operatorname{Rep}_K(\underline{d})$ und der Menge der Isomorphieklassen aller $K$-Darstellungen mit Dimensionsvektor $\underline{d}$. Ist $M$ eine $K$-Darstellung, so bezeichnen wir mit $\mathcal{O}(M)$ den zugehörigen $G_K(\dimv M)$-Orbit in $\operatorname{Rep}_K(\dimv M)$. Sind zwei $K$-Darstellungen $M,N$ gegeben, so schreiben wir
\[M\preceq_{deg}N:\Leftrightarrow \dimv M=\dimv N\text{ und }\mathcal{O}(N)\subseteq\overline{\mathcal{O}(M)}.\]

Für alle $m\in\operatorname{Rep}_K(\underline{d})$ ist $\operatorname{Iso}_{G_K(\underline{d})}(m)=\operatorname{Aut}_K(\operatorname{rep}m)\subseteq \End_K(\operatorname{rep}m)$. Also gilt für alle $K$-Darstellungen $M$ die Formel
\[\dim \mathcal{O}(M)=\sum_{x\in K_0}\dimv M(x)^2-\dim_k\End_K(M).\] 

Wir benötigen folgende Aussagen über Entartungen von Köcherdarstellungen:
\begin{prop}\label{grund:var:prop2}\ 
\begin{enumerate}
\item Seien $M,N$ $K$-Darstellungen mit $M\preceq_{deg}N$. Dann gilt für alle $K$-Darstellungen $T$:
\[\dim \Hom_K(M,T)\leq \dim \Hom_K(N,T).\]
\item Seien $M,N$ und $M',N'$ $K$-Darstellungen mit $M\preceq_{deg}N$ und $M'\preceq_{deg}N'$. Dann gilt auch $M\oplus M'\preceq_{deg} N\oplus N'$.
\item Sei 
\[0\longrightarrow M'\longrightarrow M\longrightarrow M''\longrightarrow 0\]
eine exakte Sequenz von $K$-Darstellungen. Dann gilt $M\preceq_{deg}M'\oplus M''$.
\end{enumerate}
\end{prop}

Man hat folgende Charakterisierung von offenen und abgeschlossenen Orbiten:
\begin{prop}\label{grund:var:prop2.5} Sei $M$ eine $K$-Darstellung, und sei $X:=\operatorname{Rep}_K(\dimv M)$. Dann gilt:
\begin{enumerate}
\item $\mathcal{O}(M)$ ist offen in $X$ $\Leftrightarrow$ $\Ext^1_K(M,M)=0$.
\item $\mathcal{O}(M)$ ist abgeschlossen in $X$ $\Leftrightarrow$ $M$ ist halbeinfach.
\end{enumerate}
\end{prop}
\newpage
Wir zeichnen einige $G_K(\underline{d})$-stabile Teilmengen von $\operatorname{Rep}_K(\underline{d})$ aus.
\begin{defn}\label{grund:var:def2}\ 
\begin{enumerate}
\item Für alle $i\geq0$ setze
\[
\operatorname{Rep}_K^{(i)}(\underline{d}):=\bigl\{ m\in\operatorname{Rep}_K(\underline{d})\bigm| \dim\End_K(\operatorname{rep}m)=i\bigr\}.
\]
\item Für alle $i\geq 0$ sei 
\[
U_K^{(i)}(\underline{d}):=
\bigl\{ m\in\operatorname{Rep}_K^{(i)}(\underline{d})\bigm| \operatorname{rep}m \text{ ist unzerlegbar}\bigr\}.
\]
\end{enumerate}
\end{defn}Die eben definierten Teilmengen erfüllen folgende Eigenschaften (vgl. etwa \cite[Abschnitt 2.1 und 2.2]{Bong2}:
\begin{prop}\label{grund:var:prop3}\ 
\begin{enumerate}
\item Für alle $i\geq0$ ist die Menge
\[\operatorname{Rep}_K^{(\leq i)}(\underline{d}):=\bigcup_{j\leq i}\operatorname{Rep}_K^{(j)}(\underline{d})\] 
offen in $\operatorname{Rep}_K(\underline{d})$. Die offene Schicht von $\operatorname{Rep}_K(\underline{d})$ besteht aus allen $m\in\operatorname{Rep}_K(\underline{d})$, so daß $\dim \End_K(\operatorname{rep}m)$ minimal ist.
\item Für alle $i\geq0$ ist $U_K^{(i)}(\underline{d})$ abgeschlossen in $\operatorname{Rep}_K^{(i)}(\underline{d})$.
\end{enumerate}
\end{prop}
\section{Erweiterungsgruppen von Darstellungen}\label{grund:ext}
\setcounter{zaehler}{0}
\numberwithin{zaehler}{section}
Sei $k$ ein Körper. Sei $A$ eine $k$-Algebra und $K$ ein endlicher Köcher. Wir betrachten Darstellungen beliebiger Dimension.

Seien $M,N$ zwei $A$-Moduln. Dann wird $\Hom_k(M,N)$ ein $A$-$A$-Bimodul via \[(a\cdot\phi\cdot b)(m):=a\cdot\phi(b\cdot m) \text{ für alle } a,b\in A, m\in M.\]
Setze 
\[Z^1_A(N,M):=\bigl\{\delta\in\Hom_k\bigl(A,\Hom_k(N,M)\bigr)\mid\delta(ab)=a\delta(b)+\delta(a)b\text{ für alle } a,b\in A\bigr\}.\]
Wir definieren einen Homomorphismus von $k$-Vektorräumen durch
\[d_{N,M}:\Hom_k(N,M)\longrightarrow Z^1_A(N,M),\  d_{N,M}(r)(a):=M(a)\circ r-r\circ N(a) \text{ für alle } a\in A.\]
Es sei weiter $B^1_A(N,M):=\Bild d_{N,M}$. $Z^1_A$ und $B^1_A$ sind Bifunktoren auf der Kategorie aller $A$-Moduln und $d$ ist eine natürliche Transformation. 

Die folgende Charakterisierung von $\Ext^1_A(M,N)$ basiert auf der Tatsache, daß $\Ext^1_A(M,N)$ mit der ersten Hochschild-Kohomologiegruppe von $A$ mit Werten in $\Hom_A(M,N)$ übereinstimmt und darauf, daß Hochschild-Kohomologiegruppen mit Hilfe der Standardauflösung von $A$ berechnet werden können (vgl. \cite[Kapitel IX] {CE}). Man kann die Aussage allerdings auch direkt und ohne Verwendung des allgemeinen Apparats verifizieren.
\begin{prop}\label{grund:ext:prop1}\ 
Es existiert eine in $N$ und $M$ funktorielle exakte Sequenz
\[0\longrightarrow \Hom_A(N,M)\longrightarrow \Hom_k(N,M)\stackrel{d_{N,M}}{\longrightarrow}Z^1_A(N,M)\stackrel{\phi_{N,M}}{\longrightarrow}\Ext^1_A(N,M)\longrightarrow0\]
von $k$-Vektorräumen. Interpretiert man $\Ext^1_A(N,M)$ als Gruppe von Äquivalenzklassen von Erweiterungen von $N$ durch $M$, so ist $\phi_{N,M}$ wie folgt gegeben:
Sei $\delta\in Z^1_A(N,M)$. Sei $X_\delta$ der $A$-Modul definiert durch
\begin{itemize}
\item[$\cdot$]$X_\delta:=M\oplus N$ als $k$-Vektorraum.
\item[$\cdot$]Für alle $a\in A$ sei
\[X_\delta(a):=\begin{bmatrix}M(a)&\delta(a)\\0&N(a)\end{bmatrix}.\]
\end{itemize}
$\phi_{N,M}(\delta)$ ist dann die Äquivalenzklasse der zugehörigen exakten Sequenz
\[0\longrightarrow M\stackrel{i}{\longrightarrow} X_\delta\stackrel{p}{\longrightarrow} N\longrightarrow 0\]
von $A$-Moduln, wobei $p$ die Projektion von $X_\delta$ auf $N$ und $i$ die Einbettung von $M$ in $X$ ist.
\end{prop}

Für Köcherdarstellungen erhält man eine noch einfachere Beschreibung (vgl. \cite[Abschnitt 2.1]{Ri5}). Seien $N,M$ Darstellungen von $K$. Setze
\begin{itemize}
\item[$\cdot$]$C^0_K(N,M):=\bigoplus_{x\in K_0}\Hom_k(N(x),M(x))$,
\item[$\cdot$]$C^1_K(N,M):=\bigoplus_{\alpha\in K_1}\Hom_k(N(n\alpha),M(s\alpha))$,
\item[$\cdot$]$d_{N,M}:C^0(N,M)\longrightarrow C^1(N,M),\ (f_x)_{x\in K_0}\mapsto \Bigl(M{(\alpha)}\circ f_{n\alpha}-f_{s\alpha}\circ N{(\alpha)}\Bigr)_{\alpha \in K_1}$.
\end{itemize}
$C^0_K$ und $C^1_K$ sind Bifunktoren auf der Kategorie aller $K$-Darstellungen, und $d$ ist eine natürliche Transformation. Es gilt:
\begin{prop}\label{grund:ext:prop2}\ 

Es existiert eine in $M$ und $N$ funktorielle exakte Sequenz
\[0\longrightarrow \Hom_{K}(N,M)\longrightarrow C^0_K(N,M)\stackrel{d_{N,M}}{\longrightarrow}C^1_K(N,M)\stackrel{\phi_{N,M}}{\longrightarrow}\Ext^1_K(N,M)\longrightarrow 0\]
von $k$-Vektorräumen. Dabei ist $\phi_{N,M}$ ist wie folgt gegeben:

Sei $\delta=(\delta_\alpha)_{\alpha\in K_1}\in C^1_K(N,M)$. Definiere eine Darstellung $X_\delta$ von $K$ durch
\begin{align*}
X_\delta(x)&:=M(x)\oplus N(x) & \text{ für alle } x\in K_0,\\
X_\delta(\alpha)&:=\left[\begin{array}{cc}M(\alpha)&\delta_\alpha\\0&N(\alpha)\end{array}\right] & \text{ für alle } \alpha\in K_1.
\end{align*}
$\phi_{N,M}(\delta)$ ist dann gegeben durch die zugehörige exakte Sequenz 
\[0\longrightarrow M\stackrel{i}{\longrightarrow} X_\delta\stackrel{p}{\longrightarrow} N\longrightarrow 0\]
von $K$-Darstellungen, wobei für alle $x\in K_0$ $p_x$ die Projektion von $X_\delta(x)$ auf $N(x)$ und $i_x$ die Einbettung von $M(x)$ in $X_\delta(x)$ ist.
\end{prop}

\section{Spiegelungsfunktoren}\label{grund:refl}
\setcounter{zaehler}{0}
\numberwithin{zaehler}{section}
Sei $K$ ein endlicher, zykelloser Köcher. Wir definieren die Spiegelungsfunktoren nach Bernstein, Gelfand und Ponomarev (vgl. \cite{BernGelPon}).

Sei $y\in K_0$ eine Senke. Sei $yK$ der Köcher, der aus $K$ durch Umkehren der Orientierung aller Pfeile mit Spitze $y$ entsteht. Der Spiegelungsfunktor $\sigma_y^+:\fdar{K}\longrightarrow\fdar{yK}$ bezüglich $y$ wird wie folgt definiert:

Sei $M$ eine Darstellung von $K$. Für $x\in K_0-\{y\}$ sei $\sigma_y^+M(x)=M(x)$. Für alle Pfeile $\alpha\in K_1$ mit $s\alpha\neq y$ sei $\sigma_y^+M(\alpha)=M(\alpha)$. Betrachte jetzt die $k$-lineare Abbildung
\[\phi:\bigoplus_{\begin{subarray}{c}\alpha\in K_1\\s\alpha=y\end{subarray}}M(n\alpha)\longrightarrow M(y),\]
induziert von den $M(\alpha):M(n\alpha)\longrightarrow M(y)$ für Pfeile $\alpha\in K_1$ mit Spitze $y$. Erhalte eine exakte Sequenz
\[0\longrightarrow\sigma_y^+M(y)\stackrel{\psi}{\longrightarrow}\bigoplus_{\begin{subarray}{c}\alpha\in K_1\\s\alpha=y\end{subarray}}M(n\alpha)\stackrel{\phi}{\longrightarrow} M(y)\]
von $k$-Vektorräumen. Sei $\alpha$ ein Pfeil in $K$ mit Spitze $y$. $\alpha$ liefert einen Pfeil $\alpha:y\longrightarrow n\alpha$ in $yK$. Sei
 \[\sigma_y^+M(\alpha):\sigma_y^+M(y)\longrightarrow M(n\alpha)=\sigma_y^+M(n\alpha)\]
die von $\psi$ induzierte lineare Abbildung.

Diese Definition liefert einen $k$-linearen Funktor. Dual wird für Quellen $z$ der Spiegelungsfunktor $\sigma_z^-:\fdar{K}\longrightarrow\fdar{zK}$ definiert.

Die Aussagen (1) bis (4) der folgenden Proposition sind Lemma 1.1, Theorem 1.1 und Lemma 3.1 aus \cite{BernGelPon} entnommen; (5) und (6) verifiziert man leicht.
\begin{prop}\label{grund:refl:prop1}
Sei $y$ eine Senke in $K$. Betrachte die Spiegelungsfunktoren
\[\sigma_y^+:\fdar{K}\longrightarrow\fdar{yK}\text{ und } \sigma_y^-:\fdar{yK}\longrightarrow\fdar{K}.\]
Es sei $s_y:\mathbb{Z}^{K_0}\longrightarrow \mathbb{Z}^{K_0}$ die Spiegelung zu $y$.
Dann gilt:
\begin{enumerate}
\item $\sigma_y^+(E_y)=0$ und $s_y(\dimv E_y) < 0 $.
\item Sei $M$ eine unzerlegbare Darstellung von $K$, $M\not\simeq E_y$. Dann ist $\sigma_y^+(M)$ unzerlegbar und es gilt:
\[s_y(\dimv M)=\dimv \sigma_y^+(M).\]
\item $(\sigma_y^-,\sigma_y^+)$ ist ein adjungiertes Paar von Funktoren.
\item $\sigma_y^+$ und $\sigma_y^-$ induzieren zueinander quasiinverse Äquivalenzen zwischen den vollen Unterkategorien aller Darstellungen von $K$ bzw. $yK$, die $E_y$ nicht als direkten Summanden haben.
\item Sei 
\[0\longrightarrow M'\longrightarrow M\longrightarrow M''\longrightarrow 0\]
eine exakte Sequenz von $K$-Darstellungen. Tritt $E_y$ nicht als direkter Summand von $M'$ und $M''$ auf, so tritt $E_y$ auch nicht als direkter Summand von $M$ auf, und die Sequenz
\[0\longrightarrow \sigma_y^+(M')\longrightarrow \sigma_y^+(M)\longrightarrow \sigma_y^+(M'')\longrightarrow 0\]
ist eine exakte Sequenz von $yK$-Darstellungen.
\item Seien $M,N$ $K$-Darstellungen, so daß $E_y$ kein direkter Summand von $M,N$ ist. Dann induziert $\sigma_y^+$ Isomorphismen
\[\End_{K}(M,N)\simeq\End_{yK}(\sigma_y^+(M),\sigma_y^+(N)) \text{ und } \Ext^1_{K}(M,N)\simeq\Ext^1_{yK}(\sigma_y^+(M),\sigma_y^+(N))\]
von $k$-Vektorräumen.
\end{enumerate}
\end{prop}

Wir definieren noch die Coxeter-Funktoren für $K$  und geben ihre Beziehung zu den Auslander-Reiten Verschiebungen  $\tau=\operatorname{DTr}$ und $\tau^{-1}=\operatorname{TrD}$ an (vgl. \cite[VII.5.8]{Sko}).
\begin{defn}\label{grund:refl:def2}\ 
\begin{enumerate}
\item Eine Aufzählung $(x_1,\ldots,x_n)$ der Punkte von $K$ heißt zulässig, falls gilt: Ist $i\leq j$, so existiert kein Weg von $x_i$ nach $x_j$ in $K$.
\item Sei $(x_1,\ldots,x_n)$ eine zulässige Aufzählung der Punkte von $K$. Die Coxeter-Funktoren von $K$ sind definiert durch
\[\Phi^+:=\sigma_{x_n}^+\circ\cdots\circ\sigma_{x_1}^+\text{ und }\Phi^-:=\sigma_{x_1}^-\circ\cdots\circ\sigma_{x_n}^-.\]
\end{enumerate}
\end{defn}
\begin{prop}\label{grund:refl:prop2}
Für alle $K$-Darstellungen $M$ ohne projektiven direkten Summanden existiert ein Isomorphismus \[\Phi^+(M)\simeq \tau (M)\] von $K$-Darstellungen. Besitzt $M$ keinen injektiven direkten Summanden, so gilt dual \[\Phi^-(M)\simeq \tau^{-1} (M).\]
\end{prop}
\section{Die universelle Überlagerung eines Köchers}\label{grund:cov}
\setcounter{zaehler}{0}
\numberwithin{zaehler}{section}
Sei $K$ ein endlicher zusammenhängender Köcher. Wir konstruieren die universelle Überlagerung $\pi:\widetilde{K}\longrightarrow K$ von $K$ und geben die für uns wichtigen Eigenschaften des zugehörigen Überlagerungsfunktors zwischen den Darstellungskategorien an.

Sei $\hat{K}$ der Köcher, der aus $K$ entsteht durch Hinzufügen der "`inversen"' Pfeile \[\alpha^{-1}:n\alpha\longrightarrow s\alpha,\ \alpha\in K_1.\]
Eine Wanderung in $K$ ist ein Weg in $\hat{K}$. Eine Wanderung in $K$ heißt reduziert, falls sie keine Teilwanderung der Form $\alpha\alpha^{-1}$ oder $\alpha^{-1}\alpha$, $\alpha\in K_1$ enthält. Für eine Wanderung $w$ in $K$ bezeichne $sw$ den Endpunkt und $nw$ den Startpunkt der Wanderung.

Fixiere einen Punkt $t_0\in K_0$. Sei $\mathcal{W}_{t_0}$ die Menge aller Wanderungen in $K$ mit Start $t_0$. Für zwei Wanderungen $w,w'\in\mathcal{W}_{t_0}$ schreibe
\[w\sim w':\Leftrightarrow w=ut, w'=uvv^{-1}t \text{ für gewisse Wanderungen $u,v,t$ in $K$}.\]
Es sei $\sim$ die von dieser Relation erzeugte Äquivalenzrelation.
\begin{bem}
Die reduzierten Wanderungen mit Start $t_0$ bilden ein Repräsentantensystem der Äquivalenzklassen von $\sim$.
\end{bem}
Ist $w\in\mathcal{W}_{t_0}$ eine Wanderung, so sei $[w]$ die zugehörige Äquivalenzklasse. Für eine Äquivalenzklasse $x=[w]$ einer Wanderung $w\in\mathcal{W}_{t_0}$ setze $sx:=sw\in K_0$. Ist $\alpha\in K_1$ ein Pfeil mit $n\alpha=sw$, so setze $\alpha x:=[\alpha w]$. Die universelle Überlagerung $\widetilde{K}=(\widetilde{K}_0,\widetilde{K}_1,\widetilde{n},\widetilde{s})$ von $K$ wird folgendermaßen definiert:
\begin{defn}\label{grund:cov:def1}\ 
\begin{itemize}
\item[$\cdot$]$\widetilde{K}_0:=\mathcal{W}_{t_0}/\sim$.
\item[$\cdot$]$\widetilde{K}_1:=\bigl\{(\alpha,x)\in K_1\times \widetilde{K}_0\mid n\alpha=sx\bigr\}$.
\item[$\cdot$] Für $(\alpha,x)\in\widetilde{K_1}$ sei $\widetilde{n}(\alpha,x):=x$ und $\widetilde{s}(\alpha,x):=\alpha x$.
\end{itemize}
\end{defn}
Wir erhalten einen Köchermorphismus $\pi=(\pi_0,\pi_1):\widetilde{K}\longrightarrow K$ mit
\[\pi_0:\widetilde{K}_0\longrightarrow K_0,\ y\mapsto sy\text{ und } \pi_1:\widetilde{K}_1\longrightarrow K_1,\ (\alpha,x)\mapsto\alpha.\]
\begin{bem}\ 
\begin{enumerate}
\item Die universelle Überlagerung $\pi$ ist eindeutig bestimmt durch folgende Eigenschaften:
\begin{enumerate}[i)]
\item $\widetilde{K}$ ist baumartig.
\item $\pi_0$ ist surjektiv.
\item Sei $x\in K_0$ und sei $\alpha\in K_1$ ein Pfeil mit Spitze $x$. Sei $\widetilde{x}\in\widetilde{K}_0$ ein Punkt über $x$ (d.h. mit $\pi_0(\widetilde{x})=x$). Dann läßt sich $\alpha$ eindeutig zu einem Pfeil $\widetilde{\alpha}\in\widetilde{K}_1$ mit Spitze $\widetilde{x}$ liften (d.h. es existiert genau ein Pfeil $\widetilde{\alpha}\in\widetilde{K}_1$ mit Spitze $\widetilde{x}$ und $\pi_1(\widetilde{\alpha})=\alpha$); analog läßt sich jeder Pfeil mit Start $x$ in $K$ eindeutig zu einem Pfeil mit Start $\widetilde{x}$ in $\widetilde{K}$ liften.
\end{enumerate}
\item Ist $K$ baumartig, so ist $\pi:\widetilde{K}\longrightarrow K$ ein Isomorphismus von Köchern.
\end{enumerate}
\enlargethispage{\baselineskip}
\end{bem}
Definiere einen Funktor $F_\lambda:\fdar{\widetilde{K}}\longrightarrow \fdar{K}$ wie folgt: Sei $M\in\fdar{\widetilde{K}}$. Für alle $x\in K_0$ setze
\[(F_\lambda M)(x):=\bigoplus_{y\in\pi_0^{-1}(x)}M(y).\]\newpage
Sei jetzt $\alpha:x\rightarrow x'\in K_1$. Für alle $y\in\pi_0^{-1}(x)$ sei $\alpha_y\in\widetilde{K}_1$ der eindeutig bestimmte Pfeil über $\alpha$ mit Start $y$. Die linearen Abbildungen
\[M(\alpha_y):M(y)\longrightarrow M\bigl(\widetilde{s}\alpha_y\bigr),\ y\in\pi_0^{-1}(x)\]
induzieren eine lineare Abbildung $(F_\lambda M)(\alpha):(F_\lambda M)(x)\longrightarrow(F_\lambda M)(x')$.

Ist $f:M\longrightarrow M'$ ein Homomorphismus in $\fdar{\widetilde{K}}$,
so induzieren für alle $x\in K_0$ die linearen Abbildungen $f_y:M(y)\longrightarrow M'(y),y\in\pi_0^{-1}(x)$ eine lineare Abbildung 
\[(F_\lambda f)_x:(F_\lambda M)(x)\longrightarrow (F_\lambda M')(x).\]
$F_\lambda f$ ist dann ein Homomorphismus von $K$-Darstellungen. Der Funktor $F_\lambda$ erfüllt folgende Eigenschaften:
\begin{prop}\label{grund:cov:prop1}\ 
\begin{enumerate}
\item $F_\lambda$ ist exakt und treu.
\item $\dim_kF_\lambda M=\dim_kM$ für alle $M\in\fdar{\widetilde{K}}$.
\item Sei $M\in\fdar{\widetilde{K}}$ unzerlegbar. Dann ist auch $F_\lambda M$ unzerlegbar.
\end{enumerate}
\end{prop}
Die ersten beiden Eigenschaften folgen aus der Definition, die dritte beruht auf einem Resultat von Gabriel im allgemeineren Rahmen der Theorie der Galois-Überlagerungen, vgl. \cite[Lemma 3.5]{Ga2}. Die dort benötigte freie Gruppenoperation auf $\widetilde{K}$ ist hier durch die Operation der Fundamentalgruppe von $K$ auf $\widetilde{K}$ gegeben.

\section{Kronecker-Moduln}\label{grund:kron}
\setcounter{zaehler}{0}
\numberwithin{zaehler}{section}
Sei $k$ ein algebraisch abgeschlossener Körper. Wir wiederholen die Klassifikation unzerlegbarer $\Lambda_2$-Moduln. Die $\Lambda_2$-Moduln nennt man auch Kronecker-Moduln. Für den Köcher $Q_2$ verwenden wir die Bezeichnungen
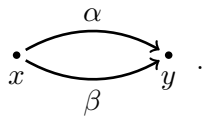
\begin{center}
\parbox[c]{2.6cm}{
\begin{tikzpicture}[line width=1pt]
	\tikzstyle{knt}=[circle, fill, inner sep=1pt]
  \node at (0,0)[knt, label=below:$x$](x){};
  \node at (2,0)[knt, label=below:$y$](y){};
  \draw[->, bend left, shorten >= 2pt, shorten <= 2pt] (x) to node[above]{$\alpha$} (y);
  \draw[->, bend right, shorten >= 2pt, shorten <= 2pt] (x) to node[below]{$\beta$} (y);
\end{tikzpicture}}.\end{center}
Wir schreiben $\underline{d}=(d_x,d_y)\in\mathbb{Z}^2$ für Dimensionsvektoren $\underline{d}\in\mathbb{Z}^{(Q_2)_0}$.
\begin{defn}\label{grund:kron:def1} Wir definieren folgende $Q_2$-Darstellungen durch Matrizen:
\begin{enumerate}
\item Für alle $n\geq0$ sei $P_n$ die Darstellung zum Dimensionsvektor $(n,n+1)$ definiert durch
\[P_n(\alpha)=\text{\small$\left[
\begin{array}{ccc}
1&        &   \\
 & \ddots &   \\
 &        & 1 \\
0 & \cdots & 0
\end{array}\right]$}\in k^{n+1 \times n} \text{ und }
P_n(\beta)=\text{\small$\left[
\begin{array}{ccc}
0 & \cdots & 0 \\
1 &        &   \\
  & \ddots &   \\
  &        & 1 
\end{array}\right]$}\in k^{n+1 \times n}.
\]
\item Für alle $n\geq0$ sei $I_n$ die Darstellung zum Dimensionsvektor $(n+1,n)$ gegeben durch
\[I_n(\alpha)=P_n(\alpha)^T\text{ und } I_n(\beta)=P_n(\beta)^T.\]
Dabei bezeichnet $A^T$ die Transponierte einer Matrix $A$.
\item Für alle $n\geq1$ und für alle $\lambda\in k$ sei $R_{n,\lambda}$ die Darstellung zum Dimensionsvektor $(n,n)$ gegeben durch
\[R_{n,\lambda}(\alpha)=\text{\small$\left[
\begin{array}{ccc}
1 &        &   \\
  & \ddots &   \\
  &        & 1 
\end{array}\right]$}\in k^{n \times n}\text{ und }
R_{n,\lambda}(\beta)=\text{\small$\left[
\begin{array}{cccc}
\lambda &        &        &   \\
1 & \ddots &        &   \\
  & \ddots & \ddots &   \\ 
  &        & 1      & \lambda \\
\end{array}\right]$}\in k^{n \times n}.
\]
\item Für alle $n\geq1$ sei $R_{n,\infty}$ die Darstellung zum Dimensionsvektor $(n,n)$ gegeben durch
\[R_{n,\infty}(\alpha)=\text{\small$\left[
\begin{array}{cccc}
0 &        &        &   \\
1 & \ddots &        &   \\
  & \ddots & \ddots &   \\ 
  &        & 1      & 0 \\
\end{array}\right]$}\in k^{n \times n}\text{ und }
R_{n,\infty}(\beta)=\text{\small$\left[
\begin{array}{ccc}
1 &        &   \\
  & \ddots &   \\
  &        & 1 
\end{array}\right]$}\in k^{n \times n}.\]
\end{enumerate}
\end{defn}
Es gilt folgender Klassifikationssatz (vgl. Kapitel VIII.7 über Kronecker-Algebren in \cite{Aus}):
\begin{satz}\label{grund:kron:satz1}
Die $Q_2$-Darstellungen aus \cref{grund:kron:def1} bilden ein Repräsentantensystem der unzerlegbaren $Q_2$-Darstellungen. Für $n\geq1$ und $\lambda\in\mathbb{P}^1=k\cup\{\infty\}$ ist $R_{\lambda,n}$ regulär, $P_n$ präprojektiv und $I_n$ präinjektiv im Sinne der Auslander-Reiten-Theorie. 
\end{satz}

\chapter{Schofield-Induktion}\label{scho}
\setcounter{zaehler}{0}
\numberwithin{zaehler}{chapter}
Sei stets $k$ ein algebraisch abgeschlossener Körper. Eine Köcherdarstellung $M$ über $k$ heißt exzeptionell, falls ihr Endomorphismenring $k$ ist und sie keine Selbsterweiterungen hat. Beispiele für exzeptionelle Köcherdarstellungen sind die präprojektiven und präinjektiven Moduln über der Köcheralgebra. Wir beweisen einen Satz von Schofield (vgl. \cref{scho:scho:satz1}) zur induktiven Beschreibung einer exzeptionellen Köcherdarstellung $M$ als Erweiterung gewisser anderer exzeptioneller Darstellungen; gleichzeitig wird $M$ als Modul einer verallgemeinerten Kronecker-Algebra gedeutet. Der Satz baut im wesentlichen auf Resultaten aus der Kipptheorie auf und ist konstruktiv bis auf die Zerlegung von endlichdimensionalen Moduln in Unzerlegbare (vgl. \cref{scho:scho:bem4}).

Wir bemerken, daß für beliebige abelsche Kategorien vermöge der Yoneda-Konstruktion alle $\Ext^i$-Gruppen existieren, und daß kurze exakte Sequenzen lange exakte Ext-Sequenzen induzieren. Ein Objekt $X$ ist projektiv genau dann, wenn $\Ext^1_\mathcal{A}(X,-)=0$ ist.
\section{Kippmoduln und universelle Erweiterungen}\label{scho:kipp}
\setcounter{zaehler}{0}
\numberwithin{zaehler}{section}
In diesem Abschnitt sei stets $A$ eine endlichdimensionale $k$-Algebra. Das Material über universelle Erweiterungen ist \cite[Abschnitt 2]{Hu} entnommen. Dort werden universelle Erweiterungen für alle $\Ext^i$ definiert, wir geben hier eine konkrete Konstruktion für für den Fall $i=1$ an.\enlargethispage{\baselineskip}
\begin{defn}\label{scho:kipp:def1}
Seien $X$ und $Y$ $A$-Moduln. Eine universelle Erweiterung von $Y$ durch Objekte in $\operatorname{add}(X)$ ist eine kurze exakte Sequenz $\xi:\ 0\longrightarrow T \longrightarrow M \longrightarrow Y\longrightarrow 0$
mit
\begin{itemize}
\item[$\cdot$] $T\in\operatorname{add}(X)$.
\item[$\cdot$] Für alle $T'\in\operatorname{add}(X)$ ist der von $\xi$ induzierte Verbindungshomomorphismus 
\[
\partial_{\:X}(\xi):\Hom_A(T, T')\longrightarrow\Ext^1_A(Y,T'),\ \phi\mapsto\phi\xi
\]
surjektiv.
\end{itemize}
$\xi$ heißt minimal, falls gilt: Ist $f\in\End_A(T)$ mit $f\xi=\xi$, so ist $f$ ein Isomorphismus.
\end{defn}

\begin{bem}\label{scho:kipp:bem1}\ 
\begin{enumerate}
\item (Minimale) universelle Erweiterungen von Objekten in $\operatorname{add}(Y)$ durch $X$ werden dual definiert. 
\item Minimale universelle Erweiterungen sind eindeutig bestimmt bis auf Isomorphie von exakten Sequenzen.
\end{enumerate}
\end{bem}

Wir skizzieren kurz den Beweis für die Existenz von minimalen universellen Erweiterungen. Zur Vermeidung größerer Rechnungen vereinbaren wir zunächst:
\begin{defn}\label{scho:kipp:def2}
Sei $\mathcal{A}$ eine $k$-lineare Kategorie und sei $F:\mathcal{A}\longrightarrow\fmod{k}$ ein $k$-linearer Funktor. Das Radikal von $F$ ist der Unterfunktor $\operatorname{Rad} F$ von $F$ definiert durch
\[\bigl(\Rad F\bigr)(X):=\sum_{T\in\mathcal{A}}\mathcal{R}_{\mathcal{A}}(T,X)F(T)\text{ für alle } X\in\mathcal{A}.\]
Dabei sei $\mathcal{R}_{\mathcal{A}}(T,X)F(T):=\ _k\bigl\langle\bigl(Ff\bigr)(x)\bigm| f\in\mathcal{R}_{\mathcal{A}}(T,X),\ x\in FT\bigr\rangle$ für alle $T,X$.
\end{defn}
\begin{bem}\label{scho:kipp:bem2}\ 
\begin{enumerate}
\item Es ist $\Rad \Hom_{\mathcal{A}}(T,-)=\mathcal{R}_{\mathcal{A}}(T,-)$ für alle $T$.
\item Sei $X$ ein $A$-Modul und sei $F:\operatorname{add}(X)\longrightarrow\fmod{k}$ ein $k$-linearer Funktor. Dann ist
\[\Rad F(T)=\mathcal{R}_{\operatorname{add}(X)}(X,T)F(X)\text{ für alle } T\in\operatorname{add}(X).\]
Insbesondere ist $(\Rad F)(X)$ das Radikal von $FX$ als $\End_A(X)$-Linksmodul. 
\end{enumerate}
\end{bem}

\begin{prop}\label{scho:kipp:prop1}
Seien $Y$ und $X$ $A$-Moduln. Dann existiert eine minimale universelle Erweiterung von $Y$ durch Objekte in $\operatorname{add}(X)$.
\end{prop}
\begin{bew}
Wir dürfen annehmen: $X=X_1\oplus\ldots\oplus X_r$ mit $X_i$ unzerlegbar, $X_i\not\simeq X_j$ für $i\neq j$. Betrachte den Funktor
\[G:\operatorname{add}(X)\longrightarrow\fmod{k},\ V\mapsto\Ext^1_A(Y,V).\]
Wähle für alle $i=1,\ldots,r$ Elemente $\zeta_{1}^{(i)},\ldots,\zeta_{n_i}^{(i)}\in \Ext^1_A(Y,X_i)=G(X_i)$, so daß die Restklassen der $\zeta_{j}^{(i)}$ eine $k$-Basis von $\bigl(G/\Rad G\bigr)(X_i)$ bilden. Setze $T:=\bigoplus_{i=1}^rX_i^{n_i}$. Für $i=1,\ldots,r$ und für $j=1,\ldots,n_i$ sei $\epsilon_j^{(i)}:X_i\longrightarrow T$ die Einbettung von $X_i$ in die $j$-te Komponente von $X_i^{n_i}$ in $T$. Sei 
\[\xi:=\sum_{i=1}^r\sum_{j=1}^{n_i}\epsilon_j^{(i)}\zeta_j^{(i)}:\ 0\longrightarrow T\longrightarrow Q\longrightarrow Y\longrightarrow 0.\]
$\xi$ ist eine minimale universelle Erweiterung: Betrachte den Funktor
\[F:\operatorname{add}(X)\longrightarrow\fmod{k},\ V\mapsto\Hom_A(T,V).\]
Für $T'\in\operatorname{add}(X)$ sei \[\partial_{\:T'}:\Hom_A(T,T')\longrightarrow\Ext^1_A(Y,T'),\ f\mapsto f\xi\]
der von $\xi$ induzierten Verbindungshomomorphismus. Wir erhalten eine natürliche Transformation $\partial:F\longrightarrow G$. $\partial$ induziert eine natürliche Transformation $\overline{\partial}:F/\Rad F\longrightarrow G/\Rad G$ auf den Deckeln. Für $i=1,\ldots,r$ folgt leicht aus der Konstruktion von $\xi$, daß $\overline{\partial}_{X_i}$ ein Isomorphismus ist. Beachte dazu, daß die vorausgesetzten Bedingungen an die $X_i$ ergeben:
\[(F/\Rad F)(X_i)=\Hom_A(T,X_i)/\mathcal{R}_A(T,X_i)\simeq k^{n_i}.\]
Also ist $\overline{\partial}$ ein Isomorphismus von Funktoren. Damit folgt, daß $\partial$ ein Epimorphismus von Funktoren ist (beachte: $\partial_{\:X}$ ist ein Homomorphismus von $\End_A(X)$-Linksmoduln und induziert einen Isomorphismus auf den Deckeln, ist also ein (wesentlicher) Epimorphismus) und das zeigt: $\xi$ ist universell. Es bleibt zu zeigen: $\xi$ ist minimal. Sei also $f\in\End_A(T)=FT$ mit $f\xi=\xi$. Da $\overline{\partial}_{\:T}$ injektiv ist, gilt $f\equiv id_T\mod\bigl(\Rad F\bigr)(T)$. Wegen $\bigl(\Rad F\bigr)(T)=J\End_A(T)$ folgt: $f$ ist ein Isomorphismus.
\end{bew}

\begin{bem}\label{scho:kipp:bem3} Sei $X=X_1\oplus\ldots\oplus X_r$, es gelte $\dim\Hom_A(X_i,X_j)=\delta_{ij}$ für alle $i,j$. In dem Fall konstruiert man die minimale universelle Erweiterung $\xi$ von $Y$ durch Objekte in $\operatorname{add}(X)$ wie folgt: Wähle für alle $i=1,\ldots,r$ eine $k$-Basis $(\zeta_{1}^{(i)},\ldots,\zeta_{n_i}^{(i)})$ von $\Ext^1_A(Y,X_i)$. Setze $T:=\bigoplus_{i=1}^rX_i^{n_i}$ und 
\[\xi:=\sum_{i=1}^r\sum_{j=1}^{n_i}\epsilon_j^{(i)}\zeta_j^{(i)}:\ 0\longrightarrow T\longrightarrow Q\longrightarrow Y\longrightarrow 0\]
(dabei seien die $\epsilon_j^{(i)}:X_i\longrightarrow T$ die Einbettungen). Für alle $T'\in\operatorname{add}(X)$ ist dann der von $\xi$ induzierte Verbindungshomomorphismus $\partial_{\:T'}:\Ext^1_A(T,T')\longrightarrow\Ext^1_A(Y,T')$ sogar ein Isomorphismus. Erfüllt umgekehrt eine exakte Sequenz diese Eigenschaft, so ist sie eine minimale universelle Erweiterung.
\end{bem}

\begin{defn}\label{scho:kipp:def3}
Sei $M$ ein $A$-Modul. $M$ heißt ein partieller Kippmodul, falls gilt:
\begin{enumerate}[i)]
\item $\Ext^1_A(M,M)=0$.
\item $\operatorname{pdim} M \leq 1$.
\end{enumerate}
Ist $M$ ein partieller Kippmodul, so heißt $M$ ein Kippmodul, falls zusätzlich gilt:
\begin{enumerate}[i)]
\setcounter{enumi}{2}
\item Es existiert eine exakte Sequenz
\[0\longrightarrow A\longrightarrow M_1\longrightarrow M_2\longrightarrow 0\]
von $A$-Moduln mit $M_1,M_2\in\operatorname{add}(M)$.
\end{enumerate}
\end{defn}

Zum Schluß dieses Abschnitts geben wir Charakterisierungen von Kippmoduln an (vgl. \cite[Teil 2]{Bong3}). \cref{scho:kipp:prop2} beweist man leicht mit Hilfe der definierenden Eigenschaften von universellen Erweiterungen; zum Beweis von \cref{scho:kipp:satz1} benötigt man Hauptresultate aus der klassischen Kipptheorie. 
\begin{prop}\label{scho:kipp:prop2}
Sei $M$ ein partieller Kippmodul. Sei
\[\xi:\ 0\longrightarrow A\longrightarrow \widetilde{M} \longrightarrow T \longrightarrow 0\]
eine universelle Erweiterung von Objekten in $\operatorname{add}(M)$ durch $A$. Dann ist $M\oplus \widetilde{M}$ ein Kippmodul.
\end{prop}

\begin{satz}\label{scho:kipp:satz1}
Sei $M$ ein partieller Kippmodul und sei $M=\bigoplus_{i=1}^r M_i^{n_i}$ eine Zerlegung von $M$ in paarweise nichtisomorphe unzerlegbare Moduln, $n_i\geq1$ für alle $i$. Dann sind äquivalent:
\begin{enumerate}[a)]
\item $M$ ist ein Kippmodul.
\item $r=\text{Anzahl der Isomorphieklassen von einfachen $A$-Moduln}$.
\end{enumerate}
\end{satz}

\section{Exzeptionelle Objekte und Orthogonale Kategorien}\label{scho:ex}
\setcounter{zaehler}{0}
\numberwithin{zaehler}{section}
\begin{defn}\label{scho:ex:def1} Sei $\mathcal{A}$ eine $k$-lineare abelsche Kategorie.
\begin{enumerate}[a)]
\item $\mathcal{A}$ heißt erblich, falls $\Ext^2_{\mathcal{A}}=0$ gilt.
\item Sei $M\in\mathcal{A}$. $M$ heißt exzeptionell, falls gilt:
\begin{itemize}
\item[$\cdot$] $\Ext^1_\mathcal{A}(M,M)=0,$
\item[$\cdot$] $\End_\mathcal{A}(M)=k.$
\end{itemize}
\item Eine exzeptionelle Folge der Länge $r$ in $\mathcal{A}$ ist eine Folge $E=(M_1,\ldots,M_r)$ von exzeptionellen Objekten, so daß gilt:
\[\Hom_\mathcal{A}(M_j,M_i)=0=\Ext^1_\mathcal{A}(M_j,M_i)\text{ für alle } i<j.\]
\item Sei $A$ eine endlichdimensionale $k$-Algebra. Eine exzeptionelle Folge der Länge $r$ in $\fmod{A}$ heißt vollständig, falls $r$ mit der Anzahl der Isomorphieklassen einfacher $A$-Moduln übereinstimmt.
\end{enumerate}
\end{defn}

Wir wiederholen einige grundlegende Fakten.
\begin{defn}\label{scho:ex:def2}
Sei $\mathcal{A}$ eine $k$-lineare abelsche Kategorie mit endlichdimensionalen Homomorphismenräumen. Ein Objekt $P\in\mathcal{A}$ heißt ein Progenerator, falls $P$ projektiv ist und falls für alle $X\in\mathcal{A}$ eine natürliche Zahl $t$ und ein Epimorphismus $P^t\longrightarrow X$ existieren.
\end{defn}
Der Beweis der folgenden Proposition ist derselbe wie der von Theorem 2.5 in \cite{Schw}.
\begin{prop}[Morita-Äquivalenz]\label{scho:ex:prop1}
Sei $\mathcal{A}$ eine $k$-lineare abelsche Kategorie mit endlichdimensionalen Homomorphismenräumen, und sei $P$ ein Progenerator von $\mathcal{A}$. Dann ist der Funktor $\Hom_{\mathcal{A}}(P,-):\mathcal{A} \longrightarrow \fmod{\End_{\mathcal{A}}(P)^{op}}$ eine Äquivalenz von $k$-linearen Kategorien.
\end{prop}
Wir benötigen außerdem folgende homologische Charakterisierung von Modulkategorien von Köcheralgebren (vgl. \cite[VII.1.4 und VII.1.7]{Sko} und beachte die Aussage über die Morita-Äquiva\-lenz):
\begin{prop}\label{scho:ex:prop2}
Ist $A$ eine endlichdimensionale Köcheralgebra, so ist $\fmod{A}$ eine erbliche Kategorie. Sei umgekehrt $\mathcal{A}$ eine erbliche abelsche Kategorie mit endlichdimensionalen Homomorphismenräumen. $\mathcal{A}$ besitze einen Progenerator $P$. Sei $K$ der Gabriel-Köcher von $\End_\mathcal{A}(P)^{op}$. Dann existiert eine Äquivalenz $\mathcal{A}\simeq\fmod{kK}$.
\end{prop}

\begin{bem}\label{scho:ex:bem1} Sei $K$ ein endlicher, zykelloser Köcher. Wir verwenden jetzt und im nächsten Beispiel, daß die quadratische Form von $K$ mit der Euler-Charakteristik von $kK$ übereinstimmt; vgl. dazu etwa \cite[VII.4]{Sko}.
\begin{enumerate}
\item Sei $M$ eine $K$-Darstellung. Nach dem Lemma von Happel-Ringel (vgl. \cite[§2]{CB1}) gilt: $M$ ist genau dann exzeptionell, wenn $M$ unzerlegbar ist und $\Ext^1_K(M,M)=0$ gilt.
\item Sei $M\in\fdar{K}$ exzeptionell. Dann ist $\dimv M$ eine relle Wurzel von $K$ (vgl. Abschnitt 1.3 in \cite{Kac} für die Definition). Sei dazu $q$ die quadratische Form des Köchers $K$. Da $M$ unzerlegbar ist, ist $\dimv M$ nach dem Satz von Kac eine Wurzel (Theorem 2 in Abschnitt 2.5 von \cite{Kac}). Es gilt $q(\dimv M)=\dim \End_K(M)-\dim\Ext^1_K(M,M)=1$. Nach \cite[Lemma 1.9] {Kac} ist $\dimv M$ reell.
Umgekehrt sind unzerlegbare Moduln zu reellen Wurzeln nicht notwendig exzeptionell: Betrachte dazu den Köcher vom Typ $\widetilde{D}_4$ mit der Orientierung
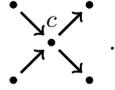
\begin{center}
\parbox[c]{1.3cm}{
\begin{tikzpicture}[line width=1pt]
	\tikzstyle{knt}=[circle, fill, inner sep=1pt]
  \node at (0,0.5)[knt](z1){};
  \node at (0,-0.5)[knt](z2){};
  
  \node at (0.5,0)[knt,label=above:\footnotesize$c$](c){};
  
  \node at (1,0.5)[knt](u1){};
  \node at (1,-0.5)[knt](u2){};
    
  \draw[->, shorten >= 2pt, shorten <= 2pt] (z1) to (c);
  \draw[->, shorten >= 2pt, shorten <= 2pt] (z2) to (c);
  \draw[->, shorten >= 2pt, shorten <= 2pt] (c) to (u1);
  \draw[->, shorten >= 2pt, shorten <= 2pt] (c) to (u2);
\end{tikzpicture}}.
\end{center}
Sei $M$ die $\widetilde{D}_4$-Darstellung
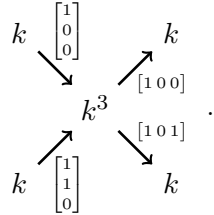
\begin{center}
\parbox[c]{2.7cm}{
\begin{tikzpicture}[line width=1pt]
  \node at (0,1)(z1){$k$};
  \node at (0,-1)(z2){$k$};
  
  \node at (1,0)(c){$k^3$};
  
  \node at (2,1)(u1){$k$};
  \node at (2,-1)(u2){$k$};
    
  \draw[->] (z1) to node[right,yshift=14pt, xshift=-6pt]{\tiny$\begin{bmatrix}1\\0\\0\end{bmatrix}$} (c);
  \draw[->] (z2) to node[right,yshift=-14pt, xshift=-6pt]{\tiny$\begin{bmatrix}1\\1\\0\end{bmatrix}$}(c);
  \draw[->] (c) to node[right, yshift=-6pt, xshift=-4pt]{\tiny$\begin{bmatrix}1\,0\,0\end{bmatrix}$} (u1);
  \draw[->] (c) to node[right, yshift=6pt, xshift=-4pt]{\tiny$\begin{bmatrix}1\,0\,1\end{bmatrix}$} (u2);
\end{tikzpicture}}.
\end{center}
Man rechnet nach: $\End M\simeq k[X]/(X^2)$. Also ist $M$ unzerlegbar, nicht exzeptionell. $\dimv M$ ist aber eine reelle Wurzel.
\item Sei $M\in\fdar{K}$ exzeptionell. Sei $N$ eine unzerlegbare $K$-Darstellung mit $\dimv M=\dimv N$. Dann gilt nach dem Satz von Kac schon $M\simeq N$ (vgl. Abschnitt 2.5, Theorem 2 in \cite{Kac}).
\item Präprojektive und präinjektive $K$-Darstellungen sind stets exzeptionell (vgl. \cite[VIII.2.7]{Sko}). Die Umkehrung ist im allgemeinen falsch: Betrachte den Köcher vom Typ $\widetilde{D}_4$ mit der Orientierung
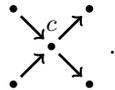
\begin{center}
\parbox[c]{1.3cm}{
\begin{tikzpicture}[line width=1pt]
	\tikzstyle{knt}=[circle, fill, inner sep=1pt]
  \node at (0,0.5)[knt](z1){};
  \node at (0,-0.5)[knt](z2){};
  
  \node at (0.5,0)[knt,label=above:\footnotesize$c$](c){};
  
  \node at (1,0.5)[knt](u1){};
  \node at (1,-0.5)[knt](u2){};
    
  \draw[->, shorten >= 2pt, shorten <= 2pt] (z1) to (c);
  \draw[->, shorten >= 2pt, shorten <= 2pt] (z2) to (c);
  \draw[->, shorten >= 2pt, shorten <= 2pt] (c) to (u1);
  \draw[->, shorten >= 2pt, shorten <= 2pt] (c) to (u2);
\end{tikzpicture}
}.
\end{center}
Der einfache Modul $E_c$ zum Punkt $c$ ist exzeptionell, hat aber Defekt $0$. Also ist $E_c$ weder präprojektiv noch präinjektiv (für genauere Informationen über den Defekt von Darstellungen von Köchern vom euklidischen Typ vgl. \cite[§7]{CB1}).
\item Sei $(x_1,\ldots,x_n)$ eine Aufzählung von $K_0$, so daß gilt: Ist $i<j$, so existiert kein Weg von $i$ nach $j$. Dann sind die Folgen
\[(E_{x_n},\ldots,E_{x_1}),(P_{x_1},\ldots,P_{x_n}),(I_{x_1},\ldots,I_{x_n})\]
exzeptionelle Folgen von $K$-Darstellungen. Beachte dazu die Gleichungen
\begin{gather*}
\dim \Hom_K(P_{x_j},P_{x_i})=\#\{\text{ Wege von $x_i$ nach $x_j$ } \}=\dim \Hom_K(I_{x_j},I_{x_i}),\\
\dim \Ext^1_K(E_{x_j},E_{x_i})=\#\{\text{ Pfeile von $x_j$ nach $x_i$ } \}\text{ für alle } i,j\in K_0.
\end{gather*}
Wir werden später sehen (vgl. \cref{scho:orko:prop1}): Jede exzeptionelle Folge von $A$-Moduln hat Länge $\leq n$.
\end{enumerate}
\end{bem}

\counterwithout{equation}{chapter}
\begin{bspl}\label{scho:ex:bsp1}Sei $n\geq0$ eine natürliche Zahl. Wir wollen alle exzeptionellen Folgen von $\Lambda_n$-Moduln klassifizieren. Sei dazu $x$ die Quelle und $y$ die Senke  von $Q_n$. Für einen $\Lambda_n$-Modul schreiben wir $\dimv M=(\dim M(x), \dim M(y))\in\mathbb{Z}^2$.
Ist $n=0$, so sind $E_x$ und $E_y$ die exzeptionelle Objekte und $(E_x,E_y)$ und $(E_y, E_x)$ die exzeptionellen Folgen. 

Ist $n=1$, so sind $P_0:=E_y$, $P_1:=P_x$ und $P_2:=E_x$ die Unzerlegbaren von $\Lambda_n$. Es ist $\dim \Hom_{\Lambda_1}(P_1,P_2)=1$, alle $P_i$ sind exzeptionell und man hat eine exakte Sequenz
\[0\longrightarrow P_0\longrightarrow P_1\longrightarrow P_2\longrightarrow 0.\]
$(P_2,P_0)$, $(P_0,P_1)$ und $(P_1,P_2)$ sind die exzeptionellen Folgen der Länge 2 in $\fmod{\Lambda_1}$.

Sei jetzt $n\geq2$. Die präprojektive Komponente des Auslander-Reiten-Köchers von $\Lambda_n$ hat die Form
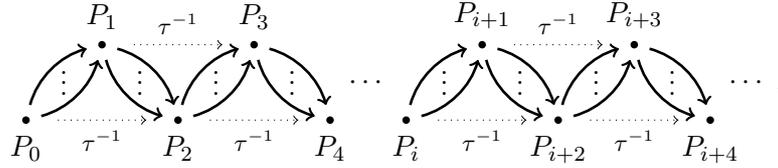
\begin{center}
\parbox[c]{10.2cm}{
\begin{tikzpicture}[line width=1pt]
	\tikzstyle{knt}=[circle, fill, inner sep=1pt]
  \node at (1,0)[knt, label=above:$P_1$](z1){};
  \node at (3,0)[knt, label=above:$P_3$](z2){};
  
  \node at (0,-1)[knt, label=below:$P_0$](u1){};
  \node at (2,-1)[knt, label=below:$P_2$](u2){};
  \node at (4,-1)[knt, label=below:$P_4$](u3){};
  
  \draw[->, bend left, shorten >= 4pt, shorten <= 4pt] (u1) to (z1);
  \draw[->, bend right, shorten >= 4pt, shorten <= 4pt] (u1) to (z1);
  \draw[->, bend left, shorten >= 4pt, shorten <= 4pt] (z1) to (u2);
  \draw[->, bend right, shorten >= 4pt, shorten <= 4pt] (z1) to (u2);
  \draw[->, bend left, shorten >= 4pt, shorten <= 4pt] (u2) to (z2);
  \draw[->, bend right, shorten >= 4pt, shorten <= 4pt] (u2) to (z2);
  \draw[->, bend left, shorten >= 4pt, shorten <= 4pt] (z2) to (u3);
  \draw[->, bend right, shorten >= 4pt, shorten <= 4pt] (z2) to (u3);
  
  \draw[->, thin, dotted, shorten >= 10pt, shorten <= 10pt] (u1) to node[below]{\footnotesize $\tau^{-1}$}(u2);
  \draw[->, thin, dotted, shorten >= 10pt, shorten <= 10pt] (u2) to node[below]{\footnotesize $\tau^{-1}$} (u3);
  \draw[->, thin, dotted, shorten >= 10pt, shorten <= 10pt] (z1) to node[above]{\footnotesize $\tau^{-1}$} (z2);
  
  \draw (4.5, -0.5) node{$\cdots$};
  
  \draw (0.5, -0.4) node{$\vdots$};
  \draw (1.5, -0.4) node{$\vdots$};
  \draw (2.5, -0.4) node{$\vdots$};
  \draw (3.5, -0.4) node{$\vdots$};
  
  \node at (6,0)[knt, label=above:$P_{i+1}$](zz1){};
  \node at (8,0)[knt, label=above:$P_{i+3}$](zz2){};
  
  \node at (5,-1)[knt, label=below:$P_{i}$](uu1){};
  \node at (7,-1)[knt, label=below:$P_{i+2}$](uu2){};
  \node at (9,-1)[knt, label=below:$P_{i+4}$](uu3){};
  
  \draw[->, bend left, shorten >= 4pt, shorten <= 4pt] (uu1) to (zz1);
  \draw[->, bend right, shorten >= 4pt, shorten <= 4pt] (uu1) to (zz1);
  \draw[->, bend left, shorten >= 4pt, shorten <= 4pt] (zz1) to (uu2);
  \draw[->, bend right, shorten >= 4pt, shorten <= 4pt] (zz1) to (uu2);
  \draw[->, bend left, shorten >= 4pt, shorten <= 4pt] (uu2) to (zz2);
  \draw[->, bend right, shorten >= 4pt, shorten <= 4pt] (uu2) to (zz2);
  \draw[->, bend left, shorten >= 4pt, shorten <= 4pt] (zz2) to (uu3);
  \draw[->, bend right, shorten >= 4pt, shorten <= 4pt] (zz2) to (uu3);
  
  \draw[->, thin, dotted, shorten >= 10pt, shorten <= 10pt] (uu1) to node[below]{\footnotesize $\tau^{-1}$} (uu2);
  \draw[->, thin, dotted, shorten >= 10pt, shorten <= 10pt] (uu2) to node[below]{\footnotesize $\tau^{-1}$}(uu3);
  \draw[->, thin, dotted, shorten >= 10pt, shorten <= 10pt] (zz1) to node[above]{\footnotesize $\tau^{-1}$} (zz2);
  
  \draw (9.5, -0.5) node{$\cdots$};
  
  \draw (5.5, -0.4) node{$\vdots$};
  \draw (6.5, -0.4) node{$\vdots$};
  \draw (7.5, -0.4) node{$\vdots$};
  \draw (8.5, -0.4) node{$\vdots$};
\end{tikzpicture}}.\end{center}
Die präinjektive Komponente hat die Form
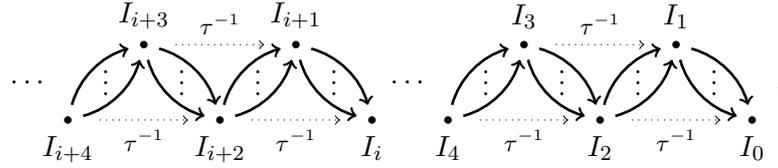
\begin{center}
\parbox[c]{10.2cm}{
\begin{tikzpicture}[line width=1pt]
	\tikzstyle{knt}=[circle, fill, inner sep=1pt]
  \node at (1,0)[knt, label=above:$I_{i+3}$](z1){};
  \node at (3,0)[knt, label=above:$I_{i+1}$](z2){};
  
  \node at (0,-1)[knt, label=below:$I_{i+4}$](u1){};
  \node at (2,-1)[knt, label=below:$I_{i+2}$](u2){};
  \node at (4,-1)[knt, label=below:$I_i$](u3){};
  
  \draw[->, bend left, shorten >= 4pt, shorten <= 4pt] (u1) to (z1);
  \draw[->, bend right, shorten >= 4pt, shorten <= 4pt] (u1) to (z1);
  \draw[->, bend left, shorten >= 4pt, shorten <= 4pt] (z1) to (u2);
  \draw[->, bend right, shorten >= 4pt, shorten <= 4pt] (z1) to (u2);
  \draw[->, bend left, shorten >= 4pt, shorten <= 4pt] (u2) to (z2);
  \draw[->, bend right, shorten >= 4pt, shorten <= 4pt] (u2) to (z2);
  \draw[->, bend left, shorten >= 4pt, shorten <= 4pt] (z2) to (u3);
  \draw[->, bend right, shorten >= 4pt, shorten <= 4pt] (z2) to (u3);
  
  \draw[->, thin, dotted, shorten >= 10pt, shorten <= 10pt] (u1) to node[below]{\footnotesize $\tau^{-1}$}(u2);
  \draw[->, thin, dotted, shorten >= 10pt, shorten <= 10pt] (u2) to node[below]{\footnotesize $\tau^{-1}$} (u3);
  \draw[->, thin, dotted, shorten >= 10pt, shorten <= 10pt] (z1) to node[above]{\footnotesize $\tau^{-1}$} (z2);
  
  \draw (4.5, -0.5) node{$\cdots$};
  
  \draw (0.5, -0.4) node{$\vdots$};
  \draw (1.5, -0.4) node{$\vdots$};
  \draw (2.5, -0.4) node{$\vdots$};
  \draw (3.5, -0.4) node{$\vdots$};
  
  \node at (6,0)[knt, label=above:$I_3$](zz1){};
  \node at (8,0)[knt, label=above:$I_1$](zz2){};
  
  \node at (5,-1)[knt, label=below:$I_4$](uu1){};
  \node at (7,-1)[knt, label=below:$I_2$](uu2){};
  \node at (9,-1)[knt, label=below:$I_0$](uu3){};
  
  \draw[->, bend left, shorten >= 4pt, shorten <= 4pt] (uu1) to (zz1);
  \draw[->, bend right, shorten >= 4pt, shorten <= 4pt] (uu1) to (zz1);
  \draw[->, bend left, shorten >= 4pt, shorten <= 4pt] (zz1) to (uu2);
  \draw[->, bend right, shorten >= 4pt, shorten <= 4pt] (zz1) to (uu2);
  \draw[->, bend left, shorten >= 4pt, shorten <= 4pt] (uu2) to (zz2);
  \draw[->, bend right, shorten >= 4pt, shorten <= 4pt] (uu2) to (zz2);
  \draw[->, bend left, shorten >= 4pt, shorten <= 4pt] (zz2) to (uu3);
  \draw[->, bend right, shorten >= 4pt, shorten <= 4pt] (zz2) to (uu3);
  
  \draw[->, thin, dotted, shorten >= 10pt, shorten <= 10pt] (uu1) to node[below]{\footnotesize $\tau^{-1}$} (uu2);
  \draw[->, thin, dotted, shorten >= 10pt, shorten <= 10pt] (uu2) to node[below]{\footnotesize $\tau^{-1}$}(uu3);
  \draw[->, thin, dotted, shorten >= 10pt, shorten <= 10pt] (zz1) to node[above]{\footnotesize $\tau^{-1}$} (zz2);
  
  \draw (-0.5, -0.5) node{$\cdots$};
  
  \draw (5.5, -0.4) node{$\vdots$};
  \draw (6.5, -0.4) node{$\vdots$};
  \draw (7.5, -0.4) node{$\vdots$};
  \draw (8.5, -0.4) node{$\vdots$};
\end{tikzpicture}}.\end{center}
Dabei ist $P_0=E_y$ der Projektive zum Punkt $y$ und $P_1=P_x$ der Projektive zum Punkt $x$. Genauso ist $I_0=E_x$ der Injektive zum Punkt $x$ und $I_1=I_y$ der Injektive zum Punkt $x$.

Definiere rekursiv eine Folge $(a_i)_{i\geq0}$ von natürlichen Zahlen durch
\[a_0:=0,\ a_1:=1\text{ und } a_{i+1}:=na_i-a_{i-1}\text{ für } i\geq 1.\]
Man zeigt per Induktion unter Verwendung der Auslander-Reiten-Sequenzen
\[0\longrightarrow P_{i-1}\longrightarrow P_i^n\longrightarrow P_{i+1}\longrightarrow 0,\]
daß $\dimv P_i=(a_i,a_{i+1})$ für $i\geq 0$ gilt. Dual ist $\dimv I_i=(a_{i+1},a_i)$ für $i\geq0$.

Sei $q:\mathbb{Z}^2\longrightarrow \mathbb{Z}$ die quadratische Form des Köchers $Q_n$, also 
\[q(v)=v_1^2+v_2^2-nv_1v_2\text{ für alle } v=(v_1,v_2)\in\mathbb{Z}^2.\]
Sei $s_x:\mathbb{Z}^2\longrightarrow\mathbb{Z}^2$ die Spiegelung zum Punkt $x\in (Q_n)_0$, also
\[s_x(v)=(nv_2-v_1,v_2)\text{ für alle } v=(v_1,v_2)\in\mathbb{Z}^2.\]
Außerdem sei die Abbildung $\omega:\mathbb{Z}^2\longrightarrow\mathbb{Z}^2$ gegeben durch $\omega(a,b):=(b,a)$. $q$ ist $s_y$- und $\omega$-invariant. Es gilt
\setcounter{equation}{0}
\numberwithin{equation}{zaehler}
\begin{equation}\label{scho:ex:bsp1:Gl1}
\dimv P_i=(\omega\circ s_x)^i(\dimv P_0) \text{ für alle } i\geq0.
\end{equation}
\newpage
Wir erhalten die folgende Klassifikation von exzeptionellen Folgen der Länge $\leq2$ in $\fmod{\Lambda_n}$:
\begin{enumerate}
\item Sei $M$ ein unzerlegbarer $\Lambda_n$-Modul. Dann sind äquivalent:
\begin{enumerate}[i)]
\item $M$ ist exzeptionell.
\item $M$ ist präprojektiv oder präinjektiv.
\end{enumerate}
\item Die Menge
\[
\bigl\{(P_{i-1},P_i)\mid (i\geq1)\bigr\}\cup\bigl\{(I_i,I_{i-1})\mid i\geq1\bigr\} \cup\bigl\{ (I_0,P_0)\bigr\}
\]
ist eine Liste der exzeptionellen Folgen der Länge 2. Dabei gelten zusätzlich für $i\geq1$:
\begin{equation*}
\begin{gathered}
\Hom_{\Lambda_n}(P_{i-1},P_i)\neq0,\ \Ext^1_{\Lambda_n}(P_{i-1},P_i)=0, \\ 
\Hom_{\Lambda_n}(I_i,I_{i-1})\neq0,\ \Ext^1_{\Lambda_n}(I_i,I_{i-1})=0, \\
\Hom_{\Lambda_n}(I_0,P_0)=0,\ \Ext^1_{\Lambda_n}(I_0,P_0)\neq0.
\end{gathered}
\end{equation*}
\end{enumerate}
Wir skizzieren die Beweise beider Aussagen.

(1.) $ii)\Rightarrow i)$: Präprojektive und präinjektive unzerlegbare Moduln über beliebigen Algebren sind exzeptionell, vgl. \cite[VIII.2.7]{Sko}.

$i)\Rightarrow ii)$: Wir zeigen zunächst
\[\bigl\{v\in\mathbb{Z}_{\geq0}^2\mid q(v)=1\bigr\}=\bigl\{\dimv P_i,\dimv I_i\mid i\geq0\bigr\}.\]
Für alle $i$ gilt
\[q(\dimv P_i)=q((\omega\circ s_x)^i(\dimv P_0))=q(\dimv P_0)=q(0,1)=1.\]
Genauso sieht man, daß $q(\dimv I_i)=1$ gilt für alle $i\geq0$. Das zeigt die Inklusion "`$\supseteq$"'.

Sei umgekehrt $v\in\mathbb{Z}_{\geq0}^2$ mit $q(v)=1$. Wir dürfen $v_1\leq v_2$ annehmen und zeigen per Induktion nach $v_2$: $v=\dimv P_i$ für ein $i$. Ist $v_2=1$, so gilt $v=\dimv P_0$. Sei also $v_2\geq2$. Sei \[w:=\bigl(\sigma_x\circ\omega\bigr)(v)=(nv_1-v_2,v_1).\]
Wir wollen die Induktionsvoraussetzung auf $w$ anwenden und müssen zeigen: $q(w)=1$, $w_1\leq w_2$ und $w_2<v_2$. $q(w)=1$ folgt aus der $\sigma_x\circ\omega$-Invarianz von $q$. Es ist $w_2=v_1<v_2$, sonst wäre $v_1=v_2$ und es folgte $1=q(v)=2v_1^2-nv_1^2=(2-n)v_1^2\leq0$ ($n\geq 2$ ist vorausgesetzt). Wegen $q(v)=1$ und $v_1\leq v_2$ erhält man außerdem
\[w_1v_2=nv_1v_2-v_2^2=v_1^2-1\leq v_1^2\leq v_1v_2=w_2v_2.\]
Wir können also die Induktionsvoraussetzung auf $w$ anwenden und erhalten $w=\dimv P_i$ für ein $i$. Es folgt $v=\bigl(\omega\circ s_x \bigr)(w)=\dimv P_{i+1}$. Damit folgt die Inklusion "`$\subseteq$"'.

Sei jetzt $M$ ein exzeptioneller $\Lambda_n$-Modul. Dann ist 
\[q(\dimv M)=\dim \End_{\Lambda_n}(M)- \dim\Ext^1_{\Lambda_n}(M,M)=1.\]
Nach dem eben Gezeigten ist $\dimv M=\dimv P_i$ oder $\dimv M=\dimv I_i$ für ein $i$. Da präprojektive und präinjektive unzerlegbare Moduln über Köcheralgebren bis auf Isomorphie eindeutig durch ihren Dimensionsvektor festliegen (vgl. \cite[VIII.2.3]{Aus}), folgt $M\simeq I_i$ oder $M\simeq P_i$.

(2.) Die Behauptung folgt unter Verwendung von (1.) aus folgenden Formeln für alle $i,j$:
\begin{equation}\label{scho:ex:bsp1:Gl2}
\begin{split}
\dim \Hom_{\Lambda_n}(P_i,P_j)=\dim \Hom_{\Lambda_n}(I_j,I_i)&=\left\{\begin{array}{cc}0&\text{ falls }i>j\\a_{j-i+1}&\text{ falls }i\leq j\end{array}\right.,\\
\dim \Hom_{\Lambda_n}(I_i,P_j)&=0,\\
\dim \Hom_{\Lambda_n}(P_i,I_j)&=a_{i+j}
\end{split}
\end{equation}
und
\begin{equation}\label{scho:ex:bsp1:Gl3}
\begin{split}
\dim \Ext^1_{\Lambda_n}(P_i,P_j)=\dim \Ext^1_{\Lambda_n}(I_j,I_i)&=\left\{\begin{array}{cc}0&\text{ falls }i<j+2\\a_{i-j-1}&\text{ falls }i\geq j+2\end{array}\right.,\\
\dim \Ext^1_{\Lambda_n}(I_i,P_j)&=a_{i+j+2},\\
\dim \Ext^1_{\Lambda_n}(P_i,I_j)&=0.
\end{split}
\end{equation}
Um die Formeln \ref{scho:ex:bsp1:Gl2} für die Dimension der Hom-Räume zu verifizieren, verwendet man folgende Bemerkungen:
\begin{itemize}
\item[$\cdot$]Sei $P$ ein präprojektiv-unzerlegbarer $\Lambda_n$-Modul, sei $M$ unzerlegbar. Ist $\Hom_{\Lambda_n}(M,P)\neq\nolinebreak[4]0$, so existiert ein Weg von $M$ nach $P$ im Auslander-Reiten-Köcher von $\Lambda_n$ (vgl. \cite[VIII.2.5]{Sko}).
\item[$\cdot$]$\Hom_{\Lambda_n}(P_0,M)=M(x)$ und $\Hom_{\Lambda_n}(P_1,M)=M(y)$ für alle $\Lambda_n$-Moduln $M$.
\item[$\cdot$]Für die Auslander-Reiten-Verschiebung $\tau$ gilt: Sind $M,N$ $\Lambda_n$-Moduln, die keinen projektiven direkten Summanden haben, so existiert ein Isomorphismus
\[\Hom_{\Lambda_n}(M,N)\simeq\Hom_{\Lambda_n}(\tau M, \tau N)\]
von $k$-Vektorräumen.
\item[$\cdot$]Für alle $n\geq0$ und für alle $j\geq0$ gelten die Formeln:
\[
\tau^{n}P_j=P_{j-2n}\text{ falls } j\geq2n\text{ und }\tau^{n}I_{j}=I_{2n+j}.
\]
\end{itemize}
Die Formeln (\ref{scho:ex:bsp1:Gl3}) für die Dimensionen der Ext-Gruppen folgen aus den Formeln (\ref{scho:ex:bsp1:Gl2}) unter Verwendung der Auslander-Reiten-Formeln (vgl. \cite[IV.2.13]{Sko}):
\[\Ext^1_{\Lambda_n}(M,N)\simeq D\Hom_{\Lambda_n}(N,\tau M)\text{ für alle } M,N\in\fmod{\Lambda_n}.\]
\end{bspl}
Das folgende Resultat (vgl. \cite[Corollary 2.3]{Ri2}) wird für den Beweis des Satzes von Schofield nützlich sein.
\begin{prop}\label{scho:ex:prop3}
Sei $K$ ein endlicher, zykelloser Köcher. Sei $M$ ein aufrichtiger $kK$-Modul ohne Selbsterweiterungen. Dann ist $M$ treu.
\end{prop}
Der Beweis beruht auf folgendem Lemma (vgl. \cite[Lemma 1.1]{Ri2}):
\begin{lemma}\label{scho:ex:lemma1} Sei $\mathcal{A}$ eine $k$-lineare erbliche abelsche Kategorie. Sei $M\in\mathcal{A}$, $\Ext^1_\mathcal{A}(M,M)=0$. Bezeichne mit $\operatorname{Sub}(M)$ die volle Unterkategorie aller Subquotienten aller Potenzen von $M$. Dann ist $\operatorname{Sub}(M)$ abgeschlossen unter Erweiterungen.
\end{lemma}
\begin{bew} Wir zeigen zunächst: Sei $\xi:\ 0\longrightarrow X'\longrightarrow X\longrightarrow X''\longrightarrow0$ eine exakte Sequenz in $\mathcal{A}$ und seien $X'\stackrel{\epsilon'}{\longrightarrow}Y'$, und $X''\stackrel{\epsilon''}{\longrightarrow}Y''$ zwei Homomorphismen, wobei $\epsilon''$ injektiv sei. Dann existiert ein kommutatives Diagramm
\[
\begin{CD}
0 @>>> X' @>>> X @>>> X'' @>>> 0 \\
@. @V\epsilon' VV @V\epsilon  VV @V{\epsilon''} VV @. \\
0 @>>>Y' @>>> Y @>>>Y'' @>>> 0
\end{CD}
\]
mit exakten Zeilen. Der injektive Homomorphismus $\epsilon'':X''\longrightarrow Y''$ induziert wegen $\Ext^2_{\mathcal{A}}=0$ einen surjektiven Homomorphismus
$\Ext^1_{\mathcal{A}}(Y'',Y')\longrightarrow\Ext^1_{\mathcal{A}}(X'',Y')$, $\zeta\mapsto\zeta\epsilon''$. Betrachte die Pushout-Sequenz $\epsilon'\xi\in\Ext^1_{\mathcal{A}}(X'',Y')$. Dann existiert eine Sequenz $\eta\in\Ext^1_{\mathcal{A}}(Y'',Y')$ mit $\eta=\epsilon'\xi\epsilon''$. Die entsprechenden Pullback- und Pushout-Diagramme liefern das gewünschte Diagramm.

Sei jetzt eine exakte Sequenz $0\longrightarrow U\longrightarrow Z\longrightarrow V\longrightarrow0$ gegeben mit $U,V\in\operatorname{Sub}(M)$. Per Definition von $\operatorname{Sub}(M)$ existieren Abbildungen $M^r\stackrel{\pi_U}{\longrightarrow} U'\stackrel{\epsilon_U}{\longleftarrow}U\text{ und }M^s\stackrel{\pi_V}{\longrightarrow} V'\stackrel{\epsilon_V}{\longleftarrow}V $ mit $\pi_U,\pi_V$ Epi und $\epsilon_U,\epsilon_V$ Mono.

Wir konstruieren einen Epimorphismus $M^r\oplus M^s\stackrel{\pi_Z}{\longrightarrow} Z'$ und einen Monomorphismus $Z'\stackrel{\epsilon_Z}{\longleftarrow}Z$. Nach der Vorbemerkung existiert ein kommutatives Diagramm
\[
\begin{CD}
0 @>>> U @>>> Z @>>> V @>>> 0 \\
@. @V\epsilon_U VV @V\epsilon_Z  VV @V{\epsilon_V} VV @. \\
0 @>>>U' @>>> Z' @>>>V' @>>> 0
\end{CD}
\]
mit exakten Zeilen. $\epsilon_Z$ ist injektiv, da $\epsilon_U$ und $\epsilon_V$ injektiv sind. Nach der dualen Version der Vorbemerkung existiert ein kommutatives Diagramm
\[
\begin{CD}
0 @>>> U' @>>> Z @>>> V' @>>> 0 \\
@. @A\pi_U AA @A\pi_Z  AA @A{\pi_V} AA @. \\
0 @>>>M^r @>>> T @>>>M^s @>>> 0
\end{CD}
\]
mit exakten Zeilen. Wegen $\Ext^1_\mathcal{A}(M,M)=0$ folgt $T\simeq M^r\oplus M^s$. $\pi_Z$ ist surjektiv, da $\pi_U$ und $\pi_V$ surjektiv sind. Es folgt $Z\in\operatorname{Sub}(M)$.
\end{bew}
\begin{bew}[von \cref{scho:ex:prop3}] Seien $E_1,\ldots,E_n$ die Kompositionsfaktoren von $M$, $E_i$ habe die Vielfachheit $\mu(i)$. Man rechnet nach: $\dimv M=\sum_{i=1}^n\mu(i)\dimv E_i$. Da $M$ aufrichtig ist, ist jeder einfache $kK$-Modul ein Kompositionsfaktor von $M$, also liegt jeder einfache $kK$-Modul in $\operatorname{Sub}(M)$. $\fmod{kK}$ ist erblich, also ist $\operatorname{Sub}(M)$ nach \cref{scho:ex:lemma1} abgeschlossen unter Erweiterungen. Man erhält $\operatorname{Sub}(M)=\fmod{kK}$. Insbesondere ist $kK$ ein Subquotient von $M$. Da $kK$ treu ist, muß auch $M$ treu sein.
\end{bew}
Wir kommen jetzt zur Definition von orthogonalen Kategorien.

\begin{defn}\label{scho:ex:def3}
Sei $\mathcal{A}$ eine $k$-lineare abelsche Kategorie. Sei $\mathcal{C}$ eine Klasse von Objekten in $\mathcal{A}$. Die orthogonalen Kategorien $\mathcal{C}^\bot$ und $^\bot\mathcal{C}$ in $\mathcal{A}$ sind definiert als die vollen Unterkategorien zu den folgenden Klassen von Objekten:
\begin{enumerate}[a)]
\item $\mathcal{C}^\bot:=\{X\in\mathcal{A} \mid \Hom_\mathcal{A}(Y,X)=0=\Ext^1_\mathcal{A}(Y,X) \text{ für alle } Y\in \mathcal{C}\}.$
\item $^\bot\mathcal{C}:=\{X\in\mathcal{A} \mid \Hom_\mathcal{A}(X,Y)=0=\Ext^1_\mathcal{A}(X,Y) \text{ für alle } Y\in \mathcal{C}\}.$
\end{enumerate}
\end{defn}

\begin{bem}\label{scho:ex:bem2} Sei $\mathcal{A}$ eine abelsche $k$-lineare Kategorie.
\begin{enumerate}
\item Wir werden häufig die folgende Aussage verwenden: Sei $\mathcal{B}$ eine volle Unterkategorie von $\mathcal{A}$, abgeschlossen unter Kernen, Kokernen und Erweiterungen. Sei $\iota:\mathcal{B}\longrightarrow\mathcal{A}$ der Vergißfunktor. Dann induziert $\iota$ für alle $M,N$ und für alle $i\geq 0$ injektive Homomorphismen $\epsilon_i:\Ext^i_{\mathcal{B}}(M,N)\longrightarrow\Ext^i_{\mathcal{A}}(\iota M,\iota N)$. Für $i\leq 1$ ist $\epsilon_i$ ein Isomorphismus. Insbesondere gilt: Ist $\mathcal{A}$ erblich, so ist auch $\mathcal{B}$ erblich.
\item Sei $\mathcal{A}$ erblich und sei $\mathcal{C}$ eine Klasse von Objekten in $\mathcal{A}$. Dann ist $\mathcal{C}^\bot$ abgeschlossen unter Kernen, Kokernen und Erweiterungen (also auch unter direkten Summen). Insbesondere gilt wegen (1.) für alle $M,N\in\mathcal{C}^\bot$: $\Ext^1_{\mathcal{C}^\bot}(M,N)=\Ext^1_\mathcal{A}(M,N)$ und $\mathcal{C}^\bot$ ist erblich.

Um etwa die Abgeschlossenheit unter Kernen und Kokernen zu zeigen, betrachte für einen Morphismus $f:X\longrightarrow Y$ die exakten Sequenzen
\[
\begin{array}{cccccccccc}
0 &\longrightarrow &\Ker(f)& \longrightarrow &X& \longrightarrow &\Bild(f)&\longrightarrow &0, \\
0 &\longrightarrow &\Bild(f)& \longrightarrow &Y& \longrightarrow &\Coker(f)&\longrightarrow &0. \\
\end{array}
\]
Wende auf die Sequenzen für $Z\in\mathcal{C}$ den Funktor $\Hom_{\mathcal{A}}(Z,-)$ an und werte die langen exakten $\Ext^1_{\mathcal{A}}$-Sequenzen unter Beachtung von $\Ext^2_\mathcal{A}=0$ aus.
\item Analoge Aussagen gelten für $^\bot\mathcal{C}$.
\end{enumerate}
\end{bem}

\section{Orthogonale Kategorien von Köcherdarstellungen}\label{scho:orko}
\setcounter{zaehler}{0}
\numberwithin{zaehler}{section}
Wir fixieren in diesem Abschnitt einen endlichen zykellosen Köcher $K$ und setzen $A:=kK$. Sei $\mathcal{C}=E^\bot$ eine orthogonale Kategorie in $\fmod{A}$. Wegen $\Ext^2_A=0$ wissen wir nach \cref{scho:ex:bem2}, daß $\mathcal{C}$ abgeschlossen ist unter Kernen, Kokernen und Erweiterungen, also insbesondere abelsch und erblich ist. Wir benutzen in diesem Kapitel häufig die Tatsache, daß der Vergißfunktor $\mathcal{C}\longrightarrow\fmod{A}$ Isomorphismen der $\Ext^1$-Gruppen induziert, ohne erneut darauf hinzuweisen.

Wir werden die folgende für Induktionsbeweise sehr nützliche Aussage zeigen (vgl. \cite[Theorem 2.3]{Scho}:
\begin{prop}\label{scho:orko:prop1} Sei $E$ eine exzeptionelle Folge von $A$-Moduln der Länge $r$. Dann gelten:
\begin{enumerate}
\item $r\leq\#K_0$.
\item Es existieren ein (bis auf Isomorphie eindeutig bestimmter) endlicher zykelloser Köcher $K(E^\bot)$ und eine Äquivalenz von Kategorien $E^\bot\stackrel{\sim}{\longrightarrow}\fmod{kK(E^\bot)}$. Dabei gilt \[\#K(E^\bot)_0=\#K_0-r.\]
\item Es existieren ein (bis auf Isomorphie eindeutig bestimmter) endlicher zykelloser Köcher $K(^\bot E)$ und eine Äquivalenz von Kategorien $^\bot E\stackrel{\sim}{\longrightarrow}\fmod{kK(^\bot E)}$. Dabei gilt
\[\#K(^\bot E)_0=\#K_0-r.\]
\end{enumerate}
\end{prop}

Der wichtigste Baustein für den Beweis der Proposition ist das folgende Lemma (vgl. \cite[Theorem 3.1]{Scho}):
\begin{lemma}\label{scho:orko:lemma1} Sei $X$ ein exzeptioneller, nicht projektiver $A$-Modul und sei 
\[\xi:\ 0\longrightarrow A \stackrel{\epsilon}{\longrightarrow} \widetilde{X}\stackrel{\pi}{\longrightarrow} X^r\longrightarrow 0\]
die minimale universelle Erweiterung von Objekten in $\operatorname{add}(X)$ durch $A$. Dann gilt:
\begin{enumerate}
\item $\widetilde{X}\in X^\bot$.
\item $\widetilde{X}$ ist ein Progenerator von $X^\bot$.
\item Es existieren genau $\#K_0-1$ Isomorphieklassen unzerlegbarer direkter Summanden von $\widetilde{X}$.
\end{enumerate}
\end{lemma}
\begin{bew} (1.) Wegen $\Ext^1_A(X,X)=0$ ist $X$ ein partieller Kippmodul (beachte: $A$ ist erblich, also ist $\operatorname{pdim} X\leq 1$ klar). Nach \cref{scho:kipp:prop2} ist $X\oplus\widetilde{X}$ ein Kippmodul, insbesondere ist $\Ext^1_A(X,\widetilde{X})=0$. Außerdem ist $\Hom_A(X,\widetilde{X})=0$: $\xi$ induziert eine exakte Sequenz
\[0\longrightarrow \Hom_A(X,A)\longrightarrow \Hom_A(X,\widetilde{X})\longrightarrow\Hom_A(X,X^r)\stackrel{\partial}{\longrightarrow}\Ext^1(X,A).\ \ \ (*)\]
Da $\xi$ eine minimale universelle Erweiterung von Objekten in $\operatorname{add}(X)$ durch $A$ ist und $X$ als exzeptionelles Objekt Endomorphismenring $k$ hat, ist $\partial$ wegen der definierenden Eigenschaft von $\xi$ ein Isomorphismus (verwende die duale Version von \cref{scho:kipp:bem3}). Außerdem ist $X$ unzerlegbar, nicht projektiv; da $A$ erblich ist, folgt also $\Hom_A(X,A)=0$. Die Exaktheit der Sequenz $(*)$ zeigt dann die Behauptung.

(2.) $\widetilde{X}$ ist ein projektives Objekt in $X^\bot$: Sei $M\in X^\bot$. $\xi$ induziert eine exakte Sequenz
\[\Ext^1_A(X^r,M)\longrightarrow\Ext^1_A(\widetilde{X},M)\longrightarrow\Ext^1_A(A,M).\]
Wegen $\Ext^1_A(X,M)=0=\Ext^1_A(A,M)$ folgt $\Ext^1_{X^\bot}(\widetilde{X},M)=\Ext^1_A(\widetilde{X},M)=0$.

$\widetilde{X}$ erzeugt $X^\bot$: Sei $M$ ein Objekt in $X^\bot$ und $\phi:A^n\longrightarrow M$ ein Epimorphismus in $\fmod{A}$ für ein $n\in\mathbb{N}$. Betrachte das Pushout-Diagramm
\[
\begin{CD}
0 @>>> A^n @>>>\widetilde{X}^n @>>> X^{nr} @>>> 0 \\
@. @V\phi VV @V\psi VV @| @. \\
0 @>>>M @>\sigma>> P @>>>X^{nr} @>>> 0.
\end{CD}
\]
Hier ist die erste Zeile einfach die $n$-fache direkte Summe von $\xi$. Wegen $\Ext^1_A(X,M)=0$ spaltet die zweite Zeile, also besitzt $\sigma$ eine Retraktion $\tau$. Nach dem Schlangenlemma ist mit $\phi$ auch $\psi$ ein Epimorphismus, und wir erhalten einen Epimorphismus $\tau\circ\psi:\widetilde{X}^n\longrightarrow M$ in $\fmod{A}$; $\tau\circ\psi$ ist erst recht ein Epimorphismus in $X^\bot$.

(3.) Wir haben bereits in (1.) gesehen, daß $X\oplus \widetilde{X}$ ein Kippmodul ist. Also hat $X\oplus \widetilde{X}$ nach \cref{scho:kipp:satz1} genau $\#K_0$ Isomorphieklassen unzerlegbarer direkter Summanden. Wegen $\Hom_A(X,\widetilde{X})=0$ haben $X$ und $\widetilde{X}$ keine zueinander isomorphen direkten Summanden, also existieren nach dem Satz von Krull-Remak-Schmidt genau $\#K_0-1$ Isomorphieklassen unzerlegbarer direkter Summanden von $\widetilde{X}$.
\end{bew}

Damit können wir \cref{scho:orko:prop1} beweisen:
\begin{bew}[von \cref{scho:orko:prop1}] Es genügt, Aussage (2.) zu zeigen. Sei $E=(X_1,\ldots,X_r)$ eine exzeptionelle Folge von $A$-Moduln der Länge $r$; außerdem sei $n:=\#K_0$. Wir verwenden Induktion nach $r$.

Im Induktionsanfang sei $r=1$. Setze $X:=X_1$.

Ist $X$ projektiv, so ist $X=P_x$ für ein $x\in K_0$. Es folgt 
\[\Ext^1_A(X,M)=0,\ \Hom_A(X,M)=e_xM=M(x)\text{ für alle }M\in\fmod{A},\]
also gilt $X^\bot = \{M\in \fmod{kK}\mid M(x)=0\}$. Sei $K(X^\bot)$ der volle Unterköcher von $K$ mit Punktmenge $K(X^\bot)_0=K_0 - \{x\}$. Dann folgt offenbar $X^\bot \simeq \fmod{kK(E^\bot)}$.

Sei jetzt $X$ nicht projektiv. Sei $\widetilde{X}$ die minimale universelle Erweiterung von Objekten in $\operatorname{add}(X)$ durch $A$. Sei $B:=\End(\widetilde{X})^{op}$, und sei $K(X^\bot)$ der Gabriel-Köcher von $B$. Nach \cref{scho:orko:lemma1} ist $\widetilde{X}$ ein Progenerator von $X^\bot$ mit genau $n-1$ Isomorphieklassen unzerlegbarer direkter Summanden. Es folgt, daß genau $n-1$ Isomorphieklassen von einfachen $B$-Moduln existieren. Also hat $K(X^\bot)$ $n-1$ Punkte. Da $A$ erblich ist, wissen wir nach \cref{scho:ex:bem2}, daß auch $X^\bot$ eine erbliche Kategorie ist. Mit \cref{scho:ex:prop2} folgt $X^\bot\simeq\fmod{kK(X^\bot)}$.

Im Induktionsschritt sei $r\geq2$. Sei $\mathcal{C}:=(X_2,\ldots,X_r)^\bot$. Dann ist $X_1\in\mathcal{C}$. Per Induktion existieren ein Köcher $K'$ mit $n-r+1$ Punkten und eine Äquivalenz $F:\mathcal{C} \stackrel{\sim}{\longrightarrow} \fmod{kK'}$.

$X_1$ ist auch ein exzeptionelles Objekt der Unterkategorie $\mathcal{C}$, also ist $FX_1$ ein exzeptionelles Objekt in $\fmod{kK'}$. Nach dem Induktionsanfang existieren also ein Köcher $K(E^\bot)$ mit 
\[\#\bigl(K(E^\bot)\bigr)_0=\#K'_0-1=n-r\]
und eine Äquivalenz $(FX_1)^{\bot}\simeq \fmod{K(E^\bot)}$. Sei \[X_1^{\bot(\mathcal{C})}=\bigl\{M\in\mathcal{C}\mid \Hom_\mathcal{C}(X_1,M)=0=\Ext^1_\mathcal{C}(X_1,M)\bigr\}.\]
Dann ist $X_1^{\bot(\mathcal{C})}=E^\bot$, und $F$ induziert offenbar eine Äquivalenz $X_1^{\bot(\mathcal{C})}\simeq (FX_1)^\bot$.
\end{bew}

\begin{defn}\label{scho:orko:def1}
Sei $E$ eine exzeptionelle Folge von $A$-Moduln. Es sei $C(E)$ die kleinste volle Unterkategorie von $\fmod{A}$, die $E$ enthält und die abgeschlossen ist unter Kernen, Kokernen und Erweiterungen.
\end{defn}

\begin{lemma}[\protect{\cite[Lemma 3]{CB2}}] \label{scho:orko:lemma2}
Sei $E$ eine exzeptionelle Folge von $A$-Moduln. Dann gilt: Ist $E$ vollständig, so folgt $C(E)=\fmod{A}$.
\end{lemma}
\begin{bew}
Induktion nach $n:=\#K_0$. Ist $n=1$, so ist $A=k$. Der einzige exzeptionelle $A$-Modul ist dann $k$ selbst; natürlich ist $C(k)=\fmod{k}$. 

Sei jetzt $n\geq2$ und sei $E=(X_1,\ldots,X_n)$. Nach \cref{scho:orko:prop1} existieren ein Köcher $K'$ mit $n-1$ Punkten und eine Äquivalenz $X_n^\bot \simeq \fmod{kK'}$. Es ist $X_i\in X_n^\bot$ für alle $i=1,\ldots,n-1$. Es gilt: 
\begin{itemize}
\item [$\cdot$]$E':=(X_1,\ldots,X_{n-1})$ ist eine exzeptionelle Folge in $X_n^\bot$.
\item [$\cdot$]$C^{(X_n^\bot)}(E')=C(E')$, wobei $C^{(X_n^\bot)}(E')$ den Abschluß von $E'$ unter Kernen, Kokernen und Erweiterungen von $E'$ in $X_n^\bot$ bezeichnet.
\end{itemize}
Wir identifizieren $X_n^\bot$ mit $\fmod{kK'}$ und erhalten per Induktion:   $C^{(X_n^\bot)}(E')=X_n^\bot$. Also gilt
\[X_n^\bot=C(E')\subseteq C(E)\text{ und } X_n\in C(E).\]
Diese Gleichungen zeigen die Behauptung:

1. Fall: $X_n$ ist projektiv, also $X_n=P_x$ für ein $x\in K_0$. Wir behaupten, daß $C(E)$ alle einfachen $A$-Moduln enthält; da $C(E)$ abgeschlossen ist unter Erweiterungen, folgt dann $C(E)=\fmod{A}$. Es ist $X_n^\bot=\{M\mid M(x)=0\}$, also ist $E_y\in X_n^\bot\subseteq C(E)$ für alle $y\in K_0$ mit $y\neq x$. Um auch $E_x\in C(E)$ zu zeigen, betrachte die exakte Sequenz
\[0\longrightarrow\Rad P_x\longrightarrow P_x\longrightarrow E_x\longrightarrow 0.\]
Wegen $\Rad P_x(x)=0$ ist $\Rad P_x\in X_n^\bot\subseteq C(E)$, außerdem ist $P_x=X_n\in C(E)$. Da $C(E)$ abgeschlossen unter Kokernen ist, folgt $E_x\in C(E)$.

2. Fall: $X_n$ ist nicht projektiv. Sei $\widetilde{X}_n$ der Mittelterm der minimalen universellen Erweiterung von Objekten in $\operatorname{add}(X_n)$ durch $A$. Nach \cref{scho:orko:lemma1} ist $\widetilde{X}_n\in X_n^\bot\subseteq C(E)$, außerdem ist $X_n\in C(E)$. Da $C(E)$ abgeschlossen unter Kernen ist, folgt $A\in C(E)$. Ist $M$ ein beliebiger $A$-Modul, so existieren natürliche Zahlen $r$ und $s$ und eine exakte Sequenz \[A^r\longrightarrow A^s\longrightarrow M\longrightarrow 0.\]
Mit $A\in C(E)$ gilt auch $A^r, A^s\in C(E)$. Da $C(E)$ abgeschlossen unter Kokernen ist, folgt $M\in C(E)$.
\end{bew}

\begin{prop}[\protect{\cite[Lemma 5]{CB2}}]\label{scho:orko:prop2}
Sei $E$ eine exzeptionelle Folge der Länge $r$. Dann existieren ein (bis auf Isomorphie eindeutig bestimmter) endlicher, zykelloser Köcher $K(E)$ mit $r$ Punkten und eine Äquivalenz $C(E)\simeq \fmod{kK(E)}$.
\end{prop}
\begin{bew}
Sei $n:=\#K_0$. Dann ist gemäß \cref{scho:orko:prop1} $^\bot E$ die Modulkategorie einer endlichdimensionalen Köcheralgebra mit $n-r$ Einfachen. Insbesondere existiert in $^\bot E$ eine (vollständige) exzeptionelle Folge $E'$ der Länge $n-r$ (die Einfachen in geeigneter Reihenfolge sind etwa eine vollständige exzeptionelle Folge). $E'$ ist dann auch eine exzeptionelle Folge in $A$-mod, und es folgt nach \cref{scho:orko:prop1}: $E'^\bot$ ist die Modulkategorie einer Köcheralgebra mit $r$ Einfachen. Nun gilt
\[E\subseteq ({}^\bot E)^\bot \subseteq E'^\bot.\]
$E$ ist also eine vollständige exzeptionelle Folge in $E'^\bot$. Wir erhalten also nach \cref{scho:orko:lemma2} $C(E)=C^{(E'^\bot)}(E)=E'^\bot$ und die Behauptung folgt.
\end{bew}

\section{Schofield-Induktion}
\setcounter{zaehler}{0}
\numberwithin{zaehler}{section}\label{scho:scho} In diesem Abschnitt sei stets $K$ ein endlicher zykelloser Köcher und $A:=kK$. Wir legen einige Notationen fest: Mit $x$ bezeichnen wir in diesem Abschnitt immer die Quelle, mit $y$ die Senke von $Q_n$. Für eine Klasse $E$ von $A$-Moduln sei wie im vorangegangenen Abschnitt $C(E)$ die kleinste volle Unterkategorie von $\fmod{A}$, die $E$ enthält und die abgeschlossen ist unter Kernen, Kokernen und Erweiterungen. Für $A$-Moduln $X,Y$ sei $\mathcal{F}(X,Y)$ die volle Unterkategorie von $\fmod{A}$ bestehend aus allen $A$-Moduln $M$, für die eine exakte Sequenz 
\[0\longrightarrow V\longrightarrow M\longrightarrow W\longrightarrow 0\]
mit $V\in\operatorname{add}(Y)$ und $W\in\operatorname{add}(X)$ existiert.

Wir formulieren das Hauptresultat dieses Abschnitts (vgl. \cite[Theorem 3.1]{Ri2}):
\begin{satz}\label{scho:scho:satz1}
Sei $M$ ein exzeptioneller, nicht-einfacher $A$-Modul. Dann existieren exzeptionelle $A$-Moduln $X$ und $Y$ mit folgenden Eigenschaften:
\begin{enumerate}[i)]
\item $\Hom_A(Y,X)=0=\Hom_A(X,Y)$.
\item Es existieren natürliche Zahlen $u\geq1$ und $v\geq1$ und eine exakte Sequenz
\[0\longrightarrow Y^v\longrightarrow M \longrightarrow X^u\longrightarrow 0.\]
\end{enumerate}
Ist $r:=\#Tr(M)$, so existieren genau  $r-1$ Isomorphieklassen von Paaren $(X,Y)$ von $A$-Moduln $X$ und $Y$ mit diesen Eigenschaften.

Zusatz: $\mathcal{F}(Y,X)$ ist abgeschlossen unter Kernen, Kokernen und Erweiterungen und ist äquivalent zur Modulkategorie einer verallgemeinerten Kronecker-Algebra $\Lambda_n$ für ein $n\geq1$. $Y$ ist dabei das projektiv-einfache und $X$ das injektiv-einfache Objekt in $\mathcal{F}(Y,X)$. Insbesondere ist $(Y,X)$ eine Liste der einfachen Objekte in $\mathcal{F}(Y,X)$.
\end{satz}
\begin{bem}\label{scho:scho:bem1}\ 
Mit Hilfe der exakten Sequenz in ii) sieht man leicht: $\Ext^1_A(Y,X)=0$, also ist $(X,Y)$ eine exzeptionelle Folge.
\end{bem}
Die folgende Proposition beweist unter anderem den Zusatz in \cref{scho:scho:satz1}.
\begin{prop}\label{scho:scho:prop1}\enlargethispage{\baselineskip}
Sei $(X,Y)$ eine exzeptionelle Folge von $A$-Moduln. Dann gilt:
\begin{enumerate}
\item Es existieren eine natürliche Zahl $n\geq0$ und eine Äquivalenz $C(X,Y)\simeq\fmod{\Lambda_n}$.
\item Es gelte zusätzlich $\Hom_A(X,Y)=0$. Sei $F:C(X,Y)\stackrel{\sim}{\longrightarrow}\fmod{\Lambda_n} $ eine Äquivalenz. Dann gilt:
\begin{enumerate}[i)]
\item $(X,Y)$ ist eine Liste der einfachen Objekte in $C(X,Y)$. Ist $n\geq1$, so gilt $FX\simeq E_x$ und $\ FY\simeq E_y$.
\item $\mathcal{F}(X,Y)=C(X,Y)$. Insbesondere ist $\mathcal{F}(X,Y)$ abgeschlossen unter Kernen, Kokernen und Erweiterungen.
\end{enumerate}
\end{enumerate}
\end{prop}
\begin{bew}
(1.) ist ein Spezialfall von \cref{scho:orko:prop2}. 

(2.) i) Offenbar sind $FX$ und $FY$ exzeptionelle $\Lambda_n$-Moduln mit \[\Hom_{\Lambda_n}(FY,FX)=0=\Hom_{\Lambda_n}(FX,FY)\text{ und } \Ext^1_{\Lambda_n}(FY,FX)=0.\]
Nach der Klassifikation aller vollständigen exzeptionellen Folgen in $\fmod{\Lambda_n}$ in \cref{scho:ex:bsp1} muß also im Fall $n\geq1$ gelten: $FY\simeq E_y$ und $FX\simeq E_x$ und $(Y,X)$ ist eine Liste der einfachen Objekte; im Fall $n=0$ ist klar, daß $X,Y$ die einzigen einfachen Objekte in $C(X,Y)$ sind.

ii) Klar ist: $\mathcal{F}(X,Y)\subseteq C(X,Y)$. Sei umgekehrt $M\in C(X,Y)$. Dann existiert eine exakte Sequenz
\[0\longrightarrow E_y^v\longrightarrow FM\longrightarrow E_x^u\longrightarrow 0\]
von $\Lambda_n$-Moduln. Da $F:C(X,Y)\simeq \fmod{\Lambda_n}$ eine Äquivalenz ist, erhält man wegen  i) eine exakte Sequenz
\[0\longrightarrow Y^v\longrightarrow M\longrightarrow X^u\longrightarrow 0.\]
Also ist $M\in\mathcal{F}(Y,X)$.
\end{bew}
\begin{bem}\label{scho:scho:bem2}\ 
\begin{enumerate}
\item Konkrete Beschreibungen von Funktoren $\fdar{Q_n}\longrightarrow \mathcal{F}(Y,X)$ in einem allgemeineren Rahmen findet man in Kapitel \ref{unz}.
\item Mit Hilfe der Proposition folgt: Die Menge aller Isomorphieklassen von Paaren von $A$-Mo\-duln $(X,Y)$, die die Bedingungen \cref{scho:scho:satz1} erfüllen, stimmt überein mit der Menge aller Isomorphieklassen exzeptioneller Folgen $(X,Y)$ von $A$-Moduln, so daß $\Hom_A(X,Y)=\nolinebreak[4]0$ gilt und $M$ nicht einfach ist in $\mathcal{F}(X,Y)$.
\end{enumerate}
\end{bem}

Wir geben zunächst eine äquivalente Formulierung von \cref{scho:scho:satz1} an (vgl. \cite[Teil 3]{Ri2} und \cite[Lemma 7]{CB2}):
\begin{satz}\label{scho:scho:satz2}
Sei $M$ ein exzeptioneller, nicht-einfacher $A$-Modul. Dann existiert ein exzeptioneller $A$-Modul $N$, so daß $(N,M)$ eine exzeptionelle Folge von $A$-Moduln ist und so daß $M$ nicht einfach in $C(N,M)$ ist.

Ist $r:=\#Tr(M)$, so existieren genau $r-1$ Ismorphieklassen von $A$-Moduln mit dieser Eigenschaft.
\end{satz}

Wir werden \cref{scho:scho:satz2} beweisen. Das folgende Lemma zeigt, daß \cref{scho:scho:satz1} und \cref{scho:scho:satz2} zueinander äquivalent sind:
\begin{lemma}\label{scho:scho:lemma1} Sei $M$ ein exzeptioneller, nicht einfacher $A$-Modul. Es sei $\mathcal{S}$ die Menge aller Isomorphieklassen von exzeptionellen Folgen $(X,Y)$ von $A$-Moduln mit $\Hom_A(X,Y)=0$ und der Eigenschaft, daß $M$ nicht einfach ist in $\mathcal{F}(Y,X)$. Sei $\mathcal{T}$ die Menge aller Isomorphieklassen von $A$-Moduln $N$, so daß $(N,M)$ eine exzeptionelle Folge ist und $M$ nicht einfach in $C(N,M)$ ist.

Dann existiert eine Bijektion $\Theta:\mathcal{T}\stackrel{\sim}{\longrightarrow}\mathcal{S}$ mit folgender Eigenschaft: Ist $N\in\mathcal{T}$, so ist \linebreak$\Theta(N)=(X,Y)$ das eindeutig bestimmte Element von $\mathcal{S}$ mit $X,Y\in C(N,M)$. Außerdem gilt $\mathcal{F}(X,Y)=C(N,M)$.
\end{lemma}
\begin{bew} Sei ein $A$-Modul $N\in \mathcal{T}$ gegeben. Wir zeigen zunächst: Es existiert genau ein Paar $(X,Y)\in\mathcal{S}$ mit $X,Y\in C(N,M)$. Dabei gilt $\mathcal{F}(X,Y)=C(N,M)$. Nach \cref{scho:scho:prop1} existiert eine Äquivalenz \[F:\fmod{\Lambda_n}\stackrel{\sim}{\longrightarrow}C(N,M)\]
für ein $n\geq0$. Da $M$ nicht einfach ist in $C(N,M)$, muß offenbar $n\geq1$ sein. Sei $Y=FE_y$ das projektiv-einfache und $X=FE_x$ das injektiv-einfache Objekt in $C(N,M)$. Dann gilt\linebreak $\mathcal{F}(X,Y)=C(N,M)$, und es folgt $(X,Y)\in\mathcal{S}$. Mit der Klassifikation der exzeptionellen Folgen von $\Lambda_n$-Moduln in \cref{scho:ex:bsp1} erhält man wegen $n\geq1$: $(E_x, E_y)$ ist die einzige exzeptionelle Folge in $\fmod{\Lambda_n}$ mit $\Hom_{\Lambda_n}(E_x,E_y)=0$. Also ist $(X,Y)$ die einzige exzeptionelle Folge in $C(N,M)$ mit $\Hom_A(X,Y)=0$ und die Behauptung folgt. 

Sei $\Theta:\mathcal{T}\longrightarrow\mathcal{S}$ die so definierte Abbildung. Wir müssen zeigen, daß $\Theta$ bijektiv ist. Sei $(X,Y)\in\mathcal{S}$ gegeben. Nach \cref{scho:scho:prop1} ist dann 
\[\mathcal{F}(X,Y)=C(X,Y)\simeq\fmod{\Lambda_n}\]
für ein $n\geq0$. $M$ ist ein exzeptionelles Objekt in $\mathcal{F}(X,Y)$, also existiert nach der Klassifikation in \cref{scho:ex:bsp1} genau ein exzeptionelles Objekt $N$ in $\mathcal{F}(X,Y)$, so daß $(N,M)$ eine exzeptionelle Folge in $\mathcal{F}(X,Y)$ ist. Wir haben in \cref{scho:orko:lemma2} gesehen, daß für jede exzeptionelle Folge $E$ der Länge 2 in $\fmod{\Lambda_n}$ schon $C(E)=\fmod{\Lambda_n}$ gilt. Also folgt $C(N,M)=\mathcal{F}(X,Y)$. Das zeigt $N\in\mathcal{T}$ und $\Theta(N)=(X,Y)$, also ist $\Theta$ surjektiv. Ist $\Theta(N)=\Theta(N')$ für ein $N'\in\mathcal{T}$, so folgt
\[\mathcal{F}(X,Y)=C(N,M)=\mathcal{F}(\Theta(N))=\mathcal{F}(\Theta(N'))=C(N',M).\]
Also ist auch $(N',M)$ eine exzeptionelle Folge in $\mathcal{F}(X,Y)$. Die Eindeutigkeit von $N$ zeigt $N\simeq N'$, also ist $\Theta$ injektiv.
\end{bew}
Wir benötigen das folgende Hilfsresultat (vgl. \cite[Lemma 3.2]{Ri2}):
\begin{lemma}\label{scho:scho:lemma2} Sei $(N,M)$ eine exzeptionelle Folge von $A$-Moduln. Dann sind äquivalent:
\begin{enumerate}
\item $M$ ist nicht einfach in $C(N,M)$.\enlargethispage{\baselineskip}
\item Es existieren eine natürliche Zahl $r$ und ein injektiver Homomorphismus $N\longrightarrow M^r$.
\end{enumerate}
\end{lemma}
\begin{bew}
Nach \cref{scho:scho:prop1} existiert eine Äquivalenz $C(N,M)\simeq\fmod{\Lambda_n}$ für ein $n\geq0$. Wir dürfen annehmen, daß schon $A=\Lambda_n$ und $C(N,M)=\fmod{\Lambda_n}$ gilt und beweisen die Äquivalenz der beiden Aussagen per Fallunterscheidung nach $n$.

Ist $n=0$, so ist weder Bedingung (1) noch Bedingung (2) erfüllt, also sind (1) und (2) zueinander äquivalent. 

Sei $n=1$. $(1)\Rightarrow(2)$: Ist $M$ nicht einfach, so muß wegen der Klassifikation von exzeptionellen Folgen von $\Lambda_1$-Moduln in \cref{scho:ex:bsp1} $M=P_x$ und $N=E_y$ gelten. Offenbar existiert dann ein injektiver Homomorphismus $N\longrightarrow M$. $(2)\Rightarrow(1)$: Angenommen, $M$ ist einfach. Dann folgt $(N,M)=(P_x,E_x)$ oder $(N,M)=(E_x, E_y)$ und in keinem der beiden Fälle ist (2.) erfüllt.

Sei also $n\geq2$. Wegen \cref{scho:ex:bsp1} folgt $(N,M)=(P_{i-1}, P_{i})$,  $(N,M)=(I_{i},I_{i-1})$ für ein $i\geq1$ oder $(N,M)=(I_0,P_0)$. Ist $(N,M)=(I_0,P_0)$, so ist weder (1) noch (2) erfüllt. Sei $(N,M)=(P_{i-1}, P_{i})$ für ein $i\geq1$. Dann ist (1) erfüllt und die Auslander-Reiten-Sequenz 
\[0\longrightarrow P_{i-1}\longrightarrow P_i^n\longrightarrow P_{i+1}\longrightarrow0\]
mit Anfang $P_{i-1}$ zeigt, daß auch (2) erfüllt ist. Sei schließlich $(N,M)=(I_{i}, I_{i-1})$ für ein $i\geq1$. $(1)\Rightarrow(2)$: Da $M$ nicht einfach ist, gilt $i\geq2$. Die Auslannder-Reiten Sequenz
\[0\longrightarrow I_{i}\longrightarrow I_{i-1}^n\longrightarrow I_{i-2}\longrightarrow0\]
mit Ende $I_{i-2}$ liefert den gewünschten Monomorphismus in $(2)$. $(2)\Rightarrow(1)$: Angenommen, $M$ ist einfach. Dann ist $N=I_1=I_y$. In dem Fall existiert keine injektive Abbildung $N\longrightarrow M^r$, denn der Untermodul $E_y$ von $N$ liegt im Kern jedes solchen Homomorphismus, also ist $(2)$ nicht erfüllt.
\end{bew}

\begin{defn}\label{scho:scho:def1} Ist $M$ ein exzeptioneller $A$-Modul, so definieren wir das Bongartz-Komple\-ment $\widetilde{M}$ von $M$ wie folgt: Ist $M$ projektiv, $M=P_s$ für ein $s\in K_0$, so setze 
\[\widetilde{M}:=\bigoplus_{\begin{subarray}{c}t\in K_0\\ t\neq s\end{subarray}}P_t.\]
Ist $M$ nicht projektiv, so sei $\widetilde{M}$ der Mittelterm der minimalen universellen Erweiterung von Objekten in $\operatorname{add}(M)$ durch $A$.
\end{defn}
\begin{bem}\label{scho:scho:bem3} Wegen Teil (3) von \cref{scho:orko:lemma1} gilt: Es existieren genau $\#K_0-1$ Isomorphieklassen von unzerlegbaren direkten Summanden von $\widetilde{M}$.
\end{bem}

Wir können die Situation in \cref{scho:scho:satz2} vollständig analysieren (vgl. \cite[Lemma 3.4]{Ri2}):
\begin{satz}\label{scho:scho:satz3}
Sei $M$ ein aufrichtiger exzeptioneller $A$-Modul. Sei $N$ ein weiterer $A$-Modul. Dann sind äquivalent:
\begin{enumerate}
\item $N$ ist ein unzerlegbarer direkter Summand von $\widetilde{M}$.
\item $(N,M)$ ist eine exzeptionelle Folge und $M$ ist nicht einfach in $C(N,M)$.
\end{enumerate}
\end{satz}
\begin{bew}
$(1)\Rightarrow(2)$ Da $M$ ein aufrichtiger $A$-Modul ohne Selbsterweiterungen ist, ist $M$ nach \cref{scho:ex:prop3} treu. Also existiert ein injektiver Homomorphismus $\phi: A\longrightarrow M^r$ für ein $r$ (wähle eine Basis $(m_1,\ldots,m_r)$ von $M$ und setze $\phi(1):=m_1+\ldots+m_r$).

Sei zunächst $M$ projektiv, also $M=P_s$ für ein $s\in K_0$. Dann ist $N=P_t$ für ein $t\in K_0-\{s\}$. Da $M$ aufrichtig ist, gilt $0\neq\Hom_A(P_t,M)=\Hom_A(P_t,P_s)=e_tAe_s$. Also existiert ein Weg von $s$ nach $t$ in $K$. Da $K$ zykellos ist folgt $e_sAe_t=0$, also $\Hom_A(M,N)=0$. $\Ext^1_A(M,N)=0$ ist klar, also ist $(N,M)$ eine exzeptionelle Folge. Da $N$ ein Untermodul von $A$ ist, erhalten wir einen Monomorphismus $N\longrightarrow M^r$. Nach \cref{scho:scho:lemma2} ist dann $M$ nicht einfach in $C(N,M)$.

Sei jetzt $M$ nicht projektiv. Dann ist $\widetilde{M}\in M^\bot$ nach \cref{scho:orko:lemma1}. Insbesondere liegt jeder direkte Summand von $\widetilde{M}$ in $M^\bot$, also ist $(N,M)$ eine exzeptionelle Folge. Es bleibt, zu zeigen, daß $M$ nicht einfach in $C(N,M)$ ist. Betrachte dazu das Pushout-Diagramm
\[
\begin{CD}
0 @>>> A @>>> \widetilde{M} @>>> M^r @>>> 0 \\
@. @V\phi VV @V\epsilon VV @| @. \\
0 @>>>M^s @>>> T @>>>M^r @>>> 0.
\end{CD}
\]
Hier ist die erste Zeile die das Bongartz-Komplement $\widetilde{M}$ definierende universelle Erweiterung. Wegen $\Ext^1_A(M,M)=0$ spaltet die untere Sequenz, also ist $T\simeq M^{r+s}$. Da $\phi$ injektiv ist, ist auch $\epsilon$ injektiv. Das liefert einen injektiven Homomorphismus $N\longrightarrow M^{r+s}$, und nach \cref{scho:scho:lemma2} ist dann $M$ nicht einfach in $C(N,M)$.

$(2)\Rightarrow(1)$ Nach \cref{scho:scho:lemma2} existiert ein injektiver Homomorphismus $\phi:N\longrightarrow M^r$.

Sei zunächst $M$ projektiv, $M=P_s$. Da $A$ erblich ist, ist $N\simeq\Bild\phi$ projektiv, also $N=P_t$ für ein $t\in K_0$. Es ist $0=\Hom_A(M,N)=\Hom_A(P_s,P_t)=e_sAe_t$. Also muß $s\neq t$ gelten, also ist $N$ ein direkter Summand von $\widetilde{M}$.

Sei jetzt $M$ nicht projektiv. Nach \cref{scho:orko:lemma1} ist $\widetilde{M}$ ein Progenerator von $M^\bot$. Um zu zeigen, daß $N$ ein direkter Summand von $\widetilde{M}$ ist, genügt es also, zu zeigen: $N$ ist projektiv in $M^\bot$. Sei dazu $N'\in M^\bot$. Da $A$ erblich ist, induziert $\phi$ einen surjektiven Homomorphismus $\Ext^1_A(M^r,N')\longrightarrow\Ext^1_A(N,N')$. Aus $\Ext^1_A(M,N')=0$ folgt $\Ext^1_{M^\bot}(N,N')=\Ext^1_A(N,N')=\nolinebreak[4]0$ (beachte \cref{scho:ex:bem2}).
\end{bew}

\cref{scho:scho:satz2} ist jetzt eine einfache Folgerung:

\begin{bew}[von \cref{scho:scho:satz2}] Sei ein exzeptioneller, nicht-einfacher $A$-Modul $M$ gegeben. Sei $\mathcal{C}$ die volle Unterkategorie von $\fmod{A}$ bestehend aus allen Moduln $L$ mit $Tr(L)\subseteq Tr(M)$. $\mathcal{C}$ ist abgeschlossen unter Unter- und Quotientendarstellungen und unter Erweiterungen. 

Sei $N$ ein $A$-Modul, so daß $M$ nicht einfach in $C(N,M)$ ist. Dann ist schon $N\in\mathcal{C}$: Wegen \cref{scho:scho:prop1} existiert eine Äquivalenz $F:\fmod{\Lambda_n}\longrightarrow C(N,M)$ für ein $n\geq0$ mit $Y:=\nolinebreak[4]FE_y$ und $X:=FE_x$. Also existiert eine exakte Sequenz 
\[0\longrightarrow Y^v\longrightarrow M \longrightarrow X^u\longrightarrow 0\]\enlargethispage{\baselineskip}
mit $u,v\geq1$. Es folgt $C(N,M)=\mathcal{F}(Y,X)\subseteq\mathcal{C}$.

Wegen dieser Bemerkungen dürfen wir nach Übergang zum Träger von $M$ annehmen: $M$ ist aufrichtig. Wir können also \cref{scho:scho:satz3} anwenden: Die $A$-Moduln $N$, so daß $M$ nicht einfach ist in $C(N,M)$, sind gerade die unzerlegbaren direkten Summanden von $\widetilde{M}$. Nach \cref{scho:scho:bem3} existieren bis auf Isomorphie genau $\#K_0-1=\#Tr(M)-1$ unzerlegbare direkte Summanden von $\widetilde{M}$. Damit folgt die Behauptung.
\end{bew}

\begin{bem}\label{scho:scho:bem4} Sei $M$ ein nicht exzeptioneller, unzerlegbarer $A$-Modul. Wir wollen skizzieren, wie man bei Kenntnis einer Zerlegung des Bongartz-Komplements von $M$ in Unzerlegbare die $A$-Moduln $X$ und $Y$ aus der Schofield-Induktion (\cref{scho:scho:satz1}) konstruieren kann. Wir dürfen ohne Einschränkung annehmen, daß $M$ aufrichtig ist. Wir verwenden die Aussagen und Notationen für exzeptionelle $\Lambda_n$-Moduln aus \cref{scho:ex:bsp1}.

Wir beginnen mit einer Vorbemerkung. Sei $n\geq1$. Betrachte die exzeptionelle Folge $(P_i,P_{i+1})$ in $\fmod{\Lambda_n}$ für $i\geq0$ (im Falle $n=1$ existieren nur die Fälle $i=0,1$). Nach der Formel \ref{scho:ex:bsp1:Gl2} gilt
\[\dim \Hom_{\Lambda_n}(P_i,P_{i+1})=n.\]
Sei $(\phi_1,\ldots,\phi_n)$ eine Basis von $\Hom_{\Lambda_n}(P_i,P_{i+1})$ und sei $\Phi_i:P_i^n\longrightarrow P_{i+1}$ der von den $\phi_j$ induzierte Homomorphismus. $\Phi_i$ erfüllt die folgende universelle Eigenschaft: Ist \[f:P_i^r\longrightarrow P_{i+1}\] ein beliebiger Homomorphismus, so faktorisiert $f$ eindeutig durch $\Phi_i$, d.h. es existiert ein eindeutig bestimmter Homomorphismus $f':P_i^r\longrightarrow P_i^n$ mit $f=\Phi_i\circ f'$ (beachte bei der Verifikation dieser Aussage $\End_{\Lambda_n}(P_i)=k$). Wir nennen $\Phi_i$ den universellen Homomorphismus von Objekten in $\operatorname{add}(P_i)$ nach $P_{i+1}$. Dann gilt:
\begin{enumerate}[i)]
\item $\Phi_0$ ist injektiv und es ist $\Coker \Phi_0=E_x$.
\item Für $i\geq1$ ist die Sequenz 
\[0\longrightarrow \Ker\Phi_i\longrightarrow P_i^n\stackrel{\Phi_i}{\longrightarrow} P_{i+1}\longrightarrow 0\]
exakt und stimmt mit der Auslander-Reiten-Sequenz mit Ende $P_{i+1}$ überein. Insbesondere ist $\Ker\Phi_i\simeq P_{i-1}$.
\end{enumerate}
i) rechnet man direkt nach. Zum Beweis von ii) betrachte die Auslander-Reiten-Sequenz
\[0\longrightarrow P_{i-1}\stackrel{\epsilon}{\longrightarrow}P_i^n\stackrel{\pi}{\longrightarrow} P_{i+1}\longrightarrow 0\]
mit Ende $P_{i+1}$. Da $\Phi_i$ universell ist, existiert ein kommutatives Diagramm
\[
\begin{CD}
0@>>> \Ker\Phi_i@>>> P_i^n@>\Phi_i>> P_{i+1}\\
@. @AiAA @AAA @| @.\\
0@>>> P_{i-1}@>\epsilon>> P_i^n@>\pi>> P_{i+1}@>>>0\\
\end{CD}
\]
mit exakten Zeilen. Also muß $\Phi_i$ surjektiv sein. Wäre die Abbildung $i$ kein Schnitt, so würde $i$ per Definition der Auslander-Reiten-Sequenz über $\epsilon$ faktorisieren. Also würde die exakte Sequenz\enlargethispage{\baselineskip}
\[0\longrightarrow \Ker\Phi_i\longrightarrow P_i^n\stackrel{\Phi_i}{\longrightarrow} P_{i+1}\longrightarrow 0\]
spalten, ein Widerspruch. Also ist $i$ ein Schnitt, aus Dimensionsgründen muß $i$ schon ein Isomorphismus sein und die Behauptung folgt.

Sei jetzt $N$ ein direkter unzerlegbarer Summand des Bongartz-Komplements $\widetilde{M}$ von $M$. Es gilt (vgl. \cref{scho:scho:satz3} und \cref{scho:scho:prop1}):
\begin{enumerate}[i)]
\item $(N,M)$ ist eine exzeptionelle Folge, und $M$ ist nicht einfach in $C(N,M)$
\item Es existiert eine Äquivalenz $F:\fmod{\Lambda_n}\stackrel{\sim}{\longrightarrow}C(N,M)$ für ein $n\geq 1$.
\end{enumerate}
Wir geben einen Algorithmus zur Bestimmung der in \cref{scho:scho:satz1} gesuchten $A$-Moduln $X$ und $Y$ und der exakten Sequenz für $M$ an.

Wir nehmen $\dim N<\dim M$ an, den Fall $\dim N>\dim M$ behandelt man völlig analog (der Fall $\dim M=\dim N$ tritt nicht auf, da $M$ nicht einfach ist). Es folgt $(N,M)=(P_i,P_{i+1})$ für ein $i\geq 0$.

Es sei $\Phi_N:N^n\longrightarrow M$ der universelle Homomorphismus von Objekten in $\operatorname{add}(N)$ nach $M$. Ist $\Phi_N$ injektiv, so ist wegen der Vorbemerkung $(N,M)=(FP_0,FP_1)=(FE_y,FP_x)$. Setze $Y:=N$ und $X:=\Coker \Phi_N=FE_x$. Offenbar haben $X,Y$ die gewünschten Eigenschaften. 

Ist $\Phi_N$ nicht injektiv, so ist $i\geq 1$. Sei $N':=\Ker\Phi$. Wegen der Vorbemerkung ist $N'=\nolinebreak[4]FP_{i-1}$. Betrachte dann den universellen Morphismus $\Phi_{N'}:N'^n\longrightarrow N$. Nach endlich vielen Schritten muß $\Phi_{N'}$ injektiv sein und dann gilt $N'=FP_0=FE_y$ und $N=FP_1=FP_x$. Wir setzen $Y:=N'=FE_y$ und $X:=\Coker\Phi_{N'}=FE_x$. Dann erfüllen $X$ und $Y$ die gewünschten Eigenschaften. Außerdem sieht man leicht mit der entsprechenden Eigenschaft von $\Lambda_n$-Moduln, daß die universellen Homomorphismen $Y^v\longrightarrow M$ von Objekten in $\operatorname{add}(Y)$ nach $M$ und $M\longrightarrow \nolinebreak[4]X^u$ von $M$ in Objekte von $\operatorname{add}(X)$ eine exakte Sequenz
\[0\longrightarrow Y^v\longrightarrow M\longrightarrow X^u\longrightarrow 0\]
liefern.
\end{bem}

\section{Exzeptionelle Objekte in erblichen Kategorien}\label{scho:erb}
\setcounter{zaehler}{0}
\numberwithin{zaehler}{section}
\cref{scho:scho:satz1} aus dem letzten Abschnitt läßt sich mit Resultaten aus \cite[Teil 2]{Ri2} auf beliebige endlichdimensionale, erbliche Längenkategorien verallgemeinern (eine $k$-lineare Kategorie heißt endlichdimensional, falls alle Homomorphismenräume endlichdimensional über $k$ sind. Sie heißt eine Längenkategorie, falls sie abelsch ist und alle Objekte endliche Länge haben). Insbesondere kann die Voraussetzung, daß der Köcher zykellos ist, fallengelassen werden.

Der Prototyp einer endlichdimensionalen Längenkategorie ist die Kategorie $\fmod{A}$ aller endlichdimensionalen Moduln einer $k$-Algebra $A$. In völliger Analogie zu Modulkategorien definiert man für endlichdimensionale Längenkategorien halbeinfache Objekte, Radikale, Sockel und Loewy-Reihen. \newpage

Sei $\mathcal{A}$ eine endlichdimensionale Längenkategorie. Wir wiederholen kurz einige wesentliche Aussagen:

\begin{enumerate}
\item Es gilt Schurs Lemma. Die Länge von Objekten in $\mathcal{A}$ ist eindeutig bestimmt.
\item Einfache Objekte haben Endomorphismenring $k$: Ist $E$ einfach, so ist nämlich $\End_\mathcal{A}(E)$ ein endlichdimensionaler Schiefkörper über dem algebraisch abgeschlossenen Körper $k$.
\item Faktor- und Unterobjekte von halbeinfachen Objekten sind wieder halbeinfach: Sei $X$ ein halbeinfaches Objekt. Man verifiziert mit Schurs Lemma, daß jede exakte Sequenz
\[0\longrightarrow X''\longrightarrow X\longrightarrow X'\longrightarrow 0\]
mit einfachem $X'$ spaltet. Anschließend beweist man per Induktion nach der Länge von $X$ und mit Hilfe des Schlangenlemmas, daß jede exakte Sequenz mit Mittelterm $X$ spaltet.
\item Zu jedem Objekt $X\in\mathcal{A}$ existiert ein eindeutig bestimmtes kleinstes Unterobjekt $\Rad X$, so daß $\Dec X:=X/\Rad X$ halbeinfach ist. Die Projektion $\pi:X\longrightarrow \Dec X$ erfüllt folgende universelle Eigenschaft: Ist $Y$ halbeinfach und $f:X\longrightarrow Y$ ein Morphismus, so existiert genau ein Morphismus $\overline{f}:\Dec X\longrightarrow Y$ mit $f=\overline{f}\circ\pi$. Der Funktor $\Dec$ von $\mathcal{A}$ in die volle Unterkategorie der halbeinfachen Objekte ist also linksadjungiert zum Vergißfunktor.

Die Existenz des Radikals $\operatorname{Rad} X$ von $X$ folgt aus der Tatsache, daß die Familie $\mathcal{F}$ aller Unterobjekte $U$ von $X$ mit $X/U$ halbeinfach abgeschlossen ist unter endlichen Durchschnitten (bzw. der Bildung des Pullbacks bezüglich der Inklusionsmorphismen).
\item Es gilt Nakayamas Lemma: Ist $\Rad M=M$, so ist $M=0$. Die Projektion \[\pi:X\longrightarrow \Dec X\] ist ein wesentlicher Epimorphismus.

Nakayamas Lemma folgt aus der Definition des Radikals, die Wesentlichkeit der Projektion zum Beispiel mit der Rechtsexaktheit des Funktors $\Dec$ und Nakayamas Lemma.
\item Die Loewy-Länge oder Höhe $\Ht(X)$ eines Objekts $X$ ist die kleinste natürliche Zahl $n$, für die eine Kette von Unterobjekten $M=M_n\supset M_{n-1}\ldots\supset M_1\supset M_0$ mit $M_i/M_{i-1}$ halbeinfach für alle $i=1,\ldots,n$ existiert. Sie stimmt überein mit der kleinsten natürlichen Zahl $h$, für die $\Rad^h M=0$ gilt.
\item Ein Objekt $X$ ist genau dann einreihig (d.h.: $X$ besitzt eine bis auf Isomorphie von Filtrierungen eindeutig bestimmte Kompositionsreihe), wenn Höhe und Länge von $X$ übereinstimmen.
\end{enumerate}

Für den Rest dieses Kapitels sei stets $\mathcal{A}$ eine endlichdimensionale, erbliche Längenkategorie. 

\begin{defn}\label{scho:erb:def1} Sei $M\in\mathcal{A}$. Die (endliche) Menge $Tr(M)$ aller Isomorphieklassen von Kompositionsfaktoren von $M$ heißt der Träger von $M$. Außerdem sei $\operatorname{Sub}(M)$ die volle Unterkategorie von $\mathcal{A}$ bestehend aus allen Subquotienten aller Potenzen von $M$.
\end{defn}
\newpage
\begin{lemma}\label{scho:erb:lemma1} Sei $M\in\mathcal{A}$ mit $\Ext^1_\mathcal{A}(M,M)=0$. Dann gilt:
\begin{enumerate}
\item $\operatorname{Sub}(M)$ ist abgeschlossen unter Unter- und Faktorobjekten und unter Erweiterungen.
\item Für alle $N,N'\in\operatorname{Sub}(M)$ gilt $\Ext^1_A(N,N')=\Ext^1_{\operatorname{Sub}(M)}(N,N')$.
\item $\operatorname{Sub}(M)$ ist eine endlichdimensionale, erbliche Längenkategorie.
\item $\operatorname{Sub}(M)$ stimmt überein mit der vollen Unterkategorie von $\mathcal{A}$ aller Objekte $X$, so daß $Tr(X)\subseteq Tr(M)$ gilt.
\end{enumerate}
\end{lemma}
\begin{bew}
(1.) Die Abgeschlossenheit unter Unter- und Faktorobjekten ist klar; die Abgeschlossenheit unter Erweiterungen wurde bereits in \cref{scho:ex:lemma1} bewiesen.

(2.), (3.) und (4.) folgen unter Beachtung von \cref{scho:ex:bem2} aus (1.).
\end{bew}
Wir wollen folgende Aussage zeigen:
\begin{satz}[\protect{\cite[Corollary 2.2]{Ri2}}]\label{scho:erb:satz1}
Sei $M\in\mathcal{A}$ mit $\Ext^1_\mathcal{A}(M,M)=0$. Dann existieren ein endlicher zykelloser Köcher $K$ und eine Äquivalenz $\operatorname{Sub}(M)\simeq\fmod{kK}$.
\end{satz}

Um unsere Notationen etwas zu vereinfachen, vereinbaren wir:
\begin{defn}\label{scho:erb:def2}  Der Köcher $K_{\mathcal{A}}$ zu $\mathcal{A}$ ist definiert wie folgt:
\begin{itemize}
\item[$\cdot$] Die Punkte von $K_{\mathcal{A}}$ sind die Isomorphieklassen von einfachen Objekten in $\mathcal{A}$.
\item[$\cdot$] Für zwei Punkten $E,E'$ existiert genau dann ein (und nur ein) Pfeil $E\rightarrow E'$, wenn $\Ext^1_\mathcal{A}(E,E')\neq0$ gilt.
\end{itemize}
\end{defn}

\begin{lemma}\label{scho:erb:lemma2}
Sei $n\geq1$ und sei $E_n\longrightarrow E_{n-1}\longrightarrow\ldots\longrightarrow E_1\longrightarrow E_0$ ein Weg in $K_{\mathcal{A}}$. Dann existiert eine Folge von Objekten $0=M_0\subset M_1\subset\ldots\subset M_{n-1}\subset M_n$, so daß gilt:
\begin{itemize}
\item[$\cdot$] $\Rad M_i = M_{i-1}$ für alle $i$.
\item[$\cdot$] $M_{i+1}/M_i=E_i$ für alle $i=1,\ldots, n-1$ .
\end{itemize}
Insbesondere sind alle $M_i$ einreihig.
\end{lemma}
\begin{bew}
Konstruktion der $M_i$ per Induktion nach $i$. Wir setzen $M_0:=0$ und $M_1:=E_0$.
Sei jetzt $2\leq i\leq n$ und seien die Moduln $M_0,\ldots,M_{i-1}$ bereits konstruiert. Dann existiert eine exakte Sequenz
\[0\longrightarrow M_{i-2}\stackrel{\epsilon}{\longrightarrow} M_{i-1}\stackrel{\pi}{\longrightarrow} E_{i-2}\longrightarrow0\]
mit $\Rad M_{i-1}=M_{i-2}$. Wegen $\Ext_\mathcal{A}^1(E_{i-1},E_{i-2})\neq0$ existiert eine nicht spaltende exakte Sequenz 
\[\xi:\ 0\longrightarrow E_{i-2}\longrightarrow M'\longrightarrow E_{i-1}\longrightarrow0.\]
Da $\mathcal{A}$ erblich ist, induziert $\pi$ eine surjektiven Homomorphismus
\[\Ext_\mathcal{A}^1(E_{i-1},M_{i-1})\longrightarrow\Ext^1_\mathcal{A}(E_{i-1},E_{i-2}),\ \eta\mapsto\pi\eta.\]
Also existiert eine exakte Sequenz \[\eta:\ 0\longrightarrow M_{i-1}\longrightarrow M_i\longrightarrow E_{i-1}\longrightarrow0\] mit $\xi=\pi\eta$.

Wir erhalten ein kommutatives Diagramm
\[
\begin{CD}
@. 0 @. 0\\
@. @VVV @VVV\\
@.  M_{i-2} @= M_{i-2}  \\
@.  @V\epsilon VV @VVV \\
\eta:\ 0 @>>>M_{i-1} @>>>M_i @>>>E_{i-1} @>>>0 \\
@. @V\pi VV @VVV @|\\
\xi:\ 0 @>>>E_{i-2} @>>>M' @>>>E_{i-1} @>>>0\\
@. @VVV @VVV\\
@. 0 @. 0
\end{CD}
\]
mit exakten Zeilen und Spalten. Per Konstruktion ist $M_{i-1}$ ein Unterobjekt von $M_i$ mit $M_i/M_{i-1}=E_{i-1}$, und wir müssen nur noch $\Rad M_i=M_{i-1}$ zeigen. Da $M_i/M_{i-1}$ einfach ist, gilt $\Rad M_i\subseteq M_{i-1}$. Weil $M_{i-1}/\Rad M_i$ als Unterobjekt von $M_i/\Rad M_i$ halbeinfach ist, folgt $\Rad M_{i-1}\subseteq \Rad M_i$. Da $\Rad M_{i-1}=M_{i-2}$ ein maximales Unterobjekt von $M_{i-1}$ ist, folgt $\Rad M_i=M_{i-2}$ oder $\Rad M_i=M_{i-1}$. Angenommen, es ist $\Rad M_i=M_{i-2}$. Es folgt $M'=\Dec M_i$, also ist $M'$ halbeinfach, also spaltet die Sequenz $\xi$, ein Widerspruch. Das beendet den Beweis.
\end{bew}

\begin{lemma}[\protect{\cite[Proposition 2.1]{Ri2}}]\label{scho:erb:lemma3}
Sei $M\in\mathcal{A}$ mit $\Ext^1_\mathcal{A}(M,M)=0$. Dann gilt:
\begin{enumerate}
\item $K_{\operatorname{Sub}(M)}$ ist zykellos.
\item Für alle $S,T\in Tr(M)$ ist $\dim\Ext_\mathcal{A}^1(S,T)<\infty$.
\end{enumerate}
\end{lemma}
\begin{bew}
(1.) Wir wissen bereits (vgl. \cref{scho:erb:lemma1}), daß $\operatorname{Sub}(M)$ eine endlichdimensionale, erbliche Längenkategorie ist. Wir dürfen also \cref{scho:erb:lemma2} auf $\operatorname{Sub}(M)$ anwenden. Die Loewy-Längen von Objekten in $\operatorname{Sub}(M)$ sind beschränkt durch $\Ht_\mathcal{A}(M)$. Ein Zykel in $K_{\operatorname{Sub}(M)}$ würde einen unendlich langen Weg in $K_{\operatorname{Sub}(M)}$ liefern. Mit der Konstruktion in \cref{scho:erb:lemma2} findet man eine Folge $(X_i)_{i\geq0}$ in $\operatorname{Sub}(M)$ mit $\Ht X_i=i$ für alle $i$, ein Widerspruch.

(2.) Sei $M=M_r\supset M_{r-1}\supset\ldots\supset M_1\supset M_0=0$ eine Kompositionsreihe von $M$. Setze\linebreak $S_i:=M_i/M_{i-1}$ für alle $i=1,\ldots,r$. Wir müssen zeigen: $\dim\Ext_\mathcal{A}^1(S_i,S_j)<\infty$ für alle $i,j$. Betrachte für alle $i$ die exakten Sequenzen
\[0\longrightarrow M_{i-1}\longrightarrow M_i\longrightarrow S_i\longrightarrow 0\]
und
\[0\longrightarrow M_i\longrightarrow M\longrightarrow M/M_i\longrightarrow 0.\]
Anwenden geeigneter Hom-Funktoren und Auswerten der induzierten langen Ext-Sequenzen unter Beachtung der Endlichdimensionalität und Erblichkeit von $\mathcal{A}$ sowie der Voraussetzung $\Ext^1_\mathcal{A}(M,M)=0$ liefert die Behauptung durch einfache, aber langwierige Rechnungen.
\end{bew}

Wir können jetzt \cref{scho:erb:satz1} beweisen:
\begin{bew}[von \cref{scho:erb:satz1}]Da $\operatorname{Sub}(M)$ erblich ist, genügt es zu zeigen: $\operatorname{Sub}(M)$ besitzt einen Progenerator (vgl. \cref{scho:ex:prop2}). Seien dazu $E_1,\ldots,E_n$ die Kompositionsfaktoren von $M$ in $\mathcal{A}$. Nach \cref{scho:erb:lemma3} ist $K_{\operatorname{Sub}(M)}$ endlich und zykellos. Die $E_i$ sind die Punkte von $K_{\operatorname{Sub}(M)}$. Nach eventueller Umnumerierung dürfen wir also annehmen: Für $i\leq j$ existiert kein Weg der Länge $\geq1$ in $K_{\operatorname{Sub}(M)}$ von $E_i$ nach $E_j$. Also gilt für $i\leq j$:
\[\Ext^1_\mathcal{A}(E_i,E_j)=\Ext^1_{\operatorname{Sub}(M)}(E_i,E_j)=0.\]
Wir behaupten: Für $i=1,\ldots,n$ existieren ein projektives Objekt $P_i\in\operatorname{Sub}(M)$ und ein Epimorphismus $P_i\longrightarrow E_i$.

Dazu bemerken wir: Sei $P\in\operatorname{Sub}(M)$. Ist $\Ext_\mathcal{A}^1(P,E_i)=0$ für alle $i=1,\ldots,n$, so ist $P$ schon projektiv in $\operatorname{Sub}(M)$ (man zeigt für alle $N\in\operatorname{Sub}(M)$ per Induktion nach der Länge von $N$: $\Ext^1_{\operatorname{Sub}(M)}(P,N)=\Ext^1_{\mathcal{A}}(P,N)=0$).

Wir konstruieren die $P_i$ per Induktion nach $i$. Setze $P_1:=E_1$. Wegen $\Ext_\mathcal{A}^1(E_1,E_i)$ für alle $i$ ist $P_1$ projektiv. Sei im Induktionsschritt $i>1$. Per Induktion existieren für $j=1,\ldots,i-1$ projektive Objekte $P_j\in\operatorname{Sub}(M)$ und Epimorphismen $P_j\longrightarrow E_j$. Betrachte die minimale universelle Erweiterung
\[\xi: 0\longrightarrow T\longrightarrow\widetilde{E}_i\longrightarrow E_i\longrightarrow0\]
von $E_i$ durch Objekte in $\operatorname{add}(E_1\oplus\ldots\oplus E_{i-1})$ (beachte: Diese universelle Erweiterung existiert, da wegen \cref{scho:erb:lemma3} die involvierten Ext-Gruppen endlichdimensional sind!). Wegen \linebreak$T\in\operatorname{add}(E_1\oplus\ldots\oplus E_{i-1})$ existiert ein Epimorphismus $\pi:U\longrightarrow T$ mit $U\in\operatorname{add}(P_1\oplus\ldots\oplus P_{i-1})$. $U$ ist insbesondere projektiv. Da $\mathcal{A}$ erblich ist, induziert $\pi$ einen Epimorphismus \[\Ext_\mathcal{A}^1(E_i,U)\longrightarrow\Ext_\mathcal{A}^1(E_i,T).\] Wir erhalten ein kommutatives Diagramm
\[
\begin{CD}
\eta: 0 @>>> U @>>> P_i @>>> E_i @>>> 0 \\
@. @V\pi VV @V\pi'VV @| @. \\
\xi: 0 @>>>T @>>> \widetilde{E}_i @>>>E_i @>>> 0
\end{CD}
\]
mit exakten Zeilen. $\pi'$ ist ein Epimorphismus, da $\pi$ ein Epimorphismus ist. $\operatorname{Sub}(M)$ ist abgeschlossen unter Erweiterungen, also gilt $P_i\in\operatorname{Sub}(M)$. Es bleibt noch zu zeigen: $P_i$ ist projektiv. Dazu genügt es, $\Ext_\mathcal{A}^1(P_i,E_j)=0$ für alle $j$ zu zeigen. 

Sei zunächst $j\geq i$. Dann induziert $\eta$ eine exakte Sequenz
\[\Ext_\mathcal{A}^1(E_i,E_j)\longrightarrow\Ext_\mathcal{A}^1(P_i,E_j)\longrightarrow\Ext_\mathcal{A}^1(U,E_j).\]
Wegen $j\geq i$ ist $\Ext_\mathcal{A}^1(E_i,E_j)=0$. Da $U$ projektiv ist, ist $\Ext_\mathcal{A}^1(U,E_j)=0$, und es folgt $\Ext_\mathcal{A}^1(P_i,E_j)=0$.

Sei jetzt $j<i$. Anwenden von $\Hom_{\mathcal{A}}(-,E_j)$ auf das kommutative Diagramm liefert ein kommutatives Diagramm von $k$-Vektorräumen:
\[
\begin{CD}
\Hom_{\mathcal{A}}(U,E_j) @>\partial_{\:\eta}>> \Ext_\mathcal{A}^1(E_i,E_j) @>>> \Ext_\mathcal{A}^1(P_i,E_j) @>>> \Ext_\mathcal{A}^1(U,E_j)\\
@A\pi^* AA @| @AAA @AAA @. \\
\Hom_{\mathcal{A}}(T,E_j) @>\partial_{\:\xi}>>\Ext_\mathcal{A}^1(E_i,E_j) @>>> \Ext_\mathcal{A}^1(\widetilde{E}_i,E_j) @>>>\Ext_\mathcal{A}^1(T,E_j).
\end{CD}
\]
$\partial_{\:\xi}$ ist surjektiv, da $\xi$ universell ist. Dann muß auch $\partial_{\:\eta}$ surjektiv sein. Da $U$ projektiv ist, folgt $\Ext_\mathcal{A}^1(U,E_j)=0$ und die Exaktheit der ersten Zeile liefert $\Ext_\mathcal{A}^1(P_i,E_j)=0$.

$P:=P_1\oplus\ldots\oplus P_n$ ist der gesuchte Progenerator: Ist nämlich $M$ ein Objekt in $\operatorname{Sub}(M)$, so existiert ein wesentlicher Epimorphismus $\pi:M\longrightarrow\Dec M$. $\Dec M$ ist halbeinfach, also existiert nach dem bereits gezeigten ein Epimorphismus $\sigma:P^r\longrightarrow\Dec M$. Da $P$ projektiv ist, faktorisiert $\sigma$ über $\pi$, und da $\pi$ wesentlich ist, ist $\sigma$ surjektiv. Das beendet den Beweis.
\end{bew}

Wir erhalten die angekündigte Verallgemeinerung der Schofield-Induction auf erbliche Kategorien:
\begin{satz}\label{scho:erb:satz2}
Sei $\mathcal{A}$ eine endlichdimensionale, erbliche Längenkategorie. Sei $M$ ein exzeptionelles, nicht-einfaches Objekt. Dann existieren exzeptionelle Objekte $X$ und $Y$ mit folgenden Eigenschaften:
\begin{enumerate}[i)]
\item $\Hom_{\mathcal{A}}(Y,X)=0=\Hom_{\mathcal{A}}(X,Y)$.
\item Es existieren natürliche Zahlen $u\geq1$ und $v\geq1$ und eine exakte Sequenz
\[0\longrightarrow Y^v\longrightarrow M \longrightarrow X^u\longrightarrow 0.\]
\end{enumerate}
Ist $r:=\# Tr(M)$, so existieren genau $r-1$ Isomorphieklassen von Paaren $(Y,X)$ von Objekten in $\mathcal{A}$ mit diesen Eigenschaften.
\end{satz}
\begin{bew}
$\operatorname{Sub}(M)$ ist abgeschlossen unter Faktor- und Unterobjekten sowie unter Erweiterungen. Wir dürfen also ohne Einschränkung $\mathcal{A}=\operatorname{Sub}(M)$ annehmen. Nach \cref{scho:erb:satz1} ist dann $\mathcal{A}$ die Modulkategorie einer endlichdimensionalen Köcheralgebra $A$ mit $r$ Einfachen und $M$ ist ein aufrichtiger $A$-Modul. Mit \cref{scho:scho:satz1} folgt die Behauptung.
\end{bew}

Wir beenden dieses Kapitel mit zwei Aussagen über exzeptionelle Darstellungen von nicht notwendig zykelfreien Köchern. Sei dazu $K$ ein beliebiger endlicher Köcher.
\begin{lemma}\label{scho:erb:lemma4}
Sei $E$ eine einfache $K$-Darstellung mit $\Ext^1_K(E,E)=0$. Dann ist $E=E_x$ für ein $x\in K_0$.
\end{lemma}
\begin{bew} Wir verwenden die in Abschnitt \ref{grund:var} eingeführten Varietäten von  $K$-Darstellungen. Der Orbit $\mathcal{O}(E)$ von $E$ ist nach \cref{grund:var:prop2.5} abgeschlossen in $\operatorname{Rep}_K(\dimv E)$. Wegen $\Ext^1_K(E,E)=0$ ist $\mathcal{O}(E)$ auch offen. Man erhält $\operatorname{Rep}_K(\dimv E)=\mathcal{O}(E)$, aber das ist nur möglich, falls $E=E_x$ für ein $x\in K_0$ ist.
\end{bew}
Wir erhalten das folgende Analogon zu \cref{scho:erb:lemma3} (vgl. \cite[Teil 3, Proposition 2]{Ri3}):
\begin{prop}\label{scho:erb:prop2}
Sei $M$ eine $K$-Darstellung mit $\Ext^1_K(M,M)=0$. Dann ist $Tr(M)$ zykellos.
\end{prop}
\begin{bew}
Seien $E_1,\ldots,E_n$ die Kompositionsfaktoren von $M$. Für alle $i$ sei $\mu(i)$ die Vielfachheit von $E_i$ als Kompositionsfaktor von $M$. Dann sind die $E_i$ die einfachen Objekte in $\operatorname{Sub}(M)$. Nach \cref{scho:erb:satz1} existieren ein endlicher, zykelloser Köcher $L$ und eine Äquivalenz $\operatorname{Sub}(M)\simeq\fmod{kL}$. Dabei gilt: Die Punkte von $L$ sind gegeben durch die $E_i$ und für alle $i,j$ ist $\dim\Ext^1_{\operatorname{Sub}(M)}(E_i,E_j)= \dim\Ext^1_{A}(E_i,E_j)$ die Anzahl der Pfeile zwischen $E_i$ und $E_j$ in $L$.
Alle einfachen $L$-Darstellungen sind exzeptionell, also sind auch die $E_i$ exzeptionell. Wegen \cref{scho:erb:lemma4} existieren also Elemente $x_1,\ldots,x_n\in K_0$ mit $E_i=E_{x_i}$ für alle $i$. Es ist $\dimv M=\sum_{i=1}^n\mu(i)\dimv E_i$, also gilt
\[Tr(M)=\bigl\{x_1,\ldots,x_n\bigr\}.\]
$\dim\Ext^1_{A}(E_{x_i},E_{x_j})$ ist die Anzahl der Pfeil zwischen $x_i$ und $x_j$ in $K$ (wie man leicht unter Verwendung der exakten Sequenz aus \cref{grund:ext:prop1} nachrechnet), also existiert ein Köcher\-isomorphismus $Tr(M)\simeq L$ und die Behauptung folgt.
\end{bew}

\chapter{Ein Funktor zur Produktion von Unzerlegbaren}\label{unz}
\setcounter{zaehler}{0}
\numberwithin{zaehler}{chapter}
Sei $k$ ein algebraisch abgeschlossener Körper und sei $A$ eine (nicht notwendig endlichdimensionale) $k$-Algebra. Alle betrachteten $A$-Moduln seien endlichdimensional. Wir geben einen Funktor an, mit dessen Hilfe man unter gewissen (starken) Voraussetzungen aus vorgegebenen unzerlegbaren $A$-Moduln neue Unzerlegbare als  Erweiterungen konstruieren kann.

Wir fixieren für das gesamte Kapitel unzerlegbare $A$-Moduln $Y_1,\ldots,Y_s$ und $X_1,\ldots,X_r$ und setzen $Y:=\bigoplus_{j=1}^sY_j$ und $X:=\bigoplus_{i=1}^rX_i$.
Dabei seien folgende Bedingungen erfüllt:
\[X_i\not\simeq X_j, Y_k\not\simeq Y_l\text{ für }i\neq j,k\neq l;\ \ \Hom_A(X,Y)=0.\]
Wir werden einen bipartiten Köcher $K$ und einen $k$-linearen Funktor $F:\fdar{K}\longrightarrow\fmod{A}$ mit folgenden Eigenschaften definieren: $F$ bildet unzerlegbare $K$-Darstellungen auf unzerlegbare $A$-Moduln ab und induziert eine injektive Abbildung zwischen der Menge der Isomorphieklassen von $K$-Darstellungen und der Menge der Isomorphieklassen von $A$-Moduln. Außerdem ist $F$ exakt und treu.

\section{Notationen}\label{unz:not}
\setcounter{zaehler}{0}
\numberwithin{zaehler}{section}
Sei $U=U_1\oplus\ldots\oplus U_n\in\fmod{A}$ mit paarweise nicht-isomorphen unzerlegbaren $A$-Moduln $U_i$. Es sei $\mathcal{C}_U:=\prod_{i=1}^n\fmod{k}$
die Produktkategorie. Fasse $\mathcal{C}_U$ als $k$-lineare Kategorie auf (die $k$-Vektorraumstrukturen auf den Homomorphismenräumen sind komponentenweise definiert). Definiere einen $k$-linearen Funktor $F_U:\mathcal{C}_U\longrightarrow\fmod{A}$ wie folgt: Sei $M=(M_1,\ldots,M_n)$ ein $n$-Tupel von $k$-Vektorräumen. Setze
\[FM:=\bigoplus_{i=1}^nU_i\otimes_k M_i.\]
Ist $f=(f_1,\ldots,f_n):M\longrightarrow M'$ ein Homomorphismus in $\mathcal{C}_U$, so sei
\[Ff:=\begin{bmatrix}id_{U_1}\otimes f_1&&\\&\ddots&\\&&id_{U_n}\otimes f_n\end{bmatrix}\in\Hom_A(FM,FM').\]\newpage

Es gilt:
\begin{itemize}
\item[$\cdot$]$F_U$ ist treu und exakt.
\item[$\cdot$]$F_U$ induziert eine Äquivalenz $\mathcal{C}_U\stackrel{\sim}{\longrightarrow}\operatorname{add}(U)/\mathcal{R}_{\operatorname{add}(U)}$.
\end{itemize}
Wir benötigen noch einige Notationen für Bifunktoren: Seien $\mathcal{A}$, $\mathcal{B}$ $k$-lineare Kategorien und sei $F:\mathcal{A}\times\mathcal{B}\longrightarrow \fmod{k}$ ein ($k$-bilinearer) Bifunktor, kontravariant in der ersten, kovariant in der zweiten Variablen. Für Homomorphismen $f:M'\longrightarrow M$ und $g:N\longrightarrow N'$ sowie für Elemente $x\in F(M,N)$ schreiben wir häufig
\[gx:=F(id_M,g)(x)\in F(M,N')\text{ und } xf:=F(f,id_M)(x)\in F(M',N).\]
Das Radikal $\operatorname{Rad}F\subseteq F$ von $F$ ist definiert durch
\[\Rad F(M,N):=\sum_{T\in\mathcal{B}}\mathcal{R}_{\mathcal{B}}(T,N)\:F(M,T)+\sum_{S\in\mathcal{A}}F(S,N)\:\mathcal{R}_{\mathcal{A}}(M,S)\]
für alle $M\in\mathcal{A},N\in\mathcal{B}$. Dabei  bezeichnet $\mathcal{R}_{\mathcal{B}}(T,N)\:F(M,T)$ den von allen $f\zeta$, $\zeta\in F(M,T)$, $f\in\mathcal{R}_{\mathcal{B}}(T,N)$ erzeugten $k$-Untervektorraum von $F(M,N)$; analog ist $F(S,N)\:\mathcal{R}_{\mathcal{A}}(M,S)$ definiert.

Ist spezieller ein Bifunktor $F:\operatorname{add}(Y)\times\operatorname{add}(X)\longrightarrow\fmod{k}$ gegeben, so gilt
\[\Rad F(M,N):=\mathcal{R}_{\operatorname{add}(X)}(X,N)\:F(M,X)+F(Y,N)\:\mathcal{R}_{\operatorname{add}(Y)}(M,Y).\]
$\Rad F(Y,X)$ ist das Radikal von $F(Y,X)$ als $\End_A(X)$-$\End_A(Y)$-Bimodul.

\section{Konstruktion des Funktors $F:\fdar{K}\longrightarrow\fmod{A}$}\label{unz:konst}
\setcounter{zaehler}{0}
\numberwithin{zaehler}{section}
Betrachte die in Abschnitt \ref{grund:ext} zur Beschreibung von $\Ext^1_A$ eingeführten Bifunktoren
\[Z^1_A,B^1_A:\fmod{A}\times\fmod{A}\longrightarrow\fmod{k}.\]
Sei $\iota:\operatorname{add}(Y)\times\operatorname{add}(X)\longrightarrow\fmod{A}\times\fmod{A}$ der Vergißfunktor. Wir schreiben zur Abkürzung $Z:=Z^1_A\circ\iota$, $B:=B^1_A\circ\iota$ und $\Ext:=\Ext^1_A\circ\iota$. Setze
\[\overline{Z}:=Z/\bigl(\Rad Z+B\bigr)\text{ und } \overline{\Ext}:=\Ext/\Rad \Ext.\]
Die natürliche Transformation $\phi:Z\longrightarrow \Ext$ aus Abschnitt \ref{grund:ext} induziert einen Isomorphismus $\overline{\phi}:\overline{Z}\stackrel{\sim}{\longrightarrow}\overline{\Ext}$.
Für alle $i=1,\ldots,r$ und für alle $j=1,\ldots,s$ wähle eine Teilmenge \[\mathcal{B}_{ij}=\bigl\{\zeta_1^{(ij)},\ldots,\zeta_{n_{ij}}^{(ij)}\bigr\}\]
von $Z(Y_j,X_i)$, so daß die Restklassen eine Basis von $\overline{Z}(Y_j,X_i)$ bilden.

\begin{defn} Wir definieren den bipartiten Köcher $K$ durch
\begin{itemize}
\item[$\cdot$] $K_0:=\{Y_1,\ldots,Y_s;X_1,\ldots,X_r\}$,
\item[$\cdot$] $K_1:=\coprod\limits_{\begin{subarray}{c}1\leq i\leq r\\1\leq j\leq s\end{subarray}}\mathcal{B}_{ij}$,
\item[$\cdot$] Für $\alpha\in\mathcal{B}_{ij}\subseteq K_1$: $n(\alpha)=Y_j$ und $s(\alpha)=X_i$.
\end{itemize}
\end{defn}
\begin{bem}\label{unz:konst:bem1}\ 
Sei $R$ das Radikal von $\Ext^1_A(Y,X)$ als $\End_A(X)$-$\End_A(Y)$-Bimodul. Seien $1\leq i\leq r$ und $1\leq j\leq s$ gegeben. Betrachte die von der Projektion $\pi_j:Y\longrightarrow Y_j$ und von der Einbettung $\epsilon_i:X_i\longrightarrow X$ induzierte Abbildung \[\tau_{ij}:\Ext^1_A(Y_j,X_i)\longrightarrow\Ext^1_A(Y,X)/R.\]
Sei $R_{ij}:=\Ker\tau_{ij}$. Dann ist $R_{ij}=\bigl(\Rad\Ext\bigr)(Y_j,X_i)$ und die kanonische Abbildung \[Z^1_A(Y_j,X_i)\longrightarrow\Ext^1_A(Y_j,X_i)\] induziert einen $k$-linearen Isomorphismus \[\overline{Z}(Y_j,X_i)\stackrel{\sim}{\longrightarrow}\Ext^1_A(Y_j,X_i)/R_{ij}.\]
\end{bem}
Wir legen weitere Notationen fest. Für eine $K$-Darstellung $M$ setze \[M_Y:=\bigl(M(Y_1),\ldots,M(Y_s)\bigr)\in\mathcal{C}_Y\ \text{ und }\  M_X:=\bigl(M(X_1),\ldots,M(X_r)\bigr)\in\mathcal{C}_X.\]
Ein Homomorphismus $f:M\longrightarrow M'$ von $K$-Darstellungen induziert Morphismen
\[f_X:=(f_{X_1},\ldots,f_{X_r}):M_X\longrightarrow M'_X\text{ und }f_Y:=(f_{Y_1},\ldots,f_{Y_s}):M_Y\longrightarrow M'_Y.\]

Umgekehrt definiere für $V\in\mathcal{C}_X$ eine $K$-Darstellung $V_K$ durch
\[
V_K(Y_j)=0 \text{ für alle } j=1\ldots,s, V_K(X_i)=V_i \text{ für alle } i=1\ldots,r.
\]
Ein Morphismus $f:V\longrightarrow V'$ in $\mathcal{C}_X$ induziert in offensichtlicher Weise einen Homomorphismus von Darstellungen, den wir ebenfalls mit $f$ bezeichnen. Analog liefert $W\in\mathcal{C}_Y$ eine $K$-Darstellung $W_K$.

Betrachte den Bifunktor $C^1_K:\fdar{K}\times\fdar{K}\longrightarrow \fmod{k}$
(vgl. Abschnitt (\ref{grund:ext})). Wir wollen im folgenden die Beziehung zwischen $C^1_K$ und $Z$ näher untersuchen.

\begin{defn}\label{unz:konst:def2}\ 
Für alle $W\in\mathcal{C}_Y$ und $V\in\mathcal{C}_X$ sei eine $k$-lineare Abbildung
\[\Theta_{W,V}:C^1_K(W_K,V_K)\longrightarrow Z(F_YW,F_XV)\]
definiert durch
\[\Theta_{W,V}\bigl((\zeta(\alpha))_{\alpha\in K_1}\bigr):=\Biggl[\ \sum_{k=1}^{n_{ij}}\alpha\otimes\zeta_{ij}(\alpha)\ \Biggr]_{\begin{subarray}{c}1\leq i\leq r\\1\leq j\leq s\end{subarray}}.\]
Beachte dabei: Die Projektionen $F_XV\longrightarrow X_i\otimes_k V_i$ und Einbettungen $Y_j\otimes_k W_j\longrightarrow F_YW$ induzieren einen kanonischen Isomorphismus
\[
Z(F_YW,F_XV)\stackrel{\sim}{\longrightarrow}\bigoplus\limits_{i,j}Z(Y_j\otimes_k W_j,X_i\otimes_k V_i),\ 
\zeta\mapsto\Bigl[\zeta_{ij}\Bigr]_{\begin{subarray}{c}1\leq i\leq r\\1\leq j\leq s\end{subarray}}.
\]
Außerdem vereinbaren wir: Sind $M,N$ $A$-Moduln, $W,V$ $k$-Vektorräume und sind $\zeta\in Z^1_A(N,M)$ und $f\in\Hom_k(W,V)$, so sei $\zeta\otimes f\in Z^1_A(N\otimes W,M\otimes V)$ definiert durch
\[
\bigl(\zeta\otimes f\bigr)(a):=\zeta(a)\otimes f\in\Hom_k(N\otimes W,M\otimes V)\text{ für alle } a\in A.
\]
\end{defn}

\begin{lemma}\label{unz:konst:lemma1}\ 
\begin{enumerate}
\item $\Theta_{(W,V)}$ verhält sich funktoriell in $W$ und $V$.
\item $\Theta_{(W,V)}$ induziert einen Isomorphismus
\[
C^1_K(W_K,V_K)\stackrel{\sim}{\longrightarrow} \overline{Z}(F_YW,F_XV).
\]
\end{enumerate}
\end{lemma}
\begin{bew}Die Funktorialität in (1.) rechnet man leicht nach. 

(2.) Für $\lambda=1,\ldots,r$ sei
\[
E_\lambda^{(X)}\in\mathcal{C}_X,\ \  E_\lambda^{(X)}(i)=\left\{\begin{array}{cl}k&\ \ \text{falls } \lambda=i\\ 0&\ \ \text{sonst}\end{array}\right..
\]
Analog sei für $\mu=1,\ldots,s$
\[
E_\mu^{(Y)}\in\mathcal{C}_Y,\ \  E_\mu^{(Y)}(j)=\left\{\begin{array}{cl}k&\ \ \text{falls } \mu=j\\ 0&\ \ \text{sonst}\end{array}\right..
\]
Seien $\lambda,\mu$ gegeben. Setze $W:=E_\mu^{(Y)}$ und $V:=E_\lambda^{(X)}$. In dem Fall ist
\[C^1_k(W,V)=\bigoplus_{k=1}^{n_{\lambda,\mu}}\Hom_k(k,k)\simeq k^{n_{\lambda,\mu}}\]
und
\[\overline{Z}(F_YW,F_XV)=\overline{Z}(Y_\mu\otimes_kk,X_\lambda\otimes_kk)\simeq \overline{Z}(Y_\mu,X_\lambda).\]
Via dieser Identifizierungen ist $\Theta_{W,V}$ gegeben durch die $k$-lineare Abbildung
\[k^{n_{\lambda,\mu}}\longrightarrow \overline{Z}(Y_\mu,X_\lambda), e_k\mapsto \overline{\zeta_k^{(\lambda\mu)}}\text{ für alle }k=1,\ldots,n_{\lambda,\mu}.\]
Diese Abbildung ist per Konstruktion der $\zeta_k^{(ij)}$ ein Isomorphismus.

$\Theta_{(-,-)}$ ist nach (1.) eine natürliche Transformation von biadditiven Bifunktoren
\[\mathcal{C}_Y\times\mathcal{C}_X\longrightarrow\fmod{k}.\]
Außerdem wird $\mathcal{C}_X$ als additive Kategorie von den $E_\lambda^{(X)}$ und $\mathcal{C}_Y$ von den $E_\mu^{(Y)}$ erzeugt. Also ist $\Theta_{(W,V)}$ für beliebige $W,V$ ein Isomorphismus.
\end{bew}
\newpage
Wir können nun den Funktor $F$ definieren.
\begin{defn}[Definition von $F$]\label{unz:konst:def3}\ 

Sei $M\in\fdar{K}$. Setze
\[\zeta^{(M)}:=\bigl(M(\alpha)\bigr)_{\alpha\in K_1}\in C^1_K(M_Y,M_X)\]
und
\[\zeta^{(FM)}:=\Theta_{M_Y,M_X}(\zeta^{(M)})\in Z(F_YM_Y,F_XM_X).\]
Es sei $FM$ der durch $\zeta^{(FM)}$ gegebene $A$-Modul, das heißt: Als $k$-Vektorraum ist 
\[FM:=F_XM_X\oplus F_YM_Y.\]
Für alle $a\in A$ ist die Linksmultiplikation mit $a$ gegeben durch
\[
FM(a):=\begin{bmatrix}F_XM_X(a)&\zeta^{(FM)}(a)\\0&F_YM_Y(a)\end{bmatrix}.
\]
Für einen Homomorphismus $f:M\longrightarrow M'$ von $K$-Darstellungen setze
\[
Ff:=\begin{bmatrix}F_X(f_X)&0\\0&F_Y(f_Y)\end{bmatrix}\in\Hom_k(FM,FM').
\]
Man verifiziert leicht mit Hilfe der Funktorialität von $\Theta_{(-,-)}$ in beiden Variablen, daß $Ff$ $A$-linear ist.

Schließlich sei 
\[(E_{FM})\ :0\longrightarrow F_XM_X\longrightarrow FM\longrightarrow F_YM_Y\longrightarrow0\]
die zugehörige exakte Sequenz von $A$-Moduln.
\end{defn}

Ist $A=kL$ eine Köcheralgebra eines endlichen Köchers $L$, so kann man den Bifunktor $Z=Z^1_A$ in der Definition durch den Bifunktor $C^1_L$ zur Beschreibung von Ext-Gruppen von Köcherdarstellungen ersetzen. Wir führen diese Konstruktion nur in einem Spezialfall aus, der in Kapitel \ref{baum} von Bedeutung ist. Sei dazu $r=s=1$, also $X=X_1$ und $Y=Y_1$. Sei $R$ das Radikal von $\Ext^1_A(Y,X)$ als $\End_A(X)$-$\End_A(Y)$-Bimodul. Betrachte die in \cref{grund:ext:prop2} definierte $k$-lineare Abbildung $\phi:C^1_L(Y,X)\longrightarrow\Ext^1_L(Y,X)$. Es seien Elemente $\zeta_1,\ldots,\zeta_n\in C^1_L(Y,X)$ gegeben, so daß die Restklassen der $\phi(\zeta_i)$ eine Basis von $\Ext^1_A(Y,X)/R$ bilden. Wir versehen $Q_n$ mit den Bezeichnungen 
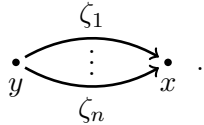
\begin{center}
\parbox[c]{2.6cm}{
\begin{tikzpicture}[line width=1pt]
	\tikzstyle{knt}=[circle, fill, inner sep=1pt]
  \node at (0,0)[knt, label=below:$y$](x){};
  \node at (2,0)[knt, label=below:$x$](y){};
  \node at (1,0.1)(z){$\vdots$};
  \draw[->, bend left, shorten >= 2pt, shorten <= 2pt] (x) to node[above]{$\zeta_1$} (y);
  \draw[->, bend right, shorten >= 2pt, shorten <= 2pt] (x) to node[below]{$\zeta_n$} (y);
\end{tikzpicture}}.
\end{center}
Wir definieren dann einen Funktor $F:\fdar{Q_n}\longrightarrow\fdar{L}$ wie folgt:
\begin{defn}\label{unz:konst:def4}
Sei $M\in\fdar{Q_n}$. Für $t\in L_0$ setze
\[FM(t):=\bigl(X(t)\otimes_kM(x)\bigr)\oplus \bigl(Y(t)\otimes_kM(y)\bigr).\] 
Für alle $\alpha\in L_1$ definiere
\[FM(\alpha):=\begin{bmatrix}X(\alpha)\otimes id_{M(x)}&\sum\limits_{i=1}^n(\zeta_i)_\alpha\otimes M(\zeta_i)\\0&Y(\alpha)\otimes id_{M(y)}\end{bmatrix}\in\Hom_k(FM(n(\alpha)),FM(s(\alpha)).\]
Ist $f=(f_y,f_x):M\longrightarrow M'$ ein Homomorphismus in $\fdar{Q_n}$, so sei für $t\in K_0$
\[(Ff)_t:=\begin{bmatrix}id_{X(t)}\otimes f_x&0\\0&id_{Y(t)}\otimes f_y\end{bmatrix}\in\Hom_k(FM(i),FM'(i)).\]
Dann ist $Ff$ ein Homomorphismus von $L$-Darstellungen.
\end{defn}

Sei $H:\fdar{L}\longrightarrow\fmod{kL}$ die kanonische Äquivalenz. Man verifiziert, daß ein Funktor 
\[G:\fdar{Q_n}\longrightarrow\fmod{kL}\]
im Sinne von \cref{unz:konst:def3} existiert, so daß $G\simeq H\circ F$ ist.
\section{Die Eigenschaften von $F$}\label{unz:eig}
\setcounter{zaehler}{0}
\numberwithin{zaehler}{section}
Wir fassen die einfachsten Eigenschaften in folgendem Lemma zusammen:
\begin{lemma}\label{unz:eig:lemma1}\ 
\begin{enumerate}
\item Für alle $i,j$ gilt: $FE_{X_i}\simeq X_i$ und $FE_{Y_j}\simeq Y_j$.
\item $F$ ist treu und exakt.
\item Sei $M\in\fdar{K}$. Betrachte den in Abschnitt (\ref{grund:ext:prop1}) definierten Homomorphismus \[\phi:Z(F_YM_Y,F_XM_X)\longrightarrow\Ext^1_A(F_YM_Y,F_XM_X).\]
Dann gilt $\phi\circ\Theta_{M_Y,M_X}(\zeta^{(M)})=E_{FM}$.
\end{enumerate}
\end{lemma}
\begin{bew}(1.) folgt sofort aus der Definition von $F$; (2.) folgt aus der Treue und Exaktheit der Funktoren $F_X$ und $F_Y$; (3.) folgt aus der Definition von $\phi$, $\Theta$ und $F$.
\end{bew}
Das Hauptresultat dieses Kapitels ist die folgende Proposition:
\begin{prop}\label{unz:eig:prop1}\ 
\begin{enumerate}
\item Ist $M$ eine unzerlegbare $K$-Darstellung, so ist $FM$ ein unzerlegbarer $A$-Modul.
\item Seien $M,M'$ $K$-Darstellungen. Gilt $FM\simeq FM'$, so schon $M\simeq M'$.
\end{enumerate}
\end{prop}

Zum Beweis benötigen wir zwei Lemmata.
\begin{lemma}\label{unz:eig:lemma2}Sei
\[
\begin{CD}
0@>>>X@>>>M@>>>Y@>>>0\\
@.@Vf'VV@VfVV@VVf''V@.\\
0@>>>X'@>>>M'@>>>Y'@>>>0\\
\end{CD}
\]
ein kommutatives Diagramm mit exakten Zeilen. Es sei $\Hom_A(X,Y)=0=\Hom_A(X',Y)$. Dann gilt: Ist $f'\in\mathcal{R}_A(X,X')$ und $f''\in\mathcal{R}_A(Y,Y')$, so ist $f\in\mathcal{R}_A(M,M')$.
\end{lemma}
\begin{bew}
Wir dürfen annehmen: Als $k$-Vektorraum ist $M=X\oplus Y$, $M'=X'\oplus Y'$ und für alle $a\in A$ gilt
\[M(a)=\begin{bmatrix}X(a)&\zeta(a)\\0&Y(a)\end{bmatrix}\text{ und } M'(a)=\begin{bmatrix} X'(a)&\zeta'(a)\\0&Y'(a)\end{bmatrix}\]
für gewisse Derivationen $\zeta\in Z^1(Y,X)$ und $\zeta'\in Z^1(Y',X')$ (vgl. Abschnitt \ref{grund:ext}). Außerdem dürfen wir wegen der Kommutativität des gegebenen Diagramms annehmen, daß $f$ von der Form
\[f=\begin{bmatrix}f'&r\\0&f''\end{bmatrix}\]
ist für eine gewisse $k$-lineare Abbildung $r:Y'\longrightarrow X$.

Wegen $\Hom_A(X,Y)=0$ gilt
\[\End_A(M)\subseteq\Bigl\{\begin{bmatrix}\phi'&r\\0&\phi''\end{bmatrix}\in\End_k(M)\Bigm|\phi'\in\End_A(X), \phi''\in\End_A(Y),r\in\Hom_k(Y,X)\Bigr\}.\]
Also ist
\[I=\Bigl\{\begin{bmatrix}\phi'&r\\0&\phi''\end{bmatrix}\in\End_A(M)\Bigm|\phi'\in J\End_A(X),\phi''\in J\End_A(Y)\Bigr\}\]
ein zweiseitiges Ideal in $\End_A(M)$. Offenbar ist $I$ nilpotent, also gilt $I\subseteq J\End_A(M)$. 

Sei jetzt $g\in\Hom_A(M',M)$ beliebig. Wegen $\Hom_A(X',Y)=0$ gilt:
\[g=\begin{bmatrix}g'&s\\0&g''\end{bmatrix}, s\in\Hom_k(Y',X), g'\in\Hom_A(X',X), g''\in\Hom_A(Y',Y)\]
Da $f'$ und $f''$ in $\mathcal{R}_A$ liegen, folgt $g\circ f\in I\subseteq J\End_A(M)$. Das zeigt $f\in\mathcal{R}_A(M,M')$.
\end{bew}

\begin{lemma}\label{unz:eig:lemma3}$F$ induziert einen vollen Funktor $\fdar{K}\longrightarrow\fmod{A}/\mathcal{R}_A$. Zusatz: Gilt außerdem
\[
\begin{gathered}
\dim\Hom_A(X_i,X_j)=\delta_{ij} \text{ und } \dim\Hom_A(Y_k,Y_l)=\delta_{kl}\text{ für alle } i,j;k,l,\\
\Hom_A(Y,X)=0=\Hom_A(X,Y),
\end{gathered}
\]
so ist schon $F$ voll. 
\end{lemma}
\begin{bew}
Seien $M,M'\in\fdar{K}$ und sei $f:M\longrightarrow M'$ ein Homomorphismus von $A$-Moduln. Wir  müssen zeigen: Es existiert ein Homomorphismus $\phi\in\Hom_{K}(M,M')$ mit $f-F\phi\in\mathcal{R}_A(FM,FM')$. Wegen $\Hom_A(X,Y)=0$ erhalten wir ein kommutatives Diagramm
\[
\begin{CD}
E_{FM}:\: 0@>>>F_XM_X@>>>FM@>>>F_YM_Y@>>>0\\
@.@Vf^{(X)}VV@VfVV@VVf^{(Y)}V@.\\
E_{FM'}: 0@>>>F_XM'_X@>>>FM'@>>>F_YM'_Y@>>>0
\end{CD}\ \ \ \ \text{ (D) }
\]
mit exakten Zeilen. Demnach ist $f$ von der Form
\[
f=\begin{bmatrix}f^{(X)}&r\\0&f^{(Y)}\end{bmatrix}
\]
für eine $k$-lineare Abbildung $r\in\Hom_k(F_YM_Y,F_XM'_X)$.

Weil der Funktor $F_X$ eine Äquivalenz $\mathcal{C}_X\longrightarrow\operatorname{add}(X)/\mathcal{R}_{\operatorname{add}(X)}$ induziert, existieren ein Homomorphismus $\phi_{X}:M_X\longrightarrow M'_X$ und eine $A$-lineare Abbildung $\widetilde{f}^{\ (X)}\in\mathcal{R}_A(F_XM_X,F_XM'_X)$ mit
\counterwithout{equation}{chapter}
\setcounter{equation}{0}
\begin{equation}\label{unz:eig:lemma3:Gl1}
f^{(X)}=F_X(\phi_X)+\widetilde{f}^{\ (X)}.
\end{equation}
Analog findet man Homomorphismen $\phi_Y:M_Y\longrightarrow M'_Y$ und $\widetilde{f}^{\ (Y)}\in\mathcal{R}_A(F_YM_Y,F_YM'_Y)$, so daß gilt:
\begin{equation}\label{unz:eig:lemma3:Gl2}
f^{(Y)}=F_Y(\phi_Y)+\widetilde{f}^{\ (Y)}.
\end{equation}

Wir zeigen: $\phi_X$ und $\phi_Y$ definieren den gesuchten Homomorphismus $\phi:M\longrightarrow M'$.

Sei dazu $d:\Hom_k(F_YM_Y,F_XM'_X)\longrightarrow Z(F_YM_Y,F_XM'_X)$ definiert durch
\[d(f)(a):=F_XM'_X(a)\circ f-f\circ F_YM_Y(a)\text{ für alle } a\in A.\]
Es ist $\Bild d=B(F_YM_Y,F_XM'_X)$ (vgl. Abschnitt (\ref{grund:ext})). Die $A$-Linearität von $f$ liefert die Formel
\begin{equation}\label{unz:eig:lemma3:Gl3}
f^{(X)}\zeta^{(FM)}-\zeta^{(FM')}f^{(Y)}=d(r).
\end{equation}
Da $\widetilde{f}^{\ (X)}$ und $\widetilde{f}^{\ (Y)}$ im Radikal von $\operatorname{add}{X}$ bzw. $\operatorname{add}{Y}$ liegen, folgt daraus mit Hilfe der Gleichungen (\ref{unz:eig:lemma3:Gl1}) und (\ref{unz:eig:lemma3:Gl2}):
\[
F_X(\phi_X)\zeta^{(FM)}-\zeta^{(FM')}F_Y(\phi_Y)\in B(F_YM_Y,F_XM'_X)+\Rad Z(F_YM_Y,F_XM'_X)=:R.
\]
Betrachte die natürliche Transformation $\Theta$ aus \cref{unz:konst:def2}. Nach \cref{unz:konst:lemma1} induziert $\Theta_{M_Y,M'_X}$ einen Isomorphismus
\[\overline{\Theta}_{M_Y,M'_X}:C^1_K(M_Y,M'_X)\stackrel{\sim}{\longrightarrow} \overline{Z}(F_YM_Y,F_XM'_X)=Z(F_YM_Y,F_XM'_X)/R.\]
Wegen der Funktorialität von $\Theta$ gilt
\[
\begin{aligned}
&\Theta_{M_Y,M'_X}(\phi_X\zeta^{(M)}-\zeta^{(M')}\phi_Y)\\
=\ &F_X(\phi_X)\Theta_{M_Y,M_X}(\zeta^{(M)})-\Theta_{M'_Y,M'_X}(\zeta^{(M')})F_Y(\phi_Y)\\
=\ &F_X(\phi_X)\zeta^{(FM)}-\zeta^{(FM')}F_Y(\phi_Y)\\
\in\ & R.
\end{aligned}
\]
Es folgt $\overline{\Theta}_{M_Y,M'_X}(\phi_X\zeta^{(M)}-\zeta^{(M')}\phi_Y)=0$. Also ist schon $\phi_X\zeta^{(M)}-\zeta^{(M')}\phi_Y=0$. Das zeigt: $\phi$ ist ein Homomorphismus von $K$-Darstellungen.

Es bleibt, zu zeigen: $f-F\phi\in\mathcal{R}_A(FM,FM')$. Das Diagramm (D) liefert ein kommutatives Diagramm
\[
\begin{CD}
E_{FM}:\: 0@>>>F_XM_X@>>>FM@>>>F_YM_Y@>>>0\\
@.@V\widetilde{f}^{\ (X)}VV@Vf-F\phi VV@VV\widetilde{f}^{\ (Y)}V@.\\
E_{FM'}: 0@>>>F_XM'_X@>>>FM'@>>>F_YM'_Y@>>>0
\end{CD}
\]
mit exakten Zeilen. In Matrizenschreibweise:
\[
f-F\phi=\left[\begin{array}{cc}\widetilde{f}^{\ (X)}&r\\0&\widetilde{f}^{\ (Y)}\end{array}\right].
\]
Da $\widetilde{f}^{\ (X)}$ und $\widetilde{f}^{\ (Y)}$ im Radikal von $\fmod{A}$ liegen, folgt die Behauptung mit \cref{unz:eig:lemma2}.

Seien jetzt die Zusatzbedingungen erfüllt. Dann gilt $\mathcal{R}_{\operatorname{add}(X)}=0=\mathcal{R}_{\operatorname{add}(Y)}$. Damit folgt $\widetilde{f}^{\ (X)}=0$ und $\widetilde{f}^{\ (Y)}=0$. Wir wissen außerdem bereits, daß $\phi$ ein Homomorphismus von $K$-Darstellungen ist. Mit den Gleichungen (\ref{unz:eig:lemma3:Gl1}), (\ref{unz:eig:lemma3:Gl2}) und (\ref{unz:eig:lemma3:Gl3}) folgt daraus
\[0=F_X(\phi_X)\zeta^{(FM)}-\zeta^{(FM')}F_Y(\phi_Y)=d(r).\]
Also ist $r\in\Ker d=\Hom_A(Y,X)=0$, also ist $r=0$, und es folgt $f=F\phi$.
\end{bew}

Wir kommen nun zum Beweis von \cref{unz:eig:prop1}.
\begin{bew}[von \cref{unz:eig:prop1}]
(1.) Sei $M$ eine unzerlegbare Darstellung von $K$. Nach \cref{unz:eig:lemma3} induziert $F$ einen Epimorphismus von $k$-Algebren
\[\End_{K}(M)\longrightarrow\End_A(FM)/J\End_A(FM).\]
Offenbar ist $FM\neq0$, also ist $\End_A(FM)/J\End_A(FM)$ als Quotient einer lokalen $k$-Algebra wieder lokal. Außerdem ist $\End_A(FM)/J\End_A(FM)$ halbeinfach. Das zeigt 
\[\End_A(FM)/J\End_A(FM)\simeq k.\]
Also ist $\End_A(FM)$ lokal, und damit ist $FM$ unzerlegbar.

(2.) Man verifiziert leicht mit Hilfe der Definition des Radikals von additiven  Kategorien folgende Aussagen:\enlargethispage{\baselineskip}
\begin{enumerate}
\item Ist $F:\mathcal{A}\longrightarrow\mathcal{B}$ ein voller Funktor zwischen additiven Kategorien, so induziert $F$ einen (vollen) Funktor $F:\mathcal{A}/\mathcal{R}_{\mathcal{A}}\longrightarrow\mathcal{B}/\mathcal{R}_{\mathcal{B}}$.
\item Sei $\mathcal{A}$ eine additive Kategorie. Dann ist $\mathcal{R}_{\mathcal{A}/\mathcal{R}_{\mathcal{A}}}=0$.
\end{enumerate}

Wegen \cref{unz:eig:lemma3} folgt damit: Der Funktor $F:\fdar{K}\longrightarrow \fmod{A}$ induziert einen vollen Funktor $\overline{F}:\fdar{K}/\mathcal{R}_{\fdar{K}}\longrightarrow \fmod{A}/\mathcal{R}_A$.

Seien jetzt $M,M'\in\fdar{K}$ mit $FM\simeq FM'$. 

1. Fall: $M$ und $M'$ sind unzerlegbar. Wegen (2.) sind auch $FM$ und $FM'$ unzerlegbar, also ist $\Hom_A(FM,FM')/\mathcal{R}_A(FM,FM')\simeq k$. Da $\overline{F}$ voll ist, existiert ein Epimorphismus von $k$-Vektorräumen
\[\Hom_{K}(M,M')/\mathcal{R}_{\fdar{K}}(M,M')\longrightarrow \Hom_A(FM,FM')/\mathcal{R}_A(FM,FM').\]
Also ist $\Hom_{K}(M,M')/\mathcal{R}_{\fdar{K}}(M,M')\neq 0$, also $M\simeq M'$.

2. Fall: $M,M'$ beliebig. Es gelte $M\simeq\bigoplus_{i=1}^tU_i^{n_i}\text{ und } M'\simeq\bigoplus_{i=1}^tU_i^{n'_i}$ mit $U_i$ unzerlegbar in $\fdar{K}$ für alle $i$, $U_i\not\simeq U_j$ für $i\neq j$ und $n_i,n'_i\geq 0$ für alle $i$. Erhalte
\[\bigoplus_{i=1}^t(FU_i)^{n_i}\simeq FM\simeq FM'\simeq\bigoplus_{i=1}^t(FU_i)^{n'_i}.\]
Nach (2.) sind alle $FU_i$ unzerlegbare $A$-Moduln; außerdem gilt nach dem oben behandelten ersten Fall $FU_i\not\simeq FU_j$ für $i\neq j$. Nach dem Satz von Krull-Remak-Schmidt folgt $n_i=n'_i$ für alle $i$, also $M\simeq M'$.
\end{bew}

Wir wollen den Funktor $F$ unter gewissen Zusatzbedingungen genauer untersuchen. Sei dazu $\mathcal{F}(Y,X)$ die volle Unterkategorie von $\fmod{A}$ aller $A$-Moduln $M$, für die eine exakte Sequenz der Form
\[0\longrightarrow V\longrightarrow M\longrightarrow W\longrightarrow 0\]
existiert mit $V\in\operatorname{add}(X)$ und $W\in\operatorname{add}(Y)$. Offenbar ist das Bild von $F$ in $\mathcal{F}(Y,X)$ enthalten. Es gilt:
\begin{prop}\label{unz:eig:prop2}
Es seien die zusätzlichen Bedingungen 
\[\dim\Hom_A(X_i,X_j)=\delta_{ij} \text{ und } \dim\Hom_A(Y_k,Y_l)=\delta_{kl}\text{ für alle } i,j;k,l\]
erfüllt. Dann gilt:
\begin{enumerate}
\item Ist $N\in\mathcal{F}(Y,X)$, so existiert eine $K$-Darstellung $M$ mit $FM\simeq N$.
\item Gilt auch noch $\Hom_A(Y,X)=0$, so ist $\mathcal{F}(Y,X)$ abgeschlossen unter Kernen und Kokernen und $F$ induziert eine Äquivalenz $\fdar{K}\stackrel{\sim}{\longrightarrow}\mathcal{F}(Y,X)$.
\end{enumerate}
\end{prop}
\begin{bew} Sei $N\in\mathcal{F}(Y,X)$ gegeben. Dann existiert eine exakte Sequenz
\[(E):\ 0\longrightarrow F_XV\longrightarrow N \longrightarrow F_YW\longrightarrow 0\]
für gewisse $V\in\mathcal{C}_X$, $W\in\mathcal{C}_Y$. Betrachte die Homomorphismen
\[\Theta:=\Theta_{W,V}:C^1_K(W_K,V_K)\longrightarrow Z(F_YW,F_XV)\text{ und }\phi:Z(F_YW,F_XV)\longrightarrow\Ext^1_A(F_YW,F_XV)\]
(vgl. \cref{unz:konst:def2} und \cref{grund:ext:prop1}). Dann ist $\phi\circ\Theta$ ein Isomorphismus: Wegen der Zusatzbedingungen gilt $\mathcal{R}_{\operatorname{add}(X)}=0=\mathcal{R}_{\operatorname{add}(Y)}$, also ist $\Rad Z=0$ und es folgt $\overline{Z}=Z/B$. Wir wissen wegen \cref{grund:ext:prop1}, daß $\phi$ einen Isomorphismus \[\overline{Z}(F_YW,F_XV)=\overline{Z}(F_YW,F_XV)/B(F_YW,F_XV)\simeq\Ext^1_A(F_YW,F_XV)\]
induziert. Nach \cref{unz:konst:lemma1} induziert $\Theta$ einen Isomorphismus $C^1_K(W_K,V_K)\simeq\overline{Z}(F_YW,F_XV)$. Also ist $\phi\circ\Theta$ ein Isomorphismus.

Wähle jetzt einen Kozykel $\zeta\in C^1_K(W_K,V_K)$ mit $\phi\circ\Theta(\zeta)=E\in\Ext^1_A(F_YW,F_XV)$. Sei $M$ die durch $\zeta$ gegebene Darstellung von $K$, also
\begin{itemize}
\item[$\cdot$]$M(X_i)=V_i$, $M(Y_j)=W_j$ für alle $i,j$.
\item[$\cdot$]$M(\alpha)=\zeta_\alpha$ für alle $\alpha\in K_1$.
\end{itemize}
Dann ist $V=M_X$, $W=M_Y$ und $\zeta=\zeta^{(M)}$. Betrachte die zugehörige kanonische exakte Sequenz
\[(E_{FM}):\ 0\longrightarrow F_XV\longrightarrow FM\longrightarrow F_YW\longrightarrow 0\]
von $A$-Moduln (vgl. \cref{unz:konst:def3}). Es gilt $\phi\circ\overline{\Theta}(\zeta)=E_{FM}$ in $\Ext^1_A(F_YW,F_XV)$  (vgl. \cref{unz:eig:lemma1}, (3.)). Also gilt $E=E_{FM}$ in $\Ext^1_A(F_YW,F_XV)$. Damit folgt $FM\simeq N$.

(2.) $F$ ist treu. Wir haben im Zusatz zu \cref{unz:eig:lemma3} gesehen, daß $F$ unter den gegebenen Voraussetzungen auch voll ist. Wegen (1.) induziert $F$ also die gewünschte Äquivalenz. Da $F$ aufgefaßt als Funktor zwischen $\fdar{K}$ und $\fmod{A}$ exakt ist, folgt, daß $\mathcal{F}(Y,X)$ abgeschlossen ist unter Kernen und Kokernen.
\end{bew}

\chapter{Baummoduln und exzeptionelle $n$-Kronecker-Darstellungen}\label{baum}
\setcounter{zaehler}{0}
\numberwithin{zaehler}{chapter}
Stets sei $k$ ein algebraisch abgeschlossener Körper. Wir interessieren uns für die Beschreibung von Köcherdarstellungen über $k$ durch ihren Koeffizientenköcher und beweisen nach Ringel \cite{Ri3} das folgende Resultat: Jede unzerlegbare endlichdimensionale Köcherdarstellung ohne Selbsterweiterungen besitzt einen baumartigen Koeffizientenköcher. Eine wichtige Rolle spielen dabei die exzeptionellen Darstellungen der Köcher $Q_n$, die man mit Hilfe von Spiegelungsfunktoren aus den beiden einfachen $Q_n$-Darstellungen konstruieren kann.
\section{Grundlegende Begriffe}\label{baum:beg}
\setcounter{zaehler}{0}
\numberwithin{zaehler}{section}
In diesem Abschnitt sei stets $K$ ein endlicher Köcher.
\begin{defn}\label{baum:beg:def1}
Sei $M$ eine endlichdimensionale Darstellung von $K$.
\begin{enumerate}[a)]
\item Eine $k$-Basis von $M$ ist ein Tupel $\mathcal{B}=(\mathcal{B}_x)_{x\in K_0}$, wobei $\mathcal{B}_x$ eine $k$-Basis von $M(x)$ ist für alle $x\in K_0$.
\item Sei $\mathcal{B}=(\mathcal{B}_x)_{x\in K_0}$ eine $k$-Basis von $M$. Für einen Pfeil $\alpha: x\rightarrow y$ in $K_1$ sei 
\[M_\mathcal{B}(\alpha):=\Bigl[\bigl(M_\mathcal{B}(\alpha)\bigr)(b', b)\Bigr]_{\begin{subarray}{l} b'\in \mathcal{B}_y\\ b\in \mathcal{B}_x\end{subarray}}\in k^{\mathcal{B}_y\times \mathcal{B}_x}\]
die Darstellungsmatrix von $M(\alpha)$ bezüglich der Basen $\mathcal{B}_x, \mathcal{B}_y$. Also ist
\[\bigl(M(\alpha)\bigr)(b)=\sum_{b'\in \mathcal{B}_{y}}\Bigl(\bigl(M_\mathcal{B}(\alpha)\bigr)(b',b)\Bigr)\cdot b'\text{ für alle } b\in\mathcal{B}_x\]
mit (eindeutig bestimmten) $\bigl(M_\mathcal{B}(\alpha)\bigr)(b,b')\in k$. Der Koeffizientenköcher $\Gamma:=\Gamma(M,\mathcal{B})$ von $M$ bezüglich $\mathcal{B}$ ist definiert durch
 \begin{itemize}
 \item[$\cdot$]$\Gamma_0:=\coprod\limits_{x\in K_0} \mathcal{B}_x$,
 \item[$\cdot$]$\Gamma_1:=\bigl\{(\alpha, b, b')\in K_1\times\Gamma_0\times\Gamma_0\mid b\in \mathcal{B}_{n(\alpha)}, b'\in \mathcal{B}_{s(\alpha)}, \bigl(M_\mathcal{B}(\alpha)\bigr)(b',b)\neq0\bigr\}$,
 \item[$\cdot$] Für $(\alpha, b, b')\in \Gamma_1$ setze $n(\alpha, b, b'):=b$ und $s(\alpha, b, b'):=b'$.
 \end{itemize}
\item $M$ heißt eine Baumdarstellung (bzw. ein Baummodul), falls es eine Basis $\mathcal{B}$ von $M$ gibt, so daß $\Gamma(M,\mathcal{B})$ baumartig ist. $\mathcal{B}$ heißt dann eine Baumbasis von $M$.
\end{enumerate}
\end{defn}

Wir wiederholen zwei elementare graphentheoretischen Resultate: Sei $G$ ein ungerichteter endlicher Graph mit Punktmenge $G_0$ und Kantenmenge $G_1$. Dann gilt:
\begin{itemize}
\item[$\cdot$]Ist $G$ zusammenhängend, so ist $\#G_1\geq\#G_0-1$.
\item[$\cdot$]$G$ ist genau dann ein Baum, wenn $G$ zusammenhängend ist und wenn $\#G_1=\#G_0-1$ gilt.

\end{itemize}
\begin{bem}\label{baum:beg:bem1}Sei $M$ eine Darstellung von $K$.
\begin{enumerate}
\item Sei $\mathcal{B}$ eine Basis von $M$. Sei $\mathcal{B^*}=(\mathcal{B}_x^*)_{x\in K_0}$ die duale Basis von $DM$. Dann ist \[\Gamma(DM,\mathcal{B}^*)=\Gamma(M,\mathcal{B})^{op}.\]
Insbesondere: Ist $\mathcal{B}$ eine Baumbasis von $M$, so ist $\mathcal{B^*}$ eine Baumbasis von $DM$. Die Klasse der Baummoduln wird also beim Dualisieren in sich überführt.
\item  (vgl. \cite[Teil 2, Property 1]{Ri3}) Folgende Aussagen sind äquivalent:
\begin{enumerate}[i)]
\item $M$ ist unzerlegbar.
\item Für alle Basen $\mathcal{B}$ von $M$ ist $\Gamma(M,\mathcal{B})$ zusammenhängend.
\end{enumerate}
\item (vgl. \cite[Teil 2, Property 2]{Ri3}) Sei $\mathcal{B}$ eine Baumbasis von $M$. Dann existieren für alle $b\in \mathcal{B}$ Elemente $\lambda_b\in k^*$, so daß bezüglich der Basis $\mathcal{B}'=(\lambda_b\cdot b)_{b\in B}$ alle Darstellungsmatrizen $M_{\mathcal{\mathcal{B}}'}(\alpha)$ Einträge in $\{0,1\}$ haben.
\item Sei $\mathcal{B}$ eine Basis von $M$. Sei $\Gamma:=\Gamma(M,\mathcal{B})$. Betrachte den Morphismus von Köchern $\phi:\Gamma\longrightarrow K$ definiert durch
\[
\begin{array}{rcccl}
\phi_0:&\Gamma_0&\longrightarrow&K_0,&\phi_0(b):=x\text{ falls } b\in\mathcal{B}_x,\\
\phi_1:&\Gamma_1&\longrightarrow&K_1, &\phi_1(\alpha,b,b'):=\alpha.
\end{array}
\]
Im allgemeinen läßt sich die Darstellung $M$ nicht aus $\phi$ zurückgewinnen; ist aber $\Gamma$ baumartig, so ist $M$ nach (3) durch $\phi$ eindeutig bestimmt.
\item Sei $\mathcal{B}$ eine Basis von $M$. Für einen Pfeil $\alpha:x\rightarrow y$ in $K_1$ sei \[Tr\bigl(M_\mathcal{B}(\alpha)\bigr):=\bigl\{(b',b)\in \mathcal{B}_y\times B_x \mid \bigl(M_\mathcal{B}(\alpha)\bigr)(b',b)\neq 0\bigr\}\] der Träger der Matrix $M_\mathcal{B}(\alpha)$. Dann gelten die Formeln
\[\#\Gamma(M,\mathcal{B})_0 = \dim M\text{ und  } \#\Gamma(M,\mathcal{B})_1 = \sum_{\alpha\in K_1}{\#Tr\bigl(M_\mathcal{B}(\alpha)\bigr)}.\]
Insbesondere gilt: Ist $M$ unzerlegbar und gilt $\dim M = \sum_{\alpha\in K_1}{\#Tr\bigl(M_\mathcal{B}(\alpha)\bigr)}+1$,
so ist $\mathcal{B}$ eine Baumbasis von $M$.
\end{enumerate}
\end{bem}

\begin{bspl}\label{baum:beg:bsp1}
Sei $K=Q_2$; bezeichne die Quelle mit $y$, die Senke mit $x$ und die zwei Pfeile mit $\alpha$ und $\beta$. Wir berechnen Koeffizientenköcher der unzerlegbaren Kronecker-Moduln und bestimmen alle unzerlegbaren Baummoduln. Wir verwenden die Notationen aus Abschnitt \ref{grund:kron}.
\begin{itemize}
\item[$\cdot$] Für alle $n\geq 0$ ist der präprojektive unzerlegbare $A$-Modul $P_n$ zum Dimensionsvektor $(n, n+1)$ ein Baummodul. Seien $(z_1,\ldots,z_n)$ und $(u_1,\ldots,u_{n+1})$ die Standardbasen von $P_n(y)$ und $P_n(x)$. Der Koeffizientenköcher $\Gamma$ von $P_n$ bezüglich dieser Basis ist gegeben durch
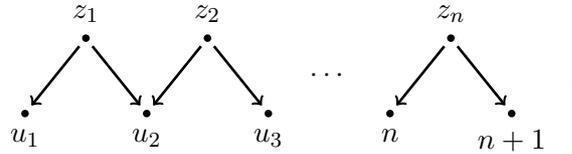
\begin{center}
\parbox[c]{7.5cm}{
\begin{tikzpicture}[line width=1pt]
	\tikzstyle{knt}=[circle, fill, inner sep=1pt]
  \node at (0.8,0)[knt, label=above:$z_1$](z1){};
  \node at (2.4,0)[knt, label=above:$z_2$](z2){};
  \node at (5.6,0)[knt, label=above:$z_n$](zn){};
  
  \node at (0,-1)[knt, label=below:$u_1$](u1){};
  \node at (1.6,-1)[knt, label=below:$u_2$](u2){};
  \node at (3.2,-1)[knt, label=below:$u_3$](u3){};
  \node at (4.8,-1)[knt, label=below:$n$](un){};
  \node at (6.4,-1)[knt, label=below:$n+1$](un+1){};
  
  \draw[->, shorten >= 2pt, shorten <= 2pt] (z1) to (u1);
  \draw[->, shorten >= 2pt, shorten <= 2pt] (z1) to (u2);
  \draw[->, shorten >= 2pt, shorten <= 2pt] (z2) to (u2);
  \draw[->, shorten >= 2pt, shorten <= 2pt] (z2) to (u3);
  \draw[->, shorten >= 2pt, shorten <= 2pt] (zn) to (un);
  \draw[->, shorten >= 2pt, shorten <= 2pt] (zn) to (un+1);
  
  \draw (4, -0.5) node{$\cdots$};
\end{tikzpicture}}.\end{center}
$\Gamma$ ist der (bis auf Isomorphie) eindeutig bestimmte baumartige Koeffizientenköcher von $P_n$: Sei dazu eine weitere Baumbasis von $P_n$ durch $k$-Basen $(z'_1,\ldots,z'_n)$ von $P_n(y)$ und $(u'_1,\ldots,u'_{n+1})$ $P_n(x)$ gegeben. Sei $\Gamma'$ der zugehörige Koeffizientenköcher. Offenbar sind alle $z'_i$ Quellen und alle $u'_i$ Senken von $\Gamma'$. Da $P_n(\alpha)$ und $P_n(\beta)$ injektiv sind, starten in jeder Quelle von $\Gamma'$ mindestens zwei Pfeile. Da $\Gamma'$ baumartig ist, hat $\Gamma'$ genau $2n$ Pfeile, also starten in jeder Quelle genau $2$ Pfeile, und man folgert leicht $\Gamma'\simeq\Gamma$.
\item[$\cdot$]Analoge Aussagen gelten für präinjektive Moduln.
\item[$\cdot$] Sei $n\geq1$. Dann gibt es (bis auf Isomorphie) genau zwei unzerlegbare Baummoduln zum Dimensionsvektor $(n, n)$, nämlich die regulären Modul $R_{n,0}$ und $R_{n,\infty}$. Die Koeffizientenköcher von $R_{n,0}$ und $R_{n,\infty}$ (bezüglich der Standardbasis) sind dann isomorph und gegeben durch
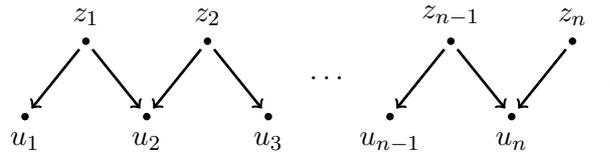
\begin{center}
\parbox[c]{8cm}{
\begin{tikzpicture}[line width=1pt]
	\tikzstyle{knt}=[circle, fill, inner sep=1pt]
  \node at (0.8,0)[knt, label=above:$z_1$](z1){};
  \node at (2.4,0)[knt, label=above:$z_2$](z2){};
  \node at (5.6,0)[knt, label=above:$z_{n-1}$](zn-1){};
  \node at (7.2,0)[knt, label=above:$z_n$](zn){};
  
  \node at (0,-1)[knt, label=below:$u_1$](u1){};
  \node at (1.6,-1)[knt, label=below:$u_2$](u2){};
  \node at (3.2,-1)[knt, label=below:$u_3$](u3){};
  \node at (4.8,-1)[knt, label=below:$u_{n-1}$](un-1){};
  \node at (6.4,-1)[knt, label=below:$u_n$](un){};
  
  \draw[->, shorten >= 2pt, shorten <= 2pt] (z1) to (u1);
  \draw[->, shorten >= 2pt, shorten <= 2pt] (z1) to (u2);
  \draw[->, shorten >= 2pt, shorten <= 2pt] (z2) to (u2);
  \draw[->, shorten >= 2pt, shorten <= 2pt] (z2) to (u3);
  \draw[->, shorten >= 2pt, shorten <= 2pt] (zn-1) to (un-1);
  \draw[->, shorten >= 2pt, shorten <= 2pt] (zn-1) to (un);
  \draw[->, shorten >= 2pt, shorten <= 2pt] (zn) to (un);
  
  \draw (4, -0.5) node{$\cdots$};
\end{tikzpicture}}.\end{center}
Mit ähnlichen Argumenten wie oben sieht man, daß diese Köcher die einzigen baumartigen Koeffizientenköcher von $R_{n,0}$ und $R_{n,\infty}$ sind.
\item[$\cdot$]Für $n\geq1$ und für $\lambda\in k^*$ ist der Koeffizientenköcher von $R_{n,\lambda}$ bezüglich der Standardbasis gegeben durch
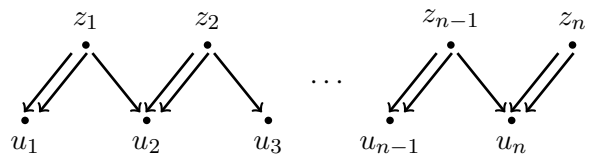
\begin{center}
\parbox[c]{8cm}{
\begin{tikzpicture}[line width=1pt]
	\tikzstyle{knt}=[circle, fill, inner sep=1pt]
  \node at (0.8,0)[knt, label=above:$z_1$](z1){};
  \node at (2.4,0)[knt, label=above:$z_2$](z2){};
  \node at (5.6,0)[knt, label=above:$z_{n-1}$](zn-1){};
  \node at (7.2,0)[knt, label=above:$z_n$](zn){};
  
  \node at (0,-1)[knt, label=below:$u_1$](u1){};
  \node at (1.6,-1)[knt, label=below:$u_2$](u2){};
  \node at (3.2,-1)[knt, label=below:$u_3$](u3){};
  \node at (4.8,-1)[knt, label=below:$u_{n-1}$](un-1){};
  \node at (6.4,-1)[knt, label=below:$u_n$](un){};
    
  \draw[transform canvas={xshift=-1mm},->, shorten >= 2pt, shorten <= 2pt] (z1) to (u1);
  \draw[transform canvas={xshift=1mm},->, shorten >= 2pt, shorten <= 2pt] (z1) to (u1);
  \draw[->, shorten >= 2pt, shorten <= 2pt] (z1) to (u2);
  \draw[transform canvas={xshift=-1mm},->, shorten >= 2pt, shorten <= 2pt] (z2) to (u2);
  \draw[transform canvas={xshift=1mm},->, shorten >= 2pt, shorten <= 2pt] (z2) to (u2);
  \draw[->, shorten >= 2pt, shorten <= 2pt] (z2) to (u3);
  \draw[transform canvas={xshift=-1mm},->, shorten >= 2pt, shorten <= 2pt] (zn-1) to (un-1);
  \draw[transform canvas={xshift=1mm},->, shorten >= 2pt, shorten <= 2pt] (zn-1) to (un-1);
  \draw[->, shorten >= 2pt, shorten <= 2pt] (zn-1) to (un);
  \draw[transform canvas={xshift=-1mm},->, shorten >= 2pt, shorten <= 2pt] (zn) to (un);
  \draw[transform canvas={xshift=1mm},->, shorten >= 2pt, shorten <= 2pt] (zn) to (un);
  
  \draw (4, -0.5) node{$\cdots$};
\end{tikzpicture}}.\end{center}
$R_{n,\lambda}$ ist kein Baummodul, denn: Jeder Koeffizientenköcher $\Gamma$ von $R_{n,\lambda}$ hat $2n$ Punkte. Wegen $\Ker R_{n,\lambda}(\alpha)=\Ker R_{n,\lambda}(\beta)=0$ hat er aber mindestens $2n$ Pfeile, ist also nicht baumartig.
\end{itemize}\enlargethispage*{\baselineskip}\enlargethispage*{\baselineskip}
Ist $M$ ein Baummodul der Dimension $n$ über einem Köcher $K$, so sind die baumartigen Koeffizientenköcher von $K$ genua diejenigen Koeffizientenköcher mit einer minimalen Anzahl an Pfeilen (und diese Anzahl ist gerade $n-1$). Im Fall von Kronecker-Moduln sind nach obigem Beispiel die baumartigen Koeffizientenköcher unzerlegbarer Baummoduln eindeutig bestimmt bis auf Isomorphie. Diese Aussage ist im allgemeinen falsch (vgl. das elementare Beispiel in \cite[Teil 2, Remark 2]{Ri3} im Falle eines Dynkin-Köchers vom Typ $D_4$); der Koeffizientenköcher eines Moduls hängt also sehr stark von der Wahl einer Basis ab.
\end{bspl}

\section{Exkurs: Spiegelungsfunktoren für Überlagerungen von Kronecker-Köchern}\label{baum:exk}
\setcounter{zaehler}{0}
\numberwithin{zaehler}{section}
Wir definieren gewisse Spiegelungsfunktoren für die universelle Überlagerung $\widetilde{Q}_n$ von $Q_n$ und untersuchen deren Beziehung zu den gewöhnlichen Spiegelungsfunktoren für $Q_n$. Die Definitionen von Spiegelungsfunktoren und universellen Überlagerungen von Köchern findet man in den Abschnitten \ref{grund:refl} und \ref{grund:cov}.

Wir legen zunächst einige Notationen fest. Für den Köcher $Q_n$ verwenden wir die Bezeichnungen
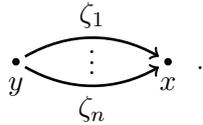
\begin{center}
\parbox[c]{2.6cm}{
\begin{tikzpicture}[line width=1pt]
	\tikzstyle{knt}=[circle, fill, inner sep=1pt]
  \node at (0,0)[knt, label=below:$y$](x){};
  \node at (2,0)[knt, label=below:$x$](y){};
  \node at (1,0.1)(z){$\vdots$};
  \draw[->, bend left, shorten >= 2pt, shorten <= 2pt] (x) to node[above]{$\zeta_1$} (y);
  \draw[->, bend right, shorten >= 2pt, shorten <= 2pt] (x) to node[below]{$\zeta_n$} (y);
\end{tikzpicture}}.\end{center}
Es sei $\pi=(\pi_0,\pi_1):\widetilde{Q}_n\longrightarrow Q_n$ die universelle Überlagerung. Setze
\[\mathcal{X}:=\pi_0^{-1}(x)\text{ und }\mathcal{Y}:=\pi_0^{-1}(y).\]
Also ist $\mathcal{Y}$ die Menge aller Quellen, $\mathcal{X}$ die Menge aller Senken in $\widetilde{Q}_n$ und es ist $(\widetilde{Q}_n)_0=\mathcal{X}\cup\mathcal{Y}$. 

\begin{lemma}\label{baum:exk:lemma1}
$\widetilde{Q}_n$ ist ein (unendlicher) baumartiger Köcher mit folgenden Eigenschaften:
\begin{enumerate}
\item Jeder Punkt in $\widetilde{Q}_n$ ist Quelle oder Senke.
\item Jeder Punkt hat genau $n$ Nachbarn.
\end{enumerate}
Ist umgekehrt $K$ ein baumartiger Köcher mit diesen Eigenschaften, so ist $K\simeq\widetilde{Q}_n$. Wir nennen $\widetilde{Q}_n$ auch den $n$-regulären, bipartiten, baumartigen Köcher.
\end{lemma}
Für alle Quellen $q\in\mathcal{Y}$ und für alle $i=1,\ldots,n$ sei \[\zeta_i^{(q)}:q\longrightarrow s_i^{(q)}\]
der eindeutig bestimmte Pfeil in $\widetilde{Q}_n$ mit Start $q$ über $\zeta_i$, analog sei für alle Senken $s\in\mathcal{X}$ und für alle $i=1,\ldots,n$ \enlargethispage{\baselineskip}\enlargethispage{\baselineskip}
\[\zeta_i^{(s)}:q_i^{(s)}\longrightarrow s\]
der eindeutig bestimmte Pfeil in $\widetilde{Q}_n$ mit Spitze $s$ über $\zeta_i$.

Der Überlagerungsfunktor $F_\lambda:\fdar{\widetilde{Q}_n}\longrightarrow\fdar{Q_n}$ bezüglich $\pi$ hat in diesen Notationen folgende Form: Sei $M\in\fdar{\widetilde{Q}_n}$. Dann ist
\[(F_\lambda M)(y)=\bigoplus_{q\in\mathcal{Y}}M(q)\text{ und }(F_\lambda M)(x)=\bigoplus_{s\in\mathcal{X}}M(s).\]
Für alle $i=1,\ldots,n$ wird die lineare Abbildung $(F_\lambda M)(\zeta_i):(F_\lambda M)(y)\longrightarrow (F_\lambda M)(x)$ induziert von den $M(\zeta_i^{(q)}):M(q)\longrightarrow M(s_i^{(q)}), q\in\mathcal{Y}$.\newpage

Wir erhalten außerdem eine universelle Überlagerung $\pi^{op}:\widetilde{Q}_n^{op}\longrightarrow Q_n^{op}$ von $Q_n^{op}$; hier ist $\pi=\pi^{op}$ als Abbildung auf den Punkten und Pfeilen. Den zugehörigen Überlagerungsfunktor bezeichnen wir mit $F_\lambda^{op}:\fdar{\widetilde{Q}_n^{op}}\longrightarrow\fdar{Q_n^{op}}$.

\begin{defn}\label{baum:exk:def1}\ 
\begin{enumerate}
\item Der Spiegelungsfunktor \[\widetilde{\sigma}_\mathcal{X}^+:\fdar{\widetilde{Q}_n}\longrightarrow\fdar{\widetilde{Q}_n^{op}}\]
bezüglich aller Senken in $\widetilde{Q}_n$ wird folgendermaßen definiert: Sei $M$ eine 
$\widetilde{Q}_n$-Darstellung. Für alle Quellen $q\in\mathcal{Y}$ sei $\widetilde{\sigma}_\mathcal{X}^+(M)(q):=M(q)$. Für eine Senke $s\in\mathcal{X}$ betrachte die $k$-lineare Abbildung
\[\phi:\bigoplus_{i=1}^n M(q_i^{(s)})\longrightarrow M(s),\  \phi:=\begin{bmatrix}M(\zeta_1^{(s)})\cdots M(\zeta_n^{(s)})\end{bmatrix}.\] 
Erhalte eine kurze exakte Sequenz
\[0\longrightarrow \widetilde{\sigma}_\mathcal{X}^+(M)(s)\stackrel{\psi}{\longrightarrow}\bigoplus_{i=1}^n M(q_i^{(s)})\stackrel{\phi}{\longrightarrow} M(s)\]
von $k$-Vektorräumen. Für $i=1,\ldots,n$ sei \[\widetilde{\sigma}_\mathcal{X}^+(M)(\zeta_i^{(s)}):\widetilde{\sigma}_\mathcal{X}^+(M)(s)\longrightarrow M(q_i^{(s)})=\widetilde{\sigma}_\mathcal{X}^+(M)(q_i^{(s)})\]
die von $\psi$ induzierte Abbildung.
\item Der Spiegelungsfunktor $\widetilde{\sigma}_\mathcal{Y}^-:\fdar{\widetilde{Q}_n}\longrightarrow\fdar{\widetilde{Q}_n^{op}}$ bezüglich aller Quellen in $\widetilde{Q}_n$ wird dual definiert.
\item Die Coxeter-Funktoren von $\widetilde{Q}_n$ sind definiert durch
\[\widetilde{\Phi}^+:=\widetilde{\sigma}_{\mathcal{Y}^{op}}^+\circ\widetilde{\sigma}_\mathcal{X}^+ \text{ und }\widetilde{\Phi}^-:\widetilde{\sigma}_{\mathcal{X}^{op}}^-\circ\widetilde{\sigma}_\mathcal{Y}^-.\]
Hier bezeichnet $\widetilde{\sigma}_{\mathcal{Y}^{op}}^+$ den auf $\fdar{\widetilde{Q}_n^{op}}$ definierten Spiegelungsfunktor bezüglich aller Senken in $\widetilde{Q}_n^{op}$, analog ist $\widetilde{\sigma}_{\mathcal{X}^{op}}^-$ definiert.
\end{enumerate}
\end{defn}

\begin{bem}\label{baum:exk:bem1} Sei $M$ eine $\widetilde{Q}_n$-Darstellung. Sei $K_M$ ein endlicher, zusammenhängender, voller Unterköcher von $\widetilde{Q}_n$ mit $Tr(M)\subseteq K_M$, so daß gilt: Ist $x\in Tr(M)$, und ist $w$ eine Wanderung in $\widetilde{Q}_n$ der Länge $\leq2$ mit Start $x$, so liegt der Zielpunkt von $w$ in $K_M$. 
\begin{enumerate} 
\item Seien $s_1,\ldots,s_n$ die Senken in $K_M$. Dann existiert ein Isomorphismus 
\[\bigl(\sigma_{s_n}^+\circ\ldots\circ\sigma_{s_1}^+\bigr)(M)\simeq \widetilde{\sigma}_{\mathcal{X}}^+(M)\]
von $K_M^{op}$-Darstellungen.
\item Sei $\widetilde{\Phi}_{K_M}^+$ der Coxeter-Funktor von $K_M$ bezüglich der Senken in $K_M$. Dann existiert ein Isomorphismus
\[\widetilde{\Phi}_{K_M}^+(M)\simeq \widetilde{\Phi}^+(M)\]
von $K_M$-Darstellungen.
\end{enumerate}
\end{bem}

Mit dieser Bemerkung können wir die üblichen Eigenschaften für Spiegelungsfunktoren auf unsere Situation übertragen (vgl. Abschnitt \ref{grund:refl}):
\begin{lemma}\label{baum:exk:lemma2}Betrachte die Funktoren
\[
\widetilde{\sigma}_\mathcal{X}^+:\fdar{\widetilde{Q}_n}\longrightarrow \fdar{\widetilde{Q}_n^{op}}\text{ und } \widetilde{\sigma}_{\mathcal{X}^{op}}^-:\fdar{\widetilde{Q}_n^{op}}\longrightarrow \fdar{\widetilde{Q}_n}.
\]
Sei $\mathcal{C}$ die volle Unterkategorie von $\fdar{\widetilde{Q}_n}$ bestehend aus allen Darstellungen, die keinen projek\-tiv-einfachen direkten Summanden haben. Analog sei $\mathcal{C}^{op}$ die volle Unterkategorie von $\fdar{\widetilde{Q}_n^{op}}$ bestehend aus allen Darstellungen, die keinen injektiv-einfachen direkten Summanden haben.
\begin{enumerate}
\item $(\widetilde{\sigma}_{\mathcal{X}^{op}}^-,\widetilde{\sigma}_\mathcal{X}^+)$ ist ein adjungiertes Paar von Funktoren.
\item $\widetilde{\sigma}_{\mathcal{X}^{op}}^-$ und $\widetilde{\sigma}_\mathcal{X}^+$ induzieren zueinander quasiinverse Äquivalenzen zwischen $\mathcal{C}$ und $\mathcal{C}^{op}$.
\item Sei 
\[0\longrightarrow M'\longrightarrow M\longrightarrow M''\longrightarrow0\]
eine exakte Sequenz von $\widetilde{Q}_n$-Darstellungen mit $M',M''\in\mathcal{C}$. Dann ist auch $M\in\mathcal{C}$, und die Sequenz
\[0\longrightarrow \widetilde{\sigma}_\mathcal{X}^+(M')\longrightarrow \widetilde{\sigma}_\mathcal{X}^+(M)\longrightarrow \widetilde{\sigma}_\mathcal{X}^+(M'')\longrightarrow0\]
ist exakt.
\item Seien $M,N\in\mathcal{C}$. Dann induziert $\widetilde{\sigma}_\mathcal{X}^+$ Isomorphismen
\begin{gather*}\Hom_{\widetilde{Q}_n}(M,N)\simeq\Hom_{\widetilde{Q}_n^{op}}(\widetilde{\sigma}_\mathcal{X}^+(M),\widetilde{\sigma}_\mathcal{X}^+(N)),\\ \Ext^1_{\widetilde{Q}_n}(M,N)\simeq\Ext^1_{\widetilde{Q}_n^{op}}(\widetilde{\sigma}_\mathcal{X}^+(M),\widetilde{\sigma}_\mathcal{X}^+(N)).\end{gather*}
\item Sei $M\in\mathcal{C}$ unzerlegbar (exzeptionell). Dann ist auch $\widetilde{\sigma}_\mathcal{X}^+(M)$ unzerlegbar (exzeptionell).
\item Seien $\widetilde{\tau}$ und $\widetilde{\tau}^{-1}$ die Auslander-Reiten-Verschiebungen für $\widetilde{Q}_n$. Ist $M$ eine Darstellung von $\widetilde{Q}_n$ ohne projektiven direkten Summanden, so existiert ein Isomorphismus \[\widetilde{\Phi}^+(M)\simeq\widetilde{\tau}(M)\]
von $\widetilde{Q}_n$-Darstellungen. Besitzt $M$ keinen injektiven direkten Summanden, so  gilt dual 
\[\widetilde{\Phi}^-(M)\simeq\widetilde{\tau}^{-1}(M).\]
\end{enumerate}
\end{lemma}
\begin{bew} Wir geben nur eine Beweisskizze für Teil (6) an. Sei $M$ unzerlegbar, nicht-projektiv. Sei $K_M$ ein endlicher, voller, zusammenhängender Unterköcher von $\widetilde{Q}_n$, der den Träger von $M$ und von $\widetilde{\tau}(M)$ enthält. Für alle Wanderungen $w$ der Länge $\leq2$ in $\widetilde{Q}_n$ mit Startpunkt in $Tr(M)$ oder $Tr(\widetilde{\tau}M)$ sei außerdem der Zielpunkt von $w$ in $K_M$ enthalten. Die Auslander-Reiten-Sequenz
\[0\longrightarrow \widetilde{\tau}(M)\longrightarrow X\longrightarrow M\longrightarrow 0\]
mit Ende $M$ in $\fdar{\widetilde{Q}_n}$ ist dann auch eine Auslander-Reiten-Sequenz in $\fdar{K_M}$, und man reduziert mit \cref{baum:exk:bem1} auf den Fall gewöhnlicher Spiegelungsfunktoren (vgl. \cref{grund:refl:prop2}). 
\end{bew}

\begin{bem}\label{baum:exk:bem2} Seien $\sigma_x^+:\fdar{Q_n}\longrightarrow \fdar{Q_n^{op}}$ und $\sigma_y^-:\fdar{Q_n}\longrightarrow \fdar{Q_n^{op}}$ die gewöhnlichen Spiegelungsfunktoren bezüglich der Senke $x$ bzw. der Quelle $y$ von $Q_n$. Man verifiziert mit Hilfe der Definitionen, daß kommutative Diagramme von Funktoren existieren:
\[
\begin{CD}
\fdar{\widetilde{Q}_n}@>\widetilde{\sigma}_\mathcal{X}^+>>\fdar{\widetilde{Q}_n^{op}}\\
@VF_\lambda VV@VVF_\lambda^{op}V\\
\fdar{Q_n}@>>\sigma_x^+>\fdar{Q_n^{op}}
\end{CD}\ \ \ \ \text{ und }\ \ \ \ 
\begin{CD}
\fdar{\widetilde{Q}_n}@>\widetilde{\sigma}_\mathcal{Y}^->>\fdar{\widetilde{Q}_n^{op}}\\
@VF_\lambda VV@VVF_\lambda^{op}V\\
\fdar{Q_n}@>>\sigma_y^->\fdar{Q_n^{op}}.
\end{CD}
\]
\end{bem}

\section{Exakte Sequenzen für exzeptionelle $n$-Kronecker-Darstellungen}\label{baum:kron}
\setcounter{zaehler}{0}
\numberwithin{zaehler}{section}
Wir verwenden die Notationen aus dem letzten Abschnitt und vereinbaren zusätzlich: Sei
\[\widetilde{\sigma}_{\mathcal{X}^{op}}^-:\fdar{\widetilde{Q}_n^{op}}\longrightarrow\fdar{\widetilde{Q}_n}\]
der Spiegelungsfunktor bezüglich aller Quellen in $\widetilde{Q}_n^{op}$, und sei
\[\sigma_{x^{op}}^-:\fdar{Q_n^{op}}\longrightarrow\fdar{Q_n}\]
der Spiegelungsfunktor bezüglich der Quelle $x$ in $Q_n^{op}$. Analog sind $\widetilde{\sigma}_{\mathcal{Y}^{op}}^+$ und $\sigma_{y^{op}}^+$ definiert. Weiter sei \[\Omega:\fdar{Q_n^{op}}\stackrel{\sim}{\longrightarrow}\fdar{Q_n}\]
die durch Vertauschen der Punkte $x$ und $y$ gegebene Äquivalenz. Es sei schließlich $G$ die Gruppe aller Automorphismen der universellen Überlagerung $\pi:\widetilde{Q}_n\longrightarrow\ Q_n$, also die Menge aller Köcherautomorphismen $\phi:\widetilde{Q}_n\longrightarrow\widetilde{Q}_n$ mit $\pi\circ\phi=\pi$. G operiert dann auf $\widetilde{Q}_n$ und man erhält eine induzierte Operation auf der Darstellungskategorie $\fdar{\widetilde{Q}_n}$. 

Sei $M$ eine exzeptionelle $Q_n$-Darstellung. Dann ist $M$ präprojektiv oder präinjektiv, und allgemeine Resultate von Gabriel in \cite{Ga2} über die Verträglichkeit von Auslander-Reiten-Sequenzen mit Push-Down-Funktoren im Rahmen der Theorie der Galois-Überlagerungen zeigen, daß sich $M$ zu einer präprojektiven oder präinjektiven $\widetilde{Q}_n$-Darstellung $\widetilde{M}$ liften läßt. Diese Liftung ist eindeutig bestimmt modulo der Operation der Gruppe $G$.

Mit Hilfe der Spiegelungsfunktoren für $Q_n$ läßt sich $M$ aus einer einfachen $Q_n$-Darstellung konstruieren. Genauso läßt sich $\widetilde{M}$ mit Hilfe der im letzten Abschnitt definierten Spiegelungsfunktoren für $\widetilde{Q}_n$ aus einer einfachen $\widetilde{Q}_n$-Darstellung konstruieren (vgl. dazu die  Diskussion in \cite[Lemma 2.3]{Ri4}).

Für $\widetilde{M}$ konstruieren wir interessante exakte Sequenzen, die ein wichtiges Hilfsmittel für einen Induktionsbeweis des Satzes von Ringel sein werden. Wir erhalten außerden die Aussage, daß die exzeptionellen $Q_n$-Darstellungen Erweiterungen eines Unzerlegbaren und eines Einfachen sind.

Wir fassen diese Bemerkungen in zwei Propositionen zusammen.
\begin{prop}\label{baum:kron:prop1} Sei $M$ eine exzeptionelle Darstellung von $Q_n$.
\begin{enumerate}
\item Es existiert eine $\widetilde{Q}_n$-Darstellung $\widetilde{M}$ mit $F_\lambda(\widetilde{M})=M$. 
\item $\widetilde{M}$ ist modulo der $G$-Operation auf $\fdar{\widetilde{Q}_n}$ eindeutig bestimmt bis auf Isomorphie, das heißt: Ist $\widetilde{N}$ eine weitere $\widetilde{Q}_n$-Darstellung mit $F_\lambda(\widetilde{N})=M$, so existieren ein Automorphismus $g\in G$ und ein Isomorphismus $\widetilde{M}\simeq\widetilde{N}^g$.
\end{enumerate}
\end{prop}
\begin{bew} (1.) Aus \cite[Theorem 3.6 b)]{Ga2} folgt, daß sich beliebige präprojektive oder präinjektive Darstellungen von endlichen, zykellosen, zusammenhängenden Köchern zu Darstellungen der universellen Überlagerung liften lassen.

(2.) Aus \cite[Lemma 3.5]{Ga2} folgt, daß Liftungen von Darstellungen (falls sie existieren) eindeutig bestimmt sind modulo der Operation der Gruppe der Decktransformationen.
\end{bew}
Wir wollen nun die Liftungen mit Hilfe der Spiegelungsfunktoren aus dem letzten Abschnitt konstruieren. Da beim Spiegeln die Orientierung umgedreht wird, werden wir die Liftungen für gerade $i$ als Darstellungen von $\widetilde{Q}_n$ und für ungerade $i$ als Darstellungen von $\widetilde{Q}_n^{op}$ deuten. Dieser Standpunkt ist auch für die konkreten Rechnungen im nächsten Abschnitt von Vorteil. 

Wir fixieren für den Rest dieses Abschnitts eine Senke $s\in\mathcal{X}$ von $\widetilde{Q}_n$. Wir definieren die Folge $(\widetilde{P}_i)_{i\geq0}$ der Liftungen wie folgt: Für alle $j$ sei $\widetilde{P}_{2j}\in\fdar{\widetilde{Q}_n}$ und $\widetilde{P}_{2j+1}\in\fdar{\widetilde{Q}_n^{op}}$. Es seien die Rekursionsformeln
\[
\widetilde{P}_0:=E_s\in\fdar{\widetilde{Q}_n}\text{ und } 
\widetilde{P}_{i}:=\left\{
\begin{array}{cc}
\widetilde{\sigma}_{\mathcal{X}^{op}}^-(\widetilde{P}_{i-1})&\text{ falls $i$ gerade}\\
\widetilde{\sigma}_{\mathcal{Y}}^-(\widetilde{P}_{i-1})&\text{ falls $i$ ungerade}
\end{array}
\right. \text{ für } i\geq 1
\]
erfüllt.
\begin{prop}\label{baum:kron:prop2}Sei $i\geq0$ eine natürliche Zahl. Dann gilt:
\begin{enumerate}
\item $\bigl(\Omega\circ\sigma_y^-\bigr)(P_i)\simeq P_{i+1}$.
\item $P_{i}\simeq\left\{
\begin{array}{cc}
F_\lambda(\widetilde{P}_{i})&\text{ falls $i$ gerade}\\
\bigl(\Omega\circ F_\lambda^{op}\bigr)(\widetilde{P}_{i}) &\text{ falls $i$ ungerade}
\end{array}
\right.$.
\item $\widetilde{P}_i$ ist exzeptionell.
\end{enumerate}
\end{prop}
\begin{bew} (1.) Verwende Induktion nach $i$ und den Zusammenhang zwischen Coxeter-Funktoren und Auslander-Reiten-Verschiebungen. Den Fall $i=0$ verifiziert man leicht. Sei also $i\geq1$. Sei $\Phi^-$ der Coxeter-Funktor bezüglich der Quelle $y$ von $Q_n$. Per Definition von $P_{i+1}$ und wegen \cref{grund:refl:prop2} gilt
\[P_{i+1}=\tau^{-1}(P_{i-1})=\Phi^-(P_{i-1}).\]
Per Induktion erhalten wir
\[\Phi^-(P_{i-1})=\sigma_{x^{op}}^-\circ\sigma_{y}^-(P_{i-1})=\sigma_{x^{op}}^-\circ\Omega^{-1}(P_i)=\Omega\circ \sigma_{y}^-(P_i).\]
(2.) Verwende wieder Induktion nach $i$. Die Fälle $i=0,1$ verifiziert man leicht. Sei im Induktionsschritt $i\geq2$. Wir nehmen an, daß $i$ gerade ist; den ungeraden Fall behandelt man analog. Sei $\widetilde{\Phi}^-$ der Coxeter-Funktor von $\widetilde{Q}_n$ bezüglich der Quellen (vgl. \cref{baum:exk:def1}). Nach einem Resultat von Gabriel (vgl. \cite[Theorem 3.6]{Ga2}) erhält der Push-Down-Funktor $F_\lambda$ Auslander-Reiten-Sequenzen, also gilt für unzerlegbare, nicht-injektive Darstellungen $M$ 
\[F_\lambda\circ\widetilde{\tau}^{-1}(M)\simeq \tau^{-1}\circ F_\lambda(M).\]
Wir wissen nach \cref{baum:exk:lemma2}, daß $\widetilde{\Phi}^-$ auf nicht-injektiven, unzerlegbaren Darstellungen mit der Auslander-Reiten-Verschiebung $\widetilde{\tau}^{-1}$ übereinstimmt. Wir erhalten mit diesen Informationen und der Induktionsvoraussetzung:
\[F_\lambda(\widetilde{P}_i)=F_\lambda\circ \widetilde{\Phi}^-(\widetilde{P}_{i-2})=
F_\lambda\circ \widetilde{\tau}^{-1}(\widetilde{P}_{i-2})=\tau^{-1}\circ F_\lambda(\widetilde{P}_{i-2})=\tau^{-1}(P_{i-2})=P_i.\]
(3.) Wegen Teil (4) von \cref{baum:exk:lemma2} überführen die Spiegelungsfunktoren $\widetilde{\sigma}_{\mathcal{X}^{op}}^-$ und $\widetilde{\sigma}_{\mathcal{Y}}^-$ exzeptionelle Darstellungen in exzeptionelle Darstellungen. Damit folgt die Behauptung per Induktion aus der Definition der $\widetilde{P}_i$.
\end{bew}

Bevor wir zur Konstruktion der exakten Sequenzen für die $\widetilde{P}_i$ kommen, benötigen wir einige einfache kombinatorische Aussagen. Wir definieren für alle $i\geq 0$ induktiv Teilmengen $L_i$ von $(\widetilde{Q}_n)_0$ durch
\[L_0:=\{s\},\ L_i:=\bigl\{ x\in(\widetilde{Q}_n)_0-\bigcup_{j=0}^{i-1}L_j\bigm| x \text{ ist ein Nachbar eines Punktes } y\in L_{i-1}\bigr\}\text{ für }i\geq1.\]
\begin{bem}\label{baum:kron:bem1} Da $\widetilde{Q}_n$ der $n$-reguläre bipartite baumartige Köcher ist, gilt:
\begin{itemize}
\item[$\cdot$]$(\widetilde{Q}_n)_0=\bigcup\limits_{i=0}^\infty L_i$.
\item[$\cdot$]Für alle $i\geq1$ und für alle $x\in L_i$: $x$ hat genau $n-1$ Nachbarn in $L_{i+1}$ und genau einen Nachbarn in $L_{i-1}$.
\item[$\cdot$]Die bipartite Orientierung von $\widetilde{Q}_n$ ist wie folgt gegeben: Ist $i$ gerade, so besteht $L_i$ aus lauter Senken von $\widetilde{Q}_n$. Ist $i$ ungerade, so besteht $L_i$ aus lauter Quellen von $\widetilde{Q}_n$.
\end{itemize}
\end{bem}

Mit diesen Vereinbarungen erhalten wir:
\begin{lemma}\label{baum:kron:lemma2} Sei $i\geq0$ gegeben. Dann gilt:
\begin{enumerate}
\item $\widetilde{P}_i(x)=0$ für alle $x\in\bigcup\limits_{j=i+1}^\infty L_j$.
\item $\widetilde{P}_i(x)\simeq k$ für alle $x\in L_{i-1}\cup L_i$ ($L_{-1}:=\emptyset$).
\item $Tr(\widetilde{P}_i)=\bigcup\limits_{j=0}^{i}L_i$.
\item Sei $i\geq1$. Dann existiert eine Senke $s$ in $Tr(\widetilde{P}_i)$, so daß genau ein Pfeil $\alpha$ in $Tr(\widetilde{P}_i)$ mit Spitze $s$ existiert und so daß gilt: \[\widetilde{P}_i(n\alpha)\simeq k\simeq \widetilde{P}_i(s),\ \widetilde{P}_i(\alpha)\neq0.\]
\end{enumerate}
\end{lemma}
\begin{bew} Wir zeigen (1.) und (2.) gleichzeitig per Induktion nach $i$. Der Induktionsanfang $i=0$ ist klar. Sei also $i\geq1$. Wir nehmen an, daß $i$ ungerade ist; für gerades $i$ geht man analog vor. $\widetilde{P}_{i-1}$ ist also eine Darstellung von $\widetilde{Q}_n$ und es ist $\widetilde{\sigma}_{\mathcal{Y}}^-(\widetilde{P}_{i-1})=\widetilde{P}_{i}$. 

Sei $t\in (\widetilde{Q}_n)_0$. Per Definition des Spiegelungsfunktors
\[
\widetilde{\sigma}_{\mathcal{Y}}^-:\fdar{\widetilde{Q}_{n}}\longrightarrow\fdar{\widetilde{Q}_{n}^{op}}
\]
bezüglich aller Quellen gilt: Ist $t$ eine Senke in $\widetilde{Q}_n$, so ist
\counterwithout{equation}{chapter}
\setcounter{equation}{0}
\begin{equation}\label{baum:kron:lemma2:Gl1}
\widetilde{P}_{i}(t)=\widetilde{P}_{i-1}(t).
\end{equation}
Ist $t$ eine Quelle in $\widetilde{Q}_n$, so existiert eine exakte Sequenz von $k$-Vektorräumen
\begin{equation}\label{baum:kron:lemma2:Gl2}
\widetilde{P}_{i-1}(t)\longrightarrow\bigoplus_{l=1}^n\widetilde{P}_{i-1}(s_l^{(t)})\longrightarrow \widetilde{P}_{i}(t)\longrightarrow 0.
\end{equation}

Wir zeigen zunächst (1.). Sei also $t\in L_j$, $j\geq i+1$. Per Induktion ist $\widetilde{P}_{i-1}(t)=0$. Ist $j$ gerade, so ist $x$ eine Senke von $\widetilde{Q}_n$  und mit Gleichung (\ref{baum:kron:lemma2:Gl1}) folgt $\widetilde{P_{i}}(t)=0$. Ist $j$ ungerade, so ist $t$ eine Quelle von $\widetilde{Q}_n$. Für alle $l=1,\ldots,n$ ist $s_l^{(t)}\in L_{j-1}\cup L_{j+1}$. Wegen $j-1\geq i$ folgt per Induktion $\widetilde{P}_{i-1}(s_l^{(t)})=0$ für alle $l=1,\ldots,n$. Die exakte Sequenz (\ref{baum:kron:lemma2:Gl2}) zeigt dann $\widetilde{P}_{i}(t)=0$ und die Behauptung folgt.

Wir zeigen (2.). Sei zunächst $t\in L_{i-1}$. Da $i$ ungerade ist, ist $t$ eine Senke in $\widetilde{Q}_n$. Per Induktion ist $\widetilde{P}_{i-1}(t)\simeq k$ und Gleichung (\ref{baum:kron:lemma2:Gl1}) zeigt $\widetilde{P}_{i}(t)\simeq k$. Sei jetzt $t\in L_{i}$. Da $i$ ungerade ist, ist $x$ eine Quelle in $\widetilde{Q}_n$. Per Induktion ist $\widetilde{P}_{i-1}(t)=0$. Die exakte Sequenz (\ref{baum:kron:lemma2:Gl2}) liefert dann einen Isomorphismus
\[\bigoplus_{l=1}^n\widetilde{P}_{i-1}(s_l^{(t)})\simeq \widetilde{P}_{i}(t).\]
Die $s_l^{(t)}$ sind per Definition gerade die Nachfolger von $t$. Nach \cref{baum:kron:bem1} besitzt $t$ genau einen Nachfolger $u$ in $L_{i-1}$, alle anderen Nachfolger liegen in $L_{i+1}$. Per Induktion ist $\widetilde{P}_{i-1}(s)=0$ für alle $s\in L_{i+1}$ und $\widetilde{P}_{i-1}(u)\simeq k$. Also gilt $\bigoplus_{l=1}^n\widetilde{P}_{i-1}(s_l^{(t)})\simeq k$, also $\widetilde{P}_{i}(t)\simeq k$.

(3.) Wegen (1.) und (2.) gelten $Tr(\widetilde{P}_i)\subseteq\bigcup_{j=0}^{i}L_i$ und $L_i\subseteq Tr(\widetilde{P}_i)$. Da $\widetilde{P}_i$ unzerlegbar ist, ist der Träger von $\widetilde{P}_i$ zusammenhängend. Weil außerdem $\widetilde{Q}_n$ ein Baum ist, muß schon gelten:\linebreak$Tr(\widetilde{P}_i)=\bigcup_{j=0}^{i}L_i$.

(4.) Für gerades $i$ ist $Tr(\widetilde{P}_i)$ ein voller Unterköcher von $\widetilde{Q}_n$, für ungerades $i$ ein voller Unterköcher von $\widetilde{Q}_n^{op}$. Wegen (2.) ist $L_i\subseteq Tr(\widetilde{P}_i)$. Außerdem sind die Elemente von $L_i$ für gerades und ungerades $i$ allesamt Senken von $Tr(\widetilde{P}_i)$. Wähle dann $s\in L_i$ beliebig; die geforderte Eigenschaft folgt unmittelbar aus (1.) und (2.).
\end{bew}

Wir fassen unsere bisherigen Ergebnisse in folgender Proposition zusammen:
\begin{prop}\label{baum:kron:prop3} Sei $M$ eine exzeptionelle $Q_n$-Darstellung und sei $\widetilde{M}\in\fdar{\widetilde{Q}_n}$ mit \linebreak $F_\lambda(\widetilde{M})=M$. Dann gilt:
\begin{enumerate}
\item $\widetilde{M}$ ist exzeptionell.
\item Sei $M$ präprojektiv. Dann existiert eine Senke $s$ in $Tr(M)$, so daß gilt: Es existiert genau ein Pfeil $\alpha:q\longrightarrow s$ in $Tr(M)$ mit Spitze $s$, und es gilt:
\[\widetilde{M}(s)\simeq k\simeq\widetilde{M}(q),\ \widetilde{M}(\alpha)\neq0.\]
\end{enumerate}
\end{prop}
\begin{bew} Es sei $p:Q_n\longrightarrow Q_n^{op}$ der Köchermorphismus mit $p_0(x)=y$, $p_0(y)=x$ und $p_1=id$. Man sieht leicht, daß ein Köcherisomorphismus $\theta:\widetilde{Q}_n\longrightarrow\widetilde{Q}_n^{op}$ existiert, so daß folgendes Diagramm kommutiert:
\[
\begin{CD}
\widetilde{Q}_n@>\theta >>\widetilde{Q}_n^{op}\\
@V\pi VV @VV\pi^{op}V\\
Q_n@>>p> Q_n^{op}.
\end{CD}
\]
Sei $\Theta:\fdar{\widetilde{Q}_n^{op}}\longrightarrow \fdar{\widetilde{Q}_n}$ der von $\theta$ induzierte Funktor. Der von $p$ induzierte Funktor ist der bereits definierte Funktor $\Omega:\fdar{Q_n^{op}}\longrightarrow\fdar{Q_n}$. Wir erhalten offenbar ein kommutatives Diagramm von Funktoren:
\[
\begin{CD}
\fdar{\widetilde{Q}_n}@<\Theta <<\fdar{\widetilde{Q}_n^{op}}\\
@VF_\lambda VV @VVF_\lambda^{op}V\\
\fdar{Q_n}@<<\Omega< \fdar{Q_n^{op}}.
\end{CD}
\]
Aufgrund dieser Bemerkungen, wegen Teil (2) von \cref{baum:kron:prop2} und wegen der Eindeutigkeit der Liftung $\widetilde{M}$ modulo der $G$-Operation (\cref{baum:kron:prop1}) dürfen wir annehmen, daß schon $\widetilde{M}=\widetilde{P}_{2j}$ oder $\widetilde{M}=\Theta(\widetilde{P}_{2j+1})$ für ein $j\geq0$ gilt. Damit folgen beide Aussagen aus \cref{baum:kron:lemma2}.
\end{bew}
\begin{bem}\label{baum:kron:bem2} Man kann für den Beweis von \cref{baum:kron:prop3} auch auf die Konstruktion der Liftungen $\widetilde{P}_i$ durch Spiegelungsfunktoren verzichten, indem man eine präprojektive Komponente der universellen Überlagerung $\widetilde{Q}_n$ über der von $Q_n$ konstruiert und die benötigten kombinatorischen Aussagen über den Träger der Liftungen aus \cref{baum:kron:lemma2} mit Hilfe der Auslander-Reiten-Sequenzen beweist. Da wir für unsere Rechnungen im nächsten Abschnitt die Beschreibung von $\widetilde{P}_{i+1}$ als gespiegelte Darstellung von $\widetilde{P}_i$ verwenden wollen, haben wir hier den etwas komplizierteren Weg gewählt.
\end{bem}

Bevor wir das Hauptresultat dieses Abschnitts beweisen, halten wir folgendes einfache Lemma fest, daß uns im nächsten Abschnitt mehrfach begegnen wird.
\begin{lemma}\label{baum:kron:lemma3}
Sei $K$ ein endlicher, zykelloser Köcher und sei $s$ eine Senke in $K$. Es gebe genau einen Pfeil $\alpha$ in $K$ mit Spitze $s$. Sei $K'$ der volle Unterköcher von $K$ mit Punktmenge $K_0-\{s\}$ und sei $F:\fdar{K}\longrightarrow\fdar{K'}$ die Einschränkung.

Betrachte die vollen Unterkategorien $\mathcal{C}^{(\text{epi},\alpha)}$ und $\mathcal{C}^{(\text{iso},\alpha)}$ gegeben durch die folgenden Klassen von Darstellungen:
\begin{align*}
\mathcal{C}^{(\text{epi},\alpha)}&:=\{M\mid M(\alpha) \text{ ist ein Epimorphismus}\},\\
\mathcal{C}^{(\text{iso},\alpha)}&:=\{M\mid M(\alpha) \text{ ist ein Isomorphismus}\}.
\end{align*}
Dann gelten folgende Aussagen:
\begin{enumerate}
\item $F$ induziert einen treuen Funktor $F:\mathcal{C}^{(\text{epi},\alpha)}\longrightarrow\fdar{K'}$. Insbesondere gilt: Ist\linebreak $M\in\mathcal{C}^{(\text{epi},\alpha)}$ eine $K$-Darstellung, und ist $FM$ unzerlegbar, so ist auch $M$ unzerlegbar.
\item $F$ induziert eine Äquivalenz $F:\mathcal{C}^{(\text{iso},\alpha)}\stackrel{\sim}{\longrightarrow}\fdar{K'}$. Für alle $K$-Darstellungen\linebreak $M,N\in\mathcal{C}^{(\text{iso},\alpha)}$ induziert $F$ Isomorphismen
\[\Hom_K(M,N)\simeq \Hom_{K'}(FM,FN)\text{ und } \Ext^1_{K}(M,N)\simeq \Ext^1_{K'}(FM,FN).\]
\end{enumerate}
\end{lemma}
\begin{bew} (1.) Man verifiziert sofort, daß $F$ treu ist. Für alle $K$-Darstellungen $M\in\mathcal{C}^{(\text{epi},\alpha)}$ ist also $\End_K(FM)$ eine Unteralgebra von $\End_K(M)$; Unteralgebren von lokalen endlichdimensionalen Algebren sind lokal. Also folgt aus der Unzerlegbarkeit von $FM$ die von $M$.

(2.) Es ist klar, daß $F$ essentiell surjektiv und volltreu, also eine Äquivalenz ist. Man beachte außerdem, daß die Kategorie $\mathcal{C}^{(\text{iso},\alpha)}$ abgeschlossen ist unter Kernen, Kokernen und Erweiterungen. Also induziert der Vergißfunktor $\mathcal{C}^{(\text{iso},\alpha)}\longrightarrow \fdar{K}$ Isomorphismen zwischen Hom-Räumen  und $\Ext^1$-Gruppen. Daraus erhält man die in (2.) geforderten Isomorphismen.
\end{bew}
Das folgende Aussage ist eine unmittelbare Folgerung aus dem letzten Lemma:
\begin{kor}\label{baum:kron:kor1}
Sei $K$ ein endlicher, zykelloser Köcher und sei $s$ eine Senke in $K$. Es gebe genau einen Pfeil $\alpha:q\longrightarrow s$ in $K$ mit Spitze $s$. Sei $M$ eine $K$-Darstellung mit
\[M(s)\simeq k\simeq M(q), M(\alpha)\neq0.\]
Dann existiert eine exakte Sequenz
\[0\longrightarrow E_s\longrightarrow M\longrightarrow C\longrightarrow 0.\]
Dabei gilt
\[\End_K(C)\simeq \End_K(M)\text{ und }\Ext^1_K(M,M)\simeq\Ext^1_K(C,C).\]
Ist insbesondere $M$ unzerlegbar (exzeptionell), so ist auch $C$ unzerlegbar (exzeptionell).
\end{kor}
\begin{bew} Sei $K'$ der volle Unterköcher von $K$ mit Punktmenge $K_0-\{s\}$. Es sei weiter\linebreak $F:\fdar{K}\longrightarrow\fdar{K'}$ die Einschränkung. Fasse $\fdar{K'}$ durch Ergänzen mit $0$ als volle Unterkategorie von $\fdar{K}$ auf. Der Vergißfunktor $\iota:\fdar{K'}\longrightarrow\fdar{K}$ induziert Isomorphismen zwischen Hom-Räumen und $\Ext^1$-Gruppen. 

Da $s$ eine Senke ist, erhält man eine exakte Sequenz wie gefordert. Es ist $C(s)=0$, also ist $C$ eine $K'$-Darstellung. In den Notationen von \cref{baum:kron:lemma3} ist $M\in\mathcal{C}^{(\text{iso},\alpha)}$ und es ist $FM=C$. Mit Aussage (2) von \cref{baum:kron:lemma3} erhält man die gewünschten Isomorphismen
\[\End_K(C)\simeq \End_K(M)\text{ und }\Ext^1_K(M,M)\simeq\Ext^1_K(C,C).\]
\end{bew}
Wir kommen nun zum Hauptresultat dieses Abschnitts.
\begin{satz}\label{baum:kron:satz1} Sei $M$ eine präprojektive, unzerlegbare, nicht einfache $Q_n$-Darstellung. Dann gilt:
\begin{enumerate}
\item Sei $\widetilde{M}\in\fdar{\widetilde{Q}_n}$ mit $F_\lambda(\widetilde{M})=M$ (beachte \cref{baum:kron:prop1} zur Existenz und Eindeutigkeit solcher Liftungen). Dann existiert eine nicht-spaltende exakte Sequenz von $\widetilde{Q}_n$-Darstellungen
\[0\longrightarrow E_s\longrightarrow\widetilde{M}\longrightarrow C\longrightarrow 0\]
mit folgenden Eigenschaften:
\begin{enumerate}[i)]
\item $C$ ist exzeptionell.
\item $\Hom_{\widetilde{Q}_n}(C,E_s)=0=\Hom_{\widetilde{Q}_n}(E_s,C)$ und $\dim\Ext^1_{\widetilde{Q}_n}(C,E_s)=1$.
\end{enumerate}
\item Es existiert eine exakte Sequenz
\[0\longrightarrow E\longrightarrow M\longrightarrow U\longrightarrow0.\]
Dabei ist $E$ die projektiv-einfache $Q_n$-Darstellung und $U$ eine unzerlegbare $Q_n$-Darstel\-lung.
\end{enumerate}
\end{satz}
\begin{bew} (1.) Wir ersetzen im Beweis $\widetilde{Q_n}$ durch den Träger $T$ von $\widetilde{M}$ und beachten: Die volle Unterkategorie von $\fdar{\widetilde{Q}_n}$ bestehend aus allen $T$-Darstellungen ist abgeschlossen unter Unter- und Faktorobjekten sowie unter Erweiterungen. Also induziert der Vergißfunktor Isomorphismen zwischen den $\Ext^1$-Gruppen von $\fdar{T}$ und denen von $\fdar{\widetilde{Q}_n}$. Der Übergang zum Träger ändert also nichts an den Aussagen. Außerdem ist $\widetilde{M}$ eine exzeptionelle $T$-Darstellung (Teil (1) von \cref{baum:kron:prop3}).

Nach Teil (2) von \cref{baum:kron:prop3} existiert eine Senke $s$ von $T$ mit folgenden Eigenschaften: 
\begin{enumerate}[a)]
\item Es existiert genau ein Pfeil $\alpha:q\longrightarrow s$ mit Spitze $s$ in $T$.
\item $\widetilde{M}(q)\simeq k \simeq \widetilde{M}(s)$ und $\widetilde{M}(\alpha)\neq0$.
\end{enumerate}

Nach \cref{baum:kron:kor1} existiert eine nicht-spaltende exakte Sequenz
\[
\xi:\ 0\longrightarrow E_s\stackrel{\phi}{\longrightarrow} \widetilde{M}\longrightarrow C\longrightarrow 0
\]
von $T$-Darstellungen, wobei $C$ exzeptionell ist.

Wir verifizieren noch die Eigenschaften in $ii)$. Wegen $C(s)=0$ ist klar:
\[\Hom_K(C,E_s)=0=\Hom_K(E_s,C).\]
Es bleibt also, $\dim\Ext^1(C,E_s)=1$ zu zeigen. Sei dazu eine weitere nicht-spaltende exakte Sequenz
\[
\xi':\ 0\longrightarrow E_s\longrightarrow N\longrightarrow C\longrightarrow 0
\]
von $T$ gegeben. Da die Sequenz nicht spaltet, ist $N(\alpha)$ ein Isomorphismus. Damit verifiziert man, daß ein kommutatives Diagramm mit exakten Zeilen existiert:
\[
\begin{CD}
0@>>>E_s@>>> \widetilde{M}@>>> C@>>> 0\\
@.@V\lambda VV @VVV @| @.\\
0@>>>E_s@>>> N@>>> C@>>> 0.
\end{CD}
\]
Es folgt $\xi'=\lambda\xi$. Das zeigt die Behauptung.

(2.) Nach \cref{grund:cov:prop1} ist $F_\lambda$ exakt und erhält Unzerlegbarkeit. Anwenden von $F_\lambda$ auf eine Sequenz aus (1.) liefert die Behauptung.
\end{bew}

\section{Zur Erreichbarkeit exzeptioneller $Q_3$-Darstellungen der Dimension $\leq 76$}\label{baum:erreich}
\setcounter{zaehler}{0}
\numberwithin{zaehler}{section}
\begin{defn}\label{baum:erreich:def2} Sei $A$ eine endlichdimensionale $k$-Algebra. Die Klasse der erreichbaren $A$-Moduln wird induktiv wie folgt definiert:
\begin{itemize}
\item[$\cdot$] Jeder einfache $A$-Modul ist erreichbar.
\item[$\cdot$] Ein $A$-Modul $M$ der Länge $l>1$ ist erreichbar, falls $M$ unzerlegbar ist und falls ein erreichbarer Unter- oder Faktormodul der Länge $l-1$ von $M$ existiert.
\end{itemize}
\end{defn}

Gemäß \cref{baum:kron:satz1} besitzt jeder exzeptionelle $\Lambda_3$-Modul $M$ einen unzerlegbaren Unter- oder Faktormodul der Länge $l(M)-1$. Wir konnten die Frage, ob die exzeptionellen $\Lambda_3$-Moduln sogar  erreichbar sind, nicht beantworten, wollen aber mit Hilfe der universellen Überlagerung von $Q_3$ nachrechnen, daß die präprojektiven $\Lambda_3$-Moduln $P_i$, $i\leq4$ erreichbar sind. Der Modul $P_4$ ist ein Beispiel für einen $76$-dimensionalen erreichbaren Modul.

Für unsere Rechnungen fixieren wir eine Senke $s$ in der universellen Überlagerung $\widetilde{Q}_3$ von $Q_3$ und setzen wie im letzten Abschnitt $L_0:=\{s\}$ und für $i\geq1$
\[L_i:=\bigl\{ x\in(\widetilde{Q}_3)_0-\bigcup_{j=0}^{i-1}L_j\bigm| x \text{ ist ein Nachbar eines Punktes } y\in L_{i-1}\bigr\}.\]
Wir visualisieren -- Fahr/Ringel folgend (vgl. \cite{FR}) -- $\widetilde{Q}_3$ als konzentrische Kreise mit Mittelpunkt $s$. Auf dem $i$-ten Kreis liegen dabei die Punkte der Menge $L_i$. In Abbildung \ref{baum:erreich:abb1} sind die ersten vier Kreise dargestellt.

\numberwithin{figure}{section}
\setcounter{figure}{\value{zaehler}}
\renewcommand{\captionformat}{~ ~}
\renewcommand{\figurename}{Abbildung}
\renewcommand{\figureformat}{\figurename~\thefigure}
\begin{figure}
\caption{}
\label{baum:erreich:abb1}
\begin{center}
\begin{tikzpicture}[line width=0.5pt]
\tikzstyle{knt}=[circle, fill, inner sep=0.8pt]
\node at (0,0)[knt,label=below:\small$s$](x_01){};

\node at (90:1cm)[knt](x_1_1){};
\node at (210:1cm)[knt](x_1_2){};
\node at (330:1cm)[knt](x_1_3){};

\node at (60:2cm)[knt](x_2_1){};
\node at (120:2cm)[knt](x_2_2){};
\node at (180:2cm)[knt](x_2_3){};
\node at (240:2cm)[knt](x_2_4){};
\node at (300:2cm)[knt](x_2_5){};
\node at (360:2cm)[knt](x_2_6){};

\node at (45:3cm)[knt](x_3_1){};
\node at (75:3cm)[knt](x_3_2){};
\node at (105:3cm)[knt](x_3_3){};
\node at (135:3cm)[knt](x_3_4){};
\node at (165:3cm)[knt](x_3_5){};
\node at (195:3cm)[knt](x_3_6){};
\node at (225:3cm)[knt](x_3_7){};
\node at (255:3cm)[knt](x_3_8){};
\node at (285:3cm)[knt](x_3_9){};
\node at (315:3cm)[knt](x_3_10){};
\node at (345:3cm)[knt](x_3_11){};
\node at (15:3cm)[knt](x_3_12){};

\node at (37.5:4cm)[knt](x_4_1){};
\node at (52.5:4cm)[knt](x_4_2){};
\node at (67.5:4cm)[knt](x_4_3){};
\node at (82.5:4cm)[knt](x_4_4){};
\node at (97.5:4cm)[knt](x_4_5){};
\node at (112.5:4cm)[knt](x_4_6){};
\node at (127.5:4cm)[knt](x_4_7){};
\node at (142.5:4cm)[knt](x_4_8){};
\node at (157.5:4cm)[knt](x_4_9){};
\node at (172.5:4cm)[knt](x_4_10){};
\node at (187.5:4cm)[knt](x_4_11){};
\node at (202.5:4cm)[knt](x_4_12){};
\node at (217.5:4cm)[knt](x_4_13){};
\node at (232.5:4cm)[knt](x_4_14){};
\node at (247.5:4cm)[knt](x_4_15){};
\node at (262.5:4cm)[knt](x_4_16){};
\node at (277.5:4cm)[knt](x_4_7){};
\node at (292.5:4cm)[knt](x_4_18){};
\node at (307.5:4cm)[knt](x_4_19){};
\node at (322.5:4cm)[knt](x_4_20){};
\node at (337.5:4cm)[knt](x_4_21){};
\node at (352.5:4cm)[knt](x_4_22){};
\node at (7.5:4cm)[knt](x_4_23){};
\node at (22.5:4cm)[knt](x_4_24){};

\draw[dotted, thin] (0,0) circle (1cm);
\draw[dotted, thin] (0,0) circle (2cm);
\draw[dotted, thin] (0,0) circle (3cm);
\draw[dotted, thin] (0,0) circle (4cm);

\draw[<-,shorten >= 3pt, shorten <= 3pt] (0,0) to (90:1cm);
\draw[<-,shorten >= 3pt, shorten <= 3pt] (0,0) to (210:1cm);
\draw[<-,shorten >= 3pt, shorten <= 3pt] (0,0) to (330:1cm);

\draw[->,shorten >= 3pt, shorten <= 3pt] (90:1cm) to (60:2cm);
\draw[->,shorten >= 3pt, shorten <= 3pt] (90:1cm) to (120:2cm);
\draw[->,shorten >= 3pt, shorten <= 3pt] (210:1cm) to (180:2cm);
\draw[->,shorten >= 3pt, shorten <= 3pt] (210:1cm) to (240:2cm);
\draw[->,shorten >= 3pt, shorten <= 3pt] (330:1cm) to (300:2cm);
\draw[->,shorten >= 3pt, shorten <= 3pt] (330:1cm) to (360:2cm);

\draw[<-,shorten >= 3pt, shorten <= 3pt] (60:2cm) to (45:3cm);
\draw[<-,shorten >= 3pt, shorten <= 3pt] (60:2cm) to (75:3cm);
\draw[<-,shorten >= 3pt, shorten <= 3pt] (120:2cm) to (105:3cm);
\draw[<-,shorten >= 3pt, shorten <= 3pt] (120:2cm) to (135:3cm);
\draw[<-,shorten >= 3pt, shorten <= 3pt] (180:2cm) to (165:3cm);
\draw[<-,shorten >= 3pt, shorten <= 3pt] (180:2cm) to (195:3cm);
\draw[<-,shorten >= 3pt, shorten <= 3pt] (240:2cm) to (225:3cm);
\draw[<-,shorten >= 3pt, shorten <= 3pt] (240:2cm) to (255:3cm);
\draw[<-,shorten >= 3pt, shorten <= 3pt] (300:2cm) to (285:3cm);
\draw[<-,shorten >= 3pt, shorten <= 3pt] (300:2cm) to (315:3cm);
\draw[<-,shorten >= 3pt, shorten <= 3pt] (360:2cm) to (345:3cm);
\draw[<-,shorten >= 3pt, shorten <= 3pt] (360:2cm) to (15:3cm);

\draw[->,shorten >= 3pt, shorten <= 3pt] (45:3cm) to (37.5:4cm);
\draw[->,shorten >= 3pt, shorten <= 3pt] (45:3cm) to (52.5:4cm);
\draw[->,shorten >= 3pt, shorten <= 3pt] (75:3cm) to (67.5:4cm);
\draw[->,shorten >= 3pt, shorten <= 3pt] (75:3cm) to (82.5:4cm);
\draw[->,shorten >= 3pt, shorten <= 3pt] (105:3cm) to (97.5:4cm);
\draw[->,shorten >= 3pt, shorten <= 3pt] (105:3cm) to (112.5:4cm);
\draw[->,shorten >= 3pt, shorten <= 3pt] (135:3cm) to (127.5:4cm);
\draw[->,shorten >= 3pt, shorten <= 3pt] (135:3cm) to (142.5:4cm);
\draw[->,shorten >= 3pt, shorten <= 3pt] (165:3cm) to (157.5:4cm);
\draw[->,shorten >= 3pt, shorten <= 3pt] (165:3cm) to (172.5:4cm);
\draw[->,shorten >= 3pt, shorten <= 3pt] (195:3cm) to (187.5:4cm);
\draw[->,shorten >= 3pt, shorten <= 3pt] (195:3cm) to (202.5:4cm);
\draw[->,shorten >= 3pt, shorten <= 3pt] (225:3cm) to (217.5:4cm);
\draw[->,shorten >= 3pt, shorten <= 3pt] (225:3cm) to (232.5:4cm);
\draw[->,shorten >= 3pt, shorten <= 3pt] (255:3cm) to (247.5:4cm);
\draw[->,shorten >= 3pt, shorten <= 3pt] (255:3cm) to (262.5:4cm);
\draw[->,shorten >= 3pt, shorten <= 3pt] (285:3cm) to (277.5:4cm);
\draw[->,shorten >= 3pt, shorten <= 3pt] (285:3cm) to (292.5:4cm);
\draw[->,shorten >= 3pt, shorten <= 3pt] (315:3cm) to (307.5:4cm);
\draw[->,shorten >= 3pt, shorten <= 3pt] (315:3cm) to (322.5:4cm);
\draw[->,shorten >= 3pt, shorten <= 3pt] (345:3cm) to (337.5:4cm);
\draw[->,shorten >= 3pt, shorten <= 3pt] (345:3cm) to (352.5:4cm);
\draw[->,shorten >= 3pt, shorten <= 3pt] (15:3cm) to (7.5:4cm);
\draw[->,shorten >= 3pt, shorten <= 3pt] (15:3cm) to (22.5:4cm);
\end{tikzpicture}
\end{center}
\end{figure}
Sei $\widetilde{\sigma}^-:\fdar{\widetilde{Q}_3}\longrightarrow \fdar{\widetilde{Q}_3^{op}}$ der in Abschnitt \ref{baum:exk} definierte Spiegelungsfunktor bezüglich aller Quellen von $\widetilde{Q_3}$. Den Spiegelungsfunktor bezüglich aller Quellen von $\fdar{\widetilde{Q}_3^{op}}$ in der anderen Richtung bezeichnen wir ebenfalls mit $\widetilde{\sigma}^-$.

Wir betrachten für $i\geq0$ die im letzten Abschnitt definierte Liftung $\widetilde{P}_i$ der präprojektiven $Q_3$-Darstellung $P_i$. Für gerades $i$ ist also $\widetilde{P}_i$ eine $\widetilde{Q}_3$-Darstellung, für ungerades $i$ eine $\widetilde{Q}_3^{op}$-Darstellung. Es gilt \[\widetilde{P}_0=E_s\text{ und } \widetilde{P}_i=\widetilde{\sigma}^-(\widetilde{P}_{i-1})\text{ für }i\geq1.\]
Anwenden des Pushdown-Funktors auf $\widetilde{P}_i$ liefert (bis auf eventuelles Vertauschen von Quelle und Senke von $Q_3$) die Darstellung $P_i$. Da der Pushdown-Funktor Exaktheit, Unzerlegbarkeit und Einfachheit erhält, genügt es also, die Erreichbarkeit der $\widetilde{P}_i$ zu zeigen.

Stets sei $T_i$ der Träger von $\widetilde{P}_i$. Wir haben in \cref{baum:kron:lemma2} gesehen, daß $T_i=\bigcup_{j=0}^iL_i$ ist und haben dort auch die Werte von $\widetilde{P}_i$ auf den beiden äußeren konzentrischen Kreise $L_i$ und $L_{i-1}$ berechnet: Für alle $x\in L_i\cup L_{i-1}$ ist $\widetilde{P}_i(x)=k$, für alle Pfeile $\alpha:x\longrightarrow y$, $y\in L_i$, $x\in L_{i-1}$ ist $\widetilde{P}_i(\alpha)=id_k$ (beachte, daß $L_i$ nur Senken von $T_i$ enthält, da $T_i$ für gerades $i$ ein Unterköcher von $\widetilde{Q}_3$, für ungerades $i$ von $\widetilde{Q}_3^{op}$ ist). 

Betrachte für $i\geq0$ den Funktor
\[F_i:=\prod_{\begin{subarray}{c}q\in T_i\\q\text{ Quelle in $T_i$}\end{subarray}}\sigma_q^-:\fdar{T_i}\longrightarrow\fdar{T_i^{op}}.\]
$F_i$ ist also der Spiegelungsfunktor bezüglich aller Quellen in $T_i$ und hängt nicht von der Wahl einer Aufzählung der Quellen ab. 

Für den Rest dieses Abschnitts sei $M_i$ die $T_i$-Darstellung, die aus $\widetilde{P}_i$ entsteht durch Annullieren der Werte auf dem äußersten Kreis $L_i$. Es ist also
\[M_i=\widetilde{P}_i/\bigoplus_{x\in L_i}E_x,\]
und wir können $M_i$ als $T_{i-1}^{op}$-Darstellung auffassen. Wir benötigen das folgende Lemma:
\begin{lemma}\label{baum:erreich:lemma1} Sei $i\geq1$. Seien $x_1,\ldots,x_t$ die Elemente von $L_i$. Dann existiert für $j=1,\ldots,t$ eine exakte Sequenz
\[0\longrightarrow E_{x_j}\longrightarrow X_j\longrightarrow X_{j+1}\longrightarrow 0\]
von $T_i$-Darstellungen mit $X_1=\widetilde{P}_{i}$ und $X_{t+1}=M_i$. Dabei gilt:
\begin{enumerate}
\item Für alle $j=1,\ldots,t$ ist $X_j$ ist exzeptionell und es gilt
\[\begin{array}{ll}
\dim X_j(t)=1&\text{für alle } t\in L_{i-1}\cup \bigl( L_i-\{x_1,\ldots,x_{j-1}\}\bigr),\\
\dim X_j(t)=\dim\widetilde{P}_i(t)\geq2&\text{für alle } t\in L_{i-2}.
\end{array}
\]
Dabei sei $L_{-1}:=\emptyset$.
\item $F_{i-1}(\widetilde{P}_{i-1})\simeq M_i$ (als $T_{i-1}^{op}$-Darstellungen).
\item Ist $M_i$ erreichbar, so ist auch $\widetilde{P}_i$ erreichbar.
\end{enumerate}
\end{lemma}
\begin{bew} Um den Leser mit unseren Notationen vertraut zu machen, geben wir einen ausführlichen Beweis dieses Lemmas. Der Leser sollte sich alle Aussagen anhand von Beispielen für kleine $i$ verdeutlichen.

(1.) Setze $X_1:=\widetilde{P}_{i}$. $X_1$ ist exzeptionell, und es gilt nach den Bemerkungen vor diesem Lemma
\[\dim X_1(t)=1\text{ für alle } t\in L_i\cup L_{i-1}.\]
Weiter ist $\dim \widetilde{P}_{i}(t)\geq2$ für alle $t\in L_{i-2}$: Wir dürfen $i\geq2$ annehmen. Es ist $\widetilde{P}_i=\widetilde{\sigma}^-(\widetilde{P}_{i-1})$, und alle Punkte von $L_{i-2}$ sind Quellen des Trägers $T_{i-1}=\bigcup_{k=0}^{i-1}L_k$ von $\widetilde{P}_{i-1}$. Sei $t\in L_{i-2}$, und seien $u_1,u_2,u_3$ die Nachfolger von $t$ in $T_{i-1}$. Per Definition des Spiegelungsfunktors $\widetilde{\sigma}^-$ bezüglich aller Quellen gilt:
\[\dim \widetilde{P}_i(t)=\dim \widetilde{\sigma}^-(\widetilde{P}_{i-1})(t)=\sum_{k=1}^3 \dim \widetilde{P}_{i-1}(u_k)-\widetilde{P}_{i-1}(t)\geq 3-1 =2.\]
$X_2$ konstruiert man folgendermaßen: Der eindeutig bestimmte Pfeil $\alpha:y\longrightarrow x_1$ mit Spitze $x_1$ wird von $X_1:=\widetilde{P}_i$ durch die Identität auf $k$ dargestellt. Mit \cref{baum:kron:kor1} folgt, daß eine exakte Sequenz
\[0\longrightarrow E_{x_1}\longrightarrow X_1\longrightarrow X_2\longrightarrow 0\]
existiert, wobei $X_2$ exzeptionell ist. Die geforderten Gleichungen in (1.) verifiziert man sofort. Die Konstruktion der $X_j$ für $j\geq3$ und der zugehörigen exakten Sequenzen sowie die Verifikation der geforderten Gleichungen erfolgen durch eine offensichtliche Induktion nach $j$. Außerdem ist klar, daß $X_{t+1}=M_i$ ist.

(2.) Wir nehmen an, daß $i$ ungerade ist. Also ist $\widetilde{P}_{i-1}\in\fdar{\widetilde{Q}_n}$. Sei $\iota:\fdar{\widetilde{Q}_n}\longrightarrow\fdar{T_{i-1}}$ die Einschränkung. Bezeichne die Einschränkung  $\fdar{\widetilde{Q}_n^{op}}\longrightarrow\fdar{T_{i-1}^{op}}$ mit $\iota'$. Es gelten die Gleichungen $M_i=\iota'(\widetilde{P}_i)$ und $F_{i-1}\circ\iota=\iota'\circ\widetilde{\sigma}^-$ (beachte, daß $T_{i-1}$ unter Nachfolgern abgeschlossen ist). Es folgt 
\[M_i=\iota'(\widetilde{P}_i)=\iota'\circ\widetilde{\sigma}^-(\widetilde{P}_{i-1})=F_{i-1}\circ\iota(\widetilde{P}_{i-1}).\]
Das zeigt die Behauptung.

(3.) Folgt sofort aus der Definition der Erreichbarkeit mit Hilfe der in (1.) konstruierten exakten Sequenzen sowie aus der Unzerlegbarkeit aller $X_i$.
\end{bew}
Wir können unsere Rechnungen weiter vereinfachen, indem wir "`an weniger Quellen spiegeln"'. Genauer gilt folgendes Reduktionslemma:
\begin{lemma}\label{baum:erreich:lemma2}\ 
\begin{enumerate}
\item $F_{i}(M_i)$ ist exzeptionell für alle $i\geq1$.
\item  Für alle $i\geq2$ gilt: Ist $F_{i}(M_i)$ erreichbar, so ist auch $\widetilde{P}_{i+1}$ erreichbar.
\end{enumerate}
\end{lemma}

Diese Aussage basiert auf einem weiteren Lemma:
\begin{lemma}\label{baum:erreich:lemma3}
Sei $K$ ein endlicher zykelloser Köcher. Sei $q$ eine Quelle und $s$ eine Senke in $K$, so daß genau einen Pfeil $\alpha:q\longrightarrow s$ von $q$ nach $s$ existiert. Sei $U$ die $K$-Darstellung mit $U(q)=U(s)=k$, $U(\alpha)=id_k$ und $U(x)=0$ falls $x\not\in\{q,s\}$. Sei eine exakte Sequenz
\[0\longrightarrow U\longrightarrow M\longrightarrow C\longrightarrow 0\]
von $K$-Darstellungen gegeben. Dabei seien die folgenden Bedingungen erfüllt:
\begin{enumerate}[i)]
\item $\Ext^1_K(M,M)=0$.
\item $\dim M(q)=1$ und $\dim M(s)\geq2$.
\end{enumerate}
Dann existieren eine $K$-Darstellung $C'$ und exakte Sequenzen
\begin{gather*}
0\longrightarrow E_s\longrightarrow M\longrightarrow C'\longrightarrow 0,\\
0\longrightarrow C\longrightarrow C'\longrightarrow E_q\longrightarrow 0,
\end{gather*}
so daß $\End_K(C')$ eine Unteralgebra von $\End_K(C)$ ist. Insbesondere gilt:
\begin{enumerate}
\item Ist $C$ unzerlegbar, so ist auch $C'$ unzerlegbar.
\item Sei $C$ erreichbar. Dann ist auch $M$ erreichbar.
\end{enumerate}
\end{lemma}
\begin{bew} Wir verwenden die Notationen aus Abschnitt \ref{grund:var} für Varietäten von Köcherdarstellungen. Sei $\dimv M=\underline{d}=(d(x))_{x\in K_0}$. Ohne Einschränkung gelte $M=\operatorname{rep} m$ für einen Punkt $m\in \operatorname{Rep}_K(\underline{d})$. Wir zeigen zunächst: Es existiert ein Punkt $n$ im Orbit $\mathcal{O}(m)$ von $m$, so daß gilt:
\begin{itemize}
\item[$\cdot$] $n_{\beta}=m_{\beta}$ für alle $\beta\in K_1-\{\alpha\}$,
\item[$\cdot$] $(m_{\alpha},n_{\alpha})$ ist linear unabhängig in $k^{d(s)\times d(q)}=k^{d(s)}$.
\end{itemize}
Betrachte dazu die Projektion
\[\pi:\operatorname{Rep}_K(\underline{d})\longrightarrow\prod_{\beta\in K_1-\{\alpha\}}k^{d(s\beta)\times d(n\beta)}.\]
Die Faser $X_m:=\pi^{-1}(\pi(m))\simeq k^{d(s)\times d(q)}$ ist irreduzibel. Die Menge
\[U_m:=\{x\in X_m\mid (x_\alpha,m_\alpha)\text{ ist linear unabhängig in } k^{d(s)\times d(q)}\}\]
ist nichtleer und offen in $X_m$ (beachte: Nach Voraussetzung ist $\dim k^{d(s)\times d(q)}\geq2$). Wegen $\Ext^1_K(M,M)=0$ ist $\mathcal{O}(m)$ offen in $\operatorname{Rep}_K(\underline{d})$, und wir erhalten aus der Irreduzibilität von $X_m$
\[(X_m\cap\mathcal{O}(m))\cap U_m\neq\emptyset.\]
Jedes $n\in X_m$ erfüllt die geforderten Bedingungen. Wähle so ein $n$ und setze $N:=\operatorname{rep} n$. Dann ist $N\simeq M$ (denn $n$ und $m$ liegen im selben Orbit). Setze
\[v:=n_\alpha\in k^{d(s)\times d(q)}=k^{d(s)}.\]
Es existiert ein kommutatives Diagramm mit exakten Zeilen:
\[
\begin{CD}
0@>>>U@>>>M@>>>C@>>>0\\
@.@|@V\text{\rotatebox{-90}{$\sim$}}VV@V\text{\rotatebox{-90}{$\sim$}}VV@.\\
0@>>>U@>>>N@>>>D@>>>0.\\
\end{CD}
\]
Die Faktordarstellung $D$ von $N$ ist aufgrund der vorausgesetzten Bedingungen an $U$ bis auf Isomorphie gegeben durch
\[D(q)=0,\ D(s)=k^{d(s)}/\langle v\rangle,\ D(x)=k^{d(x)}\text{ für alle }x\neq q,s.\]
Betrachte die Faktordarstellung $C'$ von $M$ gegeben durch
\[C'(s)=k^{d(s)}/\langle v\rangle,\ C'(x)=k^{d(x)}\text{ für alle }x\neq s.\]
Es existiert also eine exakte Sequenz
\[0\longrightarrow E_s\longrightarrow M\longrightarrow C'\longrightarrow 0.\]
Außerdem ist offenbar $D\subseteq C'$ eine Unterdarstellung von $C'$ mit einfachem Quotienten $E_q$. Wegen $D\simeq C$ erhält man also eine exakte Sequenz
\[0\longrightarrow C\longrightarrow C'\longrightarrow E_q\longrightarrow 0.\]
Der Pfeil $\alpha:q\longrightarrow s$ wird von $C'$ dargestellt durch die lineare Abbildung
\[C'(\alpha):k\stackrel{m_\alpha}{\longrightarrow}k^{d(s)}\stackrel{\pi}{\longrightarrow} k^{d(s)}/\langle v\rangle.\]
Hier ist $\pi$ die kanonische Projektion. Da $(v, m_\alpha)$ per Konstruktion von $v$ linear unabhängig in $k^{d(s)}$ ist, muß $C'(\alpha)$ injektiv sein. Damit verifiziert man leicht, daß $\End_K(C')$ eine Unteralgebra von $\End_K(C)$ ist (verwende etwa die duale Version von \cref{baum:kron:lemma3}, Aussage (1)). Es bleiben noch die Aussagen (1.) und (2.) zu zeigen.

(1.) Folgt aus der Tatsache, daß Unteralgebren von lokalen endlichdimensionalen Algebren wieder lokal sind.

(2.) Es ist nur zu zeigen, daß $M$ unzerlegbar ist. $C$ ist unzerlegbar nach Voraussetzung und die Unzerlegbarkeit von $U$ ist klar. Außerdem ist $\Hom_K(U,C)=0$. Damit folgt die Unzerlegbarkeit von $M$ zum Beispiel aus \cref{baum:erreich:lemma6}.
\end{bew}

\begin{bew}[von \cref{baum:erreich:lemma2}] (1.) Da $M_i$ nach \cref{baum:erreich:lemma1} eine exzeptionelle Darstellung von $T_{i}$ ist, muß auch $F_i(M_i)$ als gespiegelte Darstellung exzeptionell sein.

(2.) Seien $x_1,\ldots,x_t$ die Elemente von $L_{i}$. Betrachte für $j=1,\ldots,t$ die in \cref{baum:erreich:lemma1} angegebenen exakten Sequenzen
\counterwithout{equation}{chapter}
\setcounter{equation}{0}
\begin{equation}\label{baum:erreich:lemma2:gl1}
0\longrightarrow E_{x_j}\longrightarrow X_j\longrightarrow X_{j+1}\longrightarrow 0
\end{equation}
von $T_i$-Darstellungen. Es ist also $X_1=\widetilde{P}_{i}$ und $X_{t+1}=M_i$, und alle $X_j$ sind exzeptionell. Für alle $j$ sei $\alpha_j:s_j\longrightarrow x_j$ der eindeutig bestimmte Pfeil in $T_i$ mit Spitze $x_j$ und Start $s_j\in L_{i-1}$. Nach \cref{baum:erreich:lemma1} gilt
\begin{equation}\label{baum:erreich:lemma2:gl3}
\dim X_j(x_j)=\dim X_j(s_j)=1\text{ und }\dim X_j(t)=\dim\widetilde{P}_i(t)\geq2\text{ für alle } t\in L_{i-2}.
\end{equation}
Da die $X_l$ keinen injektiv-einfachen direkten Summanden haben, liefert Anwenden von $F_{i}$ auf die Sequenz (\ref{baum:erreich:lemma2:gl1}) eine exakte Sequenz
\begin{equation}\label{baum:erreich:lemma2:gl2}
0\longrightarrow F_i(E_{x_j})\longrightarrow F_i(X_j)\longrightarrow F_i(X_{j+1})\longrightarrow 0
\end{equation}
von $T_i^{op}$-Darstellungen. Sei $j\in\{1,\ldots,t\}$ gegeben. Wir zeigen, daß die Sequenz (\ref{baum:erreich:lemma2:gl2}) die Voraussetzungen von \cref{baum:erreich:lemma3} erfüllt. Per Definition von $F_i$ gilt
\[F_i(E_{x_j})(x_j)\simeq k\simeq F_i(E_{x_j})(s_j),\ F_i(E_{x_j})(y)=0\text{ falls }y\not\in\{ x_j,s_j\},\ F_i(E_{x_j})(\alpha_j)=id_k.\]
Weiter gilt:
\begin{enumerate}[i)]
\item Alle $X_{l}$ sind exzeptionelle Darstellungen. Damit sind auch alle $F_iX_{l}$ exzeptionelle Darstellungen, und es folgt $\Ext^1_{T_i^{op}}(F_i(X_j),F_i(X_j))=0$.
\item Wir zeigen $\dim F_i(X_j)(x_j)=1$ und $\dim F_i(X_j)(s_j)\geq2$ und verwenden Gleichung (\ref{baum:erreich:lemma2:gl3}). Per Definition des Spiegelungsfunktors $F_i$ gilt
\[\dim F_i(X_j)(s_j)=\sum_{\text{Nachbarn $n$ von $s_j$}}\dim X_j(n)-\dim X_j(s_j)\geq 2+1-1=2\]
und
\[\dim F_i(X_j)(x_j)=\dim X_j(x_j)=1.\]
\item $F_i(X_{j+1})$ ist unzerlegbar, da nach (i) alle $F_iX_{l}$ exzeptionell sind.
\end{enumerate}
Also sind die Voraussetzungen von \cref{baum:erreich:lemma3} erfüllt. Wir können die exakte Sequenz (\ref{baum:erreich:lemma2:gl2}) also aufspalten in zwei exakte Sequenzen
\[0\longrightarrow E_{x_j}\longrightarrow F_i(X_j)\longrightarrow C'_{j+1}\longrightarrow 0\]
und
\[0\longrightarrow F_i(X_{j+1})\longrightarrow C'_{j+1}\longrightarrow E_{s_j}\longrightarrow 0.\]
Nach Aussage (2.) in \cref{baum:erreich:lemma3} gilt: Ist $F_i(X_{j+1})$ erreichbar, so ist auch $F_i(X_{j})$ erreichbar.

Da nach Voraussetzung $F_i(M_i)=F_i(X_{t+1})$ erreichbar ist, ist auch $F_i(X_1)=F_i(\widetilde{P}_i)$ nach dem eben Gezeigten erreichbar. Nun ist nach \cref{baum:erreich:lemma1} $F_i(\widetilde{P}_i)=M_{i+1}$ und -- wieder nach \cref{baum:erreich:lemma1} -- folgt aus der Erreichbarkeit von $M_{i+1}$ die von $\widetilde{P}_{i+1}$.
\end{bew}

\textbf{$\mathbf{\widetilde{P}_0,\widetilde{P}_1,\widetilde{P}_2,\widetilde{P}_3}$ sind erreichbar}\newline
$\widetilde{P}_0$ ist erreichbar als einfache Darstellung. $\widetilde{P}_1$ ist die unzerlegbare Darstellung zum Dimensionsvektor
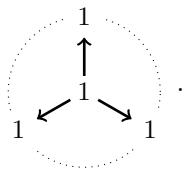
\begin{center}
\parbox[c]{2.3cm}{\begin{tikzpicture}[line width=1pt]
\node at (0,0)(x_01){\small$1$};

\draw[dotted, thin] (0,0) circle (1cm);

\draw[->,shorten >= 8pt, shorten <= 6pt] (0,0) to (90:1cm);
\draw[->,shorten >= 8pt, shorten <= 6pt] (0,0) to (210:1cm);
\draw[->,shorten >= 8pt, shorten <= 6pt] (0,0) to (330:1cm);

\node[fill=white] at (90:1cm)(x_1_1){\small$1$};
\node[fill=white] at (210:1cm)(x_1_2){\small$1$};
\node[fill=white] at (330:1cm)(x_1_3){\small$1$};

\end{tikzpicture}}.\end{center}
Diese Darstellung ist offenbar erreichbar. $\widetilde{P}_2$ ist die unzerlegbare Darstellung zum Dimensionsvektor
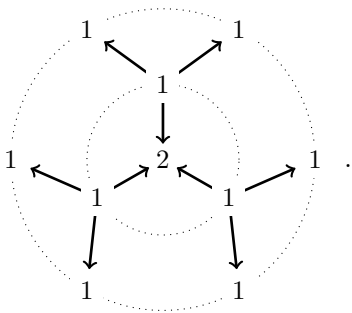
\begin{center}
\parbox[c]{4.6cm}{\begin{tikzpicture}[line width=1pt]
\draw[dotted, thin] (0,0) circle (1cm);
\draw[dotted, thin] (0,0) circle (2cm);

\draw[<-,shorten >= 2pt, shorten <= 6pt] (0,0) to (90:1cm);
\draw[<-,shorten >= 2pt, shorten <= 6pt] (0,0) to (210:1cm);
\draw[<-,shorten >= 2pt, shorten <= 6pt] (0,0) to (330:1cm);

\draw[->,shorten >= 8pt, shorten <= 2pt] (90:1cm) to (60:2cm);
\draw[->,shorten >= 8pt, shorten <= 2pt] (90:1cm) to (120:2cm);
\draw[->,shorten >= 8pt, shorten <= 2pt] (210:1cm) to (180:2cm);
\draw[->,shorten >= 8pt, shorten <= 2pt] (210:1cm) to (240:2cm);
\draw[->,shorten >= 8pt, shorten <= 2pt] (330:1cm) to (300:2cm);
\draw[->,shorten >= 8pt, shorten <= 2pt] (330:1cm) to (360:2cm);

\node at (0,0)(x_01){\small$2$};

\node[fill=white] at (90:1cm)(x_1_1){\small$1$};
\node[fill=white] at (210:1cm)(x_1_2){\small$1$};
\node[fill=white] at (330:1cm)(x_1_3){\small$1$};

\node[fill=white] at (60:2cm)(x_2_1){\small$1$};
\node[fill=white] at (120:2cm)(x_2_2){\small$1$};
\node[fill=white] at (180:2cm)(x_2_3){\small$1$};
\node[fill=white] at (240:2cm)(x_2_4){\small$1$};
\node[fill=white] at (300:2cm)(x_2_5){\small$1$};
\node[fill=white] at (360:2cm)(x_2_6){\small$1$};
\end{tikzpicture}}.\end{center}
Die sechs äußeren Einsen eliminiert man durch Einbetten von Einfachen. Die so entstehende Darstellung $M_2$ (also die Einschränkung von $\widetilde{P}_2$ auf den ersten Kreis) ist die unzerlegbare 5-dimensionale Darstellung des Dynkin-Köchers $D_4$, also erreichbar (vgl. auch \cref{baum:erreich:lemma1}).

Der Dimensionsvektor der Darstellung $\widetilde{P}_3$ ist gegeben durch
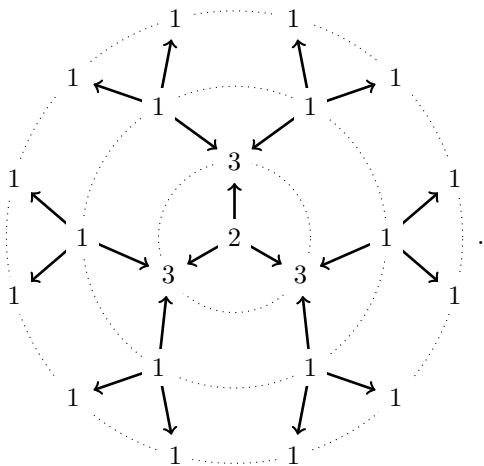
\begin{center}
\parbox[c]{6.3cm}{\begin{tikzpicture}[line width=1pt]
\draw[dotted, thin] (0,0) circle (1cm);
\draw[dotted, thin] (0,0) circle (2cm);
\draw[dotted, thin] (0,0) circle (3cm);

\draw[->,shorten >= 8pt, shorten <= 5pt] (0,0) to (90:1cm);
\draw[->,shorten >= 8pt, shorten <= 5pt] (0,0) to (210:1cm);
\draw[->,shorten >= 8pt, shorten <= 5pt] (0,0) to (330:1cm);

\draw[<-,shorten >= 3pt, shorten <= 8pt] (90:1cm) to (60:2cm);
\draw[<-,shorten >= 3pt, shorten <= 8pt] (90:1cm) to (120:2cm);
\draw[<-,shorten >= 3pt, shorten <= 8pt] (210:1cm) to (180:2cm);
\draw[<-,shorten >= 3pt, shorten <= 8pt] (210:1cm) to (240:2cm);
\draw[<-,shorten >= 3pt, shorten <= 8pt] (330:1cm) to (300:2cm);
\draw[<-,shorten >= 3pt, shorten <= 8pt] (330:1cm) to (360:2cm);

\draw[->,shorten >= 8pt, shorten <= 3pt] (60:2cm) to (45:3cm);
\draw[->,shorten >= 8pt, shorten <= 3pt] (60:2cm) to (75:3cm);
\draw[->,shorten >= 8pt, shorten <= 3pt] (120:2cm) to (105:3cm);
\draw[->,shorten >= 8pt, shorten <= 3pt] (120:2cm) to (135:3cm);
\draw[->,shorten >= 8pt, shorten <= 3pt] (180:2cm) to (165:3cm);
\draw[->,shorten >= 8pt, shorten <= 3pt] (180:2cm) to (195:3cm);
\draw[->,shorten >= 8pt, shorten <= 3pt] (240:2cm) to (225:3cm);
\draw[->,shorten >= 8pt, shorten <= 3pt] (240:2cm) to (255:3cm);
\draw[->,shorten >= 8pt, shorten <= 3pt] (300:2cm) to (285:3cm);
\draw[->,shorten >= 8pt, shorten <= 3pt] (300:2cm) to (315:3cm);
\draw[->,shorten >= 8pt, shorten <= 3pt] (360:2cm) to (345:3cm);
\draw[->,shorten >= 8pt, shorten <= 3pt] (360:2cm) to (15:3cm);

\node[fill=white] at (0,0)(x_01){\small$2$};

\node[fill=white] at (90:1cm)(x_1_1){\small$3$};
\node[fill=white] at (210:1cm)(x_1_2){\small$3$};
\node[fill=white] at (330:1cm)(x_1_3){\small$3$};

\node[fill=white] at (60:2cm)(x_2_1){\small$1$};
\node[fill=white] at (120:2cm)(x_2_2){\small$1$};
\node[fill=white] at (180:2cm)(x_2_3){\small$1$};
\node[fill=white] at (240:2cm)(x_2_4){\small$1$};
\node[fill=white] at (300:2cm)(x_2_5){\small$1$};
\node[fill=white] at (360:2cm)(x_2_6){\small$1$};
\node[fill=white] at (45:3cm)(x_3_1){\small$1$};
\node[fill=white] at (75:3cm)(x_3_2){\small$1$};
\node[fill=white] at (105:3cm)(x_3_3){\small$1$};
\node[fill=white] at (135:3cm)(x_3_4){\small$1$};
\node[fill=white] at (165:3cm)(x_3_5){\small$1$};
\node[fill=white] at (195:3cm)(x_3_6){\small$1$};
\node[fill=white] at (225:3cm)(x_3_7){\small$1$};
\node[fill=white] at (255:3cm)(x_3_8){\small$1$};
\node[fill=white] at (285:3cm)(x_3_9){\small$1$};
\node[fill=white] at (315:3cm)(x_3_10){\small$1$};
\node[fill=white] at (345:3cm)(x_3_11){\small$1$};
\node[fill=white] at (15:3cm)(x_3_12){\small$1$};
\end{tikzpicture}}.\end{center}

Zur Erreichbarkeit von $\widetilde{P}_3$ genügt es nach \cref{baum:erreich:lemma2} zu zeigen, daß $F_2(M_2)$ erreichbar ist. $M_2$ ist die unzerlegbare $T_2$-Darstellung zum Dimensionsvektor
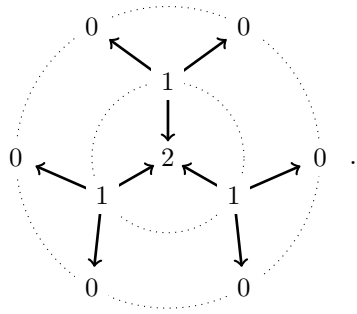
\begin{center}
\parbox[c]{4.6cm}{\begin{tikzpicture}[line width=1pt]
\draw[dotted, thin] (0,0) circle (1cm);
\draw[dotted, thin] (0,0) circle (2cm);

\draw[<-,shorten >= 2pt, shorten <= 6pt] (0,0) to (90:1cm);
\draw[<-,shorten >= 2pt, shorten <= 6pt] (0,0) to (210:1cm);
\draw[<-,shorten >= 2pt, shorten <= 6pt] (0,0) to (330:1cm);

\draw[->,shorten >= 8pt, shorten <= 2pt] (90:1cm) to (60:2cm);
\draw[->,shorten >= 8pt, shorten <= 2pt] (90:1cm) to (120:2cm);
\draw[->,shorten >= 8pt, shorten <= 2pt] (210:1cm) to (180:2cm);
\draw[->,shorten >= 8pt, shorten <= 2pt] (210:1cm) to (240:2cm);
\draw[->,shorten >= 8pt, shorten <= 2pt] (330:1cm) to (300:2cm);
\draw[->,shorten >= 8pt, shorten <= 2pt] (330:1cm) to (360:2cm);\node at (0,0)(x_01){\small$2$};

\node[fill=white] at (90:1cm)(x_1_1){\small$1$};
\node[fill=white] at (210:1cm)(x_1_2){\small$1$};
\node[fill=white] at (330:1cm)(x_1_3){\small$1$};

\node[fill=white] at (60:2cm)(x_2_1){\small$0$};
\node[fill=white] at (120:2cm)(x_2_2){\small$0$};
\node[fill=white] at (180:2cm)(x_2_3){\small$0$};
\node[fill=white] at (240:2cm)(x_2_4){\small$0$};
\node[fill=white] at (300:2cm)(x_2_5){\small$0$};
\node[fill=white] at (360:2cm)(x_2_6){\small$0$};\end{tikzpicture}}.\end{center}
Anwenden von $F_2$ bedeutet Spiegeln an allen Quellen von $T_2$, und man erhält wieder die (erreichbare) unzerlegbare $5$-dimensionale Darstellung des Dynkin-Köchers $D_4$. 

\textbf{$\mathbf{\widetilde{P}_4}$ ist erreichbar}\newline
Der Dimensionsvektor von $\widetilde{P}_4$ ist gegeben durch
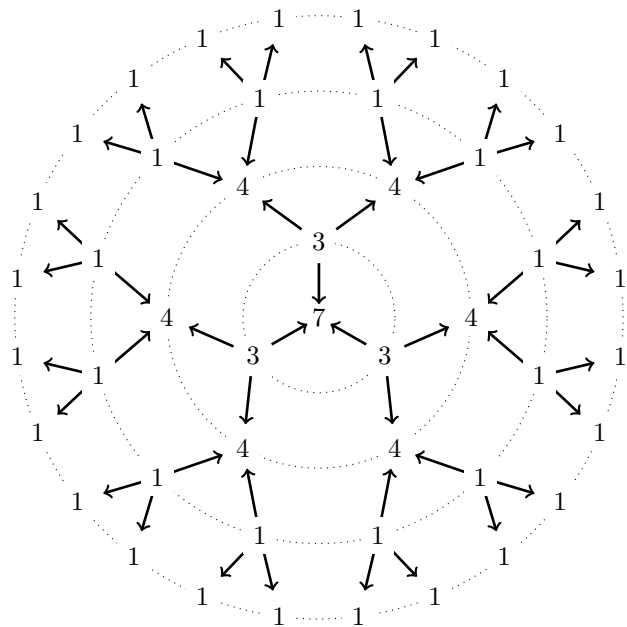
\begin{center}
\parbox[c]{9cm}{\begin{tikzpicture}[line width=1pt]
\draw[dotted, thin] (0,0) circle (1cm);
\draw[dotted, thin] (0,0) circle (2cm);
\draw[dotted, thin] (0,0) circle (3cm);
\draw[dotted, thin] (0,0) circle (4cm);

\draw[<-,shorten >= 8pt, shorten <= 5pt] (0,0) to (90:1cm);
\draw[<-,shorten >= 8pt, shorten <= 5pt] (0,0) to (210:1cm);
\draw[<-,shorten >= 8pt, shorten <= 5pt] (0,0) to (330:1cm);

\draw[->,shorten >= 9pt, shorten <= 8pt] (90:1cm) to (60:2cm);
\draw[->,shorten >= 9pt, shorten <= 8pt] (90:1cm) to (120:2cm);
\draw[->,shorten >= 9pt, shorten <= 8pt] (210:1cm) to (180:2cm);
\draw[->,shorten >= 9pt, shorten <= 8pt] (210:1cm) to (240:2cm);
\draw[->,shorten >= 9pt, shorten <= 8pt] (330:1cm) to (300:2cm);
\draw[->,shorten >= 9pt, shorten <= 8pt] (330:1cm) to (360:2cm);

\draw[<-,shorten >= 3pt, shorten <= 8pt] (60:2cm) to (45:3cm);
\draw[<-,shorten >= 3pt, shorten <= 8pt] (60:2cm) to (75:3cm);
\draw[<-,shorten >= 3pt, shorten <= 8pt] (120:2cm) to (105:3cm);
\draw[<-,shorten >= 3pt, shorten <= 8pt] (120:2cm) to (135:3cm);
\draw[<-,shorten >= 3pt, shorten <= 8pt] (180:2cm) to (165:3cm);
\draw[<-,shorten >= 3pt, shorten <= 8pt] (180:2cm) to (195:3cm);
\draw[<-,shorten >= 3pt, shorten <= 8pt] (240:2cm) to (225:3cm);
\draw[<-,shorten >= 3pt, shorten <= 8pt] (240:2cm) to (255:3cm);
\draw[<-,shorten >= 3pt, shorten <= 8pt] (300:2cm) to (285:3cm);
\draw[<-,shorten >= 3pt, shorten <= 8pt] (300:2cm) to (315:3cm);
\draw[<-,shorten >= 3pt, shorten <= 8pt] (360:2cm) to (345:3cm);
\draw[<-,shorten >= 3pt, shorten <= 8pt] (360:2cm) to (15:3cm);

\draw[->,shorten >= 10pt, shorten <= 3pt] (45:3cm) to (37.5:4cm);
\draw[->,shorten >= 10pt, shorten <= 3pt] (45:3cm) to (52.5:4cm);
\draw[->,shorten >= 10pt, shorten <= 3pt] (75:3cm) to (67.5:4cm);
\draw[->,shorten >= 10pt, shorten <= 3pt] (75:3cm) to (82.5:4cm);
\draw[->,shorten >= 10pt, shorten <= 3pt] (105:3cm) to (97.5:4cm);
\draw[->,shorten >= 10pt, shorten <= 3pt] (105:3cm) to (112.5:4cm);
\draw[->,shorten >= 10pt, shorten <= 3pt] (135:3cm) to (127.5:4cm);
\draw[->,shorten >= 10pt, shorten <= 3pt] (135:3cm) to (142.5:4cm);
\draw[->,shorten >= 10pt, shorten <= 3pt] (165:3cm) to (157.5:4cm);
\draw[->,shorten >= 10pt, shorten <= 3pt] (165:3cm) to (172.5:4cm);
\draw[->,shorten >= 10pt, shorten <= 3pt] (195:3cm) to (187.5:4cm);
\draw[->,shorten >= 10pt, shorten <= 3pt] (195:3cm) to (202.5:4cm);
\draw[->,shorten >= 10pt, shorten <= 3pt] (225:3cm) to (217.5:4cm);
\draw[->,shorten >= 10pt, shorten <= 3pt] (225:3cm) to (232.5:4cm);
\draw[->,shorten >= 10pt, shorten <= 3pt] (255:3cm) to (247.5:4cm);
\draw[->,shorten >= 10pt, shorten <= 3pt] (255:3cm) to (262.5:4cm);
\draw[->,shorten >= 10pt, shorten <= 3pt] (285:3cm) to (277.5:4cm);
\draw[->,shorten >= 10pt, shorten <= 3pt] (285:3cm) to (292.5:4cm);
\draw[->,shorten >= 10pt, shorten <= 3pt] (315:3cm) to (307.5:4cm);
\draw[->,shorten >= 10pt, shorten <= 3pt] (315:3cm) to (322.5:4cm);
\draw[->,shorten >= 10pt, shorten <= 3pt] (345:3cm) to (337.5:4cm);
\draw[->,shorten >= 10pt, shorten <= 3pt] (345:3cm) to (352.5:4cm);
\draw[->,shorten >= 10pt, shorten <= 3pt] (15:3cm) to (7.5:4cm);
\draw[->,shorten >= 10pt, shorten <= 3pt] (15:3cm) to (22.5:4cm);

\node at (0,0)(x_01){\small$7$};

\node[fill=white] at (90:1cm)(x_1_1){\small$3$};
\node[fill=white] at (210:1cm)(x_1_2){\small$3$};
\node[fill=white] at (330:1cm)(x_1_3){\small$3$};

\node[fill=white] at (60:2cm)(x_2_1){\small$4$};
\node[fill=white] at (120:2cm)(x_2_2){\small$4$};
\node[fill=white] at (180:2cm)(x_2_3){\small$4$};
\node[fill=white] at (240:2cm)(x_2_4){\small$4$};
\node[fill=white] at (300:2cm)(x_2_5){\small$4$};
\node[fill=white] at (360:2cm)(x_2_6){\small$4$};

\node[fill=white] at (45:3cm)(x_3_1){\small$1$};
\node[fill=white] at (75:3cm)(x_3_2){\small$1$};
\node[fill=white] at (105:3cm)(x_3_3){\small$1$};
\node[fill=white] at (135:3cm)(x_3_4){\small$1$};
\node[fill=white] at (165:3cm)(x_3_5){\small$1$};
\node[fill=white] at (195:3cm)(x_3_6){\small$1$};
\node[fill=white] at (225:3cm)(x_3_7){\small$1$};
\node[fill=white] at (255:3cm)(x_3_8){\small$1$};
\node[fill=white] at (285:3cm)(x_3_9){\small$1$};
\node[fill=white] at (315:3cm)(x_3_10){\small$1$};
\node[fill=white] at (345:3cm)(x_3_11){\small$1$};
\node[fill=white] at (15:3cm)(x_3_12){\small$1$};

\node[fill=white]  at (37.5:4cm)(x_4_1){\small$1$};
\node[fill=white]  at (52.5:4cm)(x_4_2){\small$1$};
\node[fill=white]  at (67.5:4cm)(x_4_3){\small$1$};
\node[fill=white]  at (82.5:4cm)(x_4_4){\small$1$};
\node[fill=white]  at (97.5:4cm)(x_4_5){\small$1$};
\node[fill=white]  at (112.5:4cm)(x_4_6){\small$1$};
\node[fill=white]  at (127.5:4cm)(x_4_7){\small$1$};
\node[fill=white]  at (142.5:4cm)(x_4_8){\small$1$};
\node[fill=white]  at (157.5:4cm)(x_4_9){\small$1$};
\node[fill=white]  at (172.5:4cm)(x_4_10){\small$1$};
\node[fill=white]  at (187.5:4cm)(x_4_11){\small$1$};
\node[fill=white]  at (202.5:4cm)(x_4_12){\small$1$};
\node[fill=white]  at (217.5:4cm)(x_4_13){\small$1$};
\node[fill=white]  at (232.5:4cm)(x_4_14){\small$1$};
\node[fill=white]  at (247.5:4cm)(x_4_15){\small$1$};
\node[fill=white]  at (262.5:4cm)(x_4_16){\small$1$};
\node[fill=white]  at (277.5:4cm)(x_4_7){\small$1$};
\node[fill=white]  at (292.5:4cm)(x_4_18){\small$1$};
\node[fill=white]  at (307.5:4cm)(x_4_19){\small$1$};
\node[fill=white]  at (322.5:4cm)(x_4_20){\small$1$};
\node[fill=white]  at (337.5:4cm)(x_4_21){\small$1$};
\node[fill=white]  at (352.5:4cm)(x_4_22){\small$1$};
\node[fill=white]  at (7.5:4cm)(x_4_23){\small$1$};
\node[fill=white]  at (22.5:4cm)(x_4_24){\small$1$};
\end{tikzpicture}}.\end{center}
$M_3$ ist die unzerlegbare $T_3$-Darstellung zum Dimensionsvektor\enlargethispage{\baselineskip}
\begin{center}
\parbox[c]{6.6cm}{\begin{tikzpicture}[line width=1pt]
\draw[dotted, thin] (0,0) circle (1cm);
\draw[dotted, thin] (0,0) circle (2cm);
\draw[dotted, thin] (0,0) circle (3cm);

\draw[->,shorten >= 8pt, shorten <= 5pt] (0,0) to (90:1cm);
\draw[->,shorten >= 8pt, shorten <= 5pt] (0,0) to (210:1cm);
\draw[->,shorten >= 8pt, shorten <= 5pt] (0,0) to (330:1cm);

\draw[<-,shorten >= 9pt, shorten <= 8pt] (90:1cm) to (60:2cm);
\draw[<-,shorten >= 9pt, shorten <= 8pt] (90:1cm) to (120:2cm);
\draw[<-,shorten >= 9pt, shorten <= 8pt] (210:1cm) to (180:2cm);
\draw[<-,shorten >= 9pt, shorten <= 8pt] (210:1cm) to (240:2cm);
\draw[<-,shorten >= 9pt, shorten <= 8pt] (330:1cm) to (300:2cm);
\draw[<-,shorten >= 9pt, shorten <= 8pt] (330:1cm) to (360:2cm);

\draw[->,shorten >= 9pt, shorten <= 3pt] (60:2cm) to (45:3cm);
\draw[->,shorten >= 9pt, shorten <= 3pt] (60:2cm) to (75:3cm);
\draw[->,shorten >= 9pt, shorten <= 3pt] (120:2cm) to (105:3cm);
\draw[->,shorten >= 9pt, shorten <= 3pt] (120:2cm) to (135:3cm);
\draw[->,shorten >= 9pt, shorten <= 3pt] (180:2cm) to (165:3cm);
\draw[->,shorten >= 9pt, shorten <= 3pt] (180:2cm) to (195:3cm);
\draw[->,shorten >= 9pt, shorten <= 3pt] (240:2cm) to (225:3cm);
\draw[->,shorten >= 9pt, shorten <= 3pt] (240:2cm) to (255:3cm);
\draw[->,shorten >= 9pt, shorten <= 3pt] (300:2cm) to (285:3cm);
\draw[->,shorten >= 9pt, shorten <= 8pt] (300:2cm) to (315:3cm);
\draw[->,shorten >= 9pt, shorten <= 8pt] (360:2cm) to (345:3cm);
\draw[->,shorten >= 9pt, shorten <= 8pt] (360:2cm) to (15:3cm);
\node at (0,0)(x_01){\small$2$};

\node[fill=white] at (90:1cm)(x_1_1){\small$3$};
\node[fill=white] at (210:1cm)(x_1_2){\small$3$};
\node[fill=white] at (330:1cm)(x_1_3){\small$3$};

\node[fill=white] at (60:2cm)(x_2_1){\small$1$};
\node[fill=white] at (120:2cm)(x_2_2){\small$1$};
\node[fill=white] at (180:2cm)(x_2_3){\small$1$};
\node[fill=white] at (240:2cm)(x_2_4){\small$1$};
\node[fill=white] at (300:2cm)(x_2_5){\small$1$};
\node[fill=white] at (360:2cm)(x_2_6){\small$1$};

\node[fill=white] at (45:3cm)(x_3_1){\small$0$};
\node[fill=white] at (75:3cm)(x_3_2){\small$0$};
\node[fill=white] at (105:3cm)(x_3_3){\small$0$};
\node[fill=white] at (135:3cm)(x_3_4){\small$0$};
\node[fill=white] at (165:3cm)(x_3_5){\small$0$};
\node[fill=white] at (195:3cm)(x_3_6){\small$0$};
\node[fill=white] at (225:3cm)(x_3_7){\small$0$};
\node[fill=white] at (255:3cm)(x_3_8){\small$0$};
\node[fill=white] at (285:3cm)(x_3_9){\small$0$};
\node[fill=white] at (315:3cm)(x_3_10){\small$0$};
\node[fill=white] at (345:3cm)(x_3_11){\small$0$};
\node[fill=white] at (15:3cm)(x_3_12){\small$0$};\end{tikzpicture}}.\end{center}
Gemäß \cref{baum:erreich:lemma2} genügt es zu zeigen, daß die Darstellung $F_3(M_3)$ unzerlegbar ist. Anwenden von $F_3$ bedeutet Spiegeln an allen Quellen von $T_3$, die Nullen auf dem äußeren Kreis bleiben erhalten. $X:=F_3(M_3)$ ist also die unzerlegbare Darstellung zum Dimensionsvektor
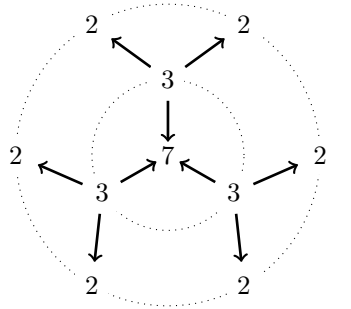
\begin{center}
\parbox[c]{4.6cm}{\begin{tikzpicture}[line width=1pt]
\draw[dotted, thin] (0,0) circle (1cm);
\draw[dotted, thin] (0,0) circle (2cm);

\draw[<-,shorten >= 8pt, shorten <= 5pt] (0,0) to (90:1cm);
\draw[<-,shorten >= 8pt, shorten <= 5pt] (0,0) to (210:1cm);
\draw[<-,shorten >= 8pt, shorten <= 5pt] (0,0) to (330:1cm);

\draw[->,shorten >= 9pt, shorten <= 8pt] (90:1cm) to (60:2cm);
\draw[->,shorten >= 9pt, shorten <= 8pt] (90:1cm) to (120:2cm);
\draw[->,shorten >= 9pt, shorten <= 8pt] (210:1cm) to (180:2cm);
\draw[->,shorten >= 9pt, shorten <= 8pt] (210:1cm) to (240:2cm);
\draw[->,shorten >= 9pt, shorten <= 8pt] (330:1cm) to (300:2cm);
\draw[->,shorten >= 9pt, shorten <= 8pt] (330:1cm) to (360:2cm);

\node at (0,0)(x_01){\small$7$};

\node[fill=white] at (90:1cm)(x_1_1){\small$3$};
\node[fill=white] at (210:1cm)(x_1_2){\small$3$};
\node[fill=white] at (330:1cm)(x_1_3){\small$3$};

\node[fill=white] at (60:2cm)(x_2_1){\small$2$};
\node[fill=white] at (120:2cm)(x_2_2){\small$2$};
\node[fill=white] at (180:2cm)(x_2_3){\small$2$};
\node[fill=white] at (240:2cm)(x_2_4){\small$2$};
\node[fill=white] at (300:2cm)(x_2_5){\small$2$};
\node[fill=white] at (360:2cm)(x_2_6){\small$2$};
\end{tikzpicture}}.\end{center}
Wir zeigen die Erreichbarkeit von $X$ mit Hilfe eines Projektionstricks und anschließender Reduktion auf den Fall eines zahmen Köchers vom Typ $\widetilde{D}_6$. Unsere Methode funktioniert in dieser Situation, weil wir genug Schur-Wurzeln "`unter"' $\dimv X$ finden und weil $\dimv X$ sogar eine reelle Schur-Wurzel ist. Genauer: Wir finden eine Folge $\dimv X=v_0,v_1,\ldots,v_n$ von Schur-Wurzeln, so daß $v_i-v_{i+1}=e_s$ die einfache Wurzel zum Punkt $s$ ist für alle $i$ und so daß $v_n$ "`klein genug"' ist. Einem Paar $(v_i,v_{i+1})$ entspricht dann eine exakte Sequenz
\[0\longrightarrow E_s\longrightarrow V_i\longrightarrow V_{i+1}\longrightarrow 0,\]
wobei $V_i$ und $V_{i+1}$ Schur-Darstellungen sind. Die Erreichbarkeit von $X$ ist damit auf die Erreichbarkeit von Schur-Darstellungen der Wurzel $v_n$ zurückgeführt. 

Wir wiederholen zunächst die in diesen Bemerkungen verwendeten Definitionen und gebrauchen für den Rest dieses Abschnitts die in Abschnitt \ref{grund:var} eingeführten Notationen für Varietäten von Köcherdarstellungen.
\begin{defn}\label{baum:erreich:def1} Sei $K$ ein endlicher Köcher ohne Schlaufen. Ein Element $\underline{d}\in \mathbb{Z}^{K_0}$ heißt eine Schur-Wurzel von $K$, falls eine $K$-Darstellung $M$ existiert mit $\End_K(M)=k$ und $\dimv M =\underline{d}$. $M$ heißt dann auch eine Schur-Darstellung von $\underline{d}$.
\end{defn}
\begin{bem}\label{baum:erreich:bem1}\ 
\begin{enumerate}
\item Die Schur-Darstellungen einer Schur-Wurzel $\underline{d}$ bilden eine offene Menge in der Varietät $\operatorname{Rep}_K(\underline{d})$.
\item Sei $x\in K_0$ eine Quelle von $K$, und sei $\underline{d}\neq e_x$ ein Dimensionsvektor. Dann gilt:
\[\underline{d} \text{ ist eine Schur-Wurzel von } K\Leftrightarrow s_x(\underline{d})\text{ ist eine Schur-Wurzel von } qK.\]
Dabei bezeichnet $s_x$ die Spiegelung zum Punkt $x$.
\item Die Menge der Schur-Wurzeln hängt im allgemeinen von der Orientierung des Köchers $K$ ab (vgl. \cite[§1.18]{Kac2}).
\item Die rellen Schur-Wurzeln sind die Dimensionsvektoren exzeptioneller $K$-Darstellungen. Sei dazu $q$ die quadratische Form des Köchers $K$. Die reellen Wurzeln sind genau die Dimensionsvektoren $\underline{d}$ unzerlegbarer Darstellungen mit $q(\underline{d})=1$ (vgl. \cite[Lemma 1.9] {Kac}). Die Formel $q(\dimv M)=\dim \End(M)-\dim\Ext^1(M,M)$ zeigt die Behauptung.
\end{enumerate}
\end{bem}

Unsere Rechnungen bestehen im Ermitteln von gewissen Schur-Wurzeln des Köchers $T_2$. Sei $e_s:=\dimv E_s$ der Dimensionsvektor des Einfachen zum Mittelpunkt $s$ von $T_2$. 
\begin{lemma}\label{baum:erreich:lemma4}\ 
\begin{enumerate}
\item Für $i=1,\ldots,4$ ist $\dimv X-i\cdot e_s$ eine Schur-Wurzel von $T_2$.
\item $\dimv X$ ist eine reelle Schur-Wurzel von $T_2$.
\item $\dimv X-5\cdot e_s$ ist eine reelle, nicht-schursche Wurzel von $T_2$.
\end{enumerate}
\end{lemma}
\begin{bew} Seien $x_1,\ldots,x_6$ die Punkte auf dem zweiten Kreis von $T_2$. Sei $K:=x_1\cdots x_6T_2$. Für alle $i\in\mathbb{N}$ sei $\underline{d}_i$ der Dimensionsvektor von $K$ gegeben durch
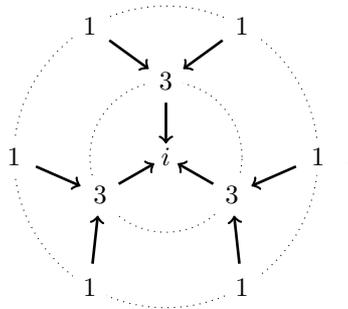
\begin{center}
\parbox[c]{4.6cm}{\begin{tikzpicture}[line width=1pt]
\draw[dotted, thin] (0,0) circle (1cm);
\draw[dotted, thin] (0,0) circle (2cm);

\draw[<-,shorten >= 8pt, shorten <= 5pt] (0,0) to (90:1cm);
\draw[<-,shorten >= 8pt, shorten <= 5pt] (0,0) to (210:1cm);
\draw[<-,shorten >= 8pt, shorten <= 5pt] (0,0) to (330:1cm);

\draw[<-,shorten >= 9pt, shorten <= 8pt] (90:1cm) to (60:2cm);
\draw[<-,shorten >= 9pt, shorten <= 8pt] (90:1cm) to (120:2cm);
\draw[<-,shorten >= 9pt, shorten <= 8pt] (210:1cm) to (180:2cm);
\draw[<-,shorten >= 9pt, shorten <= 8pt] (210:1cm) to (240:2cm);
\draw[<-,shorten >= 9pt, shorten <= 8pt] (330:1cm) to (300:2cm);
\draw[<-,shorten >= 9pt, shorten <= 8pt] (330:1cm) to (360:2cm);

\node at (0,0)(x_01){\small$i$};

\node[fill=white] at (90:1cm)(x_1_1){\small$3$};
\node[fill=white] at (210:1cm)(x_1_2){\small$3$};
\node[fill=white] at (330:1cm)(x_1_3){\small$3$};

\node[fill=white] at (60:2cm)(x_2_1){\small$1$};
\node[fill=white] at (120:2cm)(x_2_2){\small$1$};
\node[fill=white] at (180:2cm)(x_2_3){\small$1$};
\node[fill=white] at (240:2cm)(x_2_4){\small$1$};
\node[fill=white] at (300:2cm)(x_2_5){\small$1$};
\node[fill=white] at (360:2cm)(x_2_6){\small$1$};
\end{tikzpicture}}.\end{center}

(1.) Wir zeigen, daß die $\underline{d}_i$ Schur-Wurzeln von $K$ sind für $i=3,\ldots,6$. Spiegeln an allen Quellen des zweiten Kreises von $K$ zeigt dann die Behauptung für die entsprechenden Dimensionsvektoren von $T_2$. 

Sei $S$ eine $K$-Darstellung zum Dimensionsvektor $\underline{d}_i$, $i\geq 3$, so daß alle Pfeile durch injektive Abbildungen dargestellt werden. Nach Basiswahl ist $S$ gegeben durch drei dreidimensionale Unterräume $V_1,V_2,V_3\subseteq k^i$ und durch Paare $(x_j,y_j)$ von Vektoren in $V_j-\{0\}$ für $j=1,2,3$. Der Endomorphismenring von $S$ ist die Menge aller $k$-linearen Abbildungen $f:k^i\longrightarrow k^i$, so daß die $V_j$ $f$-invariante Unterräume von $k^i$ sind und alle $x_j,y_j$ Eigenvektoren von $f$ sind.

Wir konstruieren jetzt zu jedem $\underline{d}_i$ eine Schur-Darstellung $S$ von $K$, die alle Pfeile durch Monomorphismen darstellt. Wir verwenden dazu die oben angegebene Interpretation, und geben jeweils die Unterräume $V_j\subseteq k^{i}$ sowie die Paare $(x_j,y_j)$ in $V_j-\{0\}$ an.

$\mathbf{i=3}$. Es sei
\[V_1=V_2=V_3=k^i;\ (x_1,y_1)=(e_1,e_2),\ (x_2,y_2)=(e_3,\ \sum_{i=1}^3e_i),\ x_3,y_3\in V_3-\{0\}.\]
Man verifiziert sofort, daß die entsprechende $K$-Darstellung schursch ist. 

$\mathbf{i=6}$. Spiegeln am Mittelpunkt von $T_2$ liefert wieder den Dimensionsvektor $\underline{d}_3$ und man ist im Fall $i=3$. Die Umorientierung der Pfeile zum Mittelpunkt  spielt keine Rolle, da alle Pfeile mit Start im Mittelpunkt durch die Identität auf $k^3$ dargestellt werden können.

$\mathbf{i=4}$. Es sei
\[V_1=\langle e_1,e_2,e_3\rangle,\ V_2=\langle e_2,e_3,e_4\rangle,\ V_3=\langle \sum_{i=1}^4e_i,e_3,e_4\rangle\]
und 
\[(x_1,y_1)=(e_1,e_2),\ (x_2,y_2)=(e_3,e_4),\ (x_3,y_3)=(\sum_{i=1}^4e_i,e_4).\]
Die zugehörige $K$-Darstellung ist schursch, denn jeder Endomorphismus von $k^i$, der alle gegebenen Unterräume invariant läßt, hat die Eigenvektoren $e_1,e_2,e_3,e_4$ und $\sum_{i=1}^4e_i$, ist also diagonalisierbar mit einem Eigenwert.

$\mathbf{i=5}$. Sei
\[V_1=\langle e_1,e_2,e_3\rangle,\ V_2=\langle e_3,e_4,e_5\rangle,\ V_3=\langle e_5,\sum_{i=1}^5e_i,e_4\rangle\]
und 
\[(x_1,y_1)=(e_1,e_2),\ (x_2,y_2)=(e_3,e_4),\ (x_3,y_3)=(e_5,\sum_{i=1}^5e_i).\]
Mit demselben Argument wie im Falle $i=4$ folgt, daß die zugehörige $K$-Darstellung schursch ist.

(2.) Ist klar, da $X$ nach \cref{baum:erreich:lemma2} exzeptionell ist.

(3.) Betrachte den Dimensionsvektor $\underline{d}_2$ von $K$:
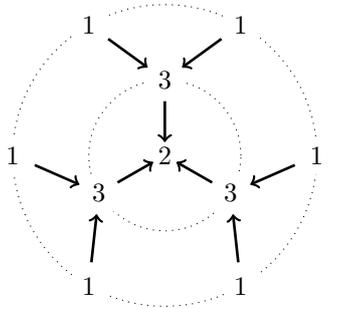
\begin{center}
\parbox[c]{4.5cm}{\begin{tikzpicture}[line width=1pt]
\draw[dotted, thin] (0,0) circle (1cm);
\draw[dotted, thin] (0,0) circle (2cm);

\draw[<-,shorten >= 8pt, shorten <= 5pt] (0,0) to (90:1cm);
\draw[<-,shorten >= 8pt, shorten <= 5pt] (0,0) to (210:1cm);
\draw[<-,shorten >= 8pt, shorten <= 5pt] (0,0) to (330:1cm);

\draw[<-,shorten >= 9pt, shorten <= 8pt] (90:1cm) to (60:2cm);
\draw[<-,shorten >= 9pt, shorten <= 8pt] (90:1cm) to (120:2cm);
\draw[<-,shorten >= 9pt, shorten <= 8pt] (210:1cm) to (180:2cm);
\draw[<-,shorten >= 9pt, shorten <= 8pt] (210:1cm) to (240:2cm);
\draw[<-,shorten >= 9pt, shorten <= 8pt] (330:1cm) to (300:2cm);
\draw[<-,shorten >= 9pt, shorten <= 8pt] (330:1cm) to (360:2cm);

\node at (0,0)(x_01){\small$2$};

\node[fill=white] at (90:1cm)(x_1_1){\small$3$};
\node[fill=white] at (210:1cm)(x_1_2){\small$3$};
\node[fill=white] at (330:1cm)(x_1_3){\small$3$};

\node[fill=white] at (60:2cm)(x_2_1){\small$1$};
\node[fill=white] at (120:2cm)(x_2_2){\small$1$};
\node[fill=white] at (180:2cm)(x_2_3){\small$1$};
\node[fill=white] at (240:2cm)(x_2_4){\small$1$};
\node[fill=white] at (300:2cm)(x_2_5){\small$1$};
\node[fill=white] at (360:2cm)(x_2_6){\small$1$};
\end{tikzpicture}}.\end{center}
Bis auf Umorientierung der äußeren Pfeile ist $\underline{d}_2$ der Dimensionsvektor der exzeptionellen $T_2^{op}$-Darstellung $M_3$. Die Menge der reellen Wurzeln hängt nicht von der Orientierung des Köchers ab, also ist $\underline{d}_2=\dimv M_3$ eine reelle Wurzel von $K$. 

Es bleibt zu zeigen, daß die zugehörige unzerlegbare Darstellung $U$ von $K$ nicht schursch ist. Sei dazu $y\in L_1$ ein Punkt auf dem ersten Kreis und seien $x_1$ und $x_2$ seine Vorgänger auf dem zweiten Kreis. Betrachte die Pfeile $\alpha_i:x_i\longrightarrow y$, $i=1,2$ mit Spitze $y$ sowie den Pfeil $\beta:y\longrightarrow s$ mit Start $y$. Wegen $\Ker U(\beta)\neq0$ existiert ein injektiver Homomorphismus $i:E_y\longrightarrow U$. Definiere eine Faktordarstellung $\overline{U}$ von $U$ durch
\[\overline{U}(t)=\left\{\begin{array}{cl}0&\text{ falls } t=x_i,i=1,2\\ U(y)/\bigl(\Bild U(\alpha_1)+\Bild U(\alpha_2)\bigr)&\text{ falls } t=y\\ U(x)&\text{ falls } t\not\in\{x_1,x_2,y\}\end{array}\right..\]
Sei $p:U\longrightarrow \overline{U}$ die kanonische Projektion. Offenbar ist $\overline{U}(y)\neq0$, also existiert ein surjektiver Homomorphismus $\phi:\overline{U}\longrightarrow E_y$. Betrachte den Endomorphismus $i\circ p\circ \phi$ von $U$. Dann ist $\phi\neq0$, aber $\phi$ ist kein Isomorphismus. Also ist $U$ nicht schursch.
\end{bew}

Das folgende Lemma faßt den angekündigten "`Projektionstrick"' zusammen:
\begin{lemma}\label{baum:erreich:lemma5}
Sei $K$ ein endlicher, zykelloser Köcher. Sei $(v_0,\ldots,v_n)$ eine Folge von Schur-Wurzeln von $K$ und sei $(t_0,\ldots,t_{n-1})$ eine Folge von Punkten in $K_0$, so daß gilt:
\begin{itemize}
\item[$\cdot$] Für alle $0\leq i\leq n-1$ ist $t_i$ eine Quelle oder eine Senke in $K$.
\item[$\cdot$] Für alle $0\leq i\leq n-1$ ist $v_{i}-v_{i+1}=e_{t_i}$. Dabei bezeichne $e_t$ die einfache Wurzel eines Punktes $t\in K_0$.
\end{itemize}
Es existiere eine offene, nichtleere Menge $U\subseteq \operatorname{Rep}_K(v_n)$ bestehend aus lauter erreichbaren (also insbesondere unzerlegbaren) Darstellungen. Dann gilt:
\begin{enumerate}
\item Es existiert auch eine offene, nichtleere Teilmenge $V\subseteq \operatorname{Rep}_K(v_0)$ bestehend aus erreichbaren Darstellungen.
\item Ist $v_0$ sogar eine reelle Schur-Wurzel und $X$ die zugehörige unzerlegbare Darstellung, so ist $X$ erreichbar.
\end{enumerate}
\end{lemma}
\begin{bew} (1.) Per Induktion dürfen wir $n=1$ annehmen. Schreibe $t=t_0$, $\underline{d}=v_0$ und $\underline{d}'=v_1$. Bis auf Dualität dürfen wir annehmen, daß $t$ eine Senke in $K$ ist. Sei $\mathcal{T}\subseteq K_1$ die Menge aller Pfeile in $K$ mit Spitze $t$ und sei $\mathcal{U}:=K_1-\mathcal{T}$. Dann gilt wegen der Voraussetzung $\underline{d}-\underline{d}'=e_t$:
\[\operatorname{Rep}_K(\underline{d}')=\prod_{\alpha\in\mathcal{U}}k^{d(s\alpha)\times d(n\alpha)}\times \prod_{\alpha\in\mathcal{T}}k^{d(t)-1\times d(n\alpha)}.\]
Betrachte die durch Eliminieren der letzten Zeile aller Matrizen $m_\alpha$, $\alpha\in\mathcal{T}$ gegebene Projektion
\[\pi:\operatorname{Rep}_K(\underline{d})\longrightarrow\operatorname{Rep}_K(\underline{d}').\]
Genauer gilt also für alle $m\in\operatorname{Rep}_K(\underline{d})$:
{\setlength{\arraycolsep}{0pt}
\[
\pi(m)_\alpha:=m_\alpha\text{ falls }\alpha\in\mathcal{U},\ \pi(m)_\alpha:=
\left[
\begin{array}{cl}
\bigl(m_\alpha\bigr)&_{1\bullet}\\
\vdots\\
\bigl(m_\alpha\bigr)&_{d(t)-1\bullet }\end{array}
\right]
\text{ falls }\alpha\in\mathcal{T}.
\]}
$A_{i\bullet}$ bezeichnet hier die $i$-te Zeile einer Matrix $A$. Da $t$ eine Senke von $K$ ist, existiert für alle $m\in\operatorname{Rep}_K(\underline{d})$ eine exakte Sequenz
\counterwithout{equation}{chapter}
\setcounter{equation}{0}
\begin{equation}\label{baum:erreich:lemma5:gl1}
0\longrightarrow E_t\longrightarrow \operatorname{rep} m\longrightarrow \operatorname{rep} \pi(m)\longrightarrow 0
\end{equation}
von $K$-Darstellungen. Sei jetzt $\operatorname{Rep}_K^{(1)}(\underline{d})$ die (offene, nichtleere) Menge aller Schur-Darstellun\-gen in $\operatorname{Rep}_K(\underline{d})$. Setze 
\[V:=\pi^{-1}(U)\cap \operatorname{Rep}_K^{(1)}(\underline{d}).\]
$V$ ist nichtleer und besteht aus lauter unzerlegbaren Darstellungen. Wegen der exakten Sequenzen in (\ref{baum:erreich:lemma5:gl1}) sind mit den Darstellungen in $U$ auch alle Darstellungen in $V$ erreichbar.

(2.) Nach (1.) existiert eine offene, nichtleere Teilmenge $V$ von $\operatorname{Rep}(v_0)$ bestehend aus erreichbaren Darstellungen. Da $v_0$ nach Voraussetzung eine reelle Schur-Wurzel ist, hat die zugehörigen Schur-Darstellung $X$ keine Selbsterweiterungen. Also ist der Orbit $\mathcal{O}(X)$ offen in $\operatorname{Rep}(v_0)$ und trifft die Menge $V$. Also ist $X$ erreichbar.
\end{bew}

\begin{prop}\label{baum:erreich:prop1}
$\widetilde{P}_4$ ist erreichbar.
\end{prop}
\begin{bew} Wir müssen zeigen, daß die zu Beginn dieses Paragraphen definierte unzerlegbare $T_2$-Darstellung $X$ erreichbar ist. Wir wissen bereits, daß $\dimv X$ eine reelle Schur-Wurzel ist (vgl. \cref{baum:erreich:lemma4}). Außerdem sind nach \cref{baum:erreich:lemma4} die $v_i:=\dimv X-i\cdot e_s$, $i=0,\ldots,4$ Schur-Wurzeln; $s$ ist überdies eine Senke von $T_2$. Nach \cref{baum:erreich:lemma5} genügt es also, eine offene Teilmenge $U\subseteq \operatorname{Rep}_{T_2}(v_4)$ zu bestimmen, die aus lauter erreichbaren Darstellungen besteht.

$v_4$ ist der Dimensionsvektor von $T_2$ gegeben durch
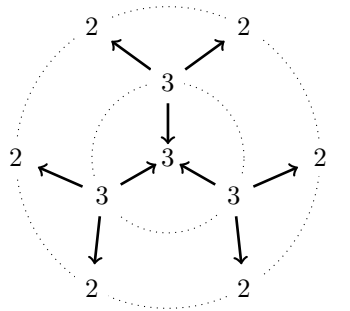
\begin{center}
\parbox[c]{4.6cm}{\begin{tikzpicture}[line width=1pt]
\draw[dotted, thin] (0,0) circle (1cm);
\draw[dotted, thin] (0,0) circle (2cm);

\draw[<-,shorten >= 8pt, shorten <= 5pt] (0,0) to (90:1cm);
\draw[<-,shorten >= 8pt, shorten <= 5pt] (0,0) to (210:1cm);
\draw[<-,shorten >= 8pt, shorten <= 5pt] (0,0) to (330:1cm);

\draw[->,shorten >= 9pt, shorten <= 8pt] (90:1cm) to (60:2cm);
\draw[->,shorten >= 9pt, shorten <= 8pt] (90:1cm) to (120:2cm);
\draw[->,shorten >= 9pt, shorten <= 8pt] (210:1cm) to (180:2cm);
\draw[->,shorten >= 9pt, shorten <= 8pt] (210:1cm) to (240:2cm);
\draw[->,shorten >= 9pt, shorten <= 8pt] (330:1cm) to (300:2cm);
\draw[->,shorten >= 9pt, shorten <= 8pt] (330:1cm) to (360:2cm);

\node at (0,0)(x_01){\small$3$};

\node[fill=white] at (90:1cm)(x_1_1){\small$3$};
\node[fill=white] at (210:1cm)(x_1_2){\small$3$};
\node[fill=white] at (330:1cm)(x_1_3){\small$3$};

\node[fill=white] at (60:2cm)(x_2_1){\small$2$};
\node[fill=white] at (120:2cm)(x_2_2){\small$2$};
\node[fill=white] at (180:2cm)(x_2_3){\small$2$};
\node[fill=white] at (240:2cm)(x_2_4){\small$2$};
\node[fill=white] at (300:2cm)(x_2_5){\small$2$};
\node[fill=white] at (360:2cm)(x_2_6){\small$2$};
\end{tikzpicture}}.\end{center}
Indem wir die oberen drei Punkte von $T_2$ eliminieren, erhalten wir einen zahmen Köcher vom Typ $\widetilde{D}_6$ als vollen Unterköcher von $T_2$. Betrachte den Dimensionsvektor $w$ von $\widetilde{D}_6$ gegeben durch
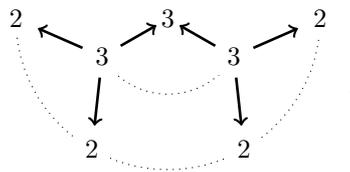
\begin{center}
\parbox[c]{4.6cm}{\begin{tikzpicture}[line width=1pt]
\draw[dotted, thin] (210:1cm) arc (210:330:1cm);
\draw[dotted, thin] (-2,0) arc (180:360:2cm);

\draw[<-,shorten >= 8pt, shorten <= 5pt] (0,0) to (210:1cm);
\draw[<-,shorten >= 8pt, shorten <= 5pt] (0,0) to (330:1cm);

\draw[->,shorten >= 9pt, shorten <= 8pt] (210:1cm) to (180:2cm);
\draw[->,shorten >= 9pt, shorten <= 8pt] (210:1cm) to (240:2cm);
\draw[->,shorten >= 9pt, shorten <= 8pt] (330:1cm) to (300:2cm);
\draw[->,shorten >= 9pt, shorten <= 8pt] (330:1cm) to (360:2cm);

\node at (0,0)(x_01){\small$3$};

\node[fill=white] at (210:1cm)(x_1_2){\small$3$};
\node[fill=white] at (330:1cm)(x_1_3){\small$3$};

\node[fill=white] at (180:2cm)(x_2_3){\small$2$};
\node[fill=white] at (240:2cm)(x_2_4){\small$2$};
\node[fill=white] at (300:2cm)(x_2_5){\small$2$};
\node[fill=white] at (360:2cm)(x_2_6){\small$2$};
\end{tikzpicture}}.\end{center}
Sei $\pi:\operatorname{Rep}_{T_2}(v_4)\longrightarrow \operatorname{Rep}_{\widetilde{D}_6}(w)$ die Projektion gegeben durch Einschränken von $T_2$-Darstel\-lungen auf den vollen Unterköcher $\widetilde{D}_6$.

Anwenden der quadratischen Form von $\widetilde{D}_6$ auf $w$ zeigt, daß $w$ eine Wurzel ist (vgl. \cite[§4]{CB1}). Durch Berechnung des Defekts von $w$ sieht man, daß $w$ eine präprojektive Wurzel von $\widetilde{D}_6$ ist (vgl. \cite[§7, Lemma 2]{CB1}). Sei $Y$ die zugehörige präprojektive unzerlegbare $\widetilde{D}_6$-Darstellung. Dann ist der Orbit $\mathcal{O}(Y)$ von $Y$ offen in $\operatorname{Rep}_{\widetilde{D}_6}(w)$. Setze \[U':=\pi^{-1}(\mathcal{O}(Y))\cap\operatorname{Rep}^{(1)}_{T_2}(v_4).\]
$\operatorname{Rep}^{(1)}_{T_2}(v_4)$ ist dabei die Menge aller Schur-Darstellungen von $T_2$ zum Dimensionsvektor $v_4$. $U'$ ist offen, nichtleer in $\operatorname{Rep}_{T_2}(v_4)$. Wir setzen
\[U:=\bigl\{u\in U'\mid \text{Alle Matrizen $u_\alpha,\ \alpha\in K_1$ haben vollen Rang}\bigr\}.\]
$U$ ist dann eine offene, nichtleere Teilmenge von $\operatorname{Rep}_{T_2}(v_4)$ bestehend aus lauter Schur-Dar\-stellungen. Wir müssen zeigen, daß alle Darstellungen in $U$ sogar erreichbar sind. 

Dazu verwenden wir die folgende Aussage $(*)$:

Sei $K$ ein endlicher, zykelloser Köcher, und sei $t$ eine Senke in $K$. Sei $\alpha:q\longrightarrow t$ en Pfeil in $K$. Sei $\mathcal{C}$ die volle Unterkategorie von $\fdar{K}$ bestehend aus allen Darstellungen $M$, so daß $M(\alpha)$ surjektiv ist, und sei $K'$ der volle Unterköcher von $K$ mit Punktmenge $K_0-\{t\}$. Betrachte den Funktor $F:\fdar{K}\longrightarrow\fdar{K'}$, der jeder $K$-Darstellung die Einschränkung zuordnet. Dann gilt für alle $M\in\mathcal{C}$: Ist $FM$ erreichbar, so ist auch $M$ erreichbar.

Sei dazu $M\in\mathcal{C}$ gegeben, $FM$ sei erreichbar. $FM$ ist insbesondere unzerlegbar; gemäß \cref{baum:kron:lemma3} ist dann auch $M$ unzerlegbar. Wir zeigen die Erreichbarkeit von $M$ per Induktion nach $n:=\dim M(t)$. Der Induktionsanfang $n=0$ ist klar. Sei also $n\geq1$. Da $t$ eine Senke ist, existiert eine exakte Sequenz
\[0\longrightarrow E_t\longrightarrow M\longrightarrow \overline{M}\longrightarrow 0\]
von $K$-Darstellungen. Offenbar ist $\overline{M}(\alpha)$ immer noch surjektiv, also ist $\overline{M}\in\mathcal{C}$. Wegen\linebreak$FM=F\overline{M}$ ist $\overline{M}$ per Induktion erreichbar, also ist auch $M$ erreichbar.

Sei jetzt $m\in U$. Per Konstruktion von $U$ ist $\pi (m)\in \mathcal{O}(Y)$, also ist $\operatorname{rep}\pi(m)\simeq Y$. Nach \cite[Theorem 6.1]{Bong2} sind unzerlegbare Darstellungen von zahmen Köchern stets erreichbar; also ist die präprojektive unzerlegbare $\widetilde{D}_6$-Darstellung $\operatorname{rep}\pi(m)\simeq Y$ erreichbar. Per Konstruktion ist weiter $M(\alpha)$ ein Epi- oder Monomorphismus für alle Pfeile $\alpha$ in $T_2$. Durch sukzessives Anwenden von $(*)$ bzw. der dualen Version von $(*)$ zeigt man, daß dann auch $\operatorname{rep}m$ erreichbar ist. Wir überlassen die Details dem Leser.
\end{bew}
\section{Exzeptionelle Moduln sind Baummoduln}\label{baum:hauptres}
\setcounter{zaehler}{0}
\numberwithin{zaehler}{section}
Wir bezeichnen in diesem Abschnitt stets mit $y$ die Quelle und mit $x$ die Senke von $Q_n$. Wir wollen den folgenden Satz von Ringel beweisen:
\begin{satz}\label{baum:hauptres:satz1}
Sei $K$ ein endlicher Köcher (nicht notwendig zykellos) und sei $M$ ein exzeptioneller $kK$-Modul. Dann ist $M$ eine Baummodul.
\end{satz}
Ringel beweist diesen Satz in \cite{Ri3} zunächst für die Köcher $Q_n$, indem er mit Hilfe von Spiegelungsfunktoren für die präprojektiven und präinjektiven $Q_n$-Darstellungen Matrixpräsentationen berechnet, die die gewünschten baumartigen Koeffizientenköcher liefern. Der allgemeine Fall folgt dann aus diesem Spezialfall per Induktion nach der Dimension der exzeptionellen Darstellung (bei festem Köcher) mit Hilfe der Schofield Induktion (\cref{scho:scho:satz1}) und der in Abschnitt (\ref{unz}) beschriebenenen Funktoren.

In \cite{Ri4} gibt Ringel einen Beweis an, der ohne die Matrix-Rechnungen im Kronecker-Fall auskommt. Dieser Beweis verwendet eine Induktion nach der Dimension der exzeptionellen Darstellung bei variierendem zugrundeliegendem Köcher $K$. Im Induktionsschritt wird zuerst der Fall einer exzeptionellen $Q_n$-Darstellung $M$ behandelt und $M$ als Darstellung der universellen Überlagerung $\widetilde{Q}_n$ gedeutet. Man findet dann eine geeignete exzeptionelle $\widetilde{Q}_n$-Unter- oder Faktordarstellung von $M$ der Dimension $\dim M-1$ und kann die Induktionsvoraussetzung anwenden. Der allgemeine Fall des Induktionsschritts wird dann wieder mit Hilfe der Schofield-Induktion behandelt. Wir geben einen Überblick über diese zweite Version des Beweises.

Wir benötigen einige Lemmata.
\begin{lemma}\label{baum:hauptres:lemma1} Sei $K$ ein endlicher zusammenhängender Köcher. Sei $\pi:\widetilde{K}\longrightarrow K$ die universelle Überlagerung und $F_\lambda$ der Überlagerungsfunktor. Sei $\widetilde{M}\in\fdar{\widetilde{Q}_n}$ eine unzerlegbare Baumdarstellung. Dann ist auch $F_\lambda(\widetilde{M})$ eine unzerlegbare Baumdarstellung.
\end{lemma}
\begin{bew} Berechne die Kardinalität der Träger geeigneter Darstellungsmatrizen der $F_\lambda(\widetilde{M})(\alpha)$. Verwende Teil (5) von \cref{baum:beg:bem1} und die Definition der universellen Überlagerung in Abschnitt \ref{grund:cov}.
\end{bew}

Wir wollen unsere Ergebnisse aus Kapitel \ref{unz} anwenden, und nicht nur Unzerlegbare sondern auch unzerlegbare Baumdarstellungen produzieren. Es sei dazu $L$ ein endlicher Köcher und es seien $X,Y$ unzerlegbare Baumdarstellungen von $L$ mit $\Hom_K(X,Y)=0$ und $X\not\simeq Y$. Setze $n:=\dim\Ext^1_L(Y,X)/R$; dabei sei $R$ das Radikal von $\Ext^1_L(Y,X)$ als $\End_L(X)$-$\End_L(Y)$-Bimodul. Sei $\mathcal{F}(Y,X)$ die volle Unterkategorie von $\fdar{L}$ aller $L$-Darstellungen $N$, für die eine exakte Sequenz der Form
\[0\longrightarrow X^u\longrightarrow N\longrightarrow Y^v\longrightarrow0\]
existiert. Wir erhalten:
\begin{lemma}\label{baum:hauptres:lemma3} Es existiert ein Funktor $F:\fdar{Q_n}\longrightarrow\fdar{L}$, der unzerlegbare Baumdarstellungen in unzerlegbare Baumdarstellungen überführt und der folgende Eigenschaften besitzt:\enlargethispage{\baselineskip}
\begin{enumerate}
\item $F$ ist treu und exakt.
\item $F$ erhält Unzerlegbarkeit.
\item $FM\simeq FM'$ impliziert $M\simeq M'$.
\item $F(E_y)=Y$ und $F(E_x)=X$.
\item Gilt zusätzlich $\Hom_L(Y,X)=0$ und $\End_L(X)=k=\End_L(Y)$, so ist $\mathcal{F}(Y,X)$ abgeschlossen unter Kernen und Cokernen und $F$ induziert eine Äquivalenz zwischen $\fdar{Q_n}$ und $\mathcal{F}(Y,X)$. 
\end{enumerate}
\end{lemma}
\begin{bew}[von \cref{baum:hauptres:lemma3}] Wähle Baumbasen $\mathcal{B}^{(X)}$ von $X$ und $\mathcal{B}^{(Y)}$ von $Y$. Betrachte den in \cref{grund:ext:prop2} definierten Epimorphismus
\[C^1_L(Y,X)\stackrel{\phi}{\longrightarrow} \Ext^1_L(Y,X).\]
Dabei ist
\[C^1_L(Y,X)=\bigoplus_{\alpha\in K_1}\Hom_k(Y(n\alpha),X(s\alpha)).\]
Wir wählen Elemente $\zeta_1,\ldots,\zeta_n\in C^1_L(Y,X)$ mit folgenden Eigenschaften:
\begin{enumerate}[i)]
\item Für alle $i$ gilt: Es existiert genau ein Pfeil $\alpha_i$ in $K$ mit $(\zeta_i)_{\alpha_i}\neq0$.
\item Für alle $i$ ist $(\zeta_i)_{\alpha_i}$ eine Elementarmatrix bezüglich der Basen $B^{(X)}_{s\alpha_i}$ von $X(s\alpha_i)$ und $B^{(Y)}_{n\alpha_i}$ von $Y(n\alpha_i)$, das heißt: Die Darstellungsmatrix von $(\zeta_i)_{\alpha_i}$ bezüglich dieser beiden Basen hat genau einen Eintrag $\lambda$ mit $\lambda=1$; alle anderen Einträge sind 0.
\item Die Restklassen der $\phi(\zeta_i)$ liefern eine Basis von $\Ext^1_L(Y,X)/R$, wobei $R$ das Radikal von $\Ext^1_L(Y,X)$ als $\End_L(X)$-$\End_L(Y)$-Bimodul bezeichnet. Insbesondere sind die $\zeta_i$ linear unabhängig.
\end{enumerate}
Für $Q_n$ verwenden wir die Bezeichnungen
\begin{center}
\parbox[c]{2.5cm}{
\begin{tikzpicture}[line width=1pt]
	\tikzstyle{knt}=[circle, fill, inner sep=1pt]
  \node at (0,0)[knt, label=below:$y$](x){};
  \node at (2,0)[knt, label=below:$x$](y){};
  \node at (1,0.1)(z){$\vdots$};
  \draw[->, bend left, shorten >= 2pt, shorten <= 2pt] (x) to node[above]{$\zeta_1$} (y);
  \draw[->, bend right, shorten >= 2pt, shorten <= 2pt] (x) to node[below]{$\zeta_n$} (y);
\end{tikzpicture}}.\end{center}
Es sei $F:\fdar{Q_n}\longrightarrow \fdar{L}$ der in \cref{unz:konst:def4} konstruierte Funktor bezüglich der Basis $(\zeta_1,\ldots,\zeta_n)$. Die Aussagen (1) bis (5) wurden bereits in Kapitel \ref{unz} bewiesen (vgl. \cref{unz:eig:lemma1}, \cref{unz:eig:prop1}, \cref{unz:eig:prop2}). Es bleibt zu zeigen, daß $F$ unzerlegbare Baumdarstellungen in unzerlegbare Baumdarstellungen überführt.

Sei also $T\in\fdar{Q_n}$ eine unzerlegbare Baumdarstellung. Setze $M:=FT$, $v:=\dim T(y)$ und $u:=\dim T(x)$. Sei $\mathcal{C}$ eine Baumbasis von $T$. Für $i=1,\ldots,n$ sei \[A^{(i)}=\Bigl[a_{kl}^{(i)}\Bigr]\in k^{u\times v}\]
die Darstellungsmatrix von $T(\zeta_i)$ bezüglich der Basen $\mathcal{C}_y$ von $T(y)$ und $\mathcal{C}_x$ von $T(x)$. Wir dürfen also annehmen: $T(x)=k^u$, $T(y)=k^v$ und $T(\zeta_i)=A_i$ für alle $i$. Gemäß Teil (5) von \cref{baum:beg:bem1} gilt die Formel
\counterwithout{equation}{chapter}
\setcounter{equation}{0}
\begin{equation}\label{baum:hauptres:lemma1:Gl1}
\sum_{i=1}^n\# Tr(A_i)=\dim T-1=u+v-1.
\end{equation}
Wir rekapitulieren die Definition von $F$ und erhalten für $q\in L_0$
\[M(q)=\bigl(X(q)\otimes_kk^u\bigr)\oplus \bigl(Y(q)\otimes_k k^v\bigr).\]
Für Pfeile $\alpha:q\longrightarrow s$ in $L$ gilt
\[
M(\alpha):=\begin{bmatrix}X(\alpha)\otimes id_{k^u}&\sum\limits_{i=1}^n(\zeta_i)_\alpha\otimes A_i\\0&Y(\alpha)\otimes id_{k^v}\end{bmatrix}.
\]
Eine kurze Rechnung zeigt: Bis auf Isomorphie von Darstellungen ist $M$ gegeben durch
\[M(q)=X(q)^u\oplus Y(q)^v\text{ für alle } q\in L_0.\]
Für Pfeile $\alpha:q\longrightarrow s$ in $L$ gilt
\begin{equation}\label{baum:hauptres:lemma1:Gl2}
M(\alpha):=\left[\begin{array}{cc}X(\alpha)^u&Z(\alpha)\\0&Y(\alpha)^v\end{array}\right],
\end{equation}
wobei $Z(\alpha)\in\Hom_k(Y(\alpha)^v,X(\alpha)^u)$ definiert ist durch die Matrix von linearen Abbildungen
\begin{equation}\label{baum:hauptres:lemma1:Gl3}
Z(\alpha)=\Biggl[\sum\limits_{i=1}^na_{kl}^{(i)}\cdot(\zeta_i)_\alpha\Biggr]_{\begin{subarray}{c}1\leq k\leq u\\1\leq l\leq v\end{subarray}}.
\end{equation}
Für alle $q\in L_0$ liefern die Basen $\mathcal{B}^{(X)}_q$ von $X(q)$ und $\mathcal{B}^{(Y)}_q$ von $Y(q)$ auf natürliche Weise eine Basis $\mathcal{B}^{(M)}_q$ von $M(q)=X(q)^u\oplus Y(q)^v$ (nimm für jeden direkten Summanden eine Kopie der jeweiligen Basis). Wir zeigen, daß $\mathcal{B}^{(M)}$ eine Baumbasis von $M$ ist.

Sei $\alpha:q\longrightarrow s$ ein Pfeil in $L$. Zur Abkürzung bezeichnen wir die Darstellungsmatrix von $X(\alpha)$ bezüglich der Basen $\mathcal{B}^{(X)}_q,\mathcal{B}^{(X)}_s$ ebenfalls mit $X(\alpha)$; genauso sei $Y(\alpha)$ die Darstellungsmatrix von $Y(\alpha)$ bezüglich der Basen $\mathcal{B}^{(Y)}_q,\mathcal{B}^{(Y)}_s$ und $(\zeta_i)_\alpha$ die Darstellungsmatrix von $(\zeta_i)_\alpha$ bezüglich der Basen $\mathcal{B}^{(Y)}_q,\mathcal{B}^{(X)}_s$. $M(\alpha)$ sei die Darstellungsmatrix von $M(\alpha)$ bezüglich der Basen $\mathcal{B}^{(M)}_q,\mathcal{B}^{(M)}_s$. Mit $t(A):=\#Tr(A)$ bezeichnen wir die Kardinalität des Trägers einer Matrix $A$. Da $\mathcal{B}^{(X)}$ und $\mathcal{B}^{(Y)}$ Baumbasen sind, erhalten wir nach Teil (5) von \cref{baum:beg:bem1} die Gleichungen
\begin{equation}\label{baum:hauptres:lemma1:Gl4}
\sum_{\alpha\in K_1}t\bigl(X(\alpha)\bigr)=\dim X-1\text{ und }\sum_{\alpha\in K_1}t\bigl(Y(\alpha)\bigr)=\dim Y-1.
\end{equation}
Wegen der in (\ref{baum:hauptres:lemma1:Gl2}) und (\ref{baum:hauptres:lemma1:Gl3}) angegebenen Form von $M(\alpha)$ gilt offenbar für jeden Pfeil $\alpha$ in $L$:
\begin{equation}\label{baum:hauptres:lemma1:Gl5}
t\bigl(M(\alpha)\bigr)=u\cdot t\bigl(X(\alpha)\bigr)+v\cdot t\bigl(Y(\alpha)\bigr)+\sum_{\begin{subarray}{c}1\leq k\leq u \\1\leq l\leq v\end{subarray}}t\Bigl(\sum_{i=1}^na_{kl}^{(i)}\cdot(\zeta_i)_\alpha\Bigr).
\end{equation}
Sei $x(s):=\dim X(s)$ und $y(s):=\dim Y(s)$. Nach unseren Vereinbarungen können wir die Kozykel $\zeta_i\in C^1_L(Y,X)$ als Tupel von Matrizen in
\[\bigoplus_{\alpha\in K_1}k^{x(s\alpha)\times y(n\alpha)}\]
auffassen. Die Bedingungen i), ii), iii) an die $\zeta_i$ ergeben:
\begin{enumerate}[a)]
\item Für alle $i=1,\ldots,n$ existiert genau ein Pfeil $\alpha_i\in K_1$ mit $(\zeta_i)_{\alpha_i}\neq 0$.
\item Für alle Pfeile $\alpha\in K_1$ sind die $(\zeta_i)_{\alpha}$ mit $\alpha_i=\alpha$ linear unabhängige Elementarmatrizen in $k^{x(s\alpha)\times y(n\alpha)}$.
\end{enumerate}
Mit Hilfe von Aussage b) erhält man leicht aus Gleichung (\ref{baum:hauptres:lemma1:Gl5}):
\begin{equation}\label{baum:hauptres:lemma1:Gl6}
t\bigl(M(\alpha)\bigr)=v\cdot t\bigl(X(\alpha)\bigr)+u\cdot t\bigl(Y(\alpha)\bigr)+\sum_{\begin{subarray}{c}1\leq i\leq n\\ \alpha=\alpha_i\end{subarray}}t(A^{(i)}).
\end{equation}
Mit Gleichung (\ref{baum:hauptres:lemma1:Gl4}) folgt daraus
\begin{equation}\label{baum:hauptres:lemma1:Gl7}
\sum_{\alpha\in K_1}t\bigl(M(\alpha)\bigr)=u\cdot(\dim X -1) + v\cdot(\dim Y-1)+
\sum_{\alpha\in K_1}\Bigl(\:\sum_{\begin{subarray}{c}1\leq i\leq n\\ \alpha=\alpha_i\end{subarray}}t(A^{(i)})\:\Bigr).
\end{equation}
Wegen Aussage a) gilt für alle $\alpha$:
\[\sum_{\alpha\in K_1}\Bigl(\:\sum_{\begin{subarray}{c}1\leq i\leq n\\ \alpha=\alpha_i\end{subarray}}t(A^{(i)})\:\Bigr)=\sum_{i=1}^nt(A^{(i)}).\]
Die Gleichungen (\ref{baum:hauptres:lemma1:Gl1}) und (\ref{baum:hauptres:lemma1:Gl7}) zeigen dann
\[
\begin{aligned}
\sum_{\alpha\in K_1}t\bigl(M(\alpha)\bigr)&=u\cdot(\dim X -1) + v\cdot(\dim Y-1)+(u+v-1)\\&=u\cdot\dim X+v\cdot\dim Y -1=\dim M-1.
\end{aligned}
\]
Nach \cref{unz:eig:prop1} ist mit $T$ auch $FT$ unzerlegbar, also ist $\mathcal{B}^{(M)}$ eine Baumbasis (vgl. Teil (5) von \cref{baum:beg:bem1}).
\end{bew}

Aus diesem Lemma erhält man folgende Aussage (vgl. die Bemerkung am Schluß von \cite{Ri4}):
\begin{kor}\label{baum:hauptres:kor1}
Seien $X$ und $Y$ unzerlegbare Baumdarstellungen eines endlichen Köchers $L$ mit
\[\Hom_L(X,Y)=0\text{ und }\dim\Ext^1_L(Y,X)=1.\]
Sei eine nicht-spaltende exakte Sequenz
\[0\longrightarrow X\longrightarrow M\longrightarrow Y\longrightarrow 0\]
von $L$-Darstellungen gegeben. Dann ist auch $M$ eine unzerlegbare Baumdarstellung.
\end{kor}
\begin{bew} Wegen $\Hom_L(X,Y)=0$ existiert nach \cref{baum:hauptres:lemma3} ein exakter Funktor
\[F:\fdar{Q_1}\longrightarrow\fdar{L},\]
der unzerlegbare Baumdarstellungen in unzerlegbare Baumdarstellungen überführt. Sei $T$ die $2$-dimensionale unzerlegbare Darstellung des Dynkin-Köchers $Q_1=A_2$. $T$ ist offenbar eine Baumdarstellung. Die exakte Sequenz
\[0\longrightarrow E_x\longrightarrow T\longrightarrow E_y\longrightarrow 0\]
von $Q_1$-Darstellungen wird von $F$ gemäß der in \cref{baum:hauptres:lemma3} aufgelisteten Eigenschaften in die exakte Sequenz
\[0\longrightarrow X\longrightarrow FT\longrightarrow Y\longrightarrow 0\]
überführt. Da $FT$ eine unzerlegbare Baumdarstellung ist, spaltet diese Sequenz nicht. Wegen $\dim\Ext^1_L(Y,X)=1$ muß also schon $FT\simeq M$ gelten und die Behauptung folgt.
\end{bew}

Wir wollen in diesem Zusammenhang noch das folgende Lemma erwähnen:
\begin{lemma}[\protect{\cite[Lemma 3.7]{Weist}}]\label{baum:erreich:lemma6}
Sei $A$ eine $k$-Algebra und sei
\[\xi:\ 0\longrightarrow M'\stackrel{\epsilon}{\longrightarrow} M\stackrel{\pi}{\longrightarrow} M''\longrightarrow 0\]
eine nicht-spaltende, exakte Sequenz von endlichdimensionalen $A$-Moduln. Dabei seien die Endterme $M'$ und $M''$ unzerlegbar und es gelte $\Hom_A(M',M'')=0$. Dann ist auch der Mittelterm $M$ unzerlegbar.
\end{lemma}
\begin{bew} Sei $M=U\oplus V$. Die Abbildungen $\epsilon$ und $\pi$ seien gegeben durch die Matrizen von $A$-linearen Abbildungen
\[\epsilon=\text{\setlength{\arraycolsep}{1pt}\small$\left[\begin{array}{c}\epsilon_U\\\epsilon_V\end{array}\right]$},\ \epsilon_U:M'\longrightarrow U,\ \epsilon_V:M'\longrightarrow V\] und
\[\pi=\text{\setlength{\arraycolsep}{1pt}\small$\left[\begin{array}{cc}\pi_U&\pi_V\end{array}\right]$},\ \pi_U:U\longrightarrow M'',\ \pi_V:V\longrightarrow M''.\]
Wegen $\Hom_A(M',M'')=0$ folgt $\pi_U\circ\epsilon_U=0=\pi_V\circ\epsilon_V$ und man erhält aus der Exaktheit der gegebenen Sequenz $\xi$
\[M'\simeq\Bild\epsilon=\Bild\epsilon_U\oplus\Bild\epsilon_V.\]
Da $M'$ unzerlegbar ist, dürfen wir $\epsilon_V=0$ annehmen. Wir erhalten
\[M''\simeq M/\Bild \epsilon \simeq U/\Bild \epsilon_U\oplus V/\Bild \epsilon_V=U/\Bild \epsilon_U\oplus V.\]
Angenommen, es ist $V\neq0$. Dann gilt wegen der Unzerlegbarkeit von $M''$ 
\[U/\Bild \epsilon_U=0,\text{ also } \Bild \epsilon_U=U.\]
Da $\epsilon_V=0$ ist folgt daraus: $\epsilon_U$ ist ein Isomorphismus. Also ist $\epsilon$ ein Schnitt und $\xi$ spaltet -- im Widerspruch zur Voraussetzung. Es folgt $V=0$ und damit die Behauptung.
\end{bew}

Wir kommen jetzt zum Beweis des Satzes von Ringel.
\begin{bew}[von \cref{baum:hauptres:satz1}]
Da $M$ exzeptionell ist, ist der Träger $Tr(M)$ von $M$ (aufgefaßt als voller Unterköcher von $K$) nach \cref{scho:erb:prop2} zykellos. $M$ ist als Darstellung von $Tr(M)$ immer noch exzeptionell. Wir müssen den Satz also nur für alle zykellosen Köcher beweisen. Wir verwenden eine einzige Induktion über die Dimension aller exzeptionellen Darstellungen aller endlichen zykellosen Köcher. Sei also $M$ eine exzeptionelle Darstellung eines endlichen zykellosen Köchers $K_M$. Der Induktionsanfang ist klar: Ist $\dim M=1$, so ist $M$ einfach, also ein Baummodul.

Sei im Induktionsschritt $\dim M>1$. Da $M$ unzerlegbar ist, hat $K_M$ mehr als zwei Punkte. Wir unterscheiden zwei Fälle.

1. Fall: $K_M$ hat zwei Punkte. In dem Fall ist $K_M=Q_n$ für eine natürliche Zahl $n\geq1$ (der Fall $n=0$ tritt nicht auf). Ist $n=1$, so ist $M$ die zweidimensionale unzerlegbare $Q_1$-Darstellung und die Behauptung folgt. Sei also $n\geq2$. Wir dürfen bis auf Dualität annehmen, daß $M$ präprojektiv unzerlegbar ist (vgl. \cref{scho:ex:bsp1}). Sei $\widetilde{Q}_n$ die universelle Überlagerung von $Q_n$ und $F_\lambda:\fdar{\widetilde{Q}_n}\longrightarrow\fdar{Q_n}$ der Überlagerungsfunktor. Gemäß \cref{baum:kron:prop1} existiert eine exzeptionelle $\widetilde{Q}_n$-Darstellung $\widetilde{M}$ mit $F_\lambda(\widetilde{M})=M$ und wir erhalten wegen \cref{baum:kron:satz1} eine exakte Sequenz
\counterwithout{equation}{chapter}
\setcounter{equation}{0}
\begin{equation}\label{baum:hauptres:satz1:gl1}
0\longrightarrow E\longrightarrow \widetilde{M}\longrightarrow C\longrightarrow 0
\end{equation}
von $\widetilde{Q}_n$-Darstellungen. Dabei ist $E$ projektiv einfach, $C$ exzeptionell und es gilt
\[\Hom_{\widetilde{Q}_n}(C,E)=0=\Hom_{\widetilde{Q}_n}(E,C)\text{ und }\dim\Ext^1_{\widetilde{Q}_n}(C,E)=1.\]
Wir dürfen $\widetilde{Q}_n$ durch den (endlichen und zykellosen) Träger $T$ von $\widetilde{M}$ ersetzen und alle involvierten Darstellungen als Darstellungen von $T$ auffassen (vgl. etwa die Bemerkungen zum Beginn des Beweises von \cref{baum:kron:satz1}). Per Induktion ist die exzeptionelle $T$-Darstellung $C$ eine Baumdarstellung. Außerdem sind die Voraussetzungen von \cref{baum:hauptres:kor1} erfüllt, also ist auch $\widetilde{M}$ eine unzerlegbare Baumdarstellung. Nach \cref{baum:hauptres:lemma1} ist dann $M=F_\lambda(\widetilde{M})$ ebenfalls eine unzerlegbare Baumdarstellung.

2. Fall: $K_M$ hat wenigstens $3$ Punkte. Wir wollen Schofield-Induktion anwenden: Gemäß Satz \ref{scho:scho:satz1} existieren exzeptionelle $K_M$-Darstellungen $X$ und $Y$ mit $\Hom_A(X,Y)=0=\Hom_A(Y,X)$, natürliche Zahlen $u,v\geq1$ und eine exakte Sequenz
\[0\longrightarrow X^u\longrightarrow M\longrightarrow Y^v\longrightarrow0.\]
Also ist $M\in\mathcal{F}(Y,X)$. Wegen $u,v\geq 1$ sind $X$ und $Y$ per Induktion unzerlegbare Baumdarstellungen. Wir dürfen außerdem annehmen: $X$ ist nicht einfach oder $Y$ ist nicht einfach. Wären nämlich sowohl $X$ als auch $Y$ einfach, so würde der Träger $T$ von $M$ nur aus zwei Punkten bestehen. $M$ wäre dann auch eine exzeptionelle $T$-Darstellung und die Behauptung folgte aus dem 1. Fall.

Sei $F:\fdar{Q_n}\longrightarrow \fdar{L}$ der in \cref{baum:hauptres:lemma3} konstruierte Funktor. Nach Teil (5) des Lemmas induziert $F$ eine Äquivalenz $\fdar{Q_n}\stackrel{\sim}{\longrightarrow} \mathcal{F}(Y,X)$. Sei also eine $Q_n$-Darstellung $T$ gegeben mit $FT\simeq M$. Da $\mathcal{F}(Y,X)$ nach Teil (5) des Lemmas abgeschlossen ist unter Kernen und Kokernen, ist $FT$ auch noch exzeptionell als Objekt in $\mathcal{F}(Y,X)$. Mithin ist $T$ eine exzeptionelle $Q_n$-Darstellung. Nach Teil (4) des Lemmas ist $F$ exakt, $F(E_x)\simeq X$ und $F(E_y)\simeq Y$. Da $F$ eine Äquivalenz ist, erhält man aus der exakten Sequenz (\ref{baum:hauptres:satz1:gl1}) von $K$-Darstellungen eine exakte Sequenz
\[0\longrightarrow E_x^u\longrightarrow T\longrightarrow E_y^v\longrightarrow0\]
von $Q_n$-Darstellungen. Also gilt $\dim T=u+v$. Da $X$ nicht einfach ist oder $Y$ nicht einfach ist, folgt
\[\dim T=u+v<u\dim X+v\dim Y=\dim M.\]
Also ist die $Q_n$-Darstellung $T$ per Induktion ein unzerlegbarer Baummodul und nach \cref{baum:hauptres:lemma3} ist dann auch $FT\simeq M$ ein Baummodul.
\end{bew}
\chapter{Unzerlegbare $Q_3$-Darstellungen in kleinen Dimensionen}\label{rech}
\setcounter{zaehler}{0}
\numberwithin{zaehler}{chapter}
Sei stets $k$ ein algebraisch abgeschlossener Körper. Wir wollen die Varietät der $Q_3$-Darstellungen zum Dimensionsvektor $(3,3)$ ein wenig näher untersuchen und mit Hilfe der bekannten Klassifikation von $Q_2$-Darstellungen die Untervarietäten $U^{(i)}_{Q_3}(3,3)$ bestimmen. Leider hat sich die Hoffnung zerschlagen, schon bei der Bestimmung der irreduziblen Komponenten der $U^{(i)}_{Q_3}(3,3)$ auf interessante Phänomene zu stoßen, wir werden nämlich zeigen, daß alle diese Varietäten irreduzibel sind. Die Rechnungen sind elementar, aber aufwendig; höherdimensionale und damit interessantere Beispiele lassen sich allein mit unseren Methoden nicht untersuchen.

\section{Notationen}\label{rech:not}
\setcounter{zaehler}{0}
\numberwithin{zaehler}{section}
Sei $Q_n$ der $n$-Kronecker Köcher mit den Bezeichnungen
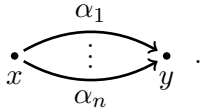
\begin{center}
\parbox[c]{2.6cm}{
\begin{tikzpicture}[line width=1pt]
	\tikzstyle{knt}=[circle, fill, inner sep=1pt]
  \node at (0,0)[knt, label=below:$x$](x){};
  \node at (2,0)[knt, label=below:$y$](y){};
  \node at (1,0.1)(z){$\vdots$};
  \draw[->, bend left, shorten >= 2pt, shorten <= 2pt] (x) to node[above]{$\alpha_1$} (y);
  \draw[->, bend right, shorten >= 2pt, shorten <= 2pt] (x) to node[below]{$\alpha_n$} (y);
\end{tikzpicture}}.\end{center}
Ist $M$ eine $Q_n$-Darstellung, so schreibe $\dimv M=(\dim M(x),\dim M(y))\in\mathbb{Z}^2$. Sei $\underline{d}=(d_x,d_y)$ ein Dimensionsvektor. Wir schreiben zur Abkürzung
\[R_n(\underline{d}):=\operatorname{Rep}_ {Q_n}(\underline{d})\text{ und } G(\underline{d}):=G_{Q_n}(\underline{d}).\]
Für alle $i\geq0$ sei $R_n^{(i)}(\underline{d}):=\operatorname{Rep}_{Q_n}^{(i)}$ die lokal abgeschlossene Untervarietät aller Darstellungen mit $i$-dimensionalem Endomorphismenring und $U_n^{(i)}(\underline{d}):=U_{Q_n}^{(i)}(\underline{d})$ die abgeschlossene Untervarietät von $R_n^{(i)}(\underline{d})$ bestehend aus allen unzerlegbaren Darstellungen (vgl. \cref{grund:var:prop3}).

Die Elemente von $R_n(\underline{d})$ sind $n$-Tupel $m=(m_1,\ldots,m_n)$ von Matrizen in $k^{d_y\times d_x}$. Die zugehörige Darstellung $\operatorname{rep} m$ ist in diesen Notationen gegeben durch
\[\bigl(\operatorname{rep}m\bigr)(x)=k^{d_x},\ \bigl(\operatorname{rep}m\bigr)(y)=k^{d_y},\ \bigl(\operatorname{rep}m\bigr)(\alpha_i)=m_i\text{ für }i=1,\ldots,n.\]
Es ist $G(\underline{d})=\operatorname{GL}_{d_x}\times\operatorname{GL}_{d_y}$. Die Operation von $G(\underline{d})$ auf $\operatorname{Rep}_{Q_n}(\underline{d})$ ist für alle $m\in R_n(\underline{d})$ und alle $g=(g_1,g_2)\in G(\underline{d})$ gegeben durch
\[g*m=(g_2\cdot m_i\cdot g_1^{-1})_{1\leq i\leq n}.\]
Wir haben noch eine Zusatzstruktur auf $R_n(\underline{d})$. Sei dazu $H_n:=\operatorname{GL}_n$. Dann operiert $H_n$ auf der Darstellungskategorie $\fdar{Q_n}$ durch Automorphismen von $k$-linearen Kategorien: Für alle $h\in H_n$ betrachte nämlich den Funktor
\[h*-:\fdar{Q_n}\stackrel{\sim}{\longrightarrow}\fdar{Q_n},\]
der wie folgt definiert ist: Für $M\in\fdar{Q_n}$ sei
\[\bigl(h*M\bigr)(x)=M(x)\text{ und } \bigl(h*M\bigr)(y)=M(y),\]
und für alle $i=1,\ldots,n$ sei
\[\bigl(h*M\bigr)(\alpha_i)=\sum_{j=1}^nh_{ij}\cdot M(\alpha_j).\]
Für Homomorphismen $f:M\longrightarrow M'$ sei außerdem $(h*f)_x=f_x$ und $(h*f)_y=f_y$. $H_n$ operiert auf völlig analoge Weise algebraisch auf dem affinen Raum $R_n(\underline{d})$; dabei gilt 
\[\operatorname{rep} (h*m) = h*(\operatorname{rep}m)\text{ für alle } h\in H_n, m\in R_n(\underline{d}).\]
Offenbar vertauschen die Operationen von $G(\underline{d})$ und $H_n$ auf $R_n(\underline{d})$, und wir erhalten eine Opera\-tion von $H_n\times G(\underline{d})$ auf $R_n(\underline{d})$.

Wir bemerken, daß die $R_n^{(i)}(\underline{d})$ und die $U_n^{(i)}(\underline{d})$ $H_n\times G(\underline{d})$-stabile Untervarietäten von $R_n(\underline{d})$ sind.

Im Fall $n=2$ hat die $H_2$-Operation auf der Darstellungskategorie eine sehr anschauliche Bedeutung. Verwende folgende Konvention: Identifizere die projektive Gerade $\mathbb{P}^1$ mit der Menge $k\cup \{\infty\}$. Ein Punkt $\lambda\in k$ entspreche dabei $[1:\lambda]\in\mathbb{P}^1$ und $\infty$ entspreche $[0:1]\in\mathbb{P}^1$. $H_2$ operiert auf $\mathbb{P}^1$ durch
\[h*[x:y]:=[h_{11}x+h_{12}y:h_{21}x+h_{22}y].\]

Es gelten (in den Notationen von \cref{grund:kron:def1}):
\begin{itemize}
\item[$\cdot$] $h*P_n\simeq P_n$ für alle $n\geq0$,
\item[$\cdot$] $h*I_n\simeq I_n$ für alle $n\geq0$,
\item[$\cdot$] $h*R_{n,\lambda}\simeq R_{n,h*\lambda}$ für alle $n\geq1$, $\lambda\in\mathbb{P}^1$.
\end{itemize}
Sei $\underline{d}\in\mathbb{Z}^2$ ein Dimensionsvektor von $Q_2$. Wir unterscheiden $Q_2$-Darstellungen nach ihrem "`Zerlegungstyp"'. Formal soll das bedeuten: Sei $r\geq0$ eine natürliche Zahl und seien für $i=1,\ldots,r$ natürliche Zahlen $\mu_1^{(i)},\ldots,\mu_{s_i}^{(i)}$ gegeben. Seien weiter $p_1,\ldots,p_t$ und $i_1,\ldots,i_u$ natürliche Zahlen. Sei
\[\underline{d}:=\sum_{k=1}^t(k,k+1)^{p_k}+\sum_{i=1}^r\sum_{n=1}^{s_i}(n,n)^{\mu_n^{(i)}}+\sum_{l=1}^u(l+1,l)^{i_l}\in \mathbb{Z}^2.\]
Die Menge $Z$ aller Punkte $m\in R_2(\underline{d})$, so daß paarweise verschieden Skalare $\lambda_1,\ldots,\lambda_r\in k$ existieren mit
\[\operatorname{rep} m\simeq \bigoplus_{k=1}^tP_k^{p_k}\oplus\Biggl(\bigoplus_{i=1}^r\Biggl(\bigoplus_{n=1}^{s_i}\bigl(R_{n,\lambda_i}\bigr)^{\mu_n^{(i)}}\Biggr)\Biggr)\oplus\bigoplus_{l=1}^uI_l^{i_l}\]
bezeichnen wir als Zerlegungsklasse von $\underline{d}$ (bezüglich der vorgegebenen Daten). Wir schreiben kurz
\[Z=\sum_{k=1}^t(k,k+1)^{p_k}+\sum_{i=1}^r\sum_{n=1}^{s_i}\bigl((n,n)_{\lambda_i}\bigr)^{\mu_n^{(i)}}+\sum_{l=1}^u(l+1,l)^{i_l}.\]
In Tabelle \ref{rech:deg:tab1} findet man eine Liste der Zerlegungsklassen von $\underline{d}=(3,3)$.

Für Dimensionsvektoren $\underline{d}\in\mathbb{Z}^2$ betrachte die Projektion
\[\pi^{(\underline{d})}:R_3(\underline{d})\longrightarrow R_2(\underline{d}),\ (m_1,m_2,m_3)\mapsto (m_1,m_2).\]
Faßt man $H_2$ via der Vorschrift $g\mapsto\Bigl[\begin{smallmatrix}g&0\\0&1\end{smallmatrix}\Bigr]$ als Untergruppe von $H_3$ auf, so wird $\pi^{(\underline{d})}$ ein $H_2\times G(\underline{d})$-invarianter Morphismus. Für einen Punkt $z\in R_2(\underline{d})$ schreiben wir zur Abkürzung:
\begin{itemize}
\item[$\cdot$] $F_z(\underline{d}):=(\pi^{(\underline{d})})^{-1}(z)$,
\item[$\cdot$] $U^{(i)}_z(\underline{d}):=F_z(\underline{d})\cap U_3^{(i)}(\underline{d})$ für alle $i\geq0$,
\item[$\cdot$] $U_z(\underline{d})=\bigcup\limits_{i\geq0} U_z^{(i)}(\underline{d})$,
\item[$\cdot$] $G_z(\underline{d}):=\operatorname{Iso}_{G(\underline{d})}(z)$.
\end{itemize}
$F_z(\underline{d})$ ist eine abgeschlossene, $G_z(\underline{d})$-stabile Untervarietät von $R_3(\underline{d})$ und die $U^{(i)}_z(\underline{d})$ sind $H_2\times G_z(\underline{d})$-stabile, lokal abgeschlossene Untervarietäten von $F_z(\underline{d})$. Wir identifizieren häufig die Operation von $G_z(\underline{d})$ auf $F_z(\underline{d})$ mit der Operation von $G_z(\underline{d})$ auf $k^{d_y\times d_x}$ via des Isomorphismus
\[F_z(\underline{d})\stackrel{\sim}{\longrightarrow}k^{d_y\times d_x},\ (z_1,z_2,A)\mapsto A,\]
sofern Mißverständnisse ausgeschlossen sind.

\section{Entartungen von Zerlegungsklassen in $\operatorname{Rep}_{Q_2}(3,3)$}\label{rech:deg}
\setcounter{zaehler}{0}
\numberwithin{zaehler}{section}
Sei in diesem Abschnit $G:=G(3,3)$. Wir betrachten die in Abschnitt \ref{rech:not} vorgestellte $H_2\times G$-Operation auf $R_2(3,3)$ und wollen das Entartungsdiagramm der $H_2\times G$-Orbiten von $R_2(3,3)$ und sämtliche Bahnendimensionen berechnen.

Es gilt folgendes Lemma:
\begin{lemma}\label{rech:deg:lemma1}
Die $H_2\times G$-Orbiten in $R_2(3,3)$ sind gerade die Zerlegungsklassen von $(3,3)$. Insbesondere besitzt $R_2(3,3)$ nur endliche viele $H_2\times G$-Orbiten.\enlargethispage{\baselineskip}
\end{lemma}
\begin{bew} Die Aussage folgt aus der Definition der Zerlegungsklassen und der im letzten Abschnitt vorgestellten Bedeutung der Operation von $H_2$ auf den $Q_2$-Darstellungen.
\end{bew}\begin{bem}\label{rech:deg:bem1} Die Aussage von \cref{rech:deg:lemma1} ist im allgemeinen falsch. Betrachte etwa die Zerlegungsklasse
\[Z:=\sum_{i=1}^4(1,1)_{\lambda_i}\subseteq R_4(4,4).\]
$Z$ ist eine $H_2\times G(4,4)$-stabile Teilmenge von $R_4(4,4)$. Sei $\overline{Z}$ der mengentheoretische Quotient der $H_2\times G(4,4)$-Operation auf $Z$. Die symmetrische Gruppe $S_4$ mit $4$ Elementen operiert auf $X:=\bigr(\mathbb{P}_1\bigl)^4-\Delta$ durch Vertauschen der Einträge ($\Delta$ bezeichnet hier die Diagonale); man erhält so eine $H_2\times S_4$-Operation auf $X$. Diese Operation ist offenbar nicht transitiv. Sei $\overline{X}$ der mengentheoretische Quotient dieser Operation. Es existiert eine Bijektion $\overline{X}\simeq \overline{Z}$, also ist $Z$ keine $H_2\times G(4,4)$-Bahn.
\end{bem}
In Tabelle \ref{rech:deg:tab1} sind die $H_2\times G$-Bahnen von $R_2(3,3)$ angegeben: In der zweiten Spalte stehen die Zerlegungsklassen von $(3,3)$ in den Notationen aus dem letzten Abschnitt, in der dritten Spalte ist die Dimension der jeweiligen $H_2\times G$-Bahn aufgeführt und in der vierten die Dimension des $Q_2$-Endomorphismenrings eines (und damit aller) Repräsentanten der Bahn.

Die Berechnung der Bahnendimensionen erfolgt immer nach dem selben Muster. Wir wollen beispielsweise die Dimension der Bahn $B_8=(1,0)+(1,1)_{\lambda_1}+(1,1)_{\lambda_2}+(0,1)$ berechnen. Setze dazu $\Delta:=\{(x,x)\in k^2\mid x\in k\}$ und betrachte den Morphismus
\[f:G\times \bigl(k^2-\Delta\bigr)\longrightarrow B_8,\ \bigl(g,(\lambda_1,\lambda_2)\bigr)\mapsto g*\Bigl(\left[\begin{smallmatrix}1&0&0\\0&1&0\\0&0&0\end{smallmatrix}\right],\left[\begin{smallmatrix}\lambda_1&0&0\\0&\lambda_2&0\\0&0&0\end{smallmatrix}\right]\Bigr).\]
Dann ist 
\[f(G\times \bigl(k^2-\Delta\bigr))=\{(m_1,m_2\in B_8)\mid \operatorname{Rang} m_1\geq2\}\] offen in $B_8$. Da $B_8$ irreduzibel ist, folgt 
\[\overline{ f(G\times \bigl(k^2-\Delta\bigr)) }= B_8.\]
Für alle $x=f\bigl(g,(\lambda_1,\lambda_2)\bigr)\in \Bild f$ gilt $f^{-1}(x)=\operatorname{Iso}_G(x) \times \{(\lambda_1,\lambda_2)\}$, also ist
\[\dim f^{-1}(x)=\dim \operatorname{Iso}_G(x)=\dim_k\End_{Q_2}(x)=8.\]
Aus den Sätzen über Faserdimensionen von Morphismen (vgl. \cite[Exercise II.3.22]{Ha}) folgt
\[\dim B_8=\dim G\times \bigl(k^2-\Delta\bigr)-\dim f^{-1}(x)=20-8=12.\]

\numberwithin{table}{section}
\setcounter{table}{\value{zaehler}}
\renewcommand{\captionformat}{~ ~}
\renewcommand{\tablename}{Tabelle}
\renewcommand{\tableformat}{\tablename~\thetable}
{\renewcommand{\arraystretch}{1.3}
\begin{table}
\caption{$H_2\times G$-Orbiten von $R_2(3,3)$}
\label{rech:deg:tab1}
\centering
\begin{tabular}{|c|c|c|c|}\hline
&& $\dim_{H_2\times G}$&$\dim \End_{Q_2}$\\\hline
$B_1$&$(0,1)^3+ (1,0)^3$& 0 & 18\\\hline
$B_2$&$(0,1)^2+ (1,1)_{\lambda_1}+ (1,0)^2$& 4 & 15\\\hline
$B_3$&$(1,2)+ (0,1)+ (1,0)^2$&8&10\\\hline
$B_4$&$(0,1)^2+ (1,0)+ (2,1)$&8&10\\\hline
$B_5$&$(0,1)+(1,1)^2_{\lambda_1} +(1,0)$&9&10\\\hline
$B_6$&$(1,1)_{\lambda_1}^3$&10&9\\\hline
$B_7$&$(0,1)+(2,2)_{\lambda_1}+(1,0)$&11&8\\\hline
$B_8$&$(1,0)+(1,1)_{\lambda_1}+(1,1)_{\lambda_2}+(0,1)$&12&8\\\hline
$B_9$&$(1,2)+(1,1)_{\lambda_1}+ (1,0)$&13&6\\\hline
$B_{10}$&$(0,1)+(1,1)_{\lambda_1}+ (2,1)$&13&6\\\hline
$B_{11}$&$(1,1)_{\lambda_1}+(2,2)_{\lambda_1}$&14&5\\\hline
$B_{12}$&$(1,2)+(2,1)$&14&4\\\hline
$B_{13}$&$(2,3)+(1,0)$&14&4\\\hline
$B_{14}$&$(0,1)+(3,2)$&14&4\\\hline
$B_{15}$&$(1,1)_{\lambda_1}+(1,1)^2_{\lambda_2}$&15&5\\\hline
$B_{16}$&$(3,3)_{\lambda_1}$&16&3\\\hline
$B_{17}$&$(1,1)_{\lambda_1}+(2,2)_{\lambda_2}$&17&3\\\hline
$B_{18}$&$(1,1)_{\lambda_1}+(1,1)_{\lambda_2}+(1,1)_{\lambda_3}$&18&3\\\hline
\end{tabular}
\end{table}}
\setcounter{zaehler}{\value{table}}

In Abbildung \ref{rech:deg:abb1} ist das Hasse-Diagramm der partiellen Ordnung $\preceq_{deg}$ auf der Menge der \linebreak $H_2\times G$-Orbiten von $R_2(3,3)$ angegeben, ganz rechts stehen die Dimensionen der jeweiligen Orbiten. Ein Pfeil $B_1\longrightarrow B_2$ im Diagramm entspricht einer minimalen Entartung $B_1\prec_{deg} B_2$. Eine Entartung $B_1\prec_{deg} B_2$ heißt dabei minimal, falls kein Orbit $C$ mit $B_1\prec_{deg}C\prec_{deg} B_2$ existiert. Beachte, daß nach \cref{grund:var:prop1} echte Entartungen eines Orbits nur in Orbiten echt kleiner Dimension existieren.

\numberwithin{figure}{section}
\setcounter{figure}{\value{zaehler}}
\renewcommand{\captionformat}{~ ~}
\renewcommand{\figurename}{Abbildung}
\renewcommand{\figureformat}{\figurename~\thefigure}
\begin{figure}
\caption{Entartungs-Diagramm für die $H_2\times G$-Operation auf $R_2(3,3)$}
\label{rech:deg:abb1}
\begin{tikzpicture}[line width=1pt]
  \node at (0,0)(x_1){$(0,1)^3+ (1,0)^3$};
  \node at (0,1.5)(x_2){$(0,1)^2+ (1,1)_{\lambda_1}+ (1,0)^2$};
  \node at (-3,3)(x_3){$(0,1)^2+ (1,0)+ (2,1)$};
  \node at (3,3)(x_4){$(1,2)+ (0,1)+ (1,0)^2$}; 
  \node at (-6,4.5)(x_5){$(0,1)+(1,1)^2_{\lambda_1} +(1,0)$};
  \node at (-6,6)(x_6){$(1,1)_{\lambda_1}^3$};
  \node at (0,7.5)(x_7){$(0,1)+(2,2)_{\lambda_1}+(1,0)$};
  \node at (0,9)(x_8){$(1,0)+(1,1)_{\lambda_1}+(1,1)_{\lambda_2}+(0,1)$};
  \node at (-3,10.5)(x_9){$(1,2)+(1,1)_{\lambda_1}+ (1,0)$};  
  \node at (3,10.5)(x_10){$(0,1)+(1,1)_{\lambda_1}+ (2,1)$};
  \node at (-6,12)(x_11){$(1,1)_{\lambda_1}+(2,2)_{\lambda_1}$};
  \node at (-3,12)(x_12){$(2,3)+(1,0)$};
  \node at (0,12)(x_13){$(1,2)+(2,1)$};
  \node at (3,12)(x_14){$(3,2)+(0,1)$};
  \node at (-6,13.5)(x_15){$(1,1)_{\lambda_1}+(1,1)^2_{\lambda_2}$};  
  \node at (0,15)(x_16){$(3,3)_{\lambda_1}$}; 
  \node at (0,16.5)(x_17){$(1,1)_{\lambda_1}+(2,2)_{\lambda_2}$};  
  \node at (0,18)(x_18){$(1,1)_{\lambda_1}+(1,1)_{\lambda_2}+(1,1)_{\lambda_3}$};
  
  \node at (7,0){0};
  \node at (7,1.5){4};
  \node at (7,3){8};
  \node at (7,4.5){9};
  \node at (7,6){10};
  \node at (7,7.5){11};
  \node at (7,9){12};
  \node at (7,10.5){13};  
  \node at (7,12){14};
  \node at (7,13.5){15};  
  \node at (7,15){16}; 
  \node at (7,16.5){17};  
  \node at (7,18){18};  
 
  \draw[->] (x_2) to (x_1);
  \draw[->] (x_3) to (x_2);
  \draw[->] (x_4) to (x_2); 
  \draw[->, bend right] (x_5) to (x_2);   
  \draw[->] (x_6) to (x_5); 
  \draw[->] (x_7) to (x_3); 
  \draw[->] (x_7) to (x_4); 
  \draw[->] (x_7) to (x_5);
  \draw[->] (x_8) to (x_7); 
  \draw[->] (x_9) to (x_8); 
  \draw[->] (x_10) to (x_8); 
  \draw[->] (x_11) to (x_6); 
  
  \path[->] (x_11) edge (x_9) 
  					(x_12) edge (x_9)
  					(x_13) edge (x_9)
  					(x_13) edge (x_10)
  					(x_14) edge (x_10) 
  					(x_11) edge[-,bend right=5,line width=4pt, draw=white] (x_10)
						       edge[bend right=5] (x_10);
						       
  \draw[->] (x_15) to (x_11); 
  \draw[->] (x_16) to (x_12);  
  \draw[->] (x_16) to (x_13);
  \draw[->] (x_16) to (x_14);
  \draw[->] (x_16) to (x_11);
  \draw[->] (x_17) to (x_15);
  \draw[->] (x_17) to (x_16);
  \draw[->] (x_18) to (x_17);
\end{tikzpicture}
\end{figure}
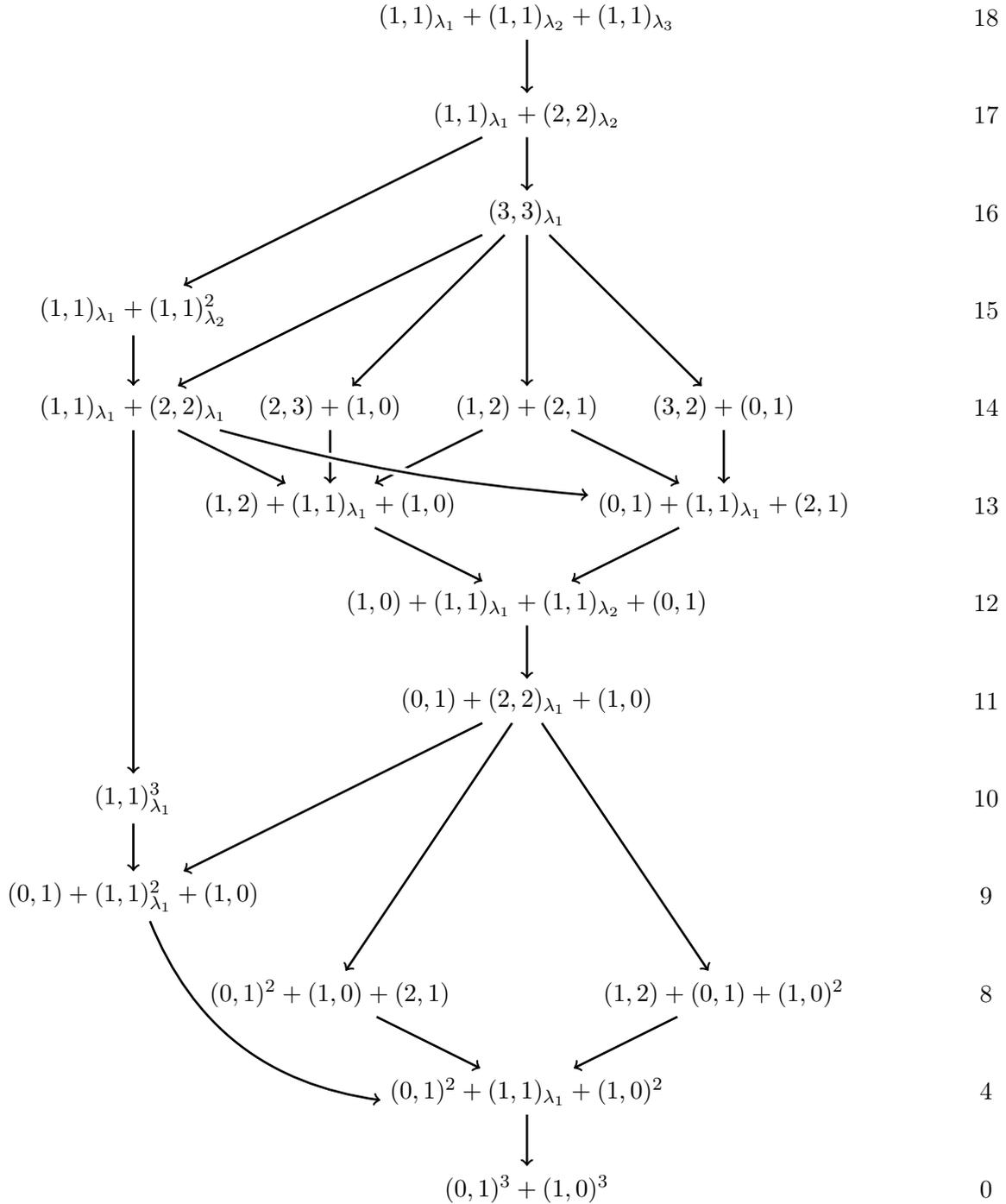
\setcounter{zaehler}{\value{figure}}

Wir wollen andeuten, wie man die Richtigkeit des Entartungsdiagramms verifizieren kann. Wir zeigen zunächst anhand zweier Beispiele, daß die Pfeile tatsächlich Entartungen von Orbiten entsprechen -- in allen anderen Situationen verwendet man genau dieselben Argumente.

1. Beispiel: Zeige, daß der Orbit $B_8=(1,0)+(1,1)_{\lambda_1}+(1,1)_{\lambda_2}+(0,1)$ in den Orbit \linebreak$B_7=(1,0)+(2,2)_{\lambda_1}+(0,1)$ entartet. Sei dazu ein Punkt $m=(m_1,m_2)$ in $B_7$  gegeben. Modulo der $H_2\times G$-Operation können wir annehmen:
\[m_1=\left[\begin{smallmatrix}1&0&0\\0&1&0\\0&0&0\end{smallmatrix}\right]\text{ und } m_2=\left[\begin{smallmatrix}0&0&0\\1&0&0\\0&0&0\end{smallmatrix}\right].\]
Betrachte für alle $\epsilon\in k$ den Punkt $t_\epsilon\in R_2(3,3)$ gegeben durch die zwei Matrizen
\[(t_\epsilon)_1=\left[\begin{smallmatrix}1&0&0\\0&1&0\\0&0&0\end{smallmatrix}\right]\text{ und } (t_\epsilon)_2=\left[\begin{smallmatrix}\epsilon&0&0\\1&0&0\\0&0&0\end{smallmatrix}\right].\]
Dann ist $t_\epsilon\in B_8$, falls $\epsilon\neq0$, also ist $m=t_0\in\overline{B_8}$.

2. Beispiel: Wir verwenden die Tatsache, daß exakte Sequenzen Entartungen liefern und zeigen, daß der der Orbit $B_{11}=(1,1)_{\lambda_1}+(2,2)_{\lambda_1}$ in den Orbit $B_9=(1,2)+(1,1)_{\lambda_1}+(1,0)$ entartet. Sei dazu ein Punkt $m\in B_9$ gegeben. Dann ist $\operatorname{rep} m\simeq P_1\oplus R_{1,\mu}\oplus I_0$ für ein $\mu\in\mathbb{P}^1$. Es existiert eine exakte Sequenz 
\[0\longrightarrow P_1\longrightarrow R_{2,\mu}\longrightarrow I_0\longrightarrow 0\]
von $Q_2$-Darstellungen (eine Liste von Mitteltermen von exakten Sequenzen von Kronecker-Moduln, so daß beide Endterme unzerlegbar sind, findet man in \cite[Teil 5]{Bong2}) und wir erhalten nach \cref{grund:var:prop2}  folgende Entartungen von $Q_2$-Darstellungen:
\[R_{2,\mu}\preceq_{deg}P_1\oplus I_0,\text{ und damit }R_{1,\mu}\oplus R_{2,\mu}\preceq_{deg}R_{1,\mu}\oplus P_1\oplus I_0\:(\simeq\operatorname{rep}m).\]
Es folgt:
\[m\in \mathcal{O}(R_{1,\mu}\oplus P_1\oplus I_0)\subseteq\overline{\mathcal{O}(R_{1,\mu}\oplus R_{2,\mu})}\subseteq \overline{B_{11}}.\]

Mit Hilfe der folgenden fünf Aussagen verifiziert man leicht, daß das Diagramm alle minimalen Entartungen beschreibt. Wir arbeiten uns von der Senke zur Quelle vor.

I) Betrachte zunächst die Teilmenge $Y$ von $R_2(3,3)$ bestehend aus allen Punkten $m$, so daß $m_1$ und $m_2$ linear abhängig sind in $k^{3\times 3}$. $Y$ besteht gerade aus den Orbiten
\[(0,1)^3+ (1,0)^3,\ (0,1)^2+ (1,1)_{\lambda_1}+ (1,0)^2,\ (0,1)+(1,1)^2_{\lambda_1} +(1,0),\ (1,1)^3_{\lambda_1}.\]
Offenbar ist $Y$ abgeschlossen und $H_2\times G$-stabil in $R_2(3,3)$, und wir erhalten, daß die Orbiten $(0,1)+(1,1)^2_{\lambda_1} +(1,0)\subseteq Y$ und $(1,1)_{\lambda_1}^3\subseteq Y$ in keinen der Orbiten $(0,1)^2+(1,0)+(2,1)\not\subseteq Y$, $(0,1)+(1,2)+(1,0^2)\not\subseteq Y$  entarten.

II) Betrachte jetzt die Teilmenge $Y$ von $R_2(3,3)$ aller Punkte $m$, so daß $\operatorname{rep} m$ einen einfachen direkten Summanden hat. $Y$ ist wieder abgeschlossen und $H_2\times G$-stabil, also existiert kein Orbit $B\subseteq Y$, so daß $B$ in $(1,1)_{\lambda_1}^3$ entartet.

III) $(1,2)+(2,1)$ entartet nicht in $(1,1)_{\lambda_1}^3$, denn sonst gäbe es Entartungen $P_1\oplus I_1 \preceq_{deg}(R_{1,\mu})^3$ von $Q_2$-Darstellungen für beliebige $\mu\in\mathbb{P}^1$; es ist aber 
\[0=\dim \Hom_{Q_2}((R_{1,\mu})^3,P_1)<\dim \Hom_{Q_2}(P_1\oplus I_1,P_1),\]
also ist das notwendige Entartungskriterium aus \cref{grund:var:prop3} nicht erfüllt.

IV) $B_{15}=(1,1)_{\lambda_1}+(1,1)_{\lambda_2}^2$ entartet nicht in die Orbiten $B_{14}:=(0,1)+(3,2)$, $B_{13}:=(1,0)+(2,3)$, $B_{12}:=(1,2)+(2,1)$: $B_{15}$ liegt nach Tabelle \ref{rech:deg:tab1} in der abgeschlossenen, $H_2\times G$-invarianten Untervarietät $Y\subseteq R_2(3,3)$ aller $m$ mit $\dim\End_{Q_2}(\operatorname{rep} m)\geq5$. Die Orbiten $B_{14}$, $B_{13}$ und $B_{12}$ sind hingegen nicht in $Y$ enthalten.

V) Schließlich entartet $B:=B_{16}=(3,3)_\lambda$ nicht in $C:=B_{15}=(1,1)_{\lambda_1}+(1,1)_{\lambda_2}^2$: Sei $Y$ die offene, $G$-stabile Untervarietät von $R_2(3,3)$ bestehend aus allen Punkten $m=(m_1,m_2)$, so daß $m_1$ invertierbar ist. Betrachte die $G$-stabilen Teilmengen $C':=C\cap Y$ und $B':=B\cap Y$ in $Y$. Wir zeigen: $C'\not\subseteq\overline{B'}$ (in $Y$). Da $B$ irreduzibel ist, liegt $B'$ dicht in $B$, und man verifiziert dann leicht, daß auch $C\not\subseteq\overline{B}$ (in $R_2(3,3)$) gilt. Betrachte die offene Abbildung
\[\pi:Y\longrightarrow k^{3\times 3},\ (m_1,m_2)\mapsto m_2\cdot(m_1)^{-1}.\]
Auf $k^{3\times 3}$ operiert die Untergruppe 
\[G':=\Delta_G=\{(g,g)\in G\mid g\in \operatorname{GL}_3\}\simeq\operatorname{GL}_3\]
von $G$ durch Konjugation und $\pi$ ist $G'$-äquivariant. Für $\lambda\in k$ betrachte die Matrix
\[J_{\lambda}:=\text{\small$\left[\begin{array}{ccc}\lambda&0&0\\1&\lambda&0\\0&1&\lambda\end{array}\right]$}\in k^{3\times3}.\]
Es sei $S:=\bigcup\limits_{\lambda\in k}G'*J_\lambda$. Dann ist $B'=\pi^{-1}(S)$. Da $\pi$ offen ist, folgt
\[\overline{B'}=\overline{\pi^{-1}(S)}=\pi^{-1}(\overline{S}).\]
Wir behaupten, daß die Menge $X$ aller Matrizen mit genau einem Eigenwert abgeschlossen ist in $k^{3\times3}$. Dann folgt $X=\overline{S}$, und man erhält
\[\overline{B'}=\overline{\pi^{-1}(S)}=\pi^{-1}(X)=\Bigl((3,3)_{\lambda_1}\cup \bigl((1,1)_{\lambda_1}+(2,2)_{\lambda_1}\bigr)\cup (1,1)_{\lambda_1}^3\Bigr)\cap Y\not\supseteq C'.\]
Um die Abgeschlosseneheit von $X$ zu zeigen, zeigen wir allgemeiner: Für alle $n\geq1$ ist die Menge 
\[X:=\{A\in k^{n\times n}\mid A \text{ hat genau einen Eigenwert}\}\]
abgeschlossen in $k^{n\times n}$. Wir identifizieren den affinen Raum $k^n$ mit der Menge aller normierten Polynome in $k[X]$ vom Grad $n$. Betrachte den Morphismus $\chi:k^{n\times n}\longrightarrow k^n$, der jeder Matrix ihr charakteristisches Polynom zuordnet. Seien $s_1,\ldots,s_n\in k[T_1,\ldots,T_n]$ die elementarsymmetrischen Polynome in $n$ Veränderlichen. Die $s_i$ sind also die Koeffizienten des Polynoms
\[\prod_{i=1}^n (X-T_i)\in \bigl(k[T_1,\ldots,T_n]\bigr)[X].\]
Der Morphismus $s:k^n\longrightarrow k^n,\ x\mapsto (s_1(x),\ldots,s_n(x))$ ist endlich, da die Ringerweiterung $k[s_1,\ldots,s_n]\subseteq k[T_1,\ldots,T_n]$ per Definition der $s_i$ ganz, also endlich ist. Insbesondere ist nach dem Going-Up-Theorem $s$ abgeschlossen. Sei jetzt $\Delta:=\{(\lambda,\ldots,\lambda)\in k^n)\mid \lambda\in k\}$. Dann ist $s(\Delta)$ abgeschlossen in $k^n$, und es ist
\[s(\Delta)=\bigl\{P\in k[X]\mid P\text{ normiert},\ \operatorname{deg} P=n,\ P \text{ hat genau eine Nullstelle}\bigr\}.\]
Es folgt $X=\chi^{-1}(s(\Delta))$, also ist $X$ abgeschlossen in $k^{n\times n}$.
\section{Bestimmung der $U_{Q_3}^{(i)}(3,3)$}\label{rech:unz}
\setcounter{zaehler}{0}
\numberwithin{zaehler}{section}
In diesem Abschnitt schreiben wir zur Abkürzung
\[G:=G(3,3),\ F_z:=F_z(3,3),\ G_z:=G_z(3,3),\ U_z^{(i)}:=U_z^{(i)}(3,3),\ U^{(i)}:=U_3^{(i)}(3,3).\]
Wir wollen die $H_3\times G$-stabilen Untervarietäten $U^{(i)}$ von $R_3(3,3)$ als $H_3\times G$-Saturierung gewisser Teilmengen beschreiben und dafür die Klassifikation der $H_2\times G$-Bahnen von $R_2(3,3)$ aus dem letzten Abschnitt verwenden. 

Sei $\mathcal{S}$ ein Repräsentantensystem der $H_2\times G$-Bahnen von $R_2(3,3)$ und sei $i\geq0$. Hat man für alle $z\in\mathcal{S}$ eine Teilmenge $T_z\subseteq U_z^{(i)}$ bestimmt mit $G_z*T_z=U_z^{(i)}$, so ist $U^{(i)}$ die $H_3\times G$-Saturierung von $T:=\bigcup_{z\in\mathcal{S}}T_z$. Zur Lösung unseres Problems genügt es also, derartige Mengen $T_z$ zu bestimmen. Wir konnten in den meisten Fällen die $U_z^{(i)}$ durch Gleichungen beschreiben.

Der Rest dieses Abschnitts ist der Bestimmung der $T_z$ gewidmet. Wir beschreiben zunächst das allgemeine Vorgehen und präsentieren anschließend unsere Ergebnisse, die wir durch ausgewählte Beispielrechnungen untermauern. In fast allen Fällen sind die $U_z^{(i)}$ irreduzibel; lediglich für $z\in B_{11}$ hat $U_z^{(2)}$ zwei irreduzible Komponenten.

\paragraph{Die Vorgehensweise}\ \newline
Es sei $B$ ein $H_2\times G$ Orbit in $R_2(3,3)$ und $z\in B$ ein Repräsentant. Den "`Standardfall"' (dies sind alle Bahnen bis auf $B_1,\ldots,B_6$ und $B_{15}$) bearbeitet man mit folgender Methode:

Sei $m=(z_1,z_2,m_3)\in F_z$. $A_m:=\End_{Q_3}(\operatorname{rep} m)$ ist eine Unteralgebra von $\End_{Q_2}(\operatorname{rep} z)$. $A_m$ läßt sich als Lösungsmenge eines homogenen linearen Gleichungssystems $\mathcal{L}_m$ beschreiben, dessen Koeffizienten lineare Funktionen in den Einträgen von $m_3$ sind. Man findet sofort einfache Bedingungen an die Einträge von $m_3$, so daß $A_m$ eine lokale Algebra ist. Durch Fallunterscheidung bestimmt man dann alle $m\in F_z$, so daß $A_m$ lokal ist (also $\operatorname{rep}m$ unzerlegbar ist). Geschicktes Anwenden der $G_z$-Operation kann dabei die Gleichungen vereinfachen.

Darüber hinaus kann man in einigen Fällen sofort für die Unzerlegbarkeit von $\operatorname{rep} m$ notwendige Bedingungen an $m_3$ angeben und so den Rechenaufwand etwas reduzieren. Man erhält (nach eventuellem Zurückverfolgen der $G_z$-Operation) eine Beschreibung der Varietäten $U_z^{(i)}$ als Vereinigung endlich vieler lokal abgeschlossener Untervarietäten von $F_z$. Wir werden dieses Vorgehen anhand des Beispiels der Bahn $B_{11}$ demonstrieren; dieser Fall ist auch der mit den aufwendigsten Rechnungen.

Die Bahnen $B_1$, $B_2$, $B_5$ und $B_6$ werden mit Hilfe der Klassifikation von Kronecker-Moduln (also $Q_2$-Darstellungen) abgehandelt; für die Bahnen $B_3$ und $B_4$ findet man einfache Argumente. Wir geben in diesen Fällen kurze Anleitungen. 

Den Sonderfall $B_{15}$ kann man nach geeigneten $G_z$-Operationen wie den Standardfall behandeln, wir geben hier die Reduktion an.

\paragraph{Ergebnisse}\ \newline
In Tabelle \ref{rech:unz:tab1} ist für jede $H_2\times G$-Bahn $B_i$ ein Repräsentant $z^{(i)}$ angegeben. Die Einträge unter den $U_z^{(i)}$ geben die Anzahl der irreduziblen Komponenten von $U_z^{(i)}$ an; ein Strich bedeutet dabei $U_z^{(i)}=\emptyset$. Man beachte, daß die Anzahl der irreduziblen Komponenten von $U_z^{(i)}$ nicht von der Wahl des Repräsentanten $z$ abhängt.

Wir beschreiben nun für jeden $H_2\times G$-Orbit $B\subseteq R_2(3,3)$ und den zugehörigen Repräsentanten $z$ aus Tabelle \ref{rech:unz:tab1} $U_z^{(i)}$ als $G_z$-Saturierung gewisser Teilmengen und bestimmen die irreduziblen Komponenten der $U_z^{(i)}$.

\numberwithin{table}{section}
\setcounter{table}{\value{zaehler}}
{\renewcommand{\arraystretch}{2}
\begin{table}
\caption{Übersicht über Unzerlegbare in den $F_z$}
\label{rech:unz:tab1}
\centering
{\small
\begin{tabular}{|@{}>{\centering\arraybackslash}m{0.7cm}|c@{}|c@{}|c@{}|c@{}|c@{}|}\hline
&Repräsentant $z$& $U_z^{(1)}\:$& $U_z^{(2)}\:$& $U_z^{(3)}\:$& $U_z^{(i)},i\geq 4\:$\\\hline
$\:\:B_1$&$(0,0)$&-&-&-&-\\\hline
$\:\:B_2$&
$z_1^{(2)}:=\left[\begin{smallmatrix}1&0&0\\0&0&0\\0&0&0\end{smallmatrix}\right],
z_2^{(2)}:=0$
&-&-&-&-\\\hline
$\:\:B_3$&
$z_1^{(3)}:=\left[\begin{smallmatrix}1&0&0\\0&0&0\\0&0&0\end{smallmatrix}\right],
z_2^{(3)}:=\left[\begin{smallmatrix}0&0&0\\1&0&0\\0&0&0\end{smallmatrix}\right]$
&-&-&1&-\\\hline
$\:\:B_4$&$z_1^{(4)}:=\left[\begin{smallmatrix}1&0&0\\0&0&0\\0&0&0\end{smallmatrix}\right],
z_2^{(4)}:=\left[\begin{smallmatrix}0&1&0\\0&0&0\\0&0&0\end{smallmatrix}\right]$
&-&-&1&-\\\hline
$\:\:B_5$&
$z_1^{(5)}:=\left[\begin{smallmatrix}0&0&0\\1&0&0\\0&1&0\end{smallmatrix}\right],
z_2^{(5)}:=0$
&-&-&1&-\\\hline
$\:\:B_6$&
$z_1^{(6)}:=\left[\begin{smallmatrix}1&0&0\\0&1&0\\0&0&1\end{smallmatrix}\right],
z_2^{(6)}:=0$
&-&-&1&-\\\hline
$\:\:B_7$&
$z_1^{(7)}:=\left[\begin{smallmatrix}1&0&0\\0&1&0\\0&0&0\end{smallmatrix}\right],
z_2^{(7)}:=\left[\begin{smallmatrix}0&0&0\\1&0&0\\0&0&0\end{smallmatrix}\right]$
&-&1&1&-\\\hline
$\:\:B_8$&
$z_1^{(8)}:=\left[\begin{smallmatrix}1&0&0\\0&1&0\\0&0&0\end{smallmatrix}\right],
z_2^{(8)}:=\left[\begin{smallmatrix}1&0&0\\0&0&0\\0&0&0\end{smallmatrix}\right]$
&-&1&-&-\\\hline
$\:\:B_9$&
$z_1^{(9)}:=\left[\begin{smallmatrix}1&0&0\\0&0&0\\0&1&0\end{smallmatrix}\right],
z_2^{(9)}:=\left[\begin{smallmatrix}0&0&0\\1&0&0\\0&0&0\end{smallmatrix}\right]$
&1&1&-&-\\\hline
$\:\:B_{10}$&
$z_1^{(10)}:=\left[\begin{smallmatrix}1&0&0\\0&0&1\\0&0&0\end{smallmatrix}\right],
z_2^{(10)}:=\left[\begin{smallmatrix}0&1&0\\0&0&0\\0&0&0\end{smallmatrix}\right]$
&1&1&-&-\\\hline
$\:\:B_{11}$&
$z_1^{(11)}:=\left[\begin{smallmatrix}1&0&0\\0&1&0\\0&0&1\end{smallmatrix}\right],
z_2^{(11)}:=\left[\begin{smallmatrix}0&0&0\\1&0&0\\0&0&0\end{smallmatrix}\right]$
&1&2&1&-\\\hline
$\:\:B_{12}$&
$z_1^{(12)}:=\left[\begin{smallmatrix}1&0&0\\0&0&0\\0&0&1\end{smallmatrix}\right],
z_2^{(12)}:=\left[\begin{smallmatrix}0&0&0\\1&0&0\\0&1&0\end{smallmatrix}\right]$
&1&1&1&-\\\hline
$\:\:B_{13}$&
$z_1^{(13)}:=\left[\begin{smallmatrix}1&0&0\\0&1&0\\0&0&0\end{smallmatrix}\right],
z_2^{(13)}:=\left[\begin{smallmatrix}0&0&0\\1&0&0\\0&1&0\end{smallmatrix}\right]$
&1&-&-&-\\\hline
$\:\:B_{14}$&
$z_1^{(14)}:=\left[\begin{smallmatrix}1&0&0\\0&1&0\\0&0&0\end{smallmatrix}\right],
z_2^{(14)}:=\left[\begin{smallmatrix}0&1&0\\0&0&1\\0&0&0\end{smallmatrix}\right]$
&1&-&-&-\\\hline
$\:\:B_{15}$&
$z_1^{(15)}:=\left[\begin{smallmatrix}1&0&0\\0&1&0\\0&0&1\end{smallmatrix}\right],
z_2^{(15)}:=\left[\begin{smallmatrix}0&0&0\\0&0&0\\0&0&1\end{smallmatrix}\right]$
&1&1&-&-\\\hline
$\:\:B_{16}$&
$z_1^{(16)}:=\left[\begin{smallmatrix}1&0&0\\0&1&0\\0&0&1\end{smallmatrix}\right],
z_2^{(16)}:=\left[\begin{smallmatrix}0&0&0\\1&0&0\\0&1&0\end{smallmatrix}\right]$
&1&1&1&-\\\hline
$\:\:B_{17}$&
$z_1^{(17)}:=\left[\begin{smallmatrix}1&0&0\\0&1&0\\0&0&1\end{smallmatrix}\right],
z_2^{(17)}:=\left[\begin{smallmatrix}0&0&0\\1&0&0\\0&0&1\end{smallmatrix}\right]$
&1&1&-&-\\\hline
$\:\:B_{18}$&
$z_1^{(18)}:=\left[\begin{smallmatrix}1&0&0\\0&1&0\\0&0&1\end{smallmatrix}\right],
z_2^{(18)}:=\left[\begin{smallmatrix}1&0&0\\0&0&0\\0&0&x\end{smallmatrix}\right],x\neq0,1\:$
&1&-&-&-\\\hline
\end{tabular}}
\end{table}}

\setcounter{zaehler}{\value{table}}
$\mathbf{B_3=(1,2)+ (0,1)+ (1,0)^2}$\\[4pt]
$U_z^{(3)}=\bigl\{C\bigm| c_{32}=c_{33}=0, \operatorname{Rang}C=3\bigr\}=G_z*\left[\begin{smallmatrix}0&1&0\\0&0&1\\1&0&0\end{smallmatrix}\right].$

$U_z^{(3)}$ ist irreduzibel, denn $G_z$ ist zusammenhängend als offene Teilmenge des affinen Raums $\End_{Q_2}(z)$, also ist $U_z^{(3)}$ irreduzibel als $G_z$-Orbit.

$\mathbf{B_4=(0,1)^2+ (1,0)+ (2,1)}$\\[4pt]
Zu $B_3$ dualer Fall.

$\mathbf{B_5=(0,1)+(1,1)^2_{\lambda_1} +(1,0)}$\\[4pt]
$U_z^{(3)}= G_z*\left[\begin{smallmatrix}1&0&0\\0&1&0\\0&0&1\end{smallmatrix}\right]$.

$U_z^{(3)}$ ist irreduzibel als $G_z$-Orbit.

$\mathbf{B_6=(1,1)_{\lambda_1}^3}$\\[4pt]
$U_z^{(3)}=\bigcup_{\lambda\in k}G_z*\left[\begin{smallmatrix}\lambda&0&0\\1&\lambda&0\\0&1&\lambda\end{smallmatrix}\right]$.

$U_z^{(3)}$ ist irreduzibel als $G_z$-Saturierung einer irreduziblen Teilmenge von $F_z$.

$\mathbf{B_7=(0,1)+(2,2)_{\lambda_1}+(1,0)}$\\[4pt]
$U_z^{(2)}=T_1\cup T_2$. Dabei ist
\begin{align*}
T_1&=\Bigl\{C\in k^{3\times 3}\Bigm| \left(\begin{aligned}&c_{13}\neq0=c_{33}\\
&c_{31}+\frac{c_{23}}{c_{13}}\cdot c_{32}\neq0\text{ oder } c_{32}\neq0\end{aligned}\ \right)\ \Bigr\},\\[12pt]
T_2&=\Bigl\{C\in k^{3\times 3}\Bigm| \left(\begin{aligned}&c_{33}=c_{13}=0\\
&c_{23}\cdot c_{32}\neq0\end{aligned}\ \right)\ \Bigr\}.\end{align*}
$U_z^{(3)}=\Bigl\{C\in k^{3\times 3}\Bigm| c_{13}=c_{32}=c_{33}=0\text{ und }c_{23}\cdot c_{31}\neq0\ \Bigr\}$.

$U_z^{(2)}$ und $U_z^{(3)}$ sind irreduzibel: Die Irreduzibilität von $U_z^{(3)}$ ist klar; außerdem ist $T_1$ irreduzibel und dicht in $U_z^{(2)}$.
\newpage
$\mathbf{B_8=(1,0)+(1,1)_{\lambda_1}+(1,1)_{\lambda_2}+(0,1)}$\\[4pt]
$U_z^{(2)}=\bigcup_{i=1}^3T_i$. Dabei ist 
\begin{align*}
T_1&=\Bigl\{C\in k^{3\times 3}\Bigm| \left(\begin{aligned}&c_{33}=0\\
&c_{13}\cdot c_{23}\neq 0\\
&c_{31}\neq0\text{ oder }c_{32}\neq 0\end{aligned}\ \right)\ \Bigr\},\\[12pt]
T_2&=\Bigl\{C\in k^{3\times 3}\Bigm| \left(\begin{aligned}&c_{33}=c_{23}=0\\
&c_{32}\cdot c_{13}\neq0\\
&c_{31}\neq0\text{ oder }c_{21}\neq 0\end{aligned}\ \right)\ \Bigr\},\\[12pt]
T_3&=\Bigl\{C\in k^{3\times 3}\Bigm| \left(\begin{aligned}&c_{33}=c_{13}=0\\
&c_{23}\cdot c_{31}\neq0\\
&c_{32}\neq0\text{ oder }c_{12}\neq 0\end{aligned}\ \right)\ \Bigr\}.\end{align*}
$U_z^{(2)}$ ist irreduzibel, da $T_1$ irreduzibel und dicht in $U_z^{(2)}$ ist.

$\mathbf{B_9=(1,2)+(1,1)_{\lambda_1}+ (1,0)}$\\[4pt]
$U_z^{(1)}=\bigcup_{i=1}^4T_i$. Dabei ist
\begin{align*}
T_1&=\Bigl\{C\in k^{3\times 3}\Bigm|c_{13}\neq0\text{ und }\left(\begin{aligned}
&c_{22}-\frac{c_{12}}{c_{13}}\cdot c_{23}\neq0\\
\text{ oder } &c_{31}+\frac{c_{33}}{c_{13}}\cdot(c_{32}-c_{11}-\frac{c_{12}}{c_{13}}\cdot c_{33})\neq0
\end{aligned}\right)
\Bigr\},\\[8pt]
T_2&=\Bigl\{C\in k^{3\times 3}\Bigm|
c_{13}=0\text{ und } c_{23}\cdot c_{12}\neq 0\Bigr\},\\[8pt]
T_3&=\Bigl\{C\in k^{3\times 3}\Bigm|\left(
\begin{aligned}
&c_{12}=c_{13}=0\\
&c_{23}\cdot c_{33}\neq0\\
&c_{11}-c_{32}+\frac{c_{22}}{c_{23}}\cdot c_{33}\neq0\\
\end{aligned}\right)\Bigr\},\\[8pt]
T_4&=\Bigl\{C\in k^{3\times 3}\Bigm|\left(
\begin{aligned}
&c_{13}=c_{23}=0\\
&c_{33}\neq0\\
&c_{22}\neq0\text{ oder } c_{12}\neq0\\
\end{aligned}\right)\Bigr\}.
\end{align*}
$U_z^{(2)}=\Bigl\{C\in k^{3\times 3}\Bigm|\text{$\left(\begin{aligned}
&c_{23}\neq 0\\
&c_{12}=c_{13}=c_{11}-c_{32}+\frac{c_{22}}{c_{23}}\cdot c_{33}=0\\
&c_{33}\neq 0\text{ oder } c_{31}\neq 0
\end{aligned}\:\right)$}\:\Bigr\}$.

$U_z^{(1)}$ ist irreduzibel als offene Teilmenge von $F_z=k^{3\times3}$. $U_z^{(2)}$ ist irreduzibel: Man verifiziert, daß das Polynom $X_1X_2-X_2X_3+X_4X_5\in k[X_1,\ldots,X_5]$ irreduzibel ist; also ist die Menge aller $C\in F_z$ mit $c_{11}\cdot c_{23}-c_{23}\cdot c_{32}+c_{22}\cdot c_{33}=c_{13}=c_{12}=0$ irreduzibel; $U_z^{(2)}$ ist offen in dieser Menge, also irreduzibel.

$\mathbf{B_{10}=(0,1)+(1,1)_{\lambda_1}+ (2,1)}$\\[4pt]
Zu $B_9$ dualer Fall.

$\mathbf{B_{11}=(1,1)_{\lambda_1}+(2,2)_{\lambda_1}}$\\[4pt]
$U_z^{(1)}=\bigcup_{i=1}^4T_i$. Dabei ist
\begin{align*}
T_1&=\Bigl\{C\in k^{3\times 3}\Bigm|c_{12}\neq0\text{ und }\left(\begin{aligned}
&c_{23}+\frac{c_{13}}{c_{12}}\cdot(c_{33}-c_{22}-\frac{c_{13}}{c_{12}}\cdot c_{32})\neq0\\
\text{ oder } &c_{31}+\frac{c_{32}}{c_{12}}\cdot(c_{33}-c_{11}-\frac{c_{13}}{c_{12}}\cdot c_{32})\neq0
\end{aligned}\right)
\Bigr\},\\[8pt]
T_2&=\Bigl\{C\in k^{3\times 3}\Bigm|
c_{12}=0\text{ und } c_{13}\cdot c_{32}\neq 0\Bigr\},\\[8pt]
T_3&=\Bigl\{C\in k^{3\times 3}\Bigm|\left(
\begin{aligned}
&c_{12}=c_{32}=0\\
&(c_{11}-c_{22})\cdot (c_{33}-c_{22})-c_{13}\cdot c_{31}\neq0\\
&c_{13}\neq0
\end{aligned}\right)\Bigr\},\\[8pt]
T_4&=\Bigl\{C\in k^{3\times 3}\Bigm|\left(
\begin{aligned}
&c_{12}=c_{13}=0\\
&(c_{11}-c_{22})\cdot (c_{11}-c_{33})-c_{23}\cdot c_{32}\neq0\\
&c_{32}\neq0
\end{aligned}\right)\Bigr\}.
\end{align*}
$U_z^{(2)}=U_1\cup U_2$. Dabei ist
\begin{align*}
U_1&=\Bigl\{C\in k^{3\times3}\Bigm|\left(\begin{aligned}&(c_{11}-c_{22})\cdot(c_{33}-c_{22})-c_{13}\cdot c_{31}=c_{12}=c_{32}=0\\&c_{13}\neq0\text{ oder }c_{23}\cdot(c_{11}-c_{22})\neq0\end{aligned}\right)\Bigr\},\\[12pt]
U_2&=\Bigl\{C\in k^{3\times3}\Bigm|\left(\begin{aligned}&(c_{11}-c_{22})\cdot(c_{11}-c_{33})-c_{23}\cdot c_{32}=c_{12}=c_{13}=0\\&c_{32}\neq0\text{ oder }c_{31}\cdot(c_{11}-c_{22})\neq0\end{aligned}\right)\Bigr\}.
\end{align*}
$U_z^{(3)}=\Bigl\{\left[\begin{smallmatrix}a&0&0\\x&a&y\\z&0&a\end{smallmatrix}\right]\Bigm|z\neq0\text{ oder } y\neq0\Bigr\}$.

$U_z^{(1)}$ und $U_z^{(3)}$ sind offenbar irreduzibel. $U_1$ und $U_2$ sind die irreduziblen Komponenten von $U_z^{(2)}$: Man verifiziert, daß $U_1$ und $U_2$ abgeschlossen sind in $U_z^{(2)}$. Die Irreduzibilität von $U_1$ folgt aus der Irreduzibilität des Polynoms $(X_1-X_2)(X_3-X_2)-X_4X_5$, genauso verifiziert man, daß $U_2$ irreduzible ist. Da $U_z^{(2)}=U_1\cup U_2$ gilt und die $U_i$ echte Teilmengen von $U_z^{(2)}$ sind, folgt die Behauptung.\enlargethispage{\baselineskip}

$\mathbf{B_{12}=(1,2)+(2,1)}$\\[4pt]
Für $C\in k^{3\times3}$ sei
{\[
E_C:=\left[\begin{array}{cc}c_{12}&c_{13}\\c_{22}&c_{23}\\c_{22}&c_{12}\\c_{23}&c_{13}\\c_{21}-c_{32}&c_{11}-c_{33}\end{array}\right].
\]}Dann gilt:

$U_z^{(1)}=\Bigl\{C\in k^{3\times 3}\Bigm|\operatorname{Rang}E_C=2\Bigr\}$.

$U_z^{(2)}=\Bigl\{C\in k^{3\times 3}\Bigm|\operatorname{Rang}E_C=1\text{ und }\Bigl[\begin{smallmatrix}c_{12}&c_{13}\\c_{22}&c_{23}\end{smallmatrix}\Bigr]\neq0\Bigr\}$.

$U_z^{(3)}=\Bigl\{\biggl[\begin{smallmatrix}x&0&0\\y&0&0\\z&y&x\end{smallmatrix}\biggr]\in k^{3\times 3}\Bigm| z\neq0\Bigr\}$.

Wieder ist die Irreduzibilität von $U_z^{(1)}$ und $U_z^{(3)}$ klar. Um die Irreduzibilität von $U_z^{(2)}$ zu zeigen betrachte die Teilmenge $X$ aller $C\in U_z^{(2)}$ mit $c_{13}\neq 0$. Dann besteht $X$ aus allen Matrizen der Form
{\[\left[\begin{array}{ccc}c_{11}&a&c_{13}\\c_{21}&\frac{a^2}{c_{13}}&a\\c_{31}&(c_{21}-\frac{a}{c_{13}}\cdot(c_{11}-c_{33}))&c_{33}\end{array}\right],\ a\in k,\ c_{13}\neq0.\]}$X$ ist irreduzibel als Bild eines Morphismus $k^*\times k^5\longrightarrow k^{3\times3}$. Man verifiziert außerdem, daß $X$ dicht ist in $U_z^{(2)}$, also ist auch $U_z^{(2)}$ irreduzibel.

$\mathbf{B_{13}=(2,3)+(1,0)}$\\[4pt]
$U_z^{(1)}=\Bigl\{C\in k^{3\times3}\Bigm| c_{13}\neq0\text{ oder } c_{23}\neq0\text{ oder } c_{33}\neq0\Bigr\}$.

$U_z^{(1)}$ ist irreduzibel.

$\mathbf{B_{14}=(0,1)+(3,2)}$\\[4pt]
Zu $B_{13}$ dualer Fall.

$\mathbf{B_{15}=(1,1)_{\lambda_1}+(1,1)^2_{\lambda_2}}$\\[4pt]
$U_z^{(1)}=G_z*\bigl(T_1\cup T_2\bigr)$. Dabei ist
\begin{align*}
T_1&=\Bigl\{\left[\begin{smallmatrix}\lambda&0&y_1\\0&\mu&y_2\\x_1&x_2&t\end{smallmatrix}\right]\in k^{3\times 3}\Bigm| (x_1,y_1)\neq0\neq(x_2,y_2)\text{ und }\lambda\neq\mu\:\Bigr\},\\[12pt]
T_2&=\Bigl\{\left[\begin{smallmatrix}\lambda&0&y_1\\1&\lambda&y_2\\x_1&x_2&t\end{smallmatrix}\right]\in k^{3\times 3}\Bigm|x_2\neq0\text{ oder }y_1\neq0\:\Bigr\}.\end{align*}
$U_z^{(2)}=G_z*\bigl(U_1\cup U_2\bigr)$. Dabei ist
\begin{align*}
U_1&=\Bigl\{\biggl[\begin{smallmatrix}\lambda&0&0\\0&\lambda&y\\x&0&t\end{smallmatrix}\biggr]\in k^{3\times 3}\Bigm| x\neq0\neq y\:\Bigr\},\\[12pt]
U_2&=\Bigl\{\biggl[\begin{smallmatrix}\lambda&0&0\\1&\lambda&y\\x&0&t\end{smallmatrix}\biggr]\in k^{3\times 3}\Bigm|x\neq0\text{ oder }y\neq0\:\Bigr\}.
\end{align*}
$U_z^{(1)}$ ist irreduzibel. Um zu zeigen, daß auch $U_z^{(2)}$ irreduzibel ist, verifiziere: $U_z^{(2)}=\overline{G_z*U_2}$.

$\mathbf{B_{16}=(3,3)_{\lambda_1}}$\\[4pt]
$U_z^{(3)}=\Bigl\{\left[\begin{smallmatrix}a&0&0\\b&a&0\\c&b&a\end{smallmatrix}\right]\in k^{3\times3}\Bigr\}$.

$U_z^{(2)}=\Bigl\{\left[\begin{smallmatrix}x&0&0\\u&y&0\\z&v&x\end{smallmatrix}\right]\in k^{3\times3}\bigm|x\neq y\text{ oder }u\neq v\Bigr\}$.

$U_z^{(1)}=k^{3\times 3}-\bigl(U_z^{(2)}\cup (U_z^{(3)}\bigr)$.

$U_z^{(1)},U_z^{(2)}$ und $U_z^{(3)}$ sind irreduzibel.

$\mathbf{B_{17}=(1,1)_{\lambda_1}+(2,2)_{\lambda_2}}$\\[4pt]
$U_z^{(2)}=\Bigl\{\left[\begin{smallmatrix}a&0&0\\b&a&d\\c&0&e\end{smallmatrix}\right]\in k^{3\times3}\Bigm|d\neq0\text{ oder }c\neq0\Bigr\}$.\\

$U_z^{(1)}=\Bigl\{C\in k^{3\times3}\Bigm|(c_{31},c_{32})\neq0 \text{ oder } (c_{13},c_{23})\neq0\Bigr\}-U_z^{(2)}$.

$U_z^{(1)}$ und $U_z^{(2)}$ sind irreduzibel.

$\mathbf{B_{18}=(1,1)_{\lambda_1}+(1,1)_{\lambda_2}+(1,1)_{\lambda_3}}$\\[4pt]
$U_z^{(1)}=\Bigl\{C\Bigm|\text{$ \left(\begin{aligned}&(c_{12},c_{21})\neq0\neq(c_{13},c_{31})\\
\text{ oder }&(c_{12},c_{21})\neq0\neq(c_{23},c_{32})\\
\text{ oder }&(c_{13},c_{31})\neq0\neq(c_{23},c_{32})
\end{aligned}\right)$}\Bigr\}$.

$U_z^{(1)}$ ist irreduzibel.

\paragraph{Anmerkungen und Beispielrechnungen}\ \newline
Wir wollen in diesem Abschnitt einige Anleitungen zur Verifikation der Beschreibungen der $U_z^{(i)}$ aus dem letzten Abschnitt geben.

{\bfseries Zu den Bahnen $\mathbf B_1$ bis $\mathbf B_6$}\newline
Um die Richtigkeit der Ergebnisse für $B_1,B_2,B_5,B_6$ zu überprüfen, beachte man: Ist $m=(m_1,0,m_3)\in R_3(\underline{d})$ gegeben, so ist $\operatorname{rep} m$ genau dann eine unzerlegbare $Q_3$-Darstellung, wenn $\operatorname{rep} (m_1,m_3)$ eine unzerlegbare $Q_2$-Darstellung ist. Verwende dann die Klassifikation der Kro\-necker-Moduln (\cref{grund:kron:satz1}).

Im Fall der Bahn $B_3$ geht man folgendermaßen vor: Betrachte den Repräsentanten \[z=\Bigl(\left[\begin{smallmatrix}1&0&0\\0&0&0\\0&0&0\end{smallmatrix}\right],
\left[\begin{smallmatrix}0&0&0\\1&0&0\\0&0&0\end{smallmatrix}\right]\Bigr)\in B_3.\]
Sei $A\in k^{3\times3}=F_z$. Dann gilt: Ist $a_{23}\neq0$ oder $a_{33}\neq0$, so existiert ein $g\in G_z$ mit 
\[g*A=\left[\begin{smallmatrix} *&*&0\\ *&*&0\\0&0&*\end{smallmatrix}\right].\]
Es folgt unmittelbar, daß $(z_1,z_2,A)$ zerlegbar ist. Eine Fallunterscheidung nach dem Rang von $A$ und Berechnung der $Q_3$-Endomorphismenringe liefert dann die gewünschten Ergebnisse.

{\bfseries Zu den Bahnen $\mathbf B_7$ und $\mathbf B_8$}\newline
Betrachte den in Tabelle \ref{rech:unz:tab1} angegebenen Repräsentanten
\[z=(z_1,z_2)=\Bigl(\left[\begin{smallmatrix}1&0&0\\0&1&0\\0&0&0\end{smallmatrix}\right],
\left[\begin{smallmatrix}0&0&0\\1&0&0\\0&0&0\end{smallmatrix}\right]\Bigr)\in B_7.\]
Sei $A\in F_z=k^{3\times 3}$. Ist $a_{33}\neq0$, so ist $A$ modulo der $G_z$-Operation konjugiert zu einer Matrix der Form
\[\widetilde{A}=\left[\begin{smallmatrix}{}*&*&0\\ *&*&0\\ 0&0&*\end{smallmatrix}\right].\]
$(z_1,z_2,\widetilde{A})$ ist eine zerlegbare $Q_3$-Darstellung; damit ist auch $(z_1,z_2,A)$ zerlegbar. Bei den Rechnungen darf man also $a_{33}=0$ annehmen.

{\bfseries Zur Bahn $\mathbf B_{11}$}\newline
Sei wie in der Tabbelle
\[z=(z_1,z_2)=\Bigl(\left[\begin{smallmatrix}1&0&0\\0&1&0\\0&0&1\end{smallmatrix}\right],
\left[\begin{smallmatrix}0&0&0\\1&0&0\\0&0&0\end{smallmatrix}\right]\Bigr)\in B_{11}.\]
Der $Q_2$-Endomorphismenring von $z$ ist gegeben durch
\[\End_{Q_2}(z)=\bigl\{\bigl(S(a,\ldots,e),S(a,\ldots,e)\bigm| S(a,\ldots,e)\bigr)=\left[\begin{smallmatrix}a&0&0\\b&a&e\\d&0&c\end{smallmatrix}\right]\bigr\}\subseteq k^{3\times 3}\times k^{3\times 3}.\]
Wir fassen $\End_{Q_2}(z)$ als Unteralgebra von $k^{3\times 3}$ auf. $G_z=\operatorname{Aut}(z)\subseteq \operatorname{GL}_3$ operiert auf $F_z$ durch Konjugation.

Sei $R\subseteq \End_{Q_2}(z)$ eine Unteralgebra. Dann sieht man leicht:
\[R\text{ ist lokal }\Leftrightarrow R\text{ enthält nur Matrizen der Form } \left[\begin{smallmatrix}a&0&0\\b&a&e\\d&0&a\end{smallmatrix}\right],\]
Sei $C\in F_z$. Für $S(a,\ldots,e)\in\End_{Q_2}(z)$ gilt: $S(a,\ldots,e)$ liegt genau dann in $\End_{Q_3}(z,C)$, wenn folgende Gleichungen erfüllt sind:

\begin{center}
$(\mathcal{L}_C)$
\counterwithout{equation}{chapter}
\begin{tabular}{p{4cm}p{7cm}}
{\begin{enumerate}
\item $b\cdot c_{12}+d\cdot c_{13}=0$
\item $e\cdot c_{12}+c'\cdot c_{13}=0$
\item $b\cdot c_{12}+e\cdot c_{32}=0$
\item $d\cdot c_{12}+c'\cdot c_{32}=0$
\end{enumerate}}&
{\begin{enumerate}\setcounter{enumi}{4}
\item $b\cdot (c_{11}-c_{22})+e\cdot e_{31}-d\cdot c_{23}=0$
\item $d\cdot (c_{11}-c_{33})+c'\cdot c_{31}-b\cdot c_{32}=0$
\item $e\cdot (c_{33}-c_{22})+b\cdot c_{13}-c'\cdot c_{23}=0$
\item $d\cdot c_{13}-e\cdot e_{32}=0$.
\end{enumerate}}
\end{tabular}
\end{center}
Dabei sei $c':=c-a$. Es folgt: $(z_1,z_2,C)$ ist genau dann eine unzerlegbare $Q_3$-Darstellung, wenn das homogene lineare Gleichungssystem $(\mathcal{L}_C)$ keine Lösung $(b,c',d,e)$ besitzt mit $c'\neq0$. Ist $X\subseteq k^4$ der Lösungsraum von $\mathcal{L}_C$, so gilt $\dim \End_{Q_3}(z_1,z_2,C)=\dim X +1$.

Sei jetzt $A\in F_z$ gegeben. Wir geben die zur Berechnung der $U_z^{(i)}$ nötigen Fallunterscheidungen an:

\textbf{1. Fall:} $\mathbf {a_{12}\neq0.}$ Man verifizert, daß ein $g\in G_z$ existiert mit \[g*A=\widetilde{A}=\text{\small$\left[\begin{array}{ccc}\widetilde{a}_{11}&\widetilde{a}_{12}&0\\ \widetilde{a}_{21}&\widetilde{a}_{22}&\widetilde{a}_{23}\\ \widetilde{a}_{31}&0&\widetilde{a}_{33}\end{array}\right]$}.\]
Dabei ist $\widetilde{a}_{12}=a_{12}\neq0$ und die $\widetilde{a}_{ij}$ sind rationale Funktionen in den Koeffizienten von $A$. Sei $S=S(a,\ldots,e)\in \End_{Q_2}(z_1,z_2)$. Die Gleichungen $1$ bis $8$ liefern
\[S\in\End_{Q_3}(z_1,z_2,\widetilde{A})\Leftrightarrow b=d=e=0\text{ und } c'\cdot\widetilde{a}_{31}=c'\cdot\widetilde{a}_{23}=0.\]
Gilt $\widetilde{a}_{31}=\widetilde{a}_{23}=0$, so existiert eine Lösung $(b,c',d,e)$ von $\mathcal{L}_{\widetilde{A}}$ mit $c'\neq 0$, also ist $(z_1,z_2,\widetilde{A})$ und damit auch $(z_1,z_2,A)$ zerlegbar. Ist $\widetilde{a}_{31}\neq 0$ oder $\widetilde{a}_{23}\neq0$, so ist der Lösungsraum von $L_{\widetilde{A}}$ trivial, also gilt $\widetilde{A}\in U_z^{(1)}$ und damit ist auch $A\in U_z^{(1)}$. Die Bedingungen ($\widetilde{a}_{31}\neq 0$ oder $\widetilde{a}_{23}\neq0$) und $a_{12}\neq0$ sind offene Bedingungen an $A$, wir erhalten also eine offene Teilmenge von $U_z^{(1)}$.

\textbf{2. Fall:} $\mathbf{a_{12}=0\text{\textbf{ und }}a_{13}\neq0.}$ Betrachte die Matrix $E_A$ definiert durch
\[E_A=\text{\small$\left[\begin{array}{cc}a_{11}-a_{22}&a_{31}\\a_{13}&a_{33}-a_{22}\end{array}\right]$}.\]
Man liest aus den Gleichungen $1$ bis $8$ ab: Ist $S=S(a,\ldots,e)\in\End_{Q_2}(z_1,z_2)$ gegeben, so gilt:
\[S\in\End_{Q_3}(z_1,z_2,A)\Leftrightarrow d=c'=0,\  E_A\cdot\text{\small$\left[\begin{array}{c}b\\e\end{array}\right]$}=0\text{ und }e\cdot a_{32}=b\cdot a_{32}=0.\]
Also ist $(z_1,z_2,A)$ iin diesem Fall stets unzerlegbar. Ist $a_{32}\neq 0$, so ist der Lösungsraum von $\mathcal{L}_A$ trivial, also $A\in U_z^{(1)}$. Ist $a_{32}=0$, so gilt $A\in U_z^{(1)}$ falls $\operatorname{Rang}E_A=2$; ansonsten ist $A\in U_z^{(2)}$. Wir erhalten also lokal abgeschlossene Teilmengen von $U_z^{(1)}$ und $U_z^{(2)}$.

\textbf{3. Fall:} $\mathbf{a_{12}=a_{13}=0\text{\textbf{ und }}a_{32}\neq0.}$ Betrachte die Matrix $E_A$ definiert durch
\[E_A=\text{\small$\left[\begin{array}{cc}a_{11}-a_{22}&-a_{23}\\-a_{32}&a_{11}-a_{33}\end{array}\right]$}.\]
Man liest aus den Gleichungen $1$ bis $8$ ab: Ist $S=S(a,\ldots,e)\in\End_{Q_2}(z_1,z_2)$ gegeben, so gilt:
\[S\in\End_{Q_3}(z_1,z_2,A)\Leftrightarrow e=c'=0\text{ und } E_A\cdot\text{\small$\left[\begin{array}{c}b\\d\end{array}\right]$}=0.\]
Wie im zweiten Fall ist dann $(z_1,z_2,A)$ unzerlegbar und man erhält je nach Rang von $E_A$ eine lokal abgeschlossene Teilmenge von $U_z^{(1)}$ oder $U_z^{(2)}$.

\textbf{4. Fall:} $\mathbf{a_{12}=a_{13}=a_{32}=0.}$ Betrachte die Matrix $E_A$ definiert durch
\[E_A=\text{\small$\left[\begin{array}{cccc}a_{11}-a_{22}&0&-a_{23}&a_{31}\\
0&a_{31}&a_{11}-a_{33}&0\\0&-a_{23}&0&a_{33}-a_{22}\end{array}\right]$}.\]
Die Gleichung (1) bis (8) ergeben in diesem Fall für $S=S(a,\ldots,e)\in\End_{Q_2}(z)$:
\[S\in\End_{Q_3}(z_1,z_2,A)\Leftrightarrow E_A\cdot\text{\small$\left[\begin{array}{c}b\\c'\\d\\e\end{array}\right]=0$}.\]
Dann gilt: $(z_1,z_2,A)$ ist genau dann eine unzerlegbare $Q_3$-Darstellung, wenn der Kern der Matrix $E_A$ in $\bigl\langle e_1,e_3,e_4\bigr\rangle$ enthalten ist (dabei sind die $e_i$ die Einheitsvektoren in $k^4$). Man löst dieses Gleichungssystem unter den gegebenen Nebenbedingungen durch langwierige, aber einfache Rechnungen.

Jeder betrachtete Fall liefert dann gewisse lokal abgeschlossene Teilmengen der $U_z^{(i)}$. Die im letzten Abschnitt angegebene Beschreibung entstand daraus durch weitere Vereinfachungen.

{\bfseries Zur Bahn $\mathbf B_{15}$}\newline
Sei
\[z=(z_1,z_2)=\Bigl(\left[\begin{smallmatrix}1&0&0\\0&1&0\\0&0&1\end{smallmatrix}\right],
\left[\begin{smallmatrix}0&0&0\\0&0&0\\0&0&1\end{smallmatrix}\right]\Bigr).\]
Der $Q_2$-Endomorphismenring von $z$ ist dann als Unteralgebra von $k^{3\times 3}$ gegeben durch Matrizen der Form
{\small\[\left[\begin{smallmatrix}*&*&0\\ *&*&0\\ 0&0&*\end{smallmatrix}\right].\]}$G_z\subseteq \operatorname{GL}_3$ operiert auf $F_z=k^{3\times 3}$ durch Konjugation. $\operatorname{GL}_2$ läßt sich in offensichtlicher Weise als Untergruppe von $G_z$ auffassen und operiert durch Konjugation auf dem oberen linken $2\times 2$-Block der Matrizen in $k^{3\times 3}$. Man erhält daraus, daß $F_z$ die $G_z$-Saturierung der Menge $X$ aller Matrizen der folgenden Form ist:
{\small \[
\left[\begin{smallmatrix}\lambda&0&*\\1&\lambda&*\\ *&*&*\end{smallmatrix}\right],\ 
\left[\begin{smallmatrix}\lambda&0&*\\0&\lambda&*\\ *&*&*\end{smallmatrix}\right],\ 
\left[\begin{smallmatrix}\lambda&0&*\\0&\mu&*\\ *&*&*\end{smallmatrix}\right],\ \lambda\neq \mu
.\]}Zur Berechnung der $U_z^{(i)}$ als $G_z$-Saturierung gewisser Teilmengen genügt es also, Matrizen in $X$ zu betrachten. Man fährt dannn wie im Fall $B_{11}$ mit der Standardmethode fort.
\numberwithin{zaehler}{section}
\setcounter{zaehler}{0}

\section{Die Irreduzibilität der $U_{Q_3}^{(i)}(3,3)$}\label{rech:irr}
\numberwithin{zaehler}{section}
\setcounter{zaehler}{0}
Wir behalten die Notationen vom Beginn des letzten Abschnitts bei und wollen mit Hilfe der dort erzielten Ergebnisse zeigen, daß die $U^{(i)}=U_3^{(i)}(3,3)$ irreduzibel sind. Wir formulieren unsere Vorgehensweise in zwei Lemmata. Das erste erlaubt es uns, in unserer speziellen Situation eine irreduzible Komponente $X$ von $U^{(i)}$ zu bestimmen; das zweite beschreibt das Verfahren, mit dem wir zeigen, daß schon $X=U^{(i)}$ gilt.
\begin{lemma}\label{rech:irr:lemma1}
Sei $G$ eine zusammenhängende algebraische Gruppe. Seien $X,Y$ $G$-Varietäten und sei $f:X\longrightarrow Y$ ein $G$-äquivarianter Morphismus. Sei $y\in Y$, so daß $U:=f^{-1}(Gy)$ offen in $X$ ist und so daß die Faser $f^{-1}(y)$ (nichtleer und) irreduzibel ist. Dann ist der Abschluß $\overline{U}$ von $U$ in $X$ eine irreduzible Komponente von $X$.
\end{lemma}
\begin{bew} Da $U$ offen ist, genügt es zu zeigen, daß $U$ irreduzibel ist. Betrachte dazu die von $f$ induzierte Abbildung $\widetilde{f}:U\longrightarrow Gy$. $\widetilde{f}$ ist ein surjektiver Morphismus von $G$-Varietäten, dabei operiert $G$ transitiv auf $Gy$. Resultate aus der algebraischen Geometrie über generische Flachheit implizieren, daß $\widetilde{f}$ generisch offen ist (vgl. \cite[14.2]{Eisen} für den Satz von der generischen Freiheit von Grothendieck und \cite[Exercise III.9.1]{Ha} für die Offenheit flacher Morphismen). Da die Operation von $G$ auf der Basis $Gy$ transitiv ist, folgt, daß $\widetilde{f}$ schon offen ist. Da $f^{-1}(y)$ irreduzibel ist, sind alle Fasern von $\widetilde{f}$ irreduzibel, außerdem ist $Gy$ irreduzibel, da $G$ zusammenhängend ist. Daraus folgt die Irreduzibilität von $U$.
\end{bew}

\begin{lemma}\label{rech:irr:lemma2}
Sei $G$ eine zusammenhängende algebraische Gruppe und sei $f:X\longrightarrow Y$ ein Morphismus von $G$-Varietäten. Auf $X$ operiere außerdem eine zusammenhängende algebraische Gruppe $G'$. Es sei $Z\subseteq X$ eine irreduzible Komponente von $X$.

Sei weiter $y\in Y$ und sei $f^{-1}(y)=T_1\cup\ldots\cup T_r$ eine Zerlegung in irreduzible Teilmengen. Für alle $i=1,\ldots,r$ gebe es eine Teilmenge $U_i\subseteq Z$ von $Z$, so daß $U_i$ offen ist in $X$ und so daß gilt: $(G'*T_i)\cap U_i\neq\emptyset$. Dann gilt schon $f^{-1}(Gy)\subseteq Z$.
\end{lemma}
\begin{bew}
Da $G'$ zusammenhängend und $T_i$ irreduzibel ist, ist auch $G'*T_i$ irreduzibel. Man erhält für alle $i$: 
\[T_i\subseteq G'*T_i=\overline{(G'*T_i)\cap U_i}^{(G'*T_i)}\subseteq \overline{(G'*T_i)\cap U_i}^{(X)}\subseteq Z.\]
Damit folgt $f^{-1}(y)\subseteq Z$. Da $G$ zusammenhängend ist und $Z$ eine irreduzible Komponente von $X$ ist, ist $Z$ $G$-stabil. Es folgt
\[f^{-1}(Gy)=Gf^{-1}(y)\subseteq GZ\subseteq Z.\]
\end{bew}

Wir konkretisieren unsere Strategie. Sei
\[\pi=\pi^{(3,3)}:R_3(3,3)\longrightarrow R_2(3,3)\]
die Projektion auf die ersten beiden Komponenten. Für $j\in\mathbb{N}$ sei $\pi_j:U^{(j)}\longrightarrow R_2(3,3)$ die Einschränkung von $\pi$ auf $U^{(j)}$. $\pi_j$ ist dann ein $H_2\times G$-äquivarianter Morphismus. Für jeden $H_2\times G$-Orbit $B$ in $R_2(3,3)$ und für alle $z\in B$ gilt
\[\pi_j^{-1}(z)=U_z^{(j)}\text{ und }\pi_j^{-1}(B)=\bigcup_{z\in B} U_z^{(j)}.\]
Für $i=1,\ldots,18$ sei $z^{(i)}\in B_i$ der in Tabelle \ref{rech:unz:tab1} angegebene Repräsentant des $H_2\times G$-Orbits $B_i$ von $R_2(3,3)$.

Wir zeigen die Irreduzibilität von $U^{(j)}$ in zwei Schritten:
\begin{enumerate}
\item Bestimme mit \cref{rech:irr:lemma1} eine irreduzible Komponente $X$ von $U^{(j)}$.
\item Zeige für alle $i=1,\ldots,18$: $\pi_j^{-1}(B_i)\subseteq X$ (denn dann folgt $X=U^{(j)}$ und damit die Behauptung). In fast allen Fällen (bis auf den Sonderfall $j=3$, $i=12$) verwenden wir folgende Argumentation: Bestimme für jede irreduzible Komponente $T$ von $\pi_j^{-1}(z^{(i)})$ eine Teilmenge $U\subseteq X$, so daß $U$ offen ist in $U^{(j)}$ und so daß gilt:
\[(H_3\times G)*\pi_j^{-1}(z^{(i)})\cap U\neq \emptyset.\]
Nach \cref{rech:irr:lemma2} folgt dann  $\pi_j^{-1}(B_i)\subseteq X$.
\end{enumerate}

$U^{(1)}$ ist irreduzibel, denn $U^{(1)}$ besteht aus allen Punkten $m\in R_3(3,3)$ mit $\End m=k$, ist also offen in $R_3(3,3)$ nach \cref{grund:var:prop3}. 

Gemäß Tabelle \ref{rech:unz:tab1} ist weiter $U_z^{(j)}=\emptyset$ für alle $z\in R_2(3,3)$ und für alle $j\geq4$; demzufolge ist $U^{(j)}=\emptyset$ für alle $j\geq4$. 

Es bleibt also zu zeigen, daß $U^{(2)}$ und $U^{(3)}$ irreduzibel sind.

\textbf{$\mathbf{U^{(2)}}$ ist irreduzibel}\newline
Sei $U_0:=\pi_2^{-1}(B_{17})$. Gemäß Tabelle \ref{rech:unz:tab1} ist $\pi_2^{-1}(B_{18})=\emptyset$, also gilt 
\[\pi_2^{-1}(B_{17})=\pi_2^{-1}(B_{17}\cup B_{18}).\]
$B_{17}\cup B_{18}$ ist die Vereinigung aller $H_2\times G$-Orbiten der Dimension $\geq 17$ in $R_2(3,3)$, also offen in $R_2(3,3)$ (vgl. \cref{grund:var:prop1}). Damit ist auch $U_0$ offen in $U^{(2)}$. Da außerdem $\pi_2^{-1}(z^{(17)})$ irreduzibel ist (vgl. wieder Tabelle \ref{rech:unz:tab1}), ist nach \cref{rech:irr:lemma1} $X:=\overline{U_0}$ eine irreduzible Komoponente von $U^{(2)}$.

Ein Blick auf Tabelle \ref{rech:unz:tab1} verrät, daß $\pi_2^{-1}(B_i)\subseteq X$ nur zu zeigen ist in den Fällen 
\[i=16,15,12,11,10,9,8,7.\]

$\mathbf{i=16.}$ $\pi_2^{-1}(z^{(16)})$ ist irreduzibel. Sei $U:=U_0$. $U$ ist offen in $U^{(2)}$ und eine Teilmenge von $X$. Mit Hilfe der Beschreibung von $U_{z^{(16)}}^{(2)}$ im letzten Abschnitt sieht man:
\[t=(t_1,t_2,t_3)=\Bigl(\left[\begin{smallmatrix}
1&0&0\\0&1&0\\0&0&1
\end{smallmatrix}\right]
\left[\begin{smallmatrix}
0&0&0\\1&0&0\\0&1&0
\end{smallmatrix}\right]
\left[\begin{smallmatrix}
0&0&0\\0&1&0\\0&0&0
\end{smallmatrix}\right]
\Bigr)\in\pi_2^{-1}(z^{(16)}).\]
Es existiert ein Element $h\in H_3$ mit
\[h*t=(t_1,\left[\begin{smallmatrix}
0&0&0\\1&1&0\\0&1&0
\end{smallmatrix}\right],t_3).\]
Die Jordansche Normalform von $\left[\begin{smallmatrix}
0&0&0\\1&1&0\\0&1&0
\end{smallmatrix}\right]$ ist $\left[\begin{smallmatrix}
0&0&0\\1&0&0\\0&0&1
\end{smallmatrix}\right]$, also folgt 
\[\pi_2(h*t)\in (1,1)_{\lambda_1}+(2,2)_{\lambda_2}=B_{17}.\]
Erhalte $h*t\in\pi_2^{-1}(B_{17})=U$, also $(H_3\times G)*\pi_2^{-1}(z^{(16)})\cap U\neq\emptyset$ wie gewünnscht.

$\mathbf{i=15,12,11.}$ Nach dem Fall $i=16$ wissen wir, daß 
\[U:=\pi_2^{-1}(B_{17}\cup B_{16})=\pi_2^{-1}(B_{18}\cup B_{17}\cup B_{16})\]
eine Teilmenge von $X$ ist. $U$ ist offen in $U^{(2)}$, da $B_{18}\cup B_{17}\cup B_{16}$ die Vereinigung aller $H_2\times G$-Bahnen der Dimension $\geq16$ in $R_2(3,3)$ ist (vgl. \cref{grund:var:prop1}). Geeignete $t=(t_1,t_2,t_3)\in\pi_2^{-1}(z^{(i)})$ und $h\in H_3$ liefern dann $h*t\in U$. Beachte noch, daß der Fall $i=11$ der einzige ist, bei dem $\pi_2^{-1}(z^{(i)})$ nicht irreduzibel ist.

$\mathbf{i=10,9,8,7.}$ Die Fälle $i\geq 11$ zeigen: 
\[U:=\{m\in U^{(2)}\mid \operatorname{Rang} m_i=3\text{ für ein $i$ }\}\subseteq X.\]
Setze dazu
\[T:=\pi_2^{-1}\Bigl(\bigcup_{i=11}^{18}B_i\Bigr)\subseteq X.\]
Aus den Listen \ref{rech:deg:tab1} und \ref{rech:unz:tab1} liest man ab: Ist $m\in U^{(2)}$ gegeben, und hat die $Q_2$-Darstellung $\pi_2(m)$ keinen nichtregulären direkten Summanden, so muß schon $m\in T$ gelten. Offenbar hat jede Darstellung $(m_1,m_2)\in R_2(3,3)$ mit $\operatorname{Rang} m_1=3$ oder $\operatorname{Rang} m_2=3$ keinen direkten nichtregulären Summanden. Also gilt:
\[U':=\{m\in U^{(2)}\mid \text{$\operatorname{Rang} m_1=3$ oder $\operatorname{Rang} m_2=3$}\}\subseteq T\subseteq X.\]
Es folgt $U\subseteq H_3 U'\subseteq U^{(2)}$. $U$ ist offen in $U^{(2)}$; geeignete $t=(t_1,t_2,t_3)\in\pi_2^{-1}(z^{(i)})$ und $h\in H_3$ liefern wieder $h*t\in U$.

Beachte noch: Die Fälle $i=9,10$ verhalten sich dual zueinander. Betrachte dazu den Automorphismus
\[D:R_2(3,3)\longrightarrow R_2(3,3),\ (A,B,C)\mapsto (A^T,B^T,C^T).\]
Dann gilt $D(\pi_2^{-1}(B_9))=\pi_2^{-1}(B_{10})$. Außerdem gilt offenbar $D(U)=U$. Da $U$ dicht in $X$ liegt, folgt $D(X)=X$. Aus $\pi_2^{-1}(B_9)\subseteq X$ folgt also sofort $\pi_2^{-1}(B_{10})\subseteq X$

\textbf{$\mathbf{U^{(3)}}$ ist irreduzibel}\newline
Sei $U_0:=\pi_3^{-1}(B_{16})=\pi_3^{-1}((3,3)_{\lambda_1})$ und sei $X:=\overline{U_0}$. Wie im Fall $U^{(2)}$ zeigt man mit Hilfe von \cref{rech:irr:lemma2}, daß $X$ eine irreduzible Komponente von $U^{(3)}$ ist. Gemäß Tabelle \ref{rech:unz:tab1} ist $\pi_3^{-1}(B_i)\subseteq X$ zu zeigen in den Fällen
\[i=12,11,7,6,5,4,3.\]
Beachte: Nach Tabelle \ref{rech:unz:tab1} sind alle $\pi_3^{-1}(z^{(i)})$ sind irreduzibel.

$\mathbf{i=12.}$ Hier müssen wir von unserer Standardargumentation abweichen. Betrachte für $\epsilon\in k$ den Punkt
\[t_\epsilon=\Bigl(\left[\begin{smallmatrix}
1&0&0\\0&\epsilon&0\\0&0&1
\end{smallmatrix}\right]
\left[\begin{smallmatrix}
0&0&0\\1&0&0\\0&1&0
\end{smallmatrix}\right]
\left[\begin{smallmatrix}
0&0&0\\0&0&0\\1&0&0
\end{smallmatrix}\right]
\Bigr)\in R_3(3,3).\]
Für $\epsilon\neq0$ existiert ein $g\in G$ mit
\[g*t=\Bigl(\left[\begin{smallmatrix}
1&0&0\\0&1&0\\0&0&1
\end{smallmatrix}\right]
\left[\begin{smallmatrix}
0&0&0\\1&0&0\\0&1&0
\end{smallmatrix}\right]
\left[\begin{smallmatrix}
0&0&0\\0&0&0\\\epsilon&0&0
\end{smallmatrix}\right]
\Bigr)\in \pi_3^{-1}(B_{16})\subseteq X.\]
Außerdem ist $t_0\in \pi_3^{-1}(B_{12})$. Da $X$ abgeschlossen ist in $U^{(3)}$ folgt schon $t_0\in X$. Man verifiziert sofort mit der Beschreibung von $U_{z^{(12)}}^{(3)}$ in Abschnitt \ref{rech:unz}, daß $\pi_3^{-1}(z^{(12)})$ im $H_3$-Orbit von $t_0$ enthalten ist. Damit folgt
\[\pi_3^{-1}(B_{12})\subseteq (H_3\times G)*\pi_3^{-1}(z^{(12)})\subseteq (H_3\times G)*t_0\subseteq X.\]

$\mathbf{i=11,6.}$ Sei $U:=U_0$. Man findet wieder $t\in\pi_3^{-1}(z^{(i)})$ und $h\in H_3$ mit $h*t\in U$.

$\mathbf{i=7,5,4,3}$ Die Fälle $i= 11,6$ zeigen: 
\[U:=\{m\in U_3^{(3)}\mid \operatorname{Rang} m_i=3\text{ für ein $i$ }\}\subseteq X.\]
(die Argumentation ist dieselbe wie für $U^{(2)}$ im Fall $i=10,9,8,7$). Man findet dann wie üblich $t\in\pi^{-1}(z^{(i)})$ und $h\in H_3$ mit $h*t\in U$. 

Wieder verhalten sich die Fälle $i=4,3$ dual zueinander.

\paragraph{Zwei Bemerkungen zu unzerlegbaren $Q_3$-Darstellungen der Dimension $\mathbf{(4,4)}$}\ \newline
Wir haben gesehen, daß der Endomorphismenring aller unzerlegbaren $Q_3$-Darstellungen zum Dimensionsvektor $(3,3)$ höchstens $3$-dimensional ist, also die Dimension des Endomorphismenrings der unzerlegbaren regulären $Q_2$-Darstellungen $R_{3,\lambda}$ zum Dimensionsvektor $(3,3)$ nicht übersteigt. Diese Aussage ist für den Dimensionsvektor $(4,4)$ falsch; betrachte etwa die $Q_3$-Darstellung gegeben durch die Matrizen
{\small\[\left[\begin{array}{cccc}1&0&0&0\\0&1&0&0\\0&0&1&0\\0&0&0&1\end{array}\right],
\left[\begin{array}{cc|cc}0&0&1&0\\0&0&0&1\\\hline0&0&0&0\\0&0&0&0\end{array}\right],
\left[\begin{array}{cc|cc}0&0&0&0\\0&0&1&0\\\hline0&0&0&0\\0&0&0&0\end{array}\right].\]}Man verifiziert, daß der Endomorphismenring dieser Darstellung lokal und $6$-dimensional ist.

Wir wollen noch zeigen, daß im allgemeinen die Varietäten unzerlegbarer $Q_3$-Darstellungen nicht irreduzibel sind. Genauer werden wir zeigen, daß $U_3^{(4)}(4,4)$ nicht irreduzibel ist und dazu die Beschreibung der irreduziblen Komponenten der Varietät $\operatorname{Alg}(4)$ der vierdimensionalen Algebren in \cite{Ga1} verwenden. Wie wiederholen die wesentlichen Definitionen.

Sei $A$ eine endlichdimensionale $k$-Algebra und $1=a_1,\ldots, a_n$ eine $k$-Basis von $A$. Für alle $i,j;\nu$ seien die $\alpha_{ij}^{(\nu)}\in k$ die Strukturkonstanten von $A$, es gilt also
\[a_i\cdot a_j=\sum_{\nu=1}^\nu\alpha_{ij}^{(\nu)}\cdot a_\nu \text{ für alle }i,j.\]
Mit $\operatorname{Mod}^d_A$ bezeichnen wir die Varietät der $d$-dimensionalen $A$-Moduln; sie besteht aus allen Punkten $m\in \prod_{i=1}^n k^{d\times d}$, so daß die folgenden Gleichungen erfüllt sind:
\[m_1=E_d,\ m_i\cdot m_j=\sum_{\nu=1}^n\alpha_{ij}^{(\nu)}\cdot m_\nu\text{ für alle }i,j.\]
$\operatorname{GL}_d$ operiert auf $\operatorname{Mod}^d_A$ durch Konjugation. Für einen Punkt $m\in\operatorname{Mod}^d_A$ bezeichnen wir mit $\operatorname{mod} m$ den zugehörigen Modul. Die Zuordnung $m\mapsto\operatorname{mod}m$ induziert eine Bijektion zwischen den $\operatorname{GL}_d$-Orbiten von $\operatorname{Mod}^d_A$ und den Isomorphieklassen $d$-dimensionaler $A$-Moduln. 

Betrachte die affine Varietät 
\[
\operatorname{Bil}(n):=\{b:k^n\times k^n\longrightarrow k^n\mid b \text{ ist $k$-bilinear }\}\simeq \mathbb{A}^{n^3}.
\]
Sei $\operatorname{Alg}(n)$ die Menge aller $b\in\operatorname{Bil}(n)$, so daß $b$ eine assoziative $k$-Algebrenstruktur mit $1$ auf $k^n$ definiert; $\operatorname{Alg}(n)$ ist dann eine lokal abgeschlossene, affine Untervarietät von $\operatorname{Bil}(n)$ (vgl. \cite[§1]{CB3}). $\operatorname{GL}_n$ operiert auf $\operatorname{Alg}(n)$; für $g\in \operatorname{GL}_n$ und $b\in \operatorname{Alg}(n)$ ist
\[\bigl(g*b\bigr)(v,w):=g\cdot b(g^{-1}v, g^{-1}w)\text{ für alle } v,w\in k^n.\]
Für einen Punkt $a\in\operatorname{Alg}(n)$ bezeichnen wir mit $\operatorname{alg}a$ die zugehörige Algebra. Die Zuordnung $a\mapsto\operatorname{alg}a$ induziert eine Bijektion zwischen den $\operatorname{GL}_n$-Orbiten von $\operatorname{Alg}(n)$ und den Isomorphieklassen $n$-dimensionaler Algebren. 

Wir benötigen den folgenden Spezialfall von \cite[Lemma 2.1]{Bong2}. Unsere Vektorbündel sind Einschränkungen der dort definierten Bündel.
\begin{lemma}\label{rech:irr:lemma3}
Sei $i\in\mathbb N$. Betrachte die Varietäten
\begin{align*}
X^{(i),d}&:=\{m\in \operatorname{Mod}^d\mid \dim \End_A(\operatorname{mod} m)=i\},\\
Y^{(i),d}&:=\{(m,f)\in X^{(i),d}\times k^{d\times d}\mid f\in\End_A(\operatorname{mod} m)\}.
\end{align*}
Sei $\phi:Y^{(i),d}\longrightarrow X^{(i),d}$ die Projektion. Dann ist $\phi$ ein Vektorbündel vom Rang $i$. Die Faser über $m\in X^{(i),d}$ ist der $i$-dimensionale $k$-Vektorraum $\End_A(\operatorname{mod} m)$.
\end{lemma}
Mit Hilfe dieses Lemmas können wir die Zuordnung $M\mapsto \End_A(M)$ auf den Teilstücken $X^{(i),d}$ lokal als Morphismus von Varietäten deuten (vgl. die Diskussion in \cite[Abschnitt 2]{Bong2}):
\begin{lemma}\label{rech:irr:lemma4}
Sei $m\in X^{(i),d}$. Dann existieren eine offene Umgebung $U$ von $m$ in $X^{(i),d}$ und ein Morphismus $\epsilon:U\longrightarrow \operatorname{Alg}(i)$, so daß für alle $u\in U$ gilt: $\operatorname{alg}\epsilon(u)\simeq \End_A(\operatorname{mod}u)$.
\end{lemma}
\begin{bew} Betrachte das Vektorbündel $\phi:Y^{(i),d}\longrightarrow X^{(i),d}$ aus dem letzten Lemma. Seien eine offene Umgebung $U$ von $m$ in $X^{(i),d}$ und eine lokale Trivialisierung
\[\tau:U\times k^i\stackrel{\sim}{\longrightarrow} \phi^{-1}(U)\]
gegeben. Für $u\in U$ sei $\tau_u:k^i\stackrel{\sim}{\longrightarrow} \phi^{-1}(u)=\End_A(\operatorname{mod} u)$ der von $\tau$ induzierte Isomorphismus von $k$-Vektorräumen. Sei $e_k$ der $k$-te Einheitsvektor von $k^i$. Definiere $\epsilon:U\longrightarrow \operatorname{Alg}(i)=\mathbb{A}^{n^{3}}$ durch 
\[\epsilon(u):=\Bigl(\tau_u^{-1}\bigl((\tau_u(e_k)\cdot\tau_u(e_l)\bigr)_\nu\Bigr)_{k,l;\nu=1,\ldots,n}.\]
$\epsilon$ ist ein Morphismus und die $\tau_u^{-1}\bigl((\tau_u(e_k)\cdot\tau_u(e_l)\bigr)_\nu$ sind einfach die Strukturkonstanten von $\End_A(\operatorname{mod} u)$ bezüglich der $k$-Basis gegeben durch die $\tau_u(e_k)$. Damit folgt die Behauptung.
\end{bew}
Analoge Aussagen gelten für Varietäten von Köcherdarstellungen; wir wollen hier nicht näher auf die Details der engen Verwandtschaft zwischen den Varietäten von endlichdimensionalen Moduln einer Köcheralgebra $A=kK$ und den Varietäten der Darstellungen von $K$ eingehen und verweisen den Leser auf \cite[Example 5.18]{Bong5}.

Wir führen unsere Bemerkung auf die folgende Aussage aus \cite[Abschnitt 5]{Ga1} zurück; dort sind sämtliche irreduzible Komponenten von $\operatorname{Alg}(4)$ beschrieben.
\begin{prop}\label{rech:irr:prop1}
Die kommutativen $4$-dimensionalen Algebren bilden eine irreduzible Komponente von $\operatorname{Alg}(4)$ und $\operatorname{Alg}(4)$ besitzt 5 irreduzible Komponenten.
\end{prop}
Wir können jetzt zeigen, daß $U_3(4,4)$ nicht irreduzibel ist. Betrachte die Projektion
\[\pi:R_3(4,4)\longrightarrow R_2(4,4),\ (A,B,C)\mapsto (A,B)).\]
Die regulären $Q_2$-Darstellungen zum Dimensionsvektor $(4,4)$, deren unzerlegbare direkte Summanden zu paarweise verschiedenen Elementen in $\mathbb{P}^1$ korrespondieren, bilden die offene Schicht in $R_2(4,4)$. Alle diese Darstellungen haben $4$-dimensionale kommutative Endomorphismenringe (man verifiziert diese Aussage zum Beispiel, indem man Entartungen von $G(4,4)$-Orbiten in $R_2(4,4)$ mit Hilfe von exakten Sequenzen von $Q_2$-Darstellungen betrachtet). 

Sei $U:=\pi^{-1}(R_2^{(4)}(4,4))\subseteq R_3(4,4)$ das Urbild der offenen Schicht von $R_2(4,4)$. Setze 
\[V:=U\cap U_3^{(4)}(4,4).\]
Dann ist $V$ offen und nichtleer in $U_3^{(4)}(4,4)$. Für alle $v\in V$ ist $\End_{Q_3}(\operatorname{rep}v)$ als Unteralgebra von $\End_{Q_2}(\operatorname{rep}\pi v)$ kommutativ. Angenommen, $U_3^{(4)}(4,4)$ ist irreduzibel. Wir behaupten, daß dann schon für alle $u\in U_3^{(4)}(4,4)$ gilt: $\End_{Q_3}(\operatorname{rep} u)$ ist kommutativ. Sei dazu $W$ eine offene Umgebung von $u$ in $U_3^{(4)}(4,4)$ und sei 
\[\epsilon:W\longrightarrow \operatorname{Alg}(4),\ w\mapsto \End_{Q_2}(\operatorname{rep} w)\]
der Morphismus aus \cref{rech:irr:lemma4}. Es sei $C$ die irreduzible Komponente aller kommmutativen Algebren in $\operatorname{Alg}(4)$ (vgl. \cref{rech:irr:prop1}). Da $U_3^{(4)}(4,4)$ nach Annahme irreduzibel ist, ist $W\cap V$ offen und dicht in $W$. Außerdem gilt $\epsilon (W\cap V)\subseteq C$. Es folgt:
\[\epsilon(W)=\epsilon(\overline{W\cap V})\subseteq \overline{\epsilon(W\cap V)}\subseteq C.\]
Um zu zeigen, daß $U_3^{(4)}(4,4)$ nicht irreduzibel ist, genügt es also, zu zeigen: Es existiert eine unzerlegbare $Q_3$-Darstellung $U$ mit nichtkommutativem, 4-dimensionalem Endomorphismenring. Betrachte dazu für $\lambda\in k^{*}-\{1\}$ die Darstellung $U_\lambda$ gegeben durch die Matrizen
{\small\[\left[\begin{array}{cccc}1&0&0&0\\0&1&0&0\\0&0&1&0\\0&0&0&1\end{array}\right],
\left[\begin{array}{cc|cc}0&0&0&0\\1&0&0&0\\\hline0&0&0&0\\0&0&1&0\end{array}\right],
\left[\begin{array}{cc|cc}0&0&0&0\\0&0&0&0\\\hline1&0&0&0\\0&\lambda^{-1}&0&0\end{array}\right].\]}
Dann gilt
\[\End_{Q_3}(U_\lambda)\simeq k\langle X,Y\rangle/\langle X^2,Y^2,XY-\lambda YX\rangle\]
und die Behauptung folgt.

\section{Zur Erreichbarkeit unzerlegbarer $Q_3$-Darstellungen der Dimension $\leq 6$}\label{rech:erreich}
\numberwithin{zaehler}{section}
\setcounter{zaehler}{0}
Wir wollen zeigen: Ist $M$ eine $Q_3$-Darstellungen, $\dim M\leq6$, so existiert eine exakte Sequenz
\[0\longrightarrow U_1\longrightarrow M\longrightarrow U_2\longrightarrow 0\]
von $Q_3$-Darstellungen, so daß die die $U_i$ unzerlegbar sind und ein $U_i$ einfach ist (vgl. dazu die Definition von Erreichbarkeit in \ref{baum:erreich:def2}). Die Aussage ist offenbar invariant unter der $H_3$-Operation; hat man die Aussage also für ein $h*M$, $h\in H_3$ bewiesen, so schon für $M$. 

Das wichtigste Hilfsmittel für unsere Rechnungen sind exakte Sequenzen von $Q_2$-Darstellungen. Das folgende Lemma ist der Liste von Erweiterungen von unzerlegbaren $Q_2$-Darstellung aus \cite[Teil 5]{Bong2} entnommen.\newpage
\begin{lemma}\label{rech:erreich:lemma1}
Sei $M$ eine $Q_2$-Darstellung, so daß eine der folgenden Bedingungen erfüllt ist:
\begin{itemize}
\item[$\cdot$]$M\simeq \bigoplus\limits_{i=1}^rR_{n_i,\lambda_i}, \lambda_i\neq\lambda_j$ für $i\neq j$.
\item[$\cdot$]$M\simeq P(k)\oplus P(l)$ oder $M\simeq I(k)\oplus I(l)$, $k,l\geq1$.
\item[$\cdot$]$M\simeq P(k)\oplus R_{n,\lambda}$ oder $M\simeq I(k)\oplus R_{n,\lambda}$, $k\geq1,n\geq0$; dabei $R_{0,\lambda}:=0$.
\item[$\cdot$]$M$ unzerlegbar.
\end{itemize}
Dann existiert eine exakte Sequenz von $Q_2$-Darstellungen
\[0\longrightarrow U_1\longrightarrow M\longrightarrow U_2\longrightarrow 0,\]
so daß die $U_i$ unzerlegbar sind und ein $U_i$ einfach ist.
\end{lemma}

Sei $\underline{d}\in\mathbb{Z}^2$ ein Dimensionsvektor von $Q_2$; betrachte die Projektion $\pi^{(\underline{d})}:R_3(\underline{d})\longrightarrow R_2(\underline{d})$ auf die ersten beiden Komponenten. Man verifiziert sofort folgende Aussage:
\begin{lemma}\label{rech:erreich:lemma2} Sei $m\in R_3(\underline{d})$. Es gebe eine exakte Sequenz
\[0\longrightarrow U_1\longrightarrow \operatorname{rep}\pi^{(\underline{d})}(m)\longrightarrow U_2\longrightarrow 0\]
von $Q_2$-Darstellungen, so daß die $U_i$ unzerlegbar sind und ein $U_i$ einfach ist. Dann existiert eine exakte Sequenz von $Q_3$-Darstellungen
\[0\longrightarrow V_1\longrightarrow \operatorname{rep} m\longrightarrow V_2\longrightarrow 0,\]
so daß die $V_i$ unzerlegbar sind und ein $V_i$ einfach ist.
\end{lemma}

\paragraph{Unzerlegbare $\mathbf{Q_3}$-Darstellungen zum Dimensionsvektor (3,3) sind Erweiterungen eines Einfachen und eines Unzerlegbaren}\ \newline
In diesem Paragraphen sei $\pi:=\pi^{(3,3)}:R_3(3,3)\longrightarrow R_2(3,3)$ die Projektion auf die ersten beiden Komponenten. Wir verwenden die zu Beginn des Abschnitts \ref{rech:unz} vereinbarten Notationen. Sei $m\in R_3(3,3)$ gegeben, $\operatorname{rep} m$ sei unzerlegbar. Betrachte die $H_2\times G$-Bahnen $B_1,\ldots,B_{18}$ von $R_2(3,3)$ (vgl. Tabelle \ref{rech:deg:tab1}). Es ist $\pi(m)\in B_i$ für ein $i\in\{1,\ldots,18\}$. Sei wie üblich $z^{(i)}\in B_i$ der in Tabelle \ref{rech:unz:tab1} angegebene Repräsentant von $B_i$. Wir merken noch an, daß $U_{z^{(i)}}$ stabil ist unter der Operation von $G_{z^{(i)}}$ und unter der Operation der Untergruppe 
\[H':=\{h\in H_3\mid h=\left[\begin{smallmatrix}1&0&0\\ 0&1&0\\ a&b&c\end{smallmatrix}\right]\}\]
von $H_3$.

\textbf{1.Fall: $\mathbf{\pi(m)}$ ist regulär.} Wir nehmen an, daß $\pi(m)=z^{(i)}$ gilt. Also ist
\[m_1=\left[\begin{smallmatrix}1&0&0\\0&1&0\\0&0&1\end{smallmatrix}\right].\]
Da ${\pi(m)}$ regulär ist, zeigt Tabelle \ref{rech:deg:tab1}: Es ist $\pi(m)\in B_i$ mit
\[i\in\{18,17,16,15,11,6\}.\]
Die Fälle 
\begin{align*}
\pi(m)\in B_{18}&=(1,1)_{\lambda_1}+(1,1)_{\lambda_2}+(1,1)_{\lambda_3},\\ 
\pi(m)\in B_{17}&=(1,1)_{\lambda_1}+(2,2)_{\lambda_2},\\
\pi(m)\in B_{16}&=(3,3)_{\lambda_1}
\end{align*}
handelt man mit Hilfe des $2$-Kronecker-Falls ab (vgl. Lemma \ref{rech:erreich:lemma1} und \ref{rech:erreich:lemma2}). Ist $\pi(m)\in B_6=(1,1)_{\lambda_1}^3$, so ist $m_2=(z^{(6)})_2=0$, also ist $(m_1,m_3)$ unzerlegbar, und man erhält die gewünschte exakte Sequenz nach Vertauschen der zweiten und dritten Komponente mit Lemma \ref{rech:erreich:lemma2}.

Sei also $i\in \{15,11\}$. Wir dürfen $\operatorname{Rang} m_3=1$ annehmen. $(m_1,m_3)$ ist nämlich eine reguläre $Q_2$-Darstellung (denn $m_1$ ist die Einheitsmatrix); nach den bereits ausgeschlossenen Fällen $i=18,17,16,6$ können wir 
\[(m_1,m_3)\in B_{15}\cup B_{11}=\Bigl((1,1)_{\lambda_1}+(1,1)_{\lambda_2}^2\Bigr)\cup\Bigl((1,1)_{\lambda_1}+(2,2)_{\lambda_1}\Bigr)\]
annehmen. Anders ausgedrückt: Die Jordansche Normalform von $m_3$ ist vom Typ
\[\left[\begin{smallmatrix}\lambda&0&0\\1&\lambda&0\\0&0&\lambda\end{smallmatrix}\right]\text{ oder }\left[\begin{smallmatrix}\lambda&0&0\\0&\lambda&0\\0&0&\mu\end{smallmatrix}\right],\ \lambda\neq \mu.\]
Ersetze $m_3$ durch $m_3-\lambda m_1$; modulo der $H'$-Operation auf $U_{z^{(i)}}$ hat also $m_3$ Rang 1.

Wir müssen also die Behauptung nur noch zeigen für alle $m\in U_{z^{(i)}}$ mit $\operatorname{Rang} m_3=1$ und $i\in \{15,11\}$.

Sei zunächst $i=15$. Also ist
\[\pi(m)=z^{(15)}=\Bigl(\left[\begin{smallmatrix}1&0&0\\0&1&0\\0&0&1\end{smallmatrix}\right],
\left[\begin{smallmatrix}0&0&0\\0&0&0\\0&0&1\end{smallmatrix}\right]\Bigr).\]
Wir wiederholen die Beschreibung von $U_{z^{(15)}}$ aus Abschnitt \ref{rech:unz}: Es ist $U_{z^{(15)}}=G_{z^{(15)}}*X$; dabei ist $X$ die Menge aller Matrizen vom Typ
{\small
\begin{gather*}
\left[\begin{smallmatrix}\lambda&0&y_1\\0&\mu&y_2\\x_1&x_2&t\end{smallmatrix}\right],\  (x_1,y_1)\neq0\neq(x_2,y_2)\text{ und }\lambda\neq\mu;\ 
\left[\begin{smallmatrix}\lambda&0&y_1\\1&\lambda&y_2\\x_1&x_2&t\end{smallmatrix}\right],\ x_2\neq0\text{ oder }y_1\neq0;\\
\biggl[\begin{smallmatrix}\lambda&0&0\\0&\lambda&y\\x&0&t\end{smallmatrix}\biggr],\  x\neq0\neq y;\ 
\biggl[\begin{smallmatrix}\lambda&0&0\\1&\lambda&y\\x&0&t\end{smallmatrix}\biggr],\ x\neq0\text{ oder }y\neq0.
\end{gather*}
}Da die Operation von $G_{z^{(15)}}$ Ränge erhält, dürfen wir $m_3\in X$ annehmen. Man verifiziert, daß wegen $\operatorname{Rang} m_3=1$ schon
\[m_3=\biggl[\begin{smallmatrix}\lambda&0&0\\1&\lambda&y\\x&0&t\end{smallmatrix}\biggr],\ x\neq0\text{ oder }y\neq0\]
gelten muß. Sei zunächst $x\neq0=y$. Dann läßt sich $m_3$ mit der $H'$-Operation auf $U_{z^{(15)}}$ auf die Form
\[m_3=\biggl[\begin{smallmatrix}0&0&0\\1&0&0\\x&0&1\end{smallmatrix}\biggr]\]
bringen. Berechnen der Jordanschen Normalform von $m_3$ zeigt dann $(m_1,m_3)\in B_{17}$; in diesem Fall wurde die Behauptung schon bewiesen. Ist $x\neq0\neq y$, so reduziert man genauso auf den Fall $B_{16}$; ist $y\neq0=x$, so reduziert man wieder auf den Fall $B_{17}$.

Sei jetzt $i=11$. Also ist
\[\pi(m)=z^{(11)}=\Bigl(\left[\begin{smallmatrix}1&0&0\\0&1&0\\0&0&1\end{smallmatrix}\right],
\left[\begin{smallmatrix}0&0&0\\1&0&0\\0&0&0\end{smallmatrix}\right]\Bigr).\]
Schreibe $m_3$ in der Form
\[m_3=\left[\begin{smallmatrix}*&*&y_1\\ *&*&y_2\\x_1&x_2&*\end{smallmatrix}\right].\]
Da $\operatorname{rep} m$ unzerlegbar ist, muß $(x_1,x_2)\neq0$ oder $(y_1,y_2)\neq0$ sein. Wir dürfen annehmen, daß schon $(x_1,x_2)\neq0$ ist (man kann $(x_1,x_2)\neq0$ erreichen, in dem man zur dualen Darstellung übergeht und anschließend so mit einem geeigneten Element $g\in G$ konjugiert, daß die ersten beiden Komponenten der konjugierten Darstellungen wieder $z^{(11)}$ ergeben). "`Eliminieren der dritte Spalte von $m$"' liefert eine Unterdarstellung $U=\operatorname{rep} u$ von $\operatorname{rep} m$ gegeben durch die Matrizen
\[u=(u_1,u_2,u_3)=\Bigl(\left[\begin{smallmatrix}1&0\\ 0&1\\0&0\end{smallmatrix}\right],\left[\begin{smallmatrix}0&0\\ 1&0\\0&0\end{smallmatrix}\right],\left[\begin{smallmatrix}*&*\\  *&*\\x_1&x_2\end{smallmatrix}\right]\Bigr).\]
Man verifiziert leicht, daß aus der Voraussetzung $x_1\neq0$ oder $x_2\neq0$ die Unzerlegbarkeit von $U$ folgt.

\textbf{2.Fall: $\mathbf{\pi(m)}$ ist nicht regulär.} Wir bemerken zunächst: Ist $\operatorname{Rang} m_j=3$ für ein $j$, so gilt modulo der $H_3$-Operation (vertausche die erste und $j$-te Komponente): $\pi(m)$ ist regulär und man ist im 1. Fall.

Wir dürfen also $\operatorname{Rang} m_j\leq2$ für alle $j$ annehmen. Man verifiziert: Ist $\operatorname{Rang} (m_j)\leq1$ für alle $j$, so ist $m$ zerlegbar. Sei also ohne Einschränkung $\operatorname{Rang} m_1=2$. Aufgrund der Unzerlegbarkeit von $m$ gilt außerdem
\[\Bild m_1+\Bild m_2+\Bild m_3=k^3,\]
also dürfen wir ohne Einschränkung annehmen, daß schon $\Bild m_1+\Bild m_2=k^3$ gilt. Insgesamt erhalten wir: $\pi(m)$ ist nicht regulär und die projektiv-einfache $Q_2$-Darstellung $E_y$ ist kein direkter Summand von $\pi(m)$. Ein Blick auf Tabelle \ref{rech:deg:tab1} zeigt, daß nur die Fälle $\pi(m)\in B_i$ mit $i\in\{13,12,9\}$ auftreten. 

Sei $i=13$, also $B_i=(2,3)+(1,0)$. Dann erhält man aus der zu dieser Zerlegung gehörigen spaltenden Sequenz von $Q_2$-Darstellungen mit Mittelterm $\pi(m)$ die gesuchte Sequenz von $Q_3$-Darstellungen. 

Sei $i=12$. Dann ist $B_i=(1,2)+(2,1)$ und wir nehmen an:
\[\pi(m)=z^{(12)}=\Bigl(\left[\begin{smallmatrix}1&0&0\\0&0&0\\0&0&1\end{smallmatrix}\right], \left[\begin{smallmatrix}0&0&0\\1&0&0\\0&1&0\end{smallmatrix}\right]\Bigr).\]
Modulo der $H'$-Operation auf $U_{z^{(12)}}$ dürfen wir annehmen:
\[m_3=\text{\small$\left[\begin{array}{ccc}a_{11}&a_{12}&a_{13}\\ a_{21}&a_{22}&a_{23}\\a_{31}&0&0\end{array}\right]$}.\]
Betrachte das Polynom
\[P(X,Y):=\operatorname{det}(X\cdot m_1+Y\cdot m_2+m_3)\in k[X,Y].\]
Wir dürfen annehmen, daß $P=0$ ist, sonst könnten wir modulo der $H'$-Operation $\operatorname{Rang} m_3=3$ annehmen, und wir wären im 1. Fall. Wir erhalten durch Berechnen der Determinante und Koeffizientenvergleich die Gleichungen
\[a_{22}=a_{13}=a_{12}\cdot a_{21}=a_{11}\cdot a_{23}=a_{23}+a_{12}=a_{31}\cdot a_{12}\cdot a_{23}=0.\ \ \ (*)\]
Wir betrachten die folgenden Fälle:
\begin{itemize}
\item[$\cdot$]$a_{12}\neq0$. Formel $(*)$ zeigt $a_{11}=0$. "`Eliminieren der zweiten Zeile von $m$"' liefert einen Faktormodul $C=\operatorname{rep} c$ von $\operatorname{rep} m$ gegeben durch die Matrizen
\[c=(c_1,c_2,c_3)=\Bigl(\bigl[\begin{smallmatrix}1&0&0\\ 0&0&1\end{smallmatrix}\bigr],\bigl[\begin{smallmatrix}0&0&0\\ 0&1&0\end{smallmatrix}\bigr],\bigl[\begin{smallmatrix}0&a_{12}&0\\  a_{31}&0&0\end{smallmatrix}\bigr]\Bigr)\]
(beachte $a_{13}=a_{11}=0$). Man verifiziert, daß $C$ unzerlegbar ist.
\item[$\cdot$]$a_{12}=0$. Formel $(*)$ zeigt, daß dann schon
\[m_3=\text{\small$\left[\begin{array}{ccc}a_{11}&0&0\\ a_{21}&0&0\\a_{31}&0&0\end{array}\right]$}\]
gilt. Ein kurzer Blick auf die Beschreibung von $U_{z^{(12)}}$ in Abschnitt \ref{rech:unz} zeigt, daß schon
\[m_3=\text{\small$\left[\begin{array}{ccc}0&0&0\\ 0&0&0\\z&0&0\end{array}\right]$}\]
gelten muß für ein $z\neq0$. "`Eliminieren der zweiten Spalte"' von $m$ liefert den Untermodul $U=\operatorname{rep} u$ von $\operatorname{rep} m$ mit
\[u=(u_1,u_2,u_3)=\Bigl(\left[\begin{smallmatrix}1&0\\ 0&0\\0&1\end{smallmatrix}\right],\left[\begin{smallmatrix}0&0\\ 1&0\\0&0\end{smallmatrix}\right],\left[\begin{smallmatrix}0&0\\  0&0\\z&0\end{smallmatrix}\right]\Bigr).\]
Man verifiziert wieder leicht, daß $U$ unzerlegbar ist.
\end{itemize}

Sei jetzt $i=9$. Also ist $B_i=(1,2)+(1,1)_{\lambda_1}+(1,0)$ und wir dürfen annehmen:
\[\pi(m)=z^{(9)}=\Bigl(\left[\begin{smallmatrix}1&0&0\\0&0&0\\0&1&0\end{smallmatrix}\right],
\left[\begin{smallmatrix}0&0&0\\1&0&0\\0&0&0\end{smallmatrix}\right]\Bigr).\]
Modulo der $H'$-Operation auf $U_{z^{(9)}}$ dürfen wir weiter annehmen:
\[m_3=\text{\small$\left[\begin{array}{ccc}0&a_{12}&a_{13}\\ 0&a_{22}&a_{23}\\a_{31}&a_{32}&a_{33}\end{array}\right]$}.\]
Betrachte das Polynom
\[P(X,Y):=\operatorname{det}(X\cdot m_1+Y\cdot m_2+m_3)\in k[X,Y].\]
$P=0$ würde die Gleichungen
\[a_{23}=a_{13}=a_{12}\cdot a_{33}=a_{22}\cdot a_{33}=0.\]
liefern; man verifiziert aber sofort anhand der Beschreibung von $U_{z^{(9)}}$ aus  Abschnitt \ref{rech:unz}, daß dann $m\not\in U_{z^{(9)}}$ ist. Also gilt $P\neq0$, mithin existieren $\lambda,\mu\in k$, so daß $\lambda\cdot m_1+\mu\cdot m_2+m_3$ invertierbar ist, modulo der $H_3$-Operation hat also $m_3$ Rang 3 und man ist im 1. Fall.

\paragraph{Unzerlegbare $\mathbf{Q_3}$-Darstellungen zu den Dimensionsvektoren (2,2), (2,3), (2,4) sind Erweiterungen eines Einfachen und eines Unzerlegbaren}\ \newline
Wir skizzieren lediglich die nötigen Fallunterscheidungen.

Sei $m\in R_3(2,2)$ gegeben, $m$ sei unzerlegbar. Ist $(m_1,m_2)$ eine reguläre $Q_2$-Darstellung, so gilt 
\[(m_1,m_2)\in (2,2)_{\lambda_1}\cup (1,1)_{\lambda_1}^2\cup \bigl((1,1)_{\lambda_1}+(1,1)_{\lambda_2}\bigr),\]
und durch Reduktion auf den $2$-Kronecker Fall (vgl. Lemma \ref{rech:erreich:lemma1} und \ref{rech:erreich:lemma2} sowie die Diskussion des regulären Falls für den Dimensionsvektor (3,3)) erhält man die gesuchten exakten  Sequenzen. Sei also $(m_1,m_2)$ nicht regulär. In Analogie zur Diskussion des nicht-regulären Falles für den Dimensionsvektor (3,3) dürfen wir annehmen, daß die projektiv-einfache $Q_2$-Darstellung kein direkter Summand von $(m_1,m_2)$ ist. Also ist $(m_1,m_2)\in (1,0)+(1,2)$, und die zugehörige spaltende Sequenz von $2$-Kronbecker-Darstellungen liefert die gewünschte Sequenz.

Sei $m\in R_3(2,3)$, $m$ sei unzerlegbar. Sei zunächst $\operatorname{Rang} m_1=2$. Wegen der Unzerlegbarkeit von $m$ dürfen wir annehmen:
\[\Bild m_1+\Bild m_2=k^3.\]
Also sind die beiden einfachen $Q_2$-Darstellungen keine direkten Summanden von $(m_1,m_2)$. Es folgt
\[(m_1,m_2)\in \bigl((1,2)+(1,1)_{\lambda_1}\bigr)\cup (2,3).\]
In beiden Fällen liefert Reduktion auf den $2$-Kronecker-Fall (vgl. Lemma \ref{rech:erreich:lemma1} und \ref{rech:erreich:lemma2}) die Behauptung. Modulo der $H_3$-Operation dürfen wir also $\operatorname{Rang}m_i\leq1$ annehmen. Man findet wegen der Unzerlegbarkeit von $m$ leicht Elemente $h\in H_3$ und $g\in G(2,3)$ mit \[\operatorname{Rang}(h*g*m)_1=2.\]
Die Behauptung folgt mit dem schon bewiesenen Fall.

Sei $m\in R_3(2,4)$, $m$ sei unzerlegbar. Wir dürfen aufgrund der Unzerlegbarkeit von $m$ nach eventuellem Vertauschen der $m_i$ annehmen, daß $\operatorname{Rang} m_1=2$ ist. Also ist $m_1$ injektiv und die injektiv-einfache $Q_2$-Darstellung $E_x$ ist kein direkter Summand von $m$.

\textbf{1. Fall: $\mathbf{\operatorname{\mathbf  {Bild}} m_1+\operatorname{\mathbf  {Bild}} m_2=k^4}$.} In dem Fall ist $(m_1,m_2)\in (1,2+(1,2)$, und die Behauptung folgt durch Reduktion auf den $2$-Kronecker-Fall.

\textbf{2. Fall: $\mathbf{\operatorname{\mathbf  {Bild}} m_1+\operatorname{\mathbf  {Bild}} m_2\neq k^4}$.} Dann dürfen wir ohne Einschränkung annehmen, daß 
\[\dim \Bild m_1+\dim\Bild m_2=3\]
gilt (wäre nämlich $\dim \Bild m_1+\dim\Bild m_2\leq2$, so folgte $\Bild m_2\subseteq\Bild m_1$; wegen der Unzerlegbarkeit von $m$ wäre dann $\Bild m_1+\Bild m_3=k^4$, und wir wären im ersten Fall). Also können nur die Fälle\enlargethispage{\baselineskip}
\[(m_1,m_2)\in (2,3)+(0,1)\text{ oder } (m_1,m_2)\in (1,1)+(1,2)+(0,1)\]
auftreten. Den ersten Fall handelt man durch Reduktion auf den $2$-Kronecker ab. Im 2. Fall erhalten wir nach Basiswahl und modulo der $H_2$-Operation auf $2$-Kronecker-Darstellungen
\[(m_1,m_2)=\Biggl(\left[\begin{smallmatrix}1&0\\0&1\\0&0\\0&0\end{smallmatrix}\right],
\left[\begin{smallmatrix}0&0\\0&0\\0&1\\0&0\end{smallmatrix}\right]\Biggr).\]
Aufgrund der Unzerlegbarkeit von $m$ muß gelten:
\[m_3=\left[\begin{smallmatrix}*&*\\ *&*\\ a&b\\ c&d\end{smallmatrix}\right],\ \neq0\text{ oder }d\neq0.\]
Wir dürfen annehmen: $a=c=0$; sonst verifiziert man, daß ein $\lambda\in k$ existiert mit
\[\Bild m_1+\Bild(\lambda m_2+m_3)=k^4,\]
und man ist im 1. Fall. Man kann mittels geeigneter Operationen von $\operatorname{Iso}_{G(2,4)}(m_1,m_2)$ und $H_3$ $m_3$ auf die Form
\[m_3=\left[\begin{smallmatrix}0&0\\ y&0\\ 0&0\\ 0&1\end{smallmatrix}\right]\]
bringen. Es ist $y\neq0$, sonst wäre $m$ zerlegbar. Nach einer weiteren $H_3$-Operation (ersetze $m_1$ durch $m_1+m_2$) folgt
\[\Bild m_1+\Bild m_3=k^4\]
und man ist erneut im 1. Fall.

Man verifiziert leicht, daß alle unzerlegbaren $Q_3$-Darstellungen $M$ mit $\dimv M=(1,i)$ für ein $i\geq 1$ (wegen der Unzerlegbarkeit muß $i\leq3$ gelten) Erweiterung eines Unzerlegbaren und eines Einfachen sind.

Wir haben also (bis auf Dualität) gezeigt, daß alle $Q_3$-Darstellungen der Dimension $\leq6$ Erweiterungen eines Unzerlegbaren und eines Einfachen sind.

\section{Baummoduln in den $U_{Q_3}^{(i)}(3,3)$}\label{rech:baum}
\numberwithin{zaehler}{section}
Wir kommen zum Ende dieser Arbeit noch einmal auf eine der Ausgangsfragen zurück: Gibt es in jeder Schicht von unzerlegbaren Darstellungen eines Köchers $K$ eine Baumdarstellung? Die Schichten von Unzerlegbaren sind dabei per Definition die irreduziblen Komponenten der Varietäten $U_K^{(i)}(\underline{d})$ bestehend aus allen unzerlegbaren $K$-Darstellungen zu einem festen Dimensionsvektor $\underline{d}$ mit konstanter Dimension $i$ des Endomorphismenrings. 

Wir begnügen uns hier damit, für jede Schicht von unzerlegbaren $Q_3$-Darstellungen zum Dimensionsvektor $(3,3)$ in $\operatorname{Rep}_{Q_3}(3,3)$ mit eine unzerlegbare Baumdarstellung anzugeben. Wir verwenden dazu die Beschreibung der $U_z^{(i)}$ aus Abschnitt \ref{rech:unz}. Mit Hilfe der dort angegebenen Gleichungen findet man noch unzählige weitere Baumdarstellungen; wir verzichten auf eine vollständigere Analyse.

\textbf{Eine Baumdarstellung in $U_{Q_3}^{(3)}(3,3)$}\newline
Aus der Beschreibung von $U_z^{(3)}$ für den Repräsentanten $z=z^{(3)}\in B_3$ in Abschnitt $\ref{rech:unz}$ liest man ab:
\[\Bigl(\left[\begin{smallmatrix}1&0&0\\ 0&0&0\\0&0&0\end{smallmatrix}\right],\left[\begin{smallmatrix}0&0&0\\ 1&0&0\\0&0&0\end{smallmatrix}\right],\left[\begin{smallmatrix}0&1&0\\  0&0&1\\1&0&0\end{smallmatrix}\right]\Bigr)\]
liegt in $U_{Q_3}^{(3)}(3,3)$. Berechnung der Träger der Matrizen zeigt, daß die zugehörige Darstellung eine Baumdarstellung ist.

\textbf{Eine Baumdarstellung in $U_{Q_3}^{(2)}(3,3)$}\newline
Wir verwenden die Beschreibung von $U_z^{(2)}$ für den Repräsentanten $z=z^{(7)}\in B_7$ und erhalten eine Baumdarstellung
\[\Bigl(\left[\begin{smallmatrix}1&0&0\\ 0&1&0\\0&0&0\end{smallmatrix}\right],\left[\begin{smallmatrix}0&0&0\\ 1&0&0\\0&0&0\end{smallmatrix}\right],\left[\begin{smallmatrix}0&0&0\\  0&0&1\\0&1&0\end{smallmatrix}\right]\Bigr).\]

\textbf{Eine Baumdarstellung in $U_{Q_3}^{(1)}(3,3)$}\newline
Betrachte die Menge $U_z^{(1)}$ für $z=z^{(13)}\in B_{13}$. Die Darstellung
\[\Bigl(\left[\begin{smallmatrix}1&0&0\\ 0&1&0\\0&0&0\end{smallmatrix}\right],\left[\begin{smallmatrix}0&0&0\\ 1&0&0\\0&1&0\end{smallmatrix}\right],\left[\begin{smallmatrix}0&0&1\\  0&0&0\\0&0&0\end{smallmatrix}\right]\Bigr)\]
liefert die gesuchte Baumdarstellung.

\end{document}